\DocumentMetadata{
    pdfversion=1.7,
    pdfstandard=a-3b,
}

\documentclass[a4paper,fleqn,preprint]{cas-dc}

\usepackage{mathtools}

\usepackage{amssymb}

\usepackage{xcolor}

\usepackage{bm}

\usepackage[range-phrase=--,range-units=single,]{siunitx}

\usepackage{ifluatex}

\usepackage{subcaption}

\usepackage{amsmath,amssymb,amsfonts}
\usepackage{graphicx}

\usepackage{textcomp}

\usepackage{microtype}

\usepackage{multicol}
\usepackage{diagbox}

\ifluatex
    \usepackage{fontspec}
    \usepackage{unicode-math}
    \NewDocumentCommand\vectorStyle{m}{\symbf{#1}}
    \NewDocumentCommand\matrixStyle{m}{\symbf{\MakeUppercase{#1}}}
\else
    \ifxetex
        \usepackage{fontspec}
        \usepackage{unicode-math}
        \NewDocumentCommand\vectorStyle{m}{\symbf{#1}}
        \NewDocumentCommand\matrixStyle{m}{\symbf{\MakeUppercase{#1}}}
    \else
        \usepackage{stix2}
        \NewDocumentCommand\vectorStyle{m}{\mathbf{#1}}
        \NewDocumentCommand\matrixStyle{m}{\mathbf{\MakeUppercase{#1}}}
    \fi
\fi

\usepackage[sorting=none,style=ieee]{biblatex}

\usepackage{cleveref}

\NewDocumentCommand\position{}{\vectorStyle{x}}

\NewDocumentCommand\positions{}{\matrixStyle{X}}

\NewDocumentCommand\currentDensity{}{\vectorStyle{j}}

\NewDocumentCommand\currentDensityMatrix{}{\matrixStyle{j}}

\NewDocumentCommand\primaryCurrentDensity{}{\currentDensity_{\mathrm{p}}}

\NewDocumentCommand\primaryCurrentDensityMatrix{}{\currentDensityMatrix_{\mathrm{p}}}

\NewDocumentCommand\secondaryCurrentDensity{}{\currentDensity_{\mathrm{s}}}

\NewDocumentCommand\conductivity{}{\vectorStyle{\sigma}}

\NewDocumentCommand\electricPotential{}{u}

\NewDocumentCommand\electricPotentialVector{}{\mathbf{u}}

\NewDocumentCommand\transferMatrix{}{\matrixStyle{T}}

\NewDocumentCommand\leadFieldMatrix{}{\matrixStyle{L}}

\NewDocumentCommand\gradient{}{\nabla}

\NewDocumentCommand\divergence{}{\gradient\cdot}

\DeclarePairedDelimiter\args{(}{)}

\DeclareMathOperator\dif{d}

\NewDocumentCommand\Hdiv{}{H(\mathrm{div})}

\NewDocumentCommand\domain{}{\Omega}

\NewDocumentCommand\integral{}{\int}

\NewDocumentCommand\inverse{m}{#1^{-1}}

\NewDocumentCommand\stiffnessMatrix{}{\matrixStyle{A}}

\DeclarePairedDelimiter{\norm}{\lVert}{\rVert}

\NewDocumentCommand\Nof{m}{N_{\text{#1}}}

\NewDocumentCommand\electrodeN{}{\Nof{e}}

\NewDocumentCommand\sourceN{}{\Nof{s}}

\NewDocumentCommand\numberSetStyle{m}{\mathbb{#1}}

\NewDocumentCommand\realNumbers{}{\numberSetStyle{R}}

\NewDocumentCommand\functionsFromTo{mm}{{#1}^{#2}}

\NewDocumentCommand\electrodePosition{}{\position_{\text{e}}}

\NewDocumentCommand\sourcePosition{}{{\position_{\text{s}}}}

\NewDocumentCommand\sourcePositions{}{{\positions_{\text{s}}}}

\NewDocumentCommand\domainPosition{}{{\position_\domain}}

\NewDocumentCommand\measurementVec{}{\vectorStyle{u}_{\mathrm{e}}}

\NewDocumentCommand\dipoleMoment{}{\vectorStyle{d}}

\NewDocumentCommand\windowFunction{}{w}

\DeclareMathOperator\opEMD{EMD}

\NewDocumentCommand\estimate{m}{{#1}^{\operatorname{est}}}

\NewDocumentCommand\ECD{}{\text{ECD}}

\NewDocumentCommand\kernelMatrix{}{\matrixStyle{K}}

\NewDocumentCommand\gramMatrix{}{\matrixStyle{G}}

\NewDocumentCommand\transpose{m}{{#1}^{\mathsf{T}}}

\DeclareMathOperator{\variance}{var}

\NewDocumentEnvironment{Proposed}{+b}{\color{blue}#1}{\relax}

\NewDocumentCommand\trialFunction{}{\psi}

\NewDocumentCommand\trialFunctionVector{}{\boldsymbol\psi}

\NewDocumentCommand\testFunction{}{\phi}

\NewDocumentCommand\trialFunctionPosition{}{\position_\trialFunction}

\DeclarePairedDelimiter\chevrons{\langle}{\rangle}

\NewDocumentCommand\innerProduct{mm}{\chevrons*{#1,#2}}

\NewDocumentCommand\bilinearForm{}{a}

\NewDocumentCommand\differentialOperator{}{D}

\NewDocumentCommand\loadFunction{}{f}

\NewDocumentCommand\loadVector{}{\mathbf\loadFunction}

\NewDocumentCommand\normalVector{}{\mathbf{n}}

\NewDocumentCommand\electrodeIndex{}{\ell}

\NewDocumentCommand\trialFunctionN{}{N_{\trialFunction}}

\NewDocumentCommand\discreteQuantity{m}{\hat{#1}}

\NewDocumentCommand\relativeDifference{}{\Delta_\mathrm{rel}}

\NewDocumentCommand\discreteElectricPotential{}{\discreteQuantity{\electricPotential}}

\begin{document}

\shorttitle{Forward--Inverse Interplay in FEM-Based EEG Source Imaging}

\shortauthors{Santtu Söderholm et~al.}

\title [mode = title]{Advanced EEG Source Models from the Perspective of FEM and Inverse Solutions}                      
\tnotemark[1]

\tnotetext[1]{This work was supported by the Research Council of Finland through the Flagship of Advanced Mathematics for Sensing, Imaging and Modelling (FAME), 2024--2031.}


\author[1]{Joonas Lahtinen}[
    style=finnish,
    orcid={0000-0002-9377-8713}
]

\fnmark[$\dagger$]

\ead{joonas.j.lahtinen@tuni.fi}

\credit{Conceptualization of this study, Methodology, Software, Writing}

\author[1]{Santtu Söderholm}[
    orcid=0009-0004-2918-3733
]

\cormark[1]

\fnmark[$\dagger$]

\ead{santtu.soderholm@tuni.fi}

\credit{Conceptualization of this study, Data curation, Software, Writing}

\affiliation[1]{
    organization={Tampere University},
    addressline={Korkeakoulunkatu 1}, 
    city={Tampere},
    postcode={33720}, 
    state={Pirkanmaa},
    country={Finland}
}

\author[1]{Sampsa Pursiainen}[%
   orcid={0000-0002-9131-9070}
]

\ead{sampsa.pursiainen@tuni.fi}

\credit{Conceptualization of this study, Data curation, Software, Writing - Original draft preparation}


\cortext[cor1]{Corresponding author}

\nonumnote{\textsuperscript{$\dagger$} S. Söderholm and J. Lahtinen contributed equally to this work.}

\begin{abstract}
In this study, we compare forward solutions computed with finite element methods in Zeffiro Interface and DUNEuro, using the different source models they provided. We compared two of the source models from DUNEuro, called Whitney basis and Local subtraction, and the divergence-conforming model of Zeffiro Interface. For source estimation, we applied sparsity-promoting standardized hierarchical adaptive L1 regression (SHAL1R), standardized Kalman filtering (SKF), classical sLORETA, and dipole scanning. Analyses include Earth Mover's Distance, depth bias scatter plots, and qualitative assessments of amplitude distribution and focality.
    Preliminary experiments with source interpolation for each method revealed that Local subtraction closely matches expectations for the local behavior of the lead field at various depths.
The main results reveal that the success of an inverse method depends strongly on the compatibility between its assumptions about the focality of the source and the chosen source model, with point-source models performing best when paired with methods designed for such sources, i.e., sparsity-promoting methods and methods that scan for a single source. Moreover, source models that admit patch sources with a wide distribution are more sensitive to additional noise.
\end{abstract}

\begin{keywords}
electroencephalography\sep
finite element method\sep
inverse problems\sep
source models
\end{keywords}

\maketitle

\section{Introduction}
Electroencephalography (EEG) source imaging aims to estimate the underlying brain activity, i.e., the source configuration, using electrical potentials measured from the scalp as observed data. The problem at hand is highly ill-posed, as different source configurations and their superpositions can yield identical sensor readings; moreover, the solutions are susceptible to noise and modelling inaccuracies \cite{kaipio2006statistical}. From a mathematical point of view, the main influencing factors on the accuracy of the source estimate are (1) the lead field matrix, denoted here as $\leadFieldMatrix\in\functionsFromTo\realNumbers{\electrodeN\times\sourceN}$, which can be viewed as the operator that maps the source activity within the brain to the potentials captured by the sensors on the scalp. And (2), the selected inversion method.

A benchmark inversion method, \emph{standardized low-resolu\-tion brain electromagnetic tomography} (sLORETA) \cite{pascual2002sloreta}, which is very meritorious in finding the source in clinical settings \cite{deGooijer2013,Coito2019}, has served as a motivator for many recent standardized methodologies \cite{SSLOFO2005Liu,lahtinen2024shalpr,lahtinen2024standardized}. The appeal of standardization stems from its theoretical ability to reconstruct sources without depth bias from noiseless data \cite{Pascual2007sqrtm,lahtinen-standardization-2024}.

Recent research has revealed how we model the forward problem can significantly shape the accuracy of source reconstruction in EEG. Notably, Söderholm \emph{et al.} demonstrated that specific preprocessing techniques --- such as anatomical "peeling", which modifies the positions where synthetic sources can be positioned within a head model --- can impact the outcomes of FEM-based EEG source analysis \cite{soderholm2024effects}. These findings underscore that decisions made during forward modeling affect not only the precision but also the spatial distribution of the estimated brain activity.

The finite element method (FEM) offers a versatile and precise approach to forward modeling by allowing the use of highly realistic head geometries and detailed tissue conductivity profiles \cite{he2020zeffiro,schrader2021duneuro}. In this context, the choice of source model plays a crucial role in shaping the inversion results. For instance, divergence-conforming $\Hdiv$ models, a quadratic refinement of the linear Whitney basis, along with partial integration and St.~Venant. Each model has its own benefits and shortcomings in terms of spatial and numerical accuracy \cite{miinalainen2019realistic}. These distinctions directly reflect on the localization error and spatial distribution of the source estimates. Recent advances in forward modeling, such as the Local subtraction method, have introduced a robust and efficient strategy for modeling singular dipolar sources in FEM \cite{holtershinken2025local}.

In this paper, we compare cutting-edge EEG source modeling techniques. By leveraging the theoretical unbiasedness of standardized estimations and dipolar fitting in Dipole Scan, we show how the $\Hdiv$ source model in Zeffiro Interface compares with two alternatives in DUNEuro: the Whitney source model (with both face-intersecting and edgewise basis functions) and Local subtraction. To deepen our insights into the source-modeling differences between the implementations, we introduce a novel patch-type forward mode constructed using a high-density Local subtraction model that captures distributed, patch-like neural activity. Our findings show that mismatches between source models and the inversion method's assumptions about the spread of neural activity can impair source estimation performance.

\section{Background}\label{sec.background}

\begin{figure}
    \centering
    \includegraphics[width=0.5\linewidth]{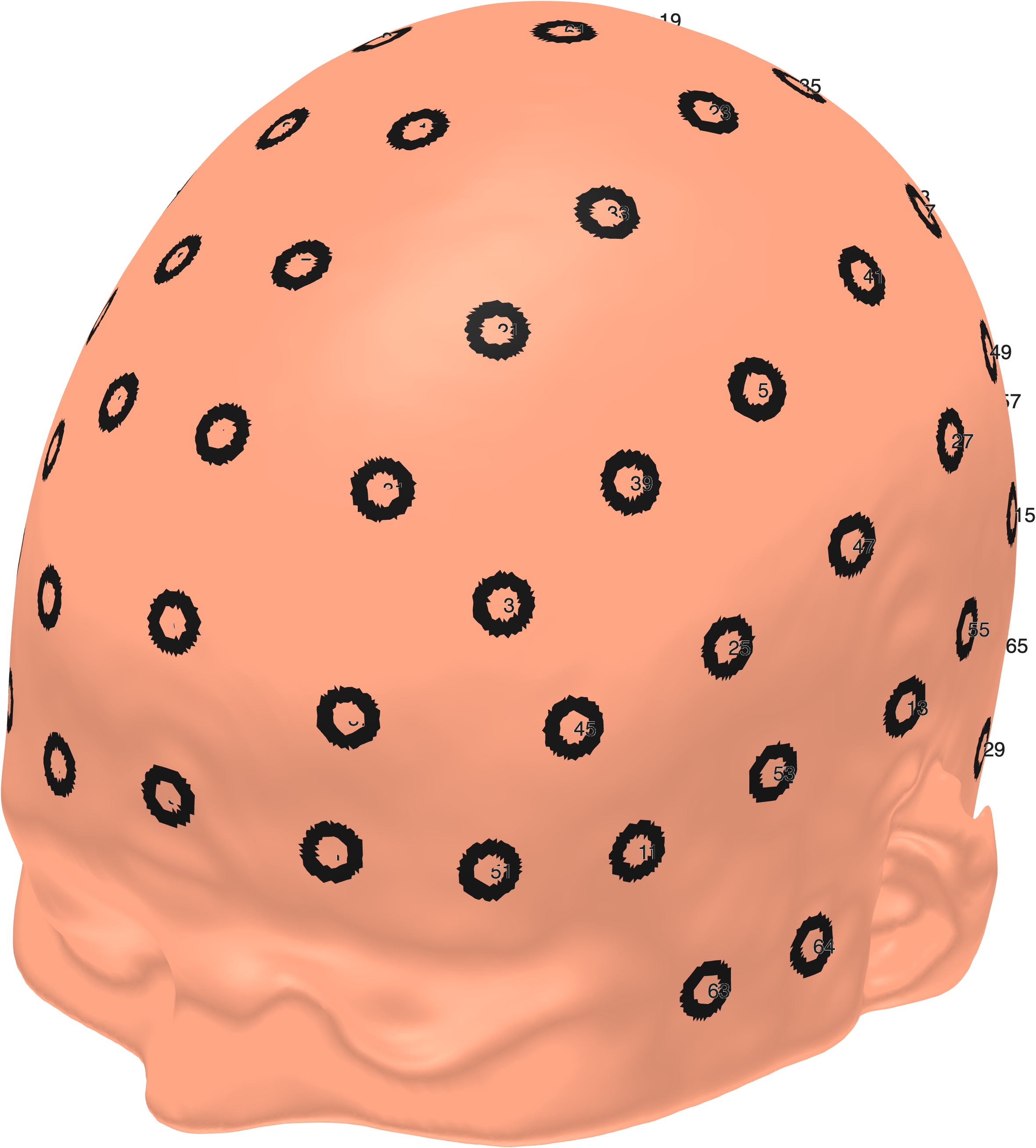}
    \caption{Electrodes positioned on the scalp of the ICBM152 2009a head model \cite{icbm-website}, according to the international 10--10 electrode positioning standard.}
    \label{fig.scalp.electrodes}
\end{figure}

The pursuit of measuring brain activity through electrical potentials has a history spanning over a century, with its origins traced to early experiments in Germany \cite{vergani-2024}. Initial approaches were inherently invasive, relying on the insertion of depth electrodes to detect neural signals given the technical limitations of the time. Although depth electrodes remain the standard for applications demanding precise localization \cite{missey-etal-2026}, contemporary EEG is predominantly conducted noninvasively, with electrodes positioned on the scalp \cite{Knosche-Haueisen-2022}, as illustrated in Figure~\ref{fig.scalp.electrodes}. The benefit of source imaging from non-invasive measurements lies in preliminary diagnostics and in guiding retrospective FLAIR-MRI or Zoomed-MRI, optimizing the location of the skull opening or the pathing to the target in neurosurgery, depth electrode placement \cite{RineyCatherine2012,Aydin2017,Mouthaan2019,Diamond2023}, since depth electrodes can detect activity only from close range \cite{Coito2019}. The steadily increasing benefits have been enabled by substantial advances in both measurement instrumentation and sophisticated signal processing algorithms \cite{hari-and-puce-2023, chaddad-etal-2023}.

The mathematical modelling of EEG sources has evolved over several decades \cite{clark-and-plonsey-1966, geselowitz-1967}. The development of the Equivalent Current Dipole (ECD) model \cite{scherg-1985,hamalainen1993magnetoencephalography}, coupled with the advent of modern computers, enabled practical solutions to the source localization problem. Forward modelling has progressed from simplistic spherical representations of two to four conductivity layers \cite{ary-etal-1981, cuffin2001spherical} to anatomically accurate domains derived from MRI segmentation \cite{cuffin2001realistic}. Computational approaches have similarly advanced, moving from boundary element methods (BEM) \cite{demunck-1992,schlitt-1995} to more sophisticated finite element methods (FEM) \cite{wolters-etal-2006,wolters-etal-2007}. While the ECD remains an effective model for highly \emph{localized} sources of brain activity, it is inadequate for describing \emph{large-scale} or patch-like activation patterns \cite{Knosche-Haueisen-2022}, and it is sensitive to modeling errors and measurement noise \cite{Lahtinen2023}. To address this limitation, alternative source models have been introduced to more accurately capture spatially extended neural processes.

One such model is the \emph{multi-pole expansion} \cite{vorwerk2019multipoles, Knosche-Haueisen-2022} or MPE approach, which bears a close resemblance to the St.~Venant forward interpolation scheme \cite{buchner1997stvenant}. Here, a truncation of the Taylor expansion of the coefficient $1 / \norm{\position - \domainPosition}$ in the volume conductor potential
\begin{equation}
    \electricPotential
    \args{\position}
    =
    -
    \frac{1}{4\pi\conductivity}
    \integral_\domain
    \frac
    {\divergence
    \primaryCurrentDensity}
    {\norm{\position - \domainPosition}}
    \dif\domain
\end{equation}
is used in place of the coefficient itself. The variable $\position$ is a  position in the conductor and $\domainPosition$ a position within the subdomain $\domain$ where a primary current density $\primaryCurrentDensity=\primaryCurrentDensity\args\domainPosition$ corresponding to brain activity can be located. The variable $\conductivity$ denotes conductivity. Where the Taylor expansion is truncated determines the order of the multipole: the first term describes a monopole, the second term a dipole, the third term a quadrupole, and so forth. The names of these different degrees come from the kinds of topographic potential patterns their use results in \cite{Knosche-Haueisen-2022}. As the order of the multipole increases, so does its decay rate as one is moving away from a source position $\sourcePosition\in\domain$. Therefore, restricting oneself to the lower-order terms corresponds well with ECDs, while increasing the number of included expansion degrees and the extent of $\domain$ can lead to a more patch-like source. However, increasing the number of terms and size of $\domain$ might lead to non-convergence of the related series, if the distances within $\domain$ extend past source--electrode distances $\norm{\electrodePosition-\sourcePosition}$ \cite{Knosche-Haueisen-2022}. The higher-order terms of the expansion are also harder to interpret physically, outside of looking at their topographic patterns.

A slightly better approach in modelling patch-like sources in terms of physical interpretability and ease of implementation in numerical solvers is to simply split the integration volume $\domain$ into further subdomains $\domain_i$, such that the extent of each $\domain_i$ is small enough to be interpreted as focal, housing a single ECD representation \cite{Knosche-Haueisen-2022}. Then the entire patch-like source could be described as a superposition or mean of each individual $\ECD_i$ dipole moment $\dipoleMoment_i$. This is referred to as \emph{distributed source modelling} (DSM) \cite{
    hämäläinen1994mne,
    lin2006dsm,
    lin2021dsm,
} or more specifically in cortically constrained 3D contexts, \emph{current density modelling} (CDM) \cite{Knosche-Haueisen-2022}. This kind of discrete summation approach by itself is rather coarse when the underlying assumption is that within $\domain$ the primary current distribution is continuous. However, one may also apply a smoothing kernel $\kernelMatrix$ to the primary current distribution $\primaryCurrentDensityMatrix$, or spatial convolution with a decreasing window function $\windowFunction$ centered at each source position over this discrete field and source positions $\sourcePosition$ to generate output which is in better conformance with actual measurements \cite{
    Knosche-Haueisen-2022,
    soderholm2026interference,
}:
\begin{equation}\label{eq.smoothing.application}
    \dipoleMoment\args\sourcePosition
    =
    \integral_\domain
    \windowFunction\args{\domainPosition - \sourcePosition}
    \cdot
    \primaryCurrentDensity\args\domainPosition
    \dif\domain
    \approx
    \kernelMatrix
    \primaryCurrentDensityMatrix
    \,.
\end{equation}
The potential measurements $\measurementVec$ at the electrodes would then be produced using a \emph{forward operator} $\leadFieldMatrix$:
\begin{equation}\label{eq.measurement.simulation}
    \measurementVec
    \approx
    \leadFieldMatrix
    \dipoleMoment
    \,.
\end{equation}

Inverse solvers originally relied on a form of dipole fitting, where a semi-analytical forward solution is compared with a measured one and the dipole model parameters are optimized such that the difference between the forward solution and the measured signal is minimized \cite{ary-etal-1981}. More recent inverse solvers either rely on the idea of regularizing a Gram matrix $\gramMatrix=\leadFieldMatrix\transpose\leadFieldMatrix$ with a linear forward operator $\leadFieldMatrix$ \cite{sekihara-etal-2008}, or use Bayesian modeling with more complicated a priori assumptions and iterative algorithms in the seeking of maximum a posteriori that maximizes the log-posterior density function \cite{kaipio2006statistical,Knosche-Haueisen-2022}.

The construction of a lead field $\leadFieldMatrix$ using FEM involves the following steps \cite{brenner2008fembook, peterson2022basisfunctions, Knosche-Haueisen-2022}:
\begin{enumerate}
    \item The governing partial differential equation (PDE) $\differentialOperator\electricPotential = \loadFunction$ with a differential operator $\differentialOperator$, a solution $\electricPotential$ and a load given by $\loadFunction$ is transformed into a \emph{weakened form}, where the equation is scaled by a \emph{test function} $\testFunction$ from a smooth-enough space, integrated over the domain $\domain$ and then reduced into a bilinear form $\bilinearForm\args{\testFunction,\electricPotential}$ with lower-order derivatives through integration by parts.
    \item The weakened PDE is \emph{discretized} using FEM, where the solution $\electricPotential$ is approximated as a finite linear combination of smooth-enough locally supported \emph{trial} or \emph{basis functions} $\trialFunctionVector_i$ and their coefficients $\electricPotentialVector_i$: $\discreteElectricPotential = \sum_i\electricPotentialVector_i\trialFunctionVector_i$.
    \item This discretized solution is substituted into the weak bilinear form $\bilinearForm\args{\testFunction,\discreteElectricPotential}$, and the linearity of the bilinear form allows one to form a system matrix $\stiffnessMatrix\in\functionsFromTo\realNumbers{\trialFunctionN\times\trialFunctionN}$ that can be inverted against a load vector $\loadVector\in\functionsFromTo\realNumbers{\trialFunctionN}$ to produce a solution for the coefficients:
    \begin{equation}\label{eq.fem.solution}
        \electricPotentialVector
        =
        \inverse\stiffnessMatrix
        \loadVector
        \,.
    \end{equation}
    \item The solution $\electricPotentialVector\in\functionsFromTo\realNumbers{\trialFunctionN}$ from \eqref{eq.fem.solution} is anchored to the trial functions $\trialFunction$, whose \emph{centers of support} $\trialFunctionPosition\in\Omega$ are attached to the discretized computational domain $\domain$, for example, the mesh nodes when the shape of $\trialFunction$ is Lagrangian, or cell facets when the shape of $\trialFunction$ follows a divergence-conforming interpolation. This means that a separate interpolation step is required if the solution is to be evaluated at arbitrary positions within $\domain$.
\end{enumerate}

There are a few specifics related to the EEG forward problem worth mentioning \cite{Knosche-Haueisen-2022}:
\begin{enumerate}
    \item The current loops inside of a human head are assumed to be closed, which allows us to write the divergence of the total current density as $\divergence\currentDensity = \divergence\args{\primaryCurrentDensity + \secondaryCurrentDensity} = 0$ or $\divergence\primaryCurrentDensity = - \divergence\secondaryCurrentDensity = - \divergence\args{\conductivity\gradient\electricPotential}$, where the primary current density $\primaryCurrentDensity = \primaryCurrentDensity\args\sourcePosition$ is the brain activity interpreted as the driving force $\loadFunction$.
    \item The EEG boundary condition states that there is no current flow out of $\domain$ in the normal direction $\normalVector$: $\innerProduct\secondaryCurrentDensity\normalVector = 0$.
    \item A forward solution could be formed based on the two assumptions above, but the naive interpretation of $\loadFunction$ is computationally too expensive because of the large number of possible $\trialFunctionPosition$ in comparison to the number of electrodes. We therefore form $\loadVector$ by relying on \emph{Helmholtz reciprocity} \cite{wolters2004transfer, gross2023reciprocity, Knosche-Haueisen-2022}: a unit current is fed in through a supposed measurement electrode contact surface and out of a ground electrode. This corresponds to setting the matching indices of $\loadVector$ such that they sum up to \num 1 and \num{-1} respectively \cite{weinstein2000reciprocity}. This is done for each electrode $\electrodeIndex$ separately and the $\electricPotentialVector$ from \eqref{eq.fem.solution} assembled into the columns $\transferMatrix_\electrodeIndex$ of a pre-computed \emph{transfer matrix} $\transferMatrix\in\functionsFromTo\realNumbers{\trialFunctionN\times\electrodeN}$, which is a solution from the electrodes to $\trialFunctionPosition$ in $\domain$. Because $\stiffnessMatrix$ is symmetric, $\transferMatrix$ can simply be transposed to produce a solution from $\trialFunctionPosition$ to the electrodes.
    \item A lead field $\leadFieldMatrix\in\functionsFromTo\realNumbers{\electrodeN\times\sourceN}$ is an interpolation of the rows of $\transpose\transferMatrix$ to a set of source positions not attached to the mesh, as in from $\trialFunctionPosition$ to $\sourcePosition$.
\end{enumerate}
When constructed using reciprocity like this, the solution $\electricPotentialVector$ should appear strongest near the electrodes. This is because electric potential falls off with the inverse distance from a source. The trial basis used for the potential and therefore for the construction of $\stiffnessMatrix$ is Lagrangian with first-order interpolation, and therefore $\trialFunctionPosition$ coincides with the mesh nodes.

\section{Methods}\label{sec.methods}

The lead field matrices used in this study are computed using the FEM-based transfer matrix approach implemented in both of the examined software packages, DUNEuro \cite{schrader2021duneuro} and Zeffiro Interface \cite{he2020zeffiro}. The inversion methods used in this study are provided by Zeffiro Interface. The models provided by DUNEuro are transferred to Zeffiro Interface via the package-connecting
    set of MATLAB functions
that the authors produced for this study.
    A tetrahedral mesh was produced by first extracting a set of triangular brain compartment surfaces from the openly available, realistic, population-based ICBM152 2009a multicompartment head model \cite{fonov-2011,fonov-2009,icbm-website}. A combination of FreeSurfer \cite{fischl2012freesurfer} and SimNIBS \cite{puonti2020simnibs} performed this extraction, after which co-registering these surfaces into the same RAS coordinate system was done based on the affine transformations describing patient orientation in the MRI files. Finally, the mesh generation and compartment subsetting algorithms of Zeffiro Interface \cite{galazprieto2023zeffiromesh} were applied to label the tetrahedra it generates into correct compartments with appropriate conductivity values $\sigma$. This mesh was then used as the computational domain of transfer-matrix-based lead field solvers of Zeffiro Interface and DUNEuro \cite{schrader2021duneuro}, through a MATLAB function written for this purpose.
The point electrode model is selected as the model of the electrodes \cite{hanke-2011,agsten2018electrodes}, which were positioned according to the international 10--10 system \cite{Knosche-Haueisen-2022}: a set of electrode points in the general shape of a 10--10 electrode cap was projected to the nearest scalp surface triangle centroids.

    The focus of the study is on different \emph{source models}, which in the context of FEM correspond to different combinations of trial $\trialFunction$ and test function $\testFunction$ spaces and their interpolations or shapes \cite{wolters2004transfer,Knosche-Haueisen-2022}, different weak forms $\bilinearForm$, together with different transfer matrix interpolation schemes after an initial FEM solution has been produced. Chosen source model implementations include the face-intersecting and edge-wise Whitney basis functions and the continuous-Galerkin-based Local subtraction (LS) approach \cite{holtershinken2025local} provided by DUNEuro, and the divergence-conforming $\Hdiv$ \cite{miinalainen2019realistic} from Zeffiro Interface.
    Sources in \num{10000} fixed positions in the shape of a hexahedral lattice but restricted within tetrahedra of active compartments, with unit orientation coordinates along the standard Cartesian axes were used to form the \emph{source space}, to which $\transferMatrix$ is interpolated via Position-Based Optimization (PBO) \cite{bauer2015whitney,Pursiainen-etal-2016} due to its ability to take source orientation better into account in comparison to alternative interpolation schemes. The source spaces are the same for each dipole-like source model.

To emulate patch-like sources using the DSM method outlined in Section~\ref{sec.background}, we first constructed a high-density lead field matrix comprising \num{160000} Cartesian sources. Here, density indicates that a large number of source positions $\sourcePosition$ are distributed evenly throughout the active gray matter areas of the head model. The rows of this dense ECD-based  $\leadFieldMatrix$ were then smoothed over $\sourcePosition$ using a Gau\ss ian window function in accordance with \eqref{eq.smoothing.application}:
\begin{equation}\label{eq.gaussian}
    \windowFunction
    \args{
        \position,
        \sourcePosition,
        \variance
    }
    \propto
    \exp\args * {
        -
        \frac 1 2
        \frac{
            \norm{
                \position
                -
                \sourcePosition
            } ^ 2
        }{
            \variance
        }
    }
    \, .
\end{equation}
Here $\variance$ denotes the squared width or variance of the window. The variance value \qty{45.5}{\milli\meter^{2}} was selected such that the sources from \qty{1}{\centi\meter} away have contribution of one third. The resulting smoothed lead field matrix was then projected to the source space of \num{10000} sources used with other models and in our experiments.

In our experiments, we consider the differences, advantages, and disadvantages of the source models from an inversion perspective. Namely, we focus on the spatial smoothness, spread, and depth bias of inversion estimations with various methods.

    In our preliminary analysis of the lead fields, we observe how the output of the utilized transfer-matrix-based lead field solvers behaves. We first visually observe the column norms of the lead fields used in experiments involving inversion, which were constructed with a source space that is in the shape of the vertices of a hexahedral lattice, but sampled at points that are located within brain compartments that are determined \emph{active}.

    We also observe how slight changes in the positioning of sources within one tetrahedral element and in the \qty{3}{\milli\meter} range can affect both the visual appearance of a lead field column norms and their relative differences
    \begin{equation}\label{eq.reldiff}
        \relativeDifference\args{
            \leadFieldMatrix_a,
            \leadFieldMatrix_b
        }
        =
        \frac{
            \norm{
                \leadFieldMatrix_a
                -
                \leadFieldMatrix_b
            }
        }{
            \norm{\leadFieldMatrix_a}
        }
    \end{equation}
    between a sub-lead field $\leadFieldMatrix_a\in\functionsFromTo\realNumbers{\electrodeN}$ acting on source position space $\positions_a$ and $\leadFieldMatrix_b\in\functionsFromTo\realNumbers{\electrodeN}$ acting on a perturbed space $\positions_b$, with $\positions_a$ being a reference set. For this perturbation experiment, we chose tetrahedra completely embedded into the active compartments, meaning all of their 4 facet-based neighbors are within the same compartment. Sources were then positioned at the first vertices $\positions_1$ of these tetrahedra, where the vertex numbering is determined by the Zeffiro Interface mesh generator \cite{galazprieto2023zeffiromesh}. The triangular faces opposing each vertex 1 were then sought, and the vertices in the neighboring tetrahedra sharing this face but opposing it were chosen as positions $\positions_2$ of a perturbed source. Sources were also placed at the midpoints $\positions_3$ of these two source positions between each pair of tetrahedra. Finally, a union of all these 3 source spaces $\positions_4$ was taken. Lead fields corresponding to all of these 4 source spaces were generated for both Whitney and Local subtraction after sampling the source spaces at regular intervals such that the sizes of each lead field $\leadFieldMatrix_{1 \to 3}$ were roughly \num{10 000} and for $\leadFieldMatrix_4$ it was \num{30 000}. These were then compared with each other.

In experiment I, we examine the spatial properties of inverse solutions in the noiseless case and under inversion crime using an anisotropic conductivity profile, where the local smoothness and relation of the neighboring sources in the forward model and source model are the only factors impacting the smoothness and spread of the source estimate. In addition, we compare the distribution of the estimate using the estimate obtained with Local subtraction as the reference. In practice, the estimated source distribution of Local subtraction is subtracted from the estimated source distribution obtained from the other models.

    Experiment II consists of positioning \num{100} randomly distributed sources with standard Cartesian orientations at \qtyrange{0}{60}{\milli\meter} relative heights from the lowest electrode $z$-coordinate at \qty{5}{\milli\meter} intervals, for a total of \num{1200} sources. A FEM-based lead field $\leadFieldMatrix$ with an anisotropic conductivity $\conductivity$ is then built for the purpose of generating simulated measurements. A corresponding isotropic $\leadFieldMatrix$ to be used by the inverse methods in question is also built to avoid an inverse crime \cite{kaipio2006statistical}. Before inversion, the synthetic measurements are generated both with and without Gau\ss ian noise, where the noisy data has a signal-to-noise ratio (SNR) of \qty{15}{\decibel}, with a mean of \qty{18}{\percent} of the signal peak to make the effects on the results obviously visible. A lower amount of noise at \qty{20}{\decibel} SNR has been shown to be insufficient to differentiate between noisy and noiseless reconstructions \cite{Lahtinen2023,lahtinen2024shalpr,lahtinen2024standardized}. Each synthetic source is then reconstructed one at a time with each inverse method.

The source models are compared using the depth-bias scatter plot (DBSP) presented in \cite{LuckaFelix2012HBif,elvetun2025depthbias} and the \emph{Earth Mover's Distance} (EMD) to measure the spread of the reconstructions. Earth Mover's Distance is a measure of similarity between two distributions, computed between the distribution of true sources and the estimate \cite{rubner-etal-1998-emd}:
\begin{equation}\label{eq.emd}
    \opEMD
    \args * {
        \sourcePositions,
        \estimate\sourcePositions
    }
    =
    \frac{
        \sum_i
        \sum_j
        \left\|\sourcePosition_{i}-\estimate\sourcePosition_{j}\right\|
        F_{i,j}
    }{
        \sum_i
        \sum_j
        F_{i,j}
    }
    \,,
\end{equation}
subject to the following constraints
\begin{align}
    F_{i,j}&\geq 0, \quad \forall i\in [1,n_a],\: j\in [1,n_b]\\
    \sum_i F_{i,j}&\leq w_{j}^{(a)},\quad \forall j\in [1,n_b]\\
    \sum_j F_{i,j}&\leq w_{i}^{(b)}, \quad \forall i\in [1,n_a]
\end{align}
Here $i$ corresponds to a source position $\sourcePosition_i\in\sourcePositions$ and $j$ to a reconstructed source position $\estimate\sourcePosition_j\in\estimate\sourcePositions$.
The quantities $w^{(a)}_j$ and $w^{(b)}_i$ are the weights of each estimated source position, i.e., the magnitude of the estimate at the source position and the magnitude of the true source, respectively. $F_{i,j}$ is a non-negative flow that minimizes the expression in the numerator of \eqref{eq.emd}. In non-technical terms, EMD measures the mechanical work required to transfer a distribution of matter (here, the distributional source estimate) to the pre-assigned "holes," which correspond to the distribution of the true activity.

For source reconstruction, we use standardized distributed methods including sLORETA \cite{pascual2002sloreta,sekihara-etal-2008}, the sparsity-promoting SHAL1R \cite{lahtinen2024shalpr}, and spatiotemporal SKF \cite{lahtinen2024standardized} with classical Dipole Scanning (DS) \cite{fuchs-1998,neugebauer-etal-2022, Knosche-Haueisen-2022}. Because DS is observed to be sensitive mainly to the difference between forward and inverse models \cite{Lahtinen2023}, we omit it in Experiment II.

The motivation for displaying DBSP for standardized methods stems from the theoretical result that the estimated source location should be unbiased, meaning we obtain perfect agreement between the true source depth and the estimated depth \cite{Pascual2007sqrtm}. We assume that greater deviation from perfect agreement indicates greater ambiguity among sources, which is undesirable. In EMD terms, a lower value indicates a smaller spread. The evolution of EMD, along with the depth of the true source, is also an important indicator.



\section{Results}\label{sec.results}
\subsection{Analysis of the lead field column norms}
\begin{figure}
    \centering
    \hfill
    \begin{subfigure}[b]{0.4\linewidth}
        \includegraphics[width=\linewidth]{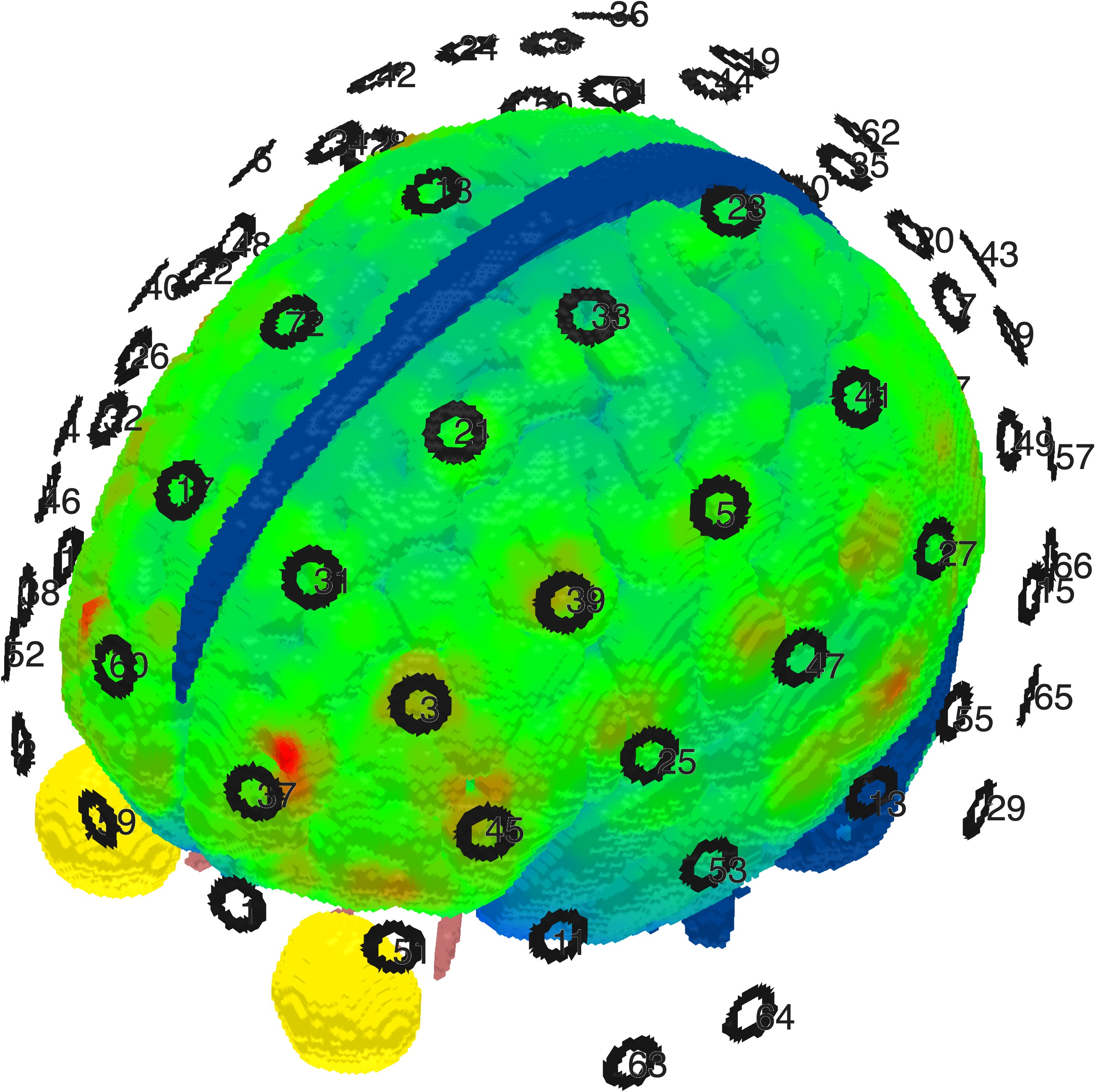}
        \caption{}
        \label{fig.duneuro.L.vecnorm.superficial}
    \end{subfigure}
    \hfill
    \begin{subfigure}[b]{0.4\linewidth}
        \includegraphics[width=\linewidth]{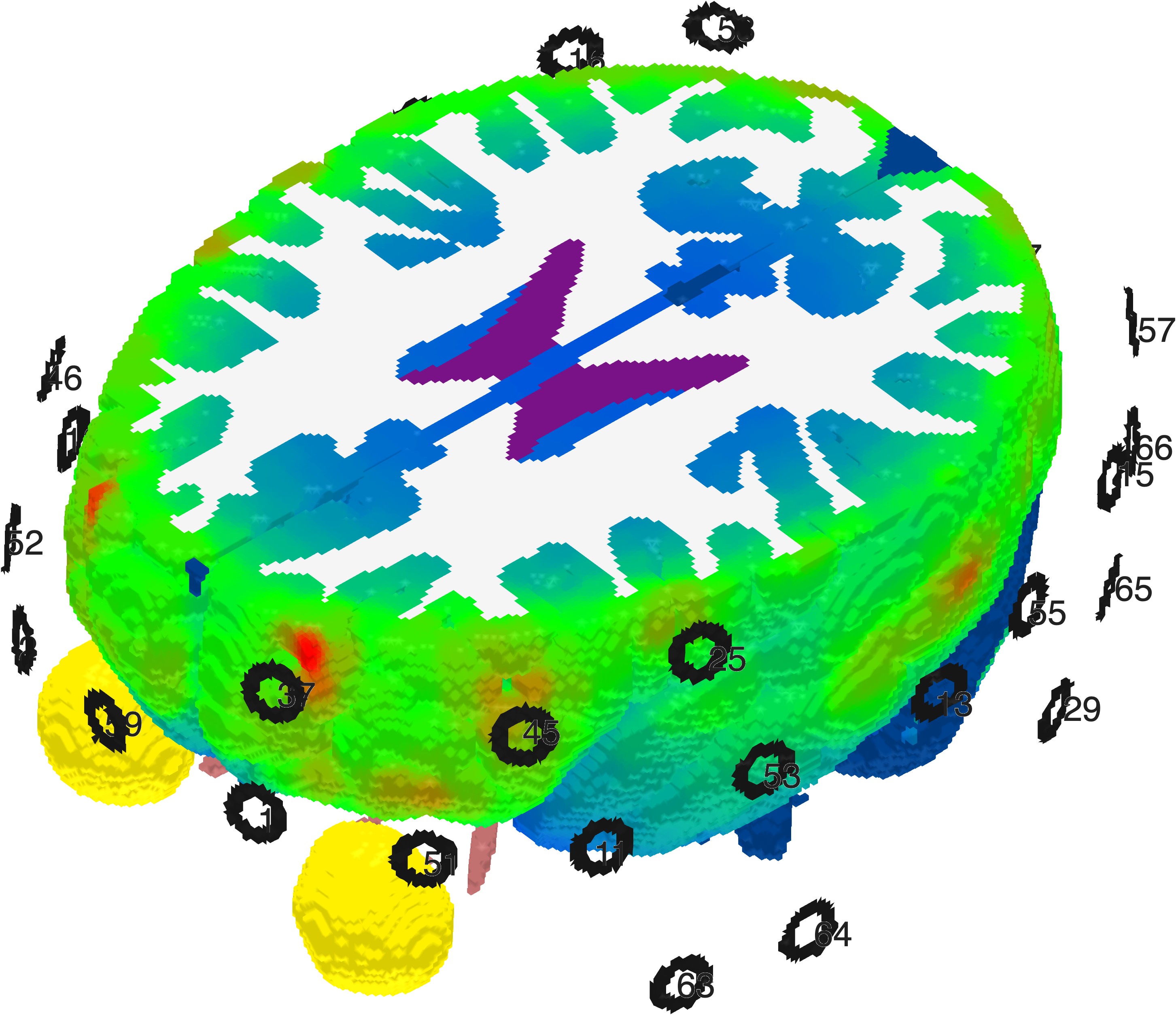}
        \caption{}
        \label{fig.duneuro.L.vecnorm.deep}
    \end{subfigure}
    \hfill
    \begin{minipage}[b]{0.1\linewidth}
        \includegraphics[height=4cm]{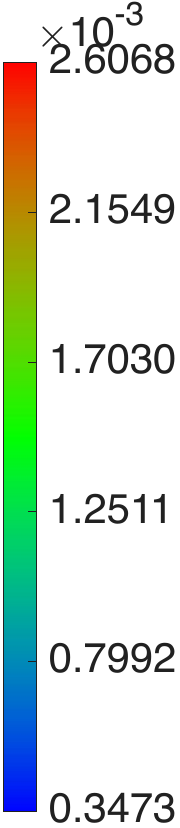}
    \end{minipage}
    \hfill
    \caption{A display of the column norms $\norm\leadFieldMatrix_2$ of a Whitney lead field $\leadFieldMatrix$ produced with DUNEuro. Subfigure~\ref{fig.duneuro.L.vecnorm.superficial} displays a cortical view while \ref{fig.duneuro.L.vecnorm.deep} shows a deeper cross-section in the RA-plane of the Right-Anterior-Superior (RAS) coordinate system, roughly at the height of the thalamus. The field is strong near the EEG electrodes depicted as black circles, and decays smoothly when moving away from them.}
    \label{fig.duneuro.L.vecnorm}
\end{figure}

\begin{figure}
    \centering
    \hfill
    \begin{subfigure}[b]{0.4\linewidth}
        \includegraphics[width=\linewidth]{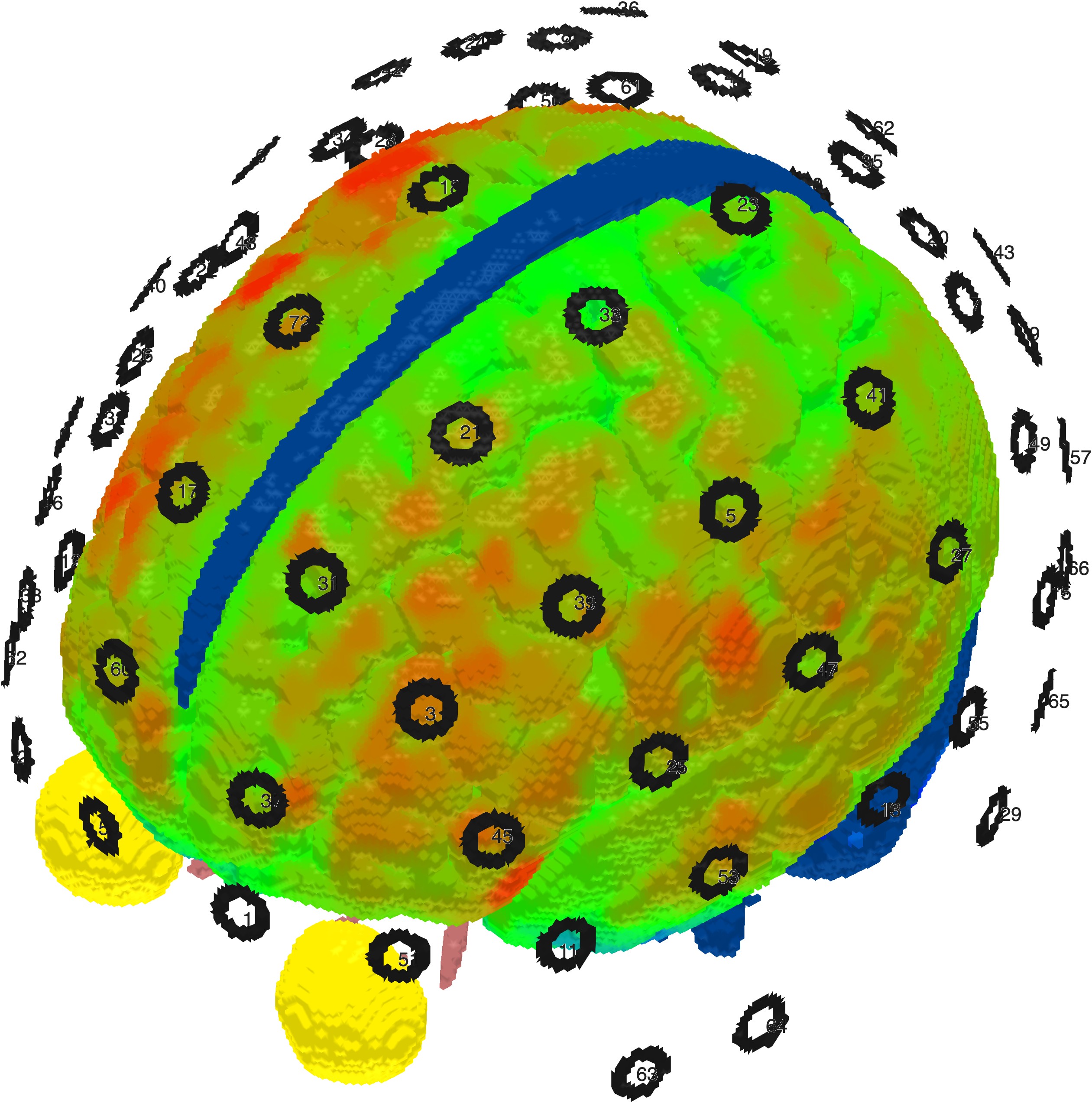}
        \caption{}
        \label{fig.zeffiro.L.vecnorm.superficial}
    \end{subfigure}
    \hfill
    \begin{subfigure}[b]{0.4\linewidth}
        \includegraphics[width=\linewidth]{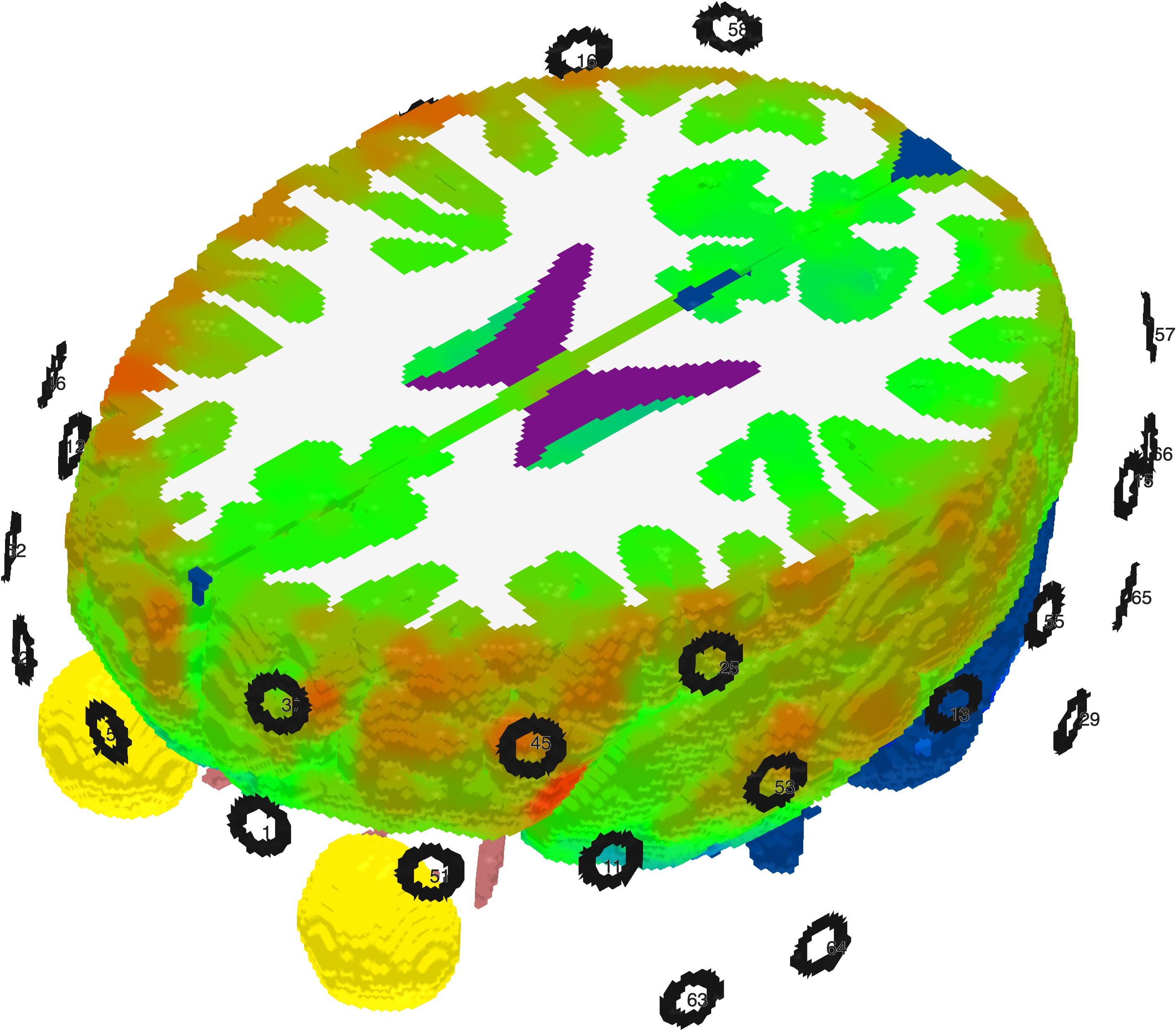}
        \caption{}
        \label{fig.zeffiro.L.vecnorm.deep}
    \end{subfigure}
    \hfill
    \begin{minipage}[b]{0.1\linewidth}
        \includegraphics[height=4cm]{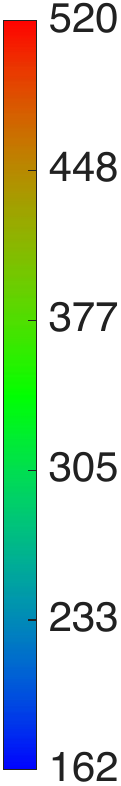}
    \end{minipage}
    \hfill
    \caption{A display of the column norms $\norm\leadFieldMatrix_2$ of a $\Hdiv$ lead field $\leadFieldMatrix$ produced with Zeffiro Interface, up to the \qty{95}{\percent} quantile of the field strenght, with higher values filtered out. The filtering was done due to an unwarranted field maximum occuring behind the yellow eye balls. The viewports in \ref{fig.zeffiro.L.vecnorm.superficial} and \ref{fig.zeffiro.L.vecnorm.deep} are the same as in \ref{fig.duneuro.L.vecnorm.superficial} and \ref{fig.duneuro.L.vecnorm.deep}, respectively. We see a higher irregularity in the field values when compared to DUNEuro, in addition to a scale of values that is \num{100000} times stronger. The field decays relatively slower when moving away from the electrodes.}
    \label{fig.zeffiro.L.vecnorm}
\end{figure}

Figure~\ref{fig.duneuro.L.vecnorm} (a--b) displays a column norm of the Whitney-basis lead field matrix on the geometric brain domain. Within this value scale, the Local subtraction produces a visually identical norm representation; hence, it is omitted. The column norm for $\Hdiv$ is presented in a similar manner in Figure~\ref{fig.zeffiro.L.vecnorm} (a--b). The column norm of the Whitney-basis lead field follows the expectations: observed electrical potential decreases with the inverse distance from the sensors.

\begin{figure}
    \centering
    \hfill
    \begin{subfigure}[b]{0.4\linewidth}
        \includegraphics[width=\linewidth]{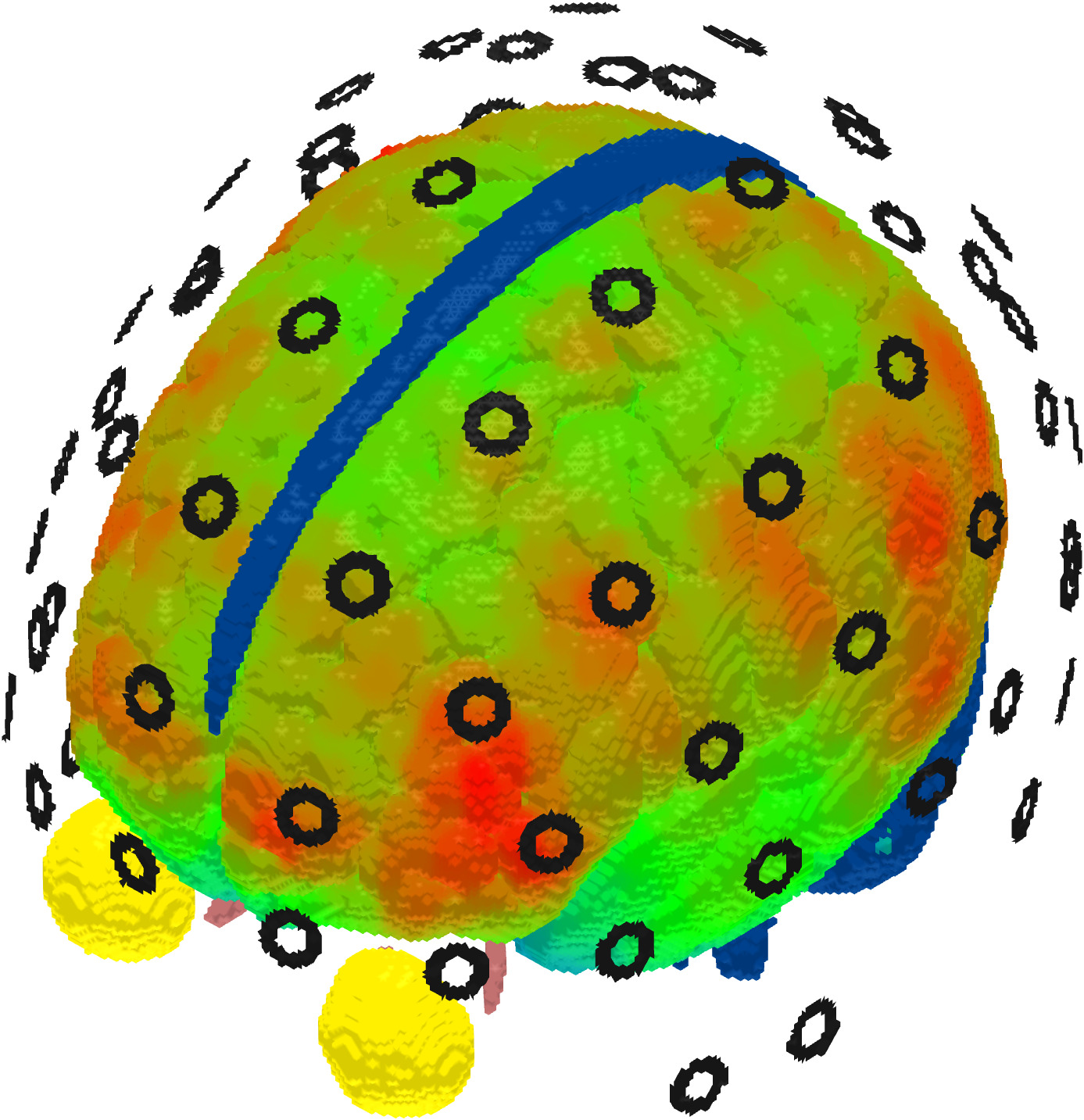}
        \caption{}
        \label{fig.patch.L.vecnorm.superficial}
    \end{subfigure}
    \hfill
    \begin{subfigure}[b]{0.4\linewidth}
        \includegraphics[width=\linewidth]{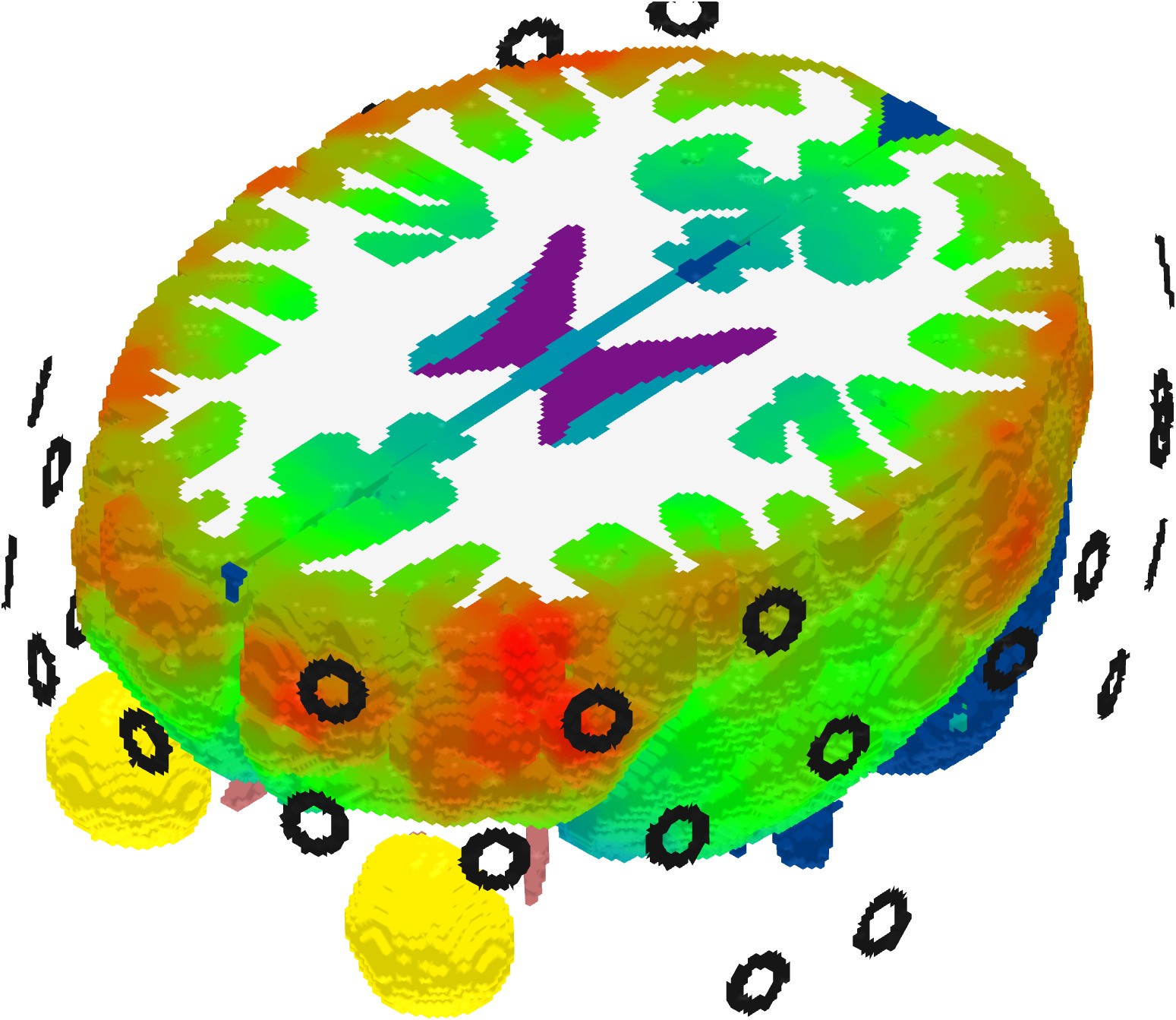}
        \caption{}
        \label{fig.patch.L.vecnorm.deep}
    \end{subfigure}
    \hfill
    \begin{minipage}[b]{0.1\linewidth}
        \includegraphics[height=4cm]{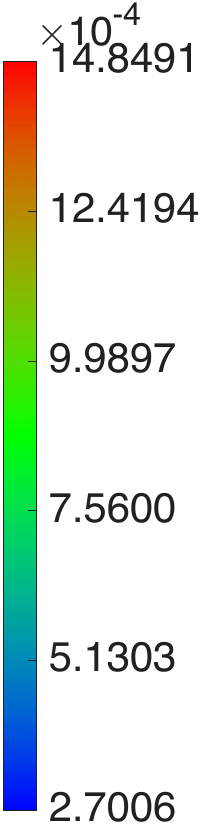}
    \end{minipage}
    \hfill
    \caption{A display of the column norms $\norm\leadFieldMatrix_2$ of the surrogate patch lead field $\leadFieldMatrix$ produced using fine resolution Local subtraction forward solution, spatial Gau\ss ian window smoothing, and projection to the common source space of other models. Subfigure~\ref{fig.patch.L.vecnorm.superficial} displays a cortical view while \ref{fig.patch.L.vecnorm.deep} shows a deeper cross-section in the RA-plane of the Right-Anterior-Superior (RAS) coordinate system, roughly at the height of the thalamus. The field is strong near the EEG electrodes, depicted as black circles, and decays smoothly as distance from them increases. Overlapping source patches produce extended regions of higher field intensity, especially near regions where surface curvature brings multiple electrodes into close proximity.}
    \label{fig.patch.L.vecnorm}
\end{figure}

\begin{figure}
    \centering
    \hfill
    \begin{subfigure}[b]{0.4\linewidth}
        \includegraphics[width=\linewidth]{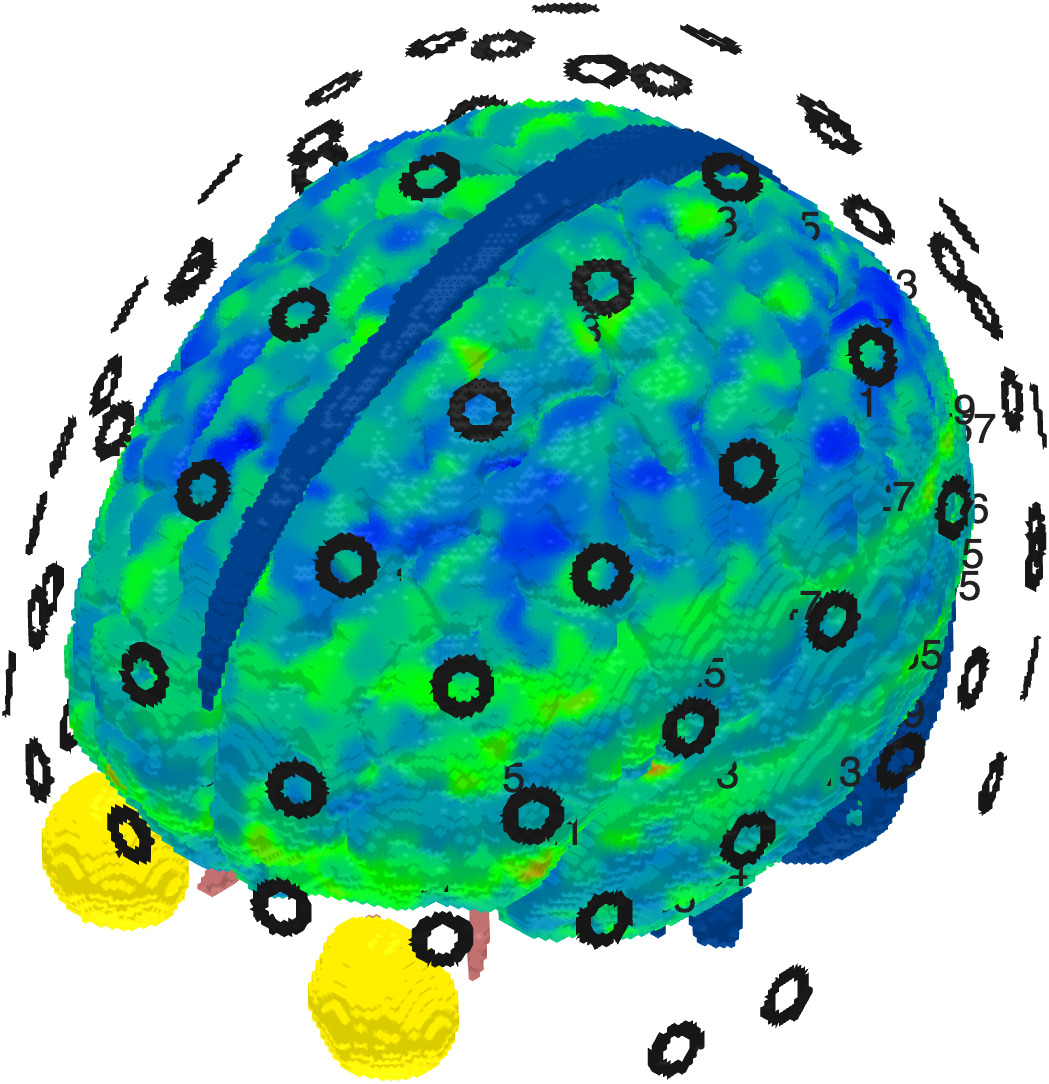}
        \caption{}
        \label{fig.perturbed.start.L.vecnorm.superficial}
    \end{subfigure}
    \hfill
    \begin{subfigure}[b]{0.4\linewidth}
        \includegraphics[width=\linewidth]{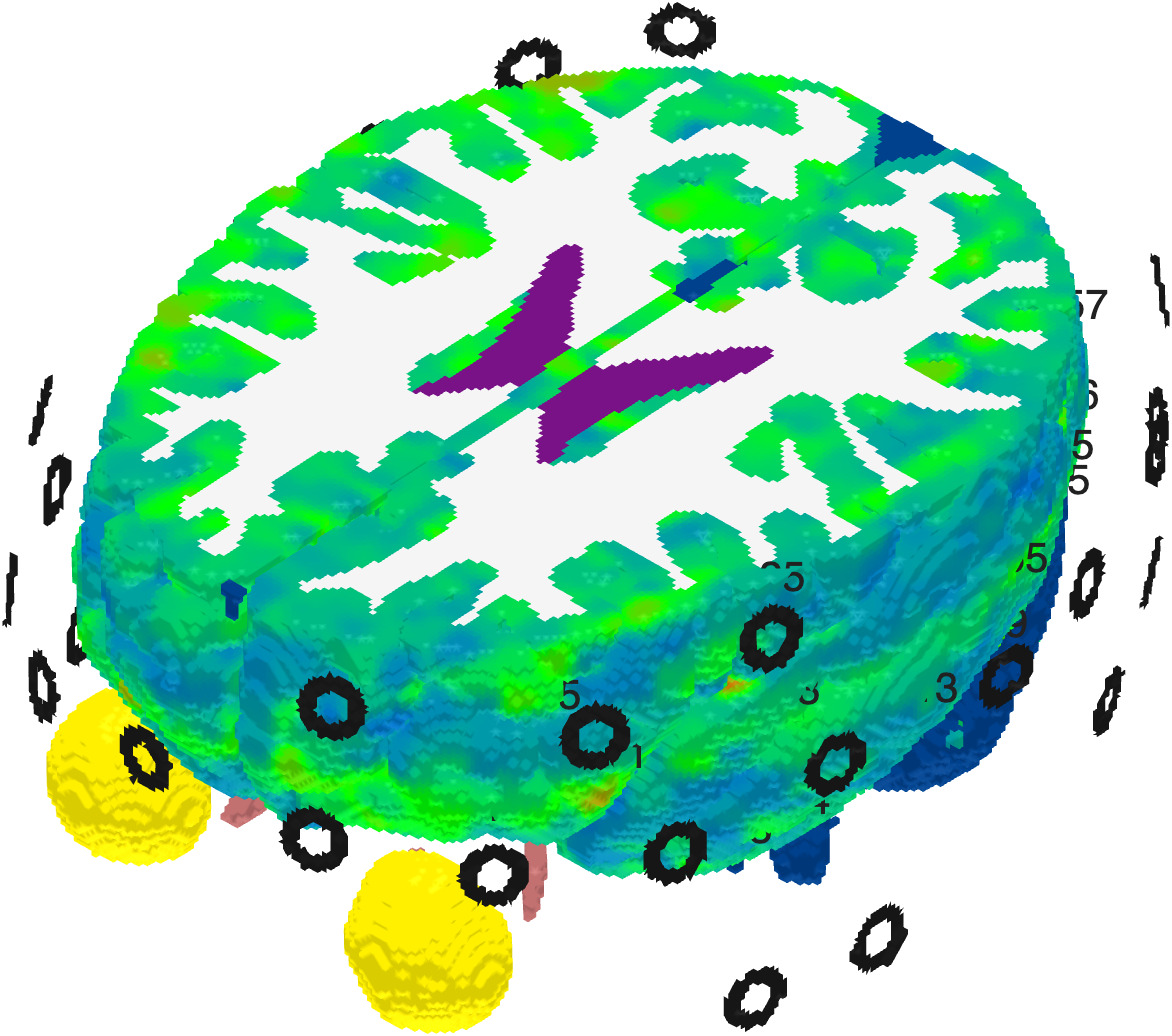}
        \caption{}
        \label{fig.perturbed.start.L.vecnorm.deep}
    \end{subfigure}
    \hfill
    \begin{minipage}[b]{0.1\linewidth}
        \includegraphics[height=4cm]{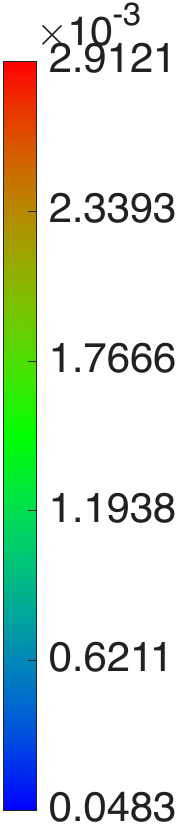}
    \end{minipage}
    \hfill
    \caption{A display of the column norms $\norm{\leadFieldMatrix_1}_2$ of the reference source space $\positions_1$ of the source perturbation experiment, with $\leadFieldMatrix_1$ computed with Local subtraction. Subfigure~\ref{fig.perturbed.start.L.vecnorm.superficial} displays a cortical view while \ref{fig.perturbed.start.L.vecnorm.deep} shows a deeper cross-section in the RA-plane of the Right-Anterior-Superior (RAS) coordinate system, roughly at the height of the thalamus. A rather drastic change in the field distribution is seen when compared to Figure~\ref{fig.duneuro.L.vecnorm}: the field is much less clearly anchored to the electrode vicinity}
    \label{fig.perturbed.start.vecnorm}
\end{figure}

\begin{figure}
    \centering
    \hfill
    \begin{subfigure}[b]{0.4\linewidth}
        \includegraphics[width=\linewidth]{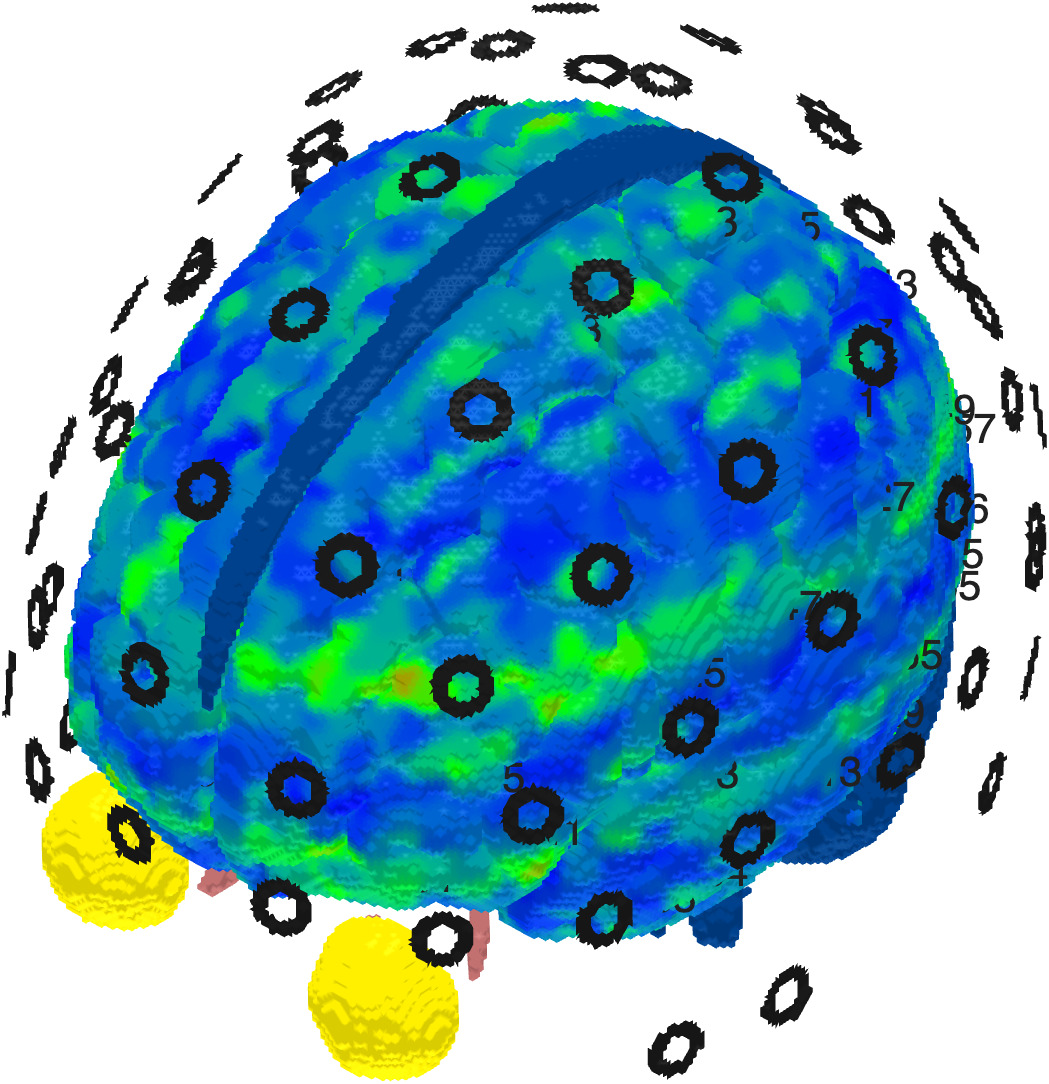}
        \caption{}
        \label{fig.perturbed.union.L.vecnorm.superficial}
    \end{subfigure}
    \hfill
    \begin{subfigure}[b]{0.4\linewidth}
        \includegraphics[width=\linewidth]{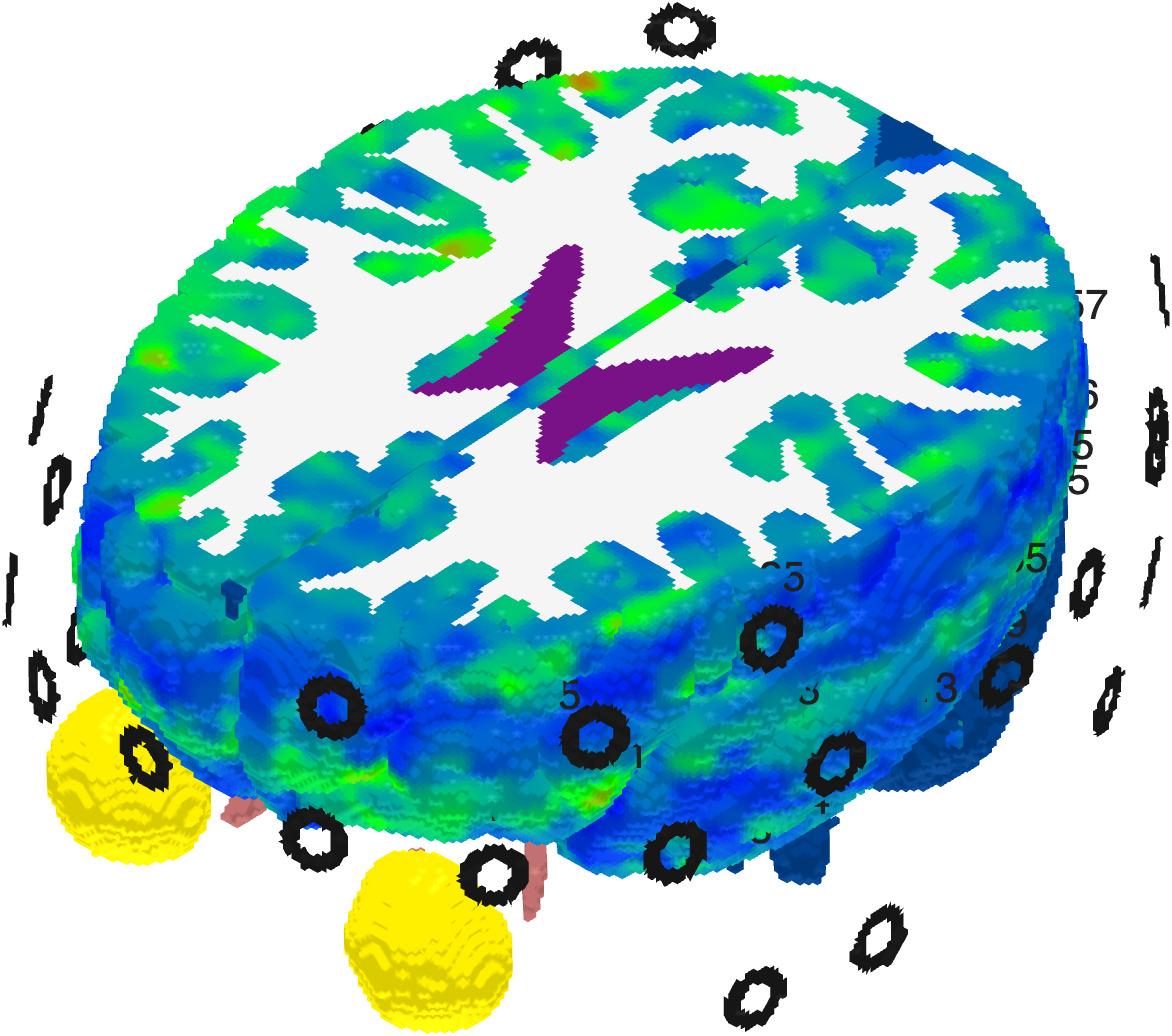}
        \caption{}
        \label{fig.perturbed.union.L.vecnorm.deep}
    \end{subfigure}
    \hfill
    \begin{minipage}[b]{0.1\linewidth}
        \includegraphics[height=4cm]{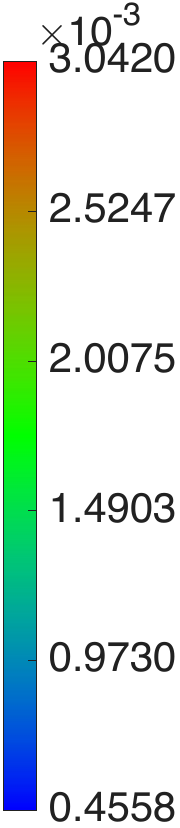}
    \end{minipage}
    \hfill
    \caption{A display of the column norms $\norm{\leadFieldMatrix_4}_2$ of the union source space $\positions_4$ of the source perturbation experiment, with $\leadFieldMatrix_4$ computed with Local subtraction. Subfigure~\ref{fig.perturbed.union.L.vecnorm.superficial} displays a cortical view while \ref{fig.perturbed.union.L.vecnorm.deep} shows a deeper cross-section in the RA-plane of the Right-Anterior-Superior (RAS) coordinate system, roughly at the height of the thalamus. Similarly to Figure~\ref{fig.perturbed.start.vecnorm}, the field is less regular in terms of its strength near the electrodes.}
    \label{fig.perturbed.union.vecnorm}
\end{figure}

 In contrast, the $\Hdiv$-based lead field built with Zeffiro Interface is presented in Figure~\ref{fig.zeffiro.L.vecnorm} (a--b). It is more irregular, and the decay of the column norm from the vicinity of the sensors is weaker as the distance from the superficial electrodes increases. The scale of the potential values is also much larger than that of the Whitney-based lead field generated with DUNEuro.

    A similar visualization of a patch modeling lead field column norm appears in Figure~\ref{fig.patch.L.vecnorm}. This bears resemblance to the $\Hdiv$ model's column norm $\norm\leadFieldMatrix_2$ presented in Figure~\ref{fig.zeffiro.L.vecnorm} in terms of field smoothness and decay rate along with surface, but again, the superficial field values seem more clearly localized near the electrodes as with the Whitney basis, and the column norm decays towards deep structures more strongly than with Zeffiro Interface's $\Hdiv$. Also, no significant column norm outliers were observed.

    For the source space perturbation experiment, the column norms of the Local-subtraction-based lead fields $\leadFieldMatrix_1$ and $\leadFieldMatrix_4$ are seen in Figures~\ref{fig.perturbed.start.vecnorm} and \ref{fig.perturbed.union.vecnorm}, respectively. Again, at this scale the visual differences between Whitney and Local subtraction are too small to warrant displaying both. We see a clear change in the behaviour of the column norms when compared to the original element-centroid-based source positioning of Whitney and Local subtraction in Figure~\ref{fig.duneuro.L.vecnorm}: the field maxima are not unambiguously near the electrodes, and the field strength also shows no clear decay as distance from the electrodes increases. Instead, it seems that the sensitivity pattern is more random across the entire volume. The column norm of $\leadFieldMatrix_1$ also seems more vertically spread out, whereas for $\leadFieldMatrix_4$ the column norms show a vertically banded region of high sensitivity at the height of the thalamus. Even if the shape of the field changes along with the source positions, the field intensities remain within a sensible \unit{\micro\volt} range for Whitney and Local subtraction produced with DUNEuro.

\subsection{Analysis of the interpolation schemes}\label{sec.interpolation.schemes}

    To observe the sensitivity of Local subtraction and Whitney models of DUNEuro and $\Hdiv$ of Zeffiro Interface to source positioning, \num 4 additional Cartesian source spaces were generated apart from the hexahedrally positioned sources used in experiments I and II. For Local subtraction and Whitney, first all tetrahedra with 4 facet-based neighbours present within the same compartment as the tetrahedron itself were selected from active mesh compartments, after which
    \begin{enumerate}
        \item sources were attached to the 1st vertex of each such tetrahedron. This is the source position set $\positions_1$.
        \item Sources were also placed at the opposing vertices $\positions_2$ of neighbouring tetrahedra sharing the facet opposing vertex 1 of the initial tetrahedra, where vertex order is determined by the Zeffiro Interface mesh generator \cite{galazprieto2023zeffiromesh}. This can be illustrated as pyramids placed such that their bases face each other and the sources are positioned at the tips of the pyramids.
        \item A third source space was positioned at the midpoints $\positions_3$ between $\positions_1$ and $\positions_2$.
        \item Finally, a union of the above 3 source spaces was created, named $\positions_4$.
    \end{enumerate}
    Corresponding lead fields $\leadFieldMatrix_1$, $\leadFieldMatrix_2$, $\leadFieldMatrix_3$ and $\leadFieldMatrix_4$ were then created based on these \num 4 source spaces.

\begin{figure*}
    \centering
    \begin{minipage}{0.03\linewidth}
        \rotatebox{90}{$\relativeDifference\args{\leadFieldMatrix_1,\leadFieldMatrix_3}$}
    \end{minipage}\begin{minipage}{0.3\linewidth}
    \begin{center}
    \qtyrange{0}{15}{\milli\meter}
    \end{center}%
    \vspace{-0.3cm}%
    
    \includegraphics[width=0.98\linewidth]{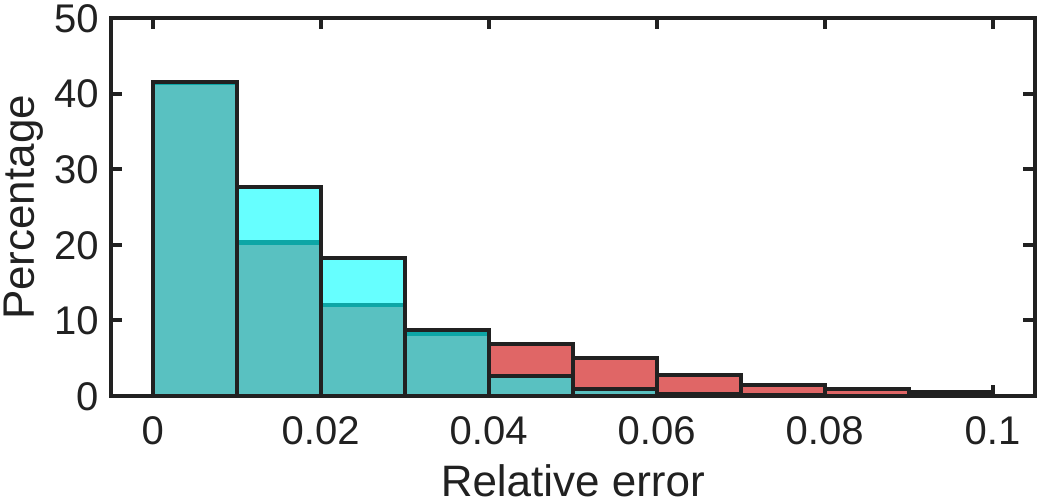}
    \end{minipage}\begin{minipage}{0.3\linewidth}
    \begin{center}
    \qtyrange{15}{30}{\milli\meter}
    \end{center}%
    \vspace{-0.3cm}
        \includegraphics[width=0.98\linewidth]{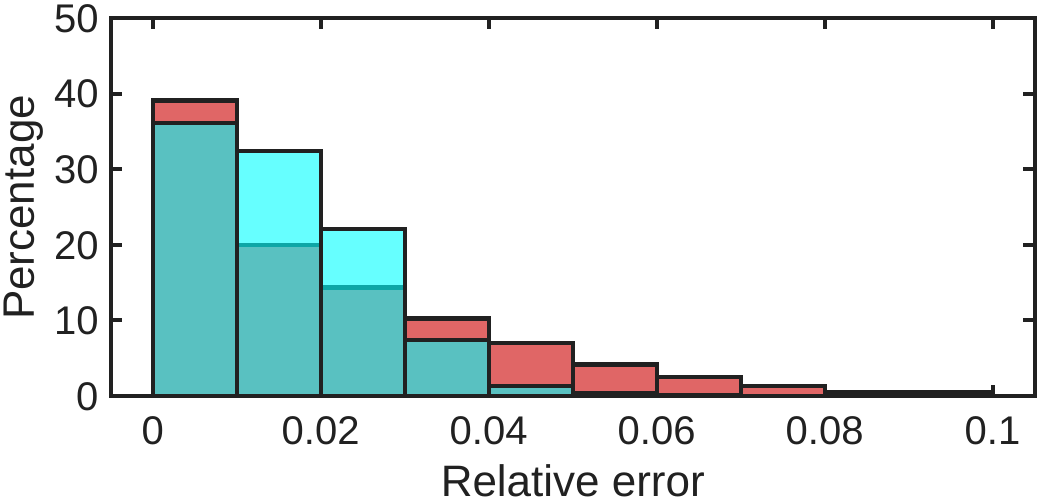}
    \end{minipage}\begin{minipage}{0.3\linewidth}
    \begin{center}
    \qtyrange{30}{60}{\milli\meter}
    \end{center}\vspace{-0.3cm}
        \includegraphics[width=0.98\linewidth]{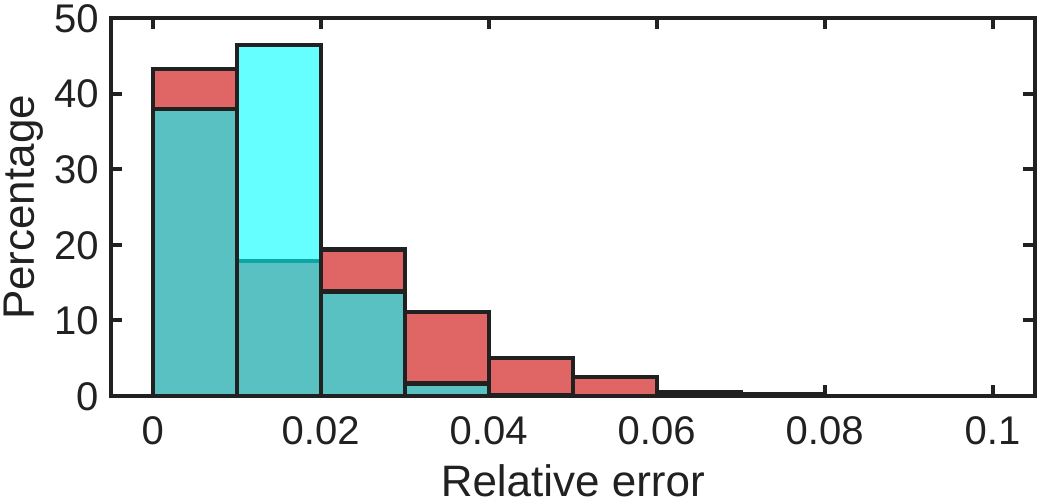}
    \end{minipage}
    \vspace{0.2 cm}
    \begin{minipage}{0.03\linewidth}
        \rotatebox{90}{$\relativeDifference\args{\leadFieldMatrix_3,\leadFieldMatrix_2}$}
    \end{minipage}\begin{minipage}{0.3\linewidth}
        \includegraphics[width=0.98\linewidth]{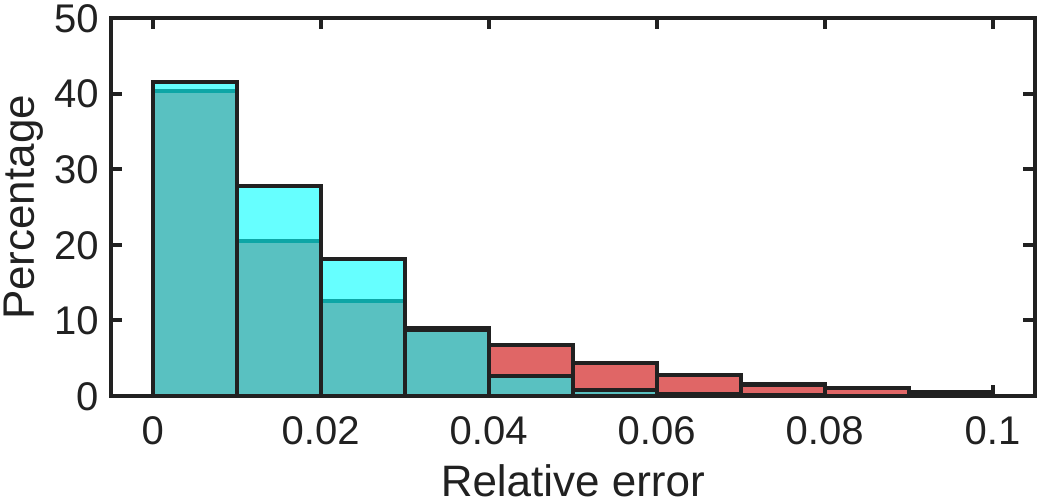}
    \end{minipage}\begin{minipage}{0.3\linewidth}
        \includegraphics[width=0.98\linewidth]{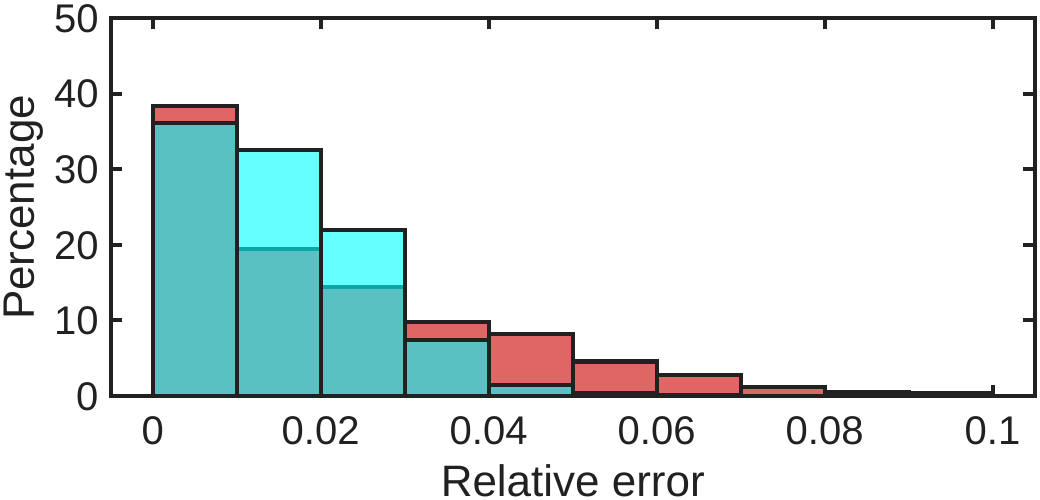}
    \end{minipage}\begin{minipage}{0.3\linewidth}
        \includegraphics[width=0.98\linewidth]{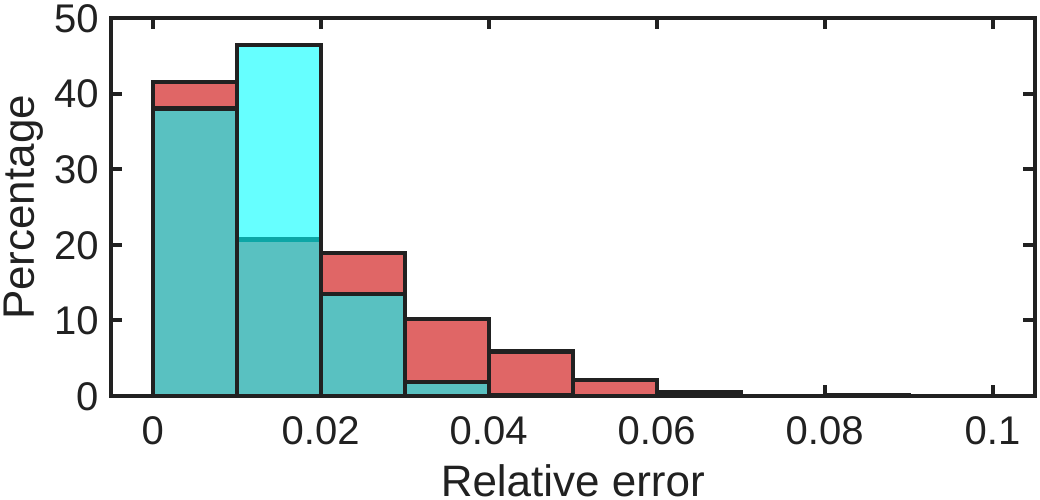}
    \end{minipage}
    \vspace{0.2 cm}
    \begin{minipage}{0.03\linewidth}
        \rotatebox{90}{$\relativeDifference\args{\leadFieldMatrix_1,\leadFieldMatrix_2}$}
    \end{minipage}\begin{minipage}{0.3\linewidth}
        \includegraphics[width=0.98\linewidth]{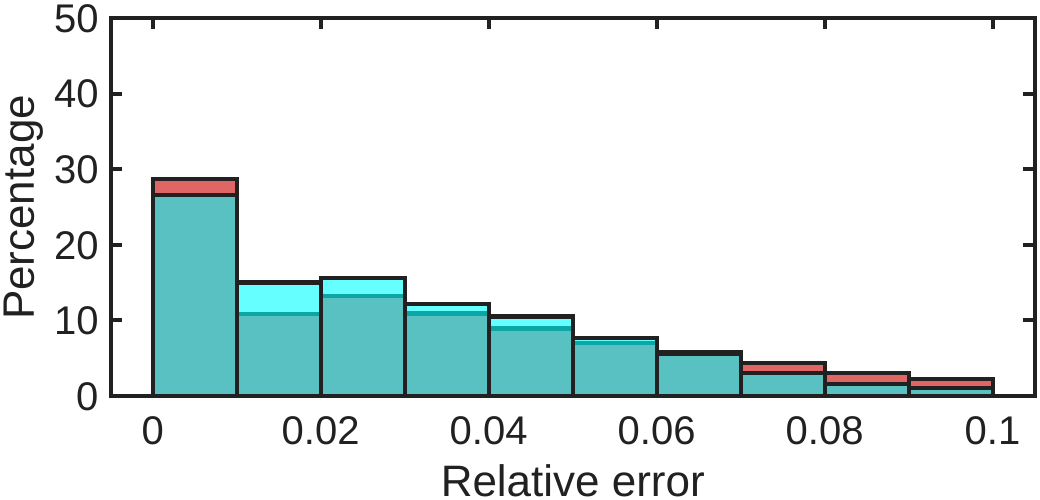}
    \end{minipage}\begin{minipage}{0.3\linewidth}
        \includegraphics[width=0.98\linewidth]{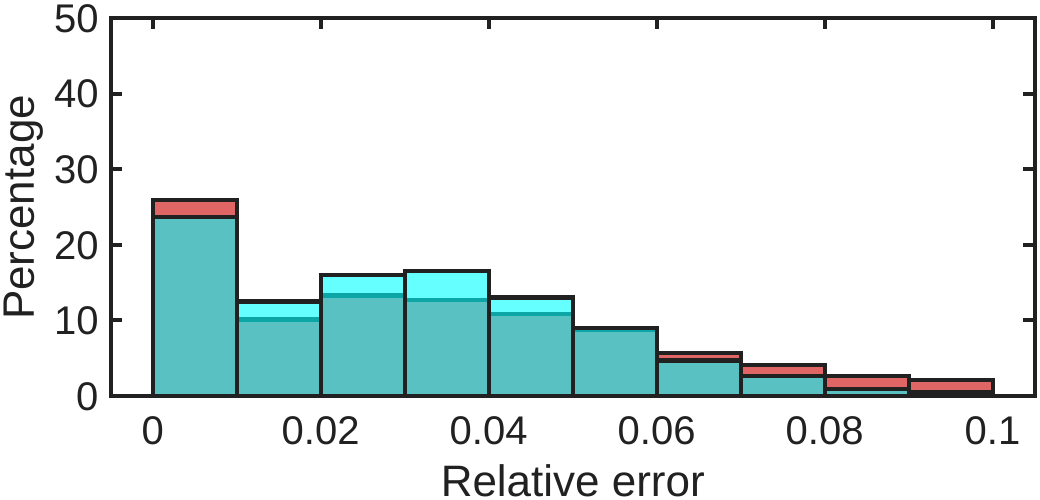}
    \end{minipage}\begin{minipage}{0.3\linewidth}
        \includegraphics[width=0.98\linewidth]{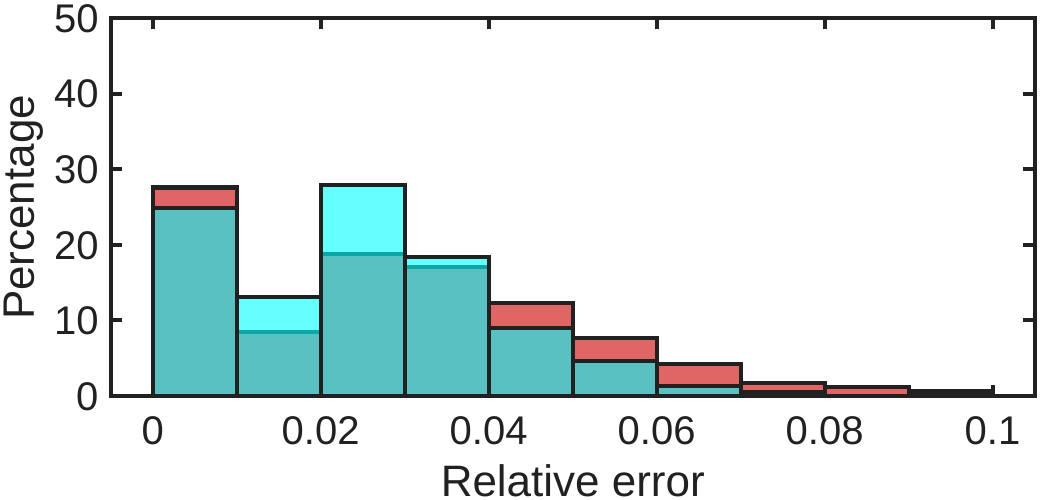}
    \end{minipage}    
    \caption{Histograms of relative differences $\relativeDifference\args{\leadFieldMatrix_1,\leadFieldMatrix_2}$, $\relativeDifference\args{\leadFieldMatrix_1,\leadFieldMatrix_3}$ and $\relativeDifference\args{\leadFieldMatrix_3,\leadFieldMatrix_2}$ at source depths \qtyrange{0}{15}{\milli\meter},  \qtyrange{15}{30}{\milli\meter} and  \qtyrange{30}{60}{\milli\meter}, where $\relativeDifference$ is defined by \eqref{eq.reldiff}. Red represents the Whitney basis, and turquoise is Local subtraction.}
    \label{fig.histograms}
\end{figure*}

    In the case of $\Hdiv$ provided by the Zeffiro Interface, the interpolation within a tetrahedron is done using a characteristic function mapping the value at each point in the tetrahedron to correspond to the value of the function at the centroid, which makes the relative difference within a tetrahedron effectively zero. For this reason, we are omitting it from the local inspection presented in Figure \ref{fig.histograms}. To get comparable results when examining source model sensitivity within \qty{3}{\milli\meter} range, a very dense \num{1 000 000} element regular source space within the active elements was created. This allows for the discovery of a good representative of the points in $\positions_4$ within the mesh by picking the best active element centroids to match with the source space $\positions_4$. The script was set to warn if a good representative centroid was not found. In the end, no warning was given.

    The lead fields $\leadFieldMatrix_1$, $\leadFieldMatrix_2$ and $\leadFieldMatrix_3$ of Local subtraction and Whitney were then compared to each other by taking the relative differences of \eqref{eq.reldiff} between them, as in computing $\relativeDifference\args{\leadFieldMatrix_1,\leadFieldMatrix_2}$, $\relativeDifference\args{\leadFieldMatrix_1,\leadFieldMatrix_3}$ and $\relativeDifference\args{\leadFieldMatrix_3,\leadFieldMatrix_2}$. Histograms of these differences at various source depths were then visualized in Figure~\ref{fig.histograms}.

    The histograms show that most of the relative differences between Whitney and Local subtraction are rather small, and remain below \qty{10}{\percent}. Altogether, Local subtraction seems to present at least \qtyrange{1}{2}{\percent} shorter tail across all comparisons, but possibly even greater at shallower source depths. The general trend for the combinations $\leadFieldMatrix_1$ vs. $\leadFieldMatrix_3$ and $\leadFieldMatrix_2$ vs. $\leadFieldMatrix_3$ is that moving to larger source depths shortens the tail of each distribution for both Whitney and Local subtraction, where a clearest jump happens between the ranges \qtyrange{15}{30}{\milli\meter} and \qtyrange{30}{60}{\milli\meter}, where the tail cutoff point moves from \qty{8}{\percent} to \qty{6}{\percent} for Whitney and from \qty{5}{\percent} to \qty{4}{\percent} for Local subtraction. This does not happen in the $\leadFieldMatrix_1$ vs. $\leadFieldMatrix_2$ comparison where Whitney maintains a nonzero tail even near the \qty{10}{\percent} relative error, even at the greatest source depth, but for Local subtraction a visible difference of roughly \qty{1}{\percent} in tail reduction is seen.

\begin{figure}
    \centering
    \hfill
    \includegraphics[width=\linewidth]{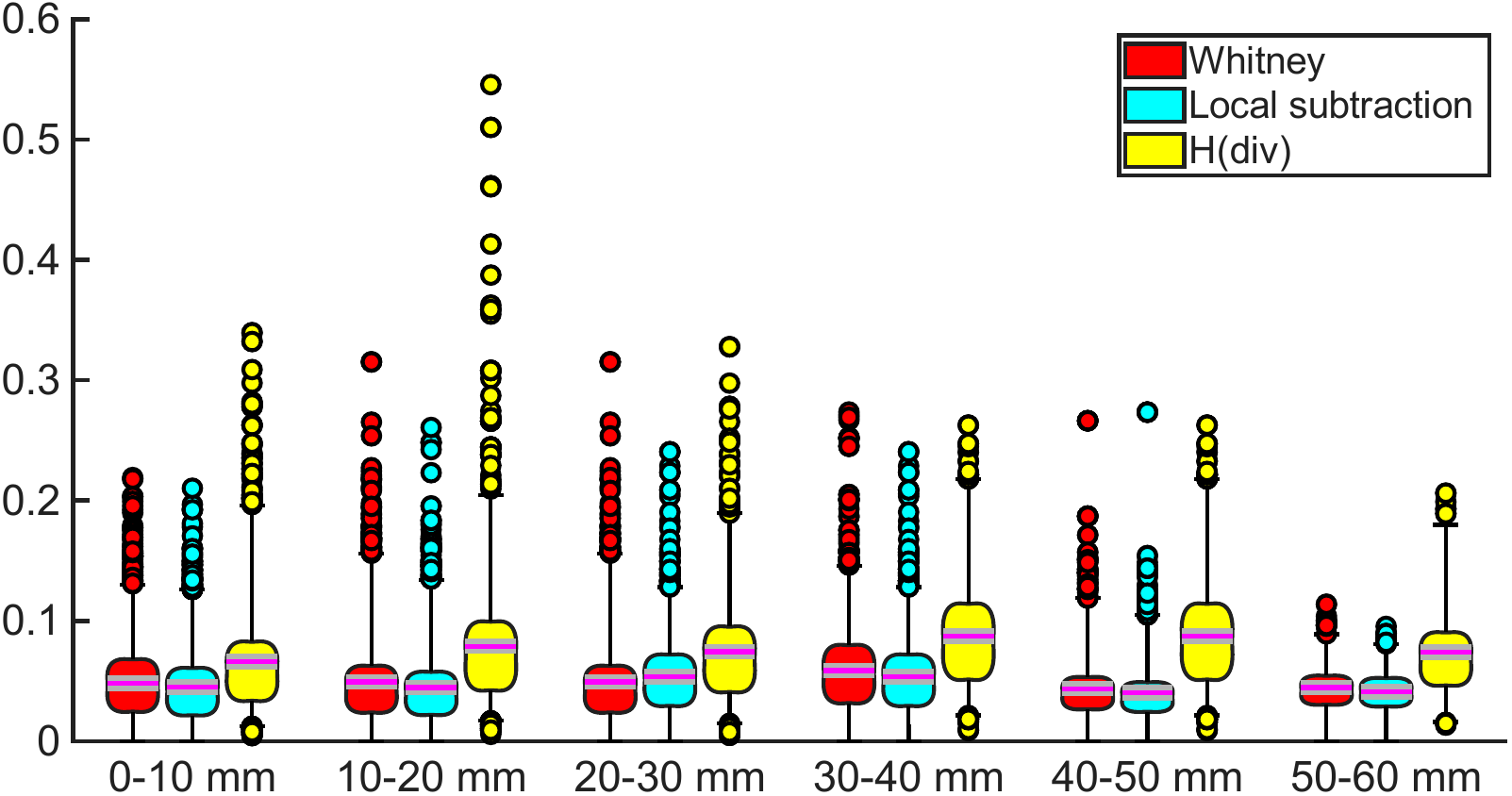}
    \hfill
    \caption{Relative differences \eqref{eq.reldiff} of lead field value perturbations within a \qty{3}{\milli\meter} range from the reference lead field element computed for various depths at intervals of \qty{10}{\milli\meter} up to \qty{60}{\milli\meter}. Red boxes represent the close-range lead field discrepancies of the Whitney basis, turquoise depicts the Local subtraction method, and yellow presents the statistics for $\Hdiv$.}
    \label{fig.boxplot}
\end{figure}

    An experiment where \eqref{eq.reldiff} was computed for $\leadFieldMatrix_1$ and $\leadFieldMatrix_4$ with $\leadFieldMatrix_1$ as a reference, with points in $\positions_4$ within \qty{3}{\milli\meter} regions of interest (ROI) around $\positions_1$. The $\Hdiv$ model was included in this experiment by choosing the closest representatives of $\positions_4$ and their corresponding columns of the corresponding lead field. These ROI-based differences were then visualized for source depths between \qtyrange{0}{60}{\milli\meter} at \qty{10}{\milli\meter} intervals using box plots with \num 1.5 inter-quantile range (IQR). The results can be seen in Figure~\ref{fig.boxplot}.

    Generally, we see a similar difference between Whitney and Local subtraction as we did when observing Figure~\ref{fig.histograms}: at slightly below \qty{5}{\percent}, the medians of relative errors are approximately \qty{1}{\percent} smaller for Local subtraction than they are for Whitney, except at the source depth \qtyrange{20}{30}{\milli\meter}, where the median of Whitney drops below that of Local subtraction. However, even in this case the tail of the distribution is shorter for Local subtraction, with that of Whitney showing an outlier above the \qty{30}{\percent} mark while Local subtraction remains at \qty{25}{\percent} at most. At all depths, $\Hdiv$ produces the weakest performance, with its median approximately in the range \qtyrange{7}{10}{\percent} and its tail seeing values near the \qty{35}{\percent} mark for shallower source depths while dropping below \qty{30}{\percent} starting from the \qtyrange{30}{40}{\milli\meter} mark. At the depth \qtyrange{10}{20}{\milli\meter} we see $\Hdiv$ outliers with as high as \qty{55}{\percent} error.

\subsection{Reconstructions of the cortical source}\label{sec.cortical.reconstructions}
\begin{figure}
    \def\imageWidth{0.21\linewidth}
    \def\inverseMethodNameWidth{0.09\linewidth}
    \def\colorbarWidth{0.075\linewidth}
    \scriptsize
    \centering
    \begin{minipage}[b]{0.9\linewidth}
        \begin{minipage}[b]{\inverseMethodNameWidth}
            \rotatebox{90}{\hspace{0.05cm} Whitney}
        \end{minipage}\begin{minipage}[b]{\imageWidth}
        \begin{center}
            sLORETA
        \end{center}
            \includegraphics[trim={6cm 2cm 6.5cm 2cm},clip,width=1\linewidth]{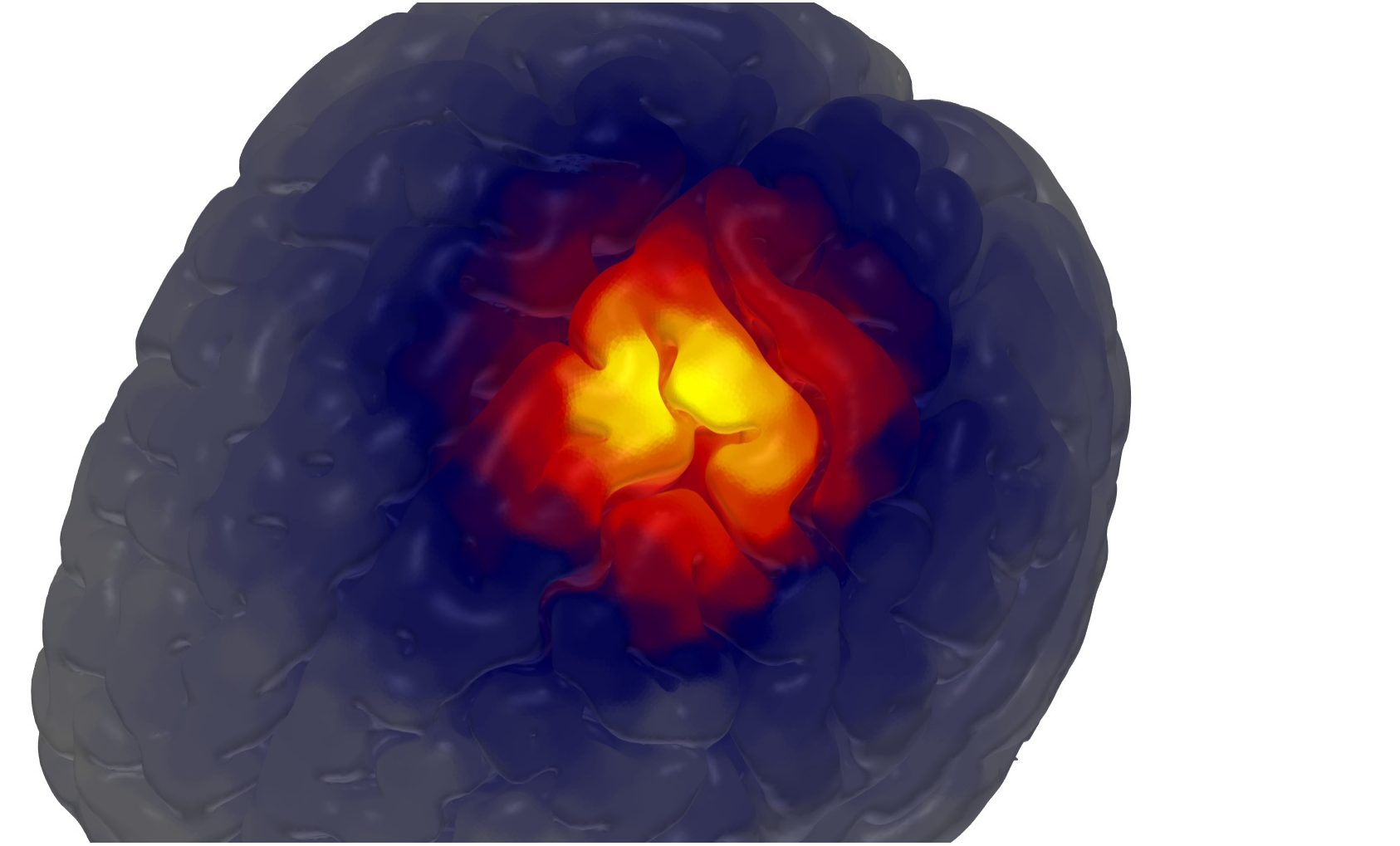}
        \end{minipage}\begin{minipage}[b]{\imageWidth}
        \begin{center}
            SHAL1R
        \end{center}
            \includegraphics[trim={6cm 2cm 6.5cm 2cm},clip,width=1\linewidth]{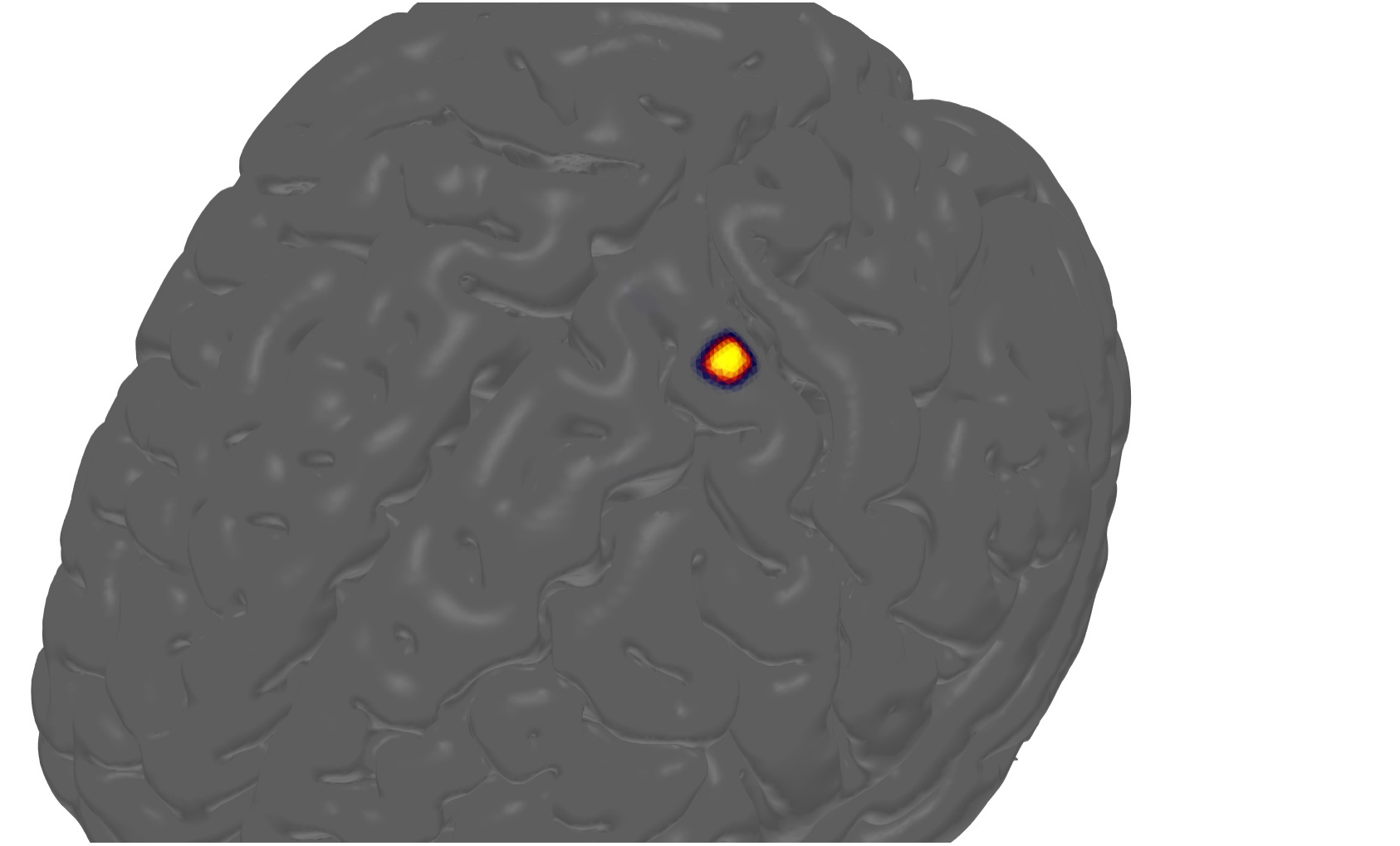}
        \end{minipage}%
        \begin{minipage}[b]{\imageWidth}
        \begin{center}
            SKF
        \end{center}
            \includegraphics[trim={6cm 2cm 6.5cm 2cm},clip,width=1\linewidth]{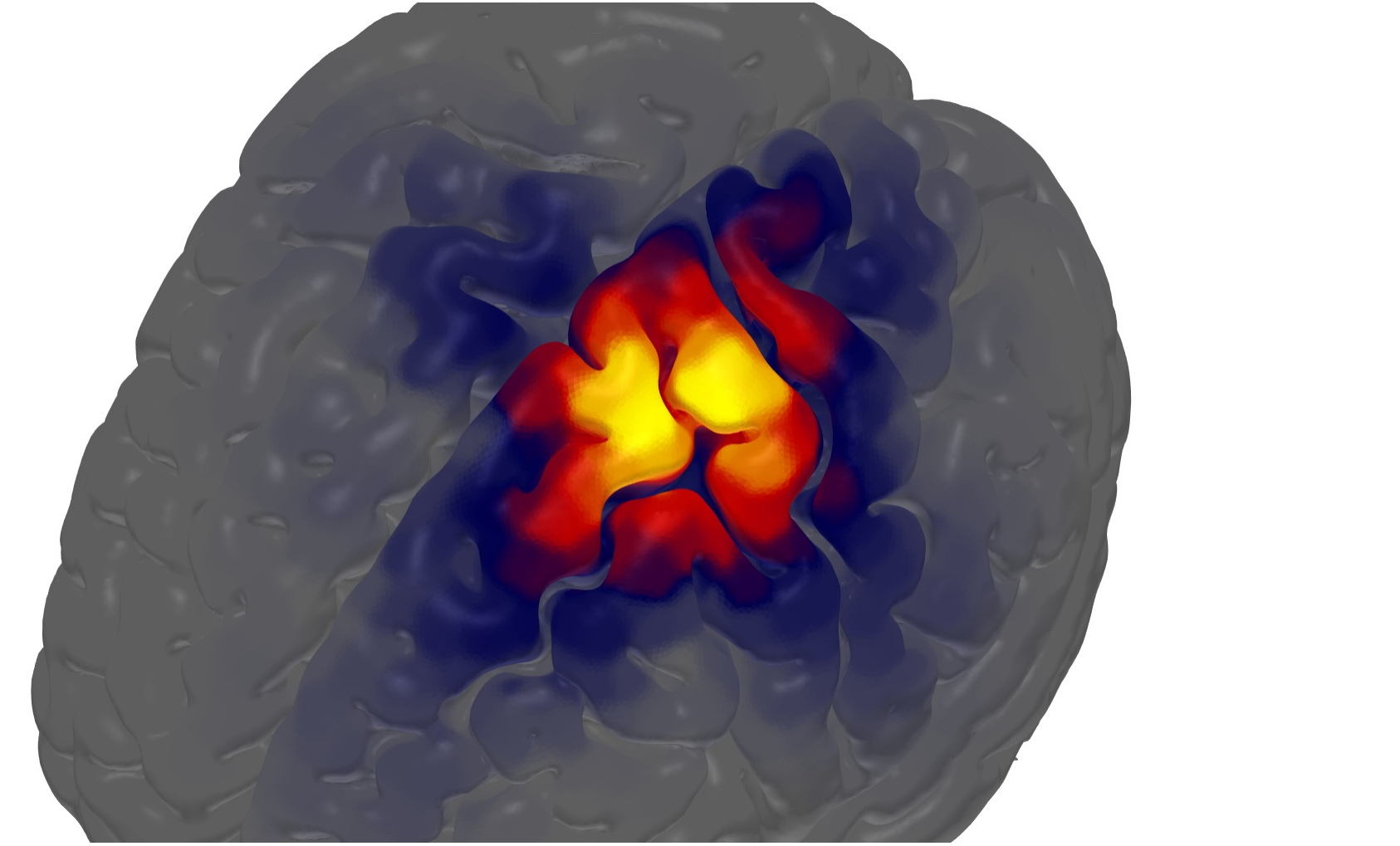}
        \end{minipage}%
        \begin{minipage}[b]{\imageWidth}
        \begin{center}
            DS
        \end{center}
            \includegraphics[trim={6cm 2cm 6.5cm 2cm},clip,width=1\linewidth]{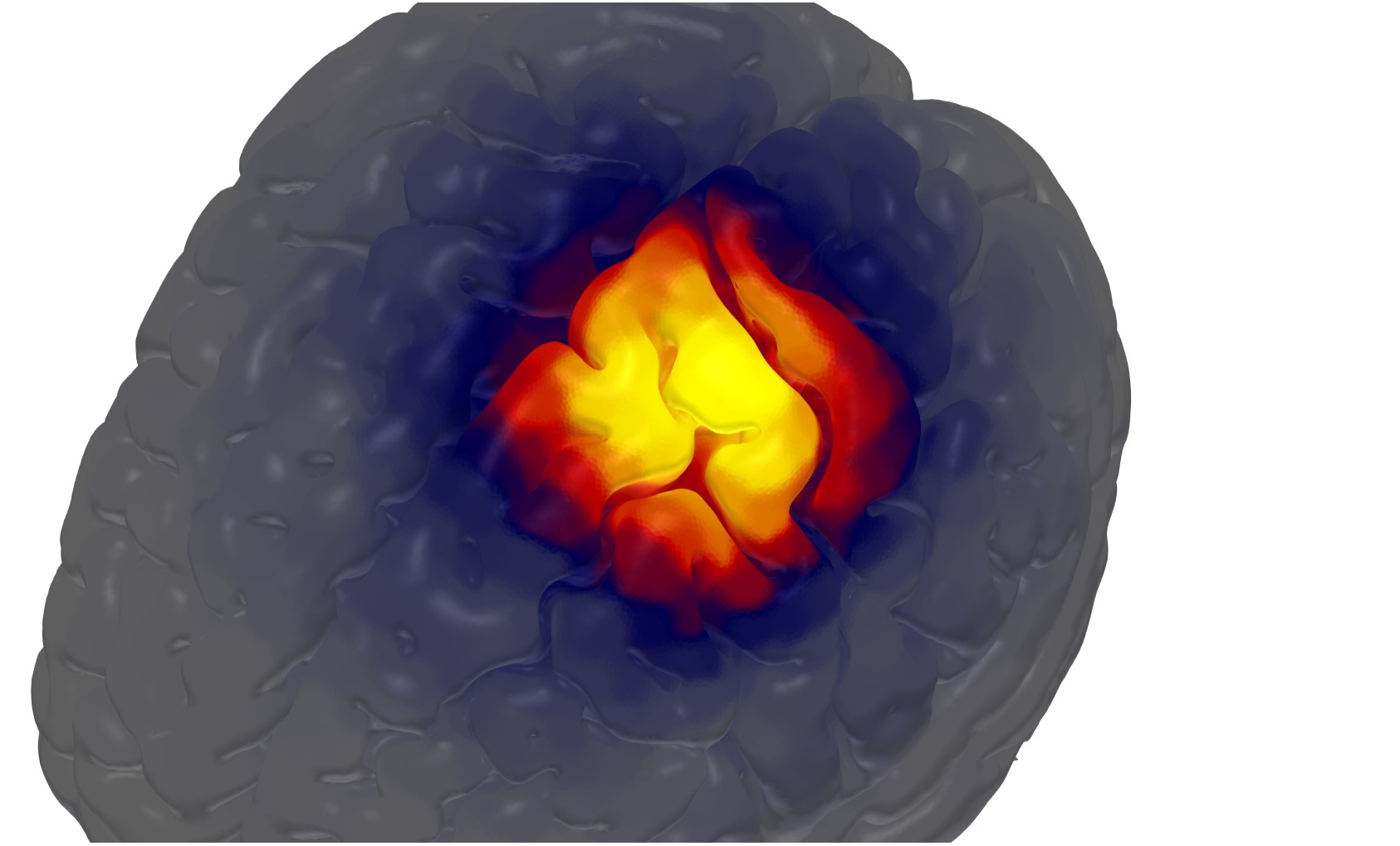}
        \end{minipage}

        \begin{minipage}[b]{\inverseMethodNameWidth}
            \rotatebox{90}{\hspace{0.3cm}Local} \rotatebox{90}{\hspace{0.32cm}subt.}
        \end{minipage}%
        \begin{minipage}[b]{\imageWidth}
            \includegraphics[trim={6cm 2cm 6.5cm 2cm},clip,width=1\linewidth]{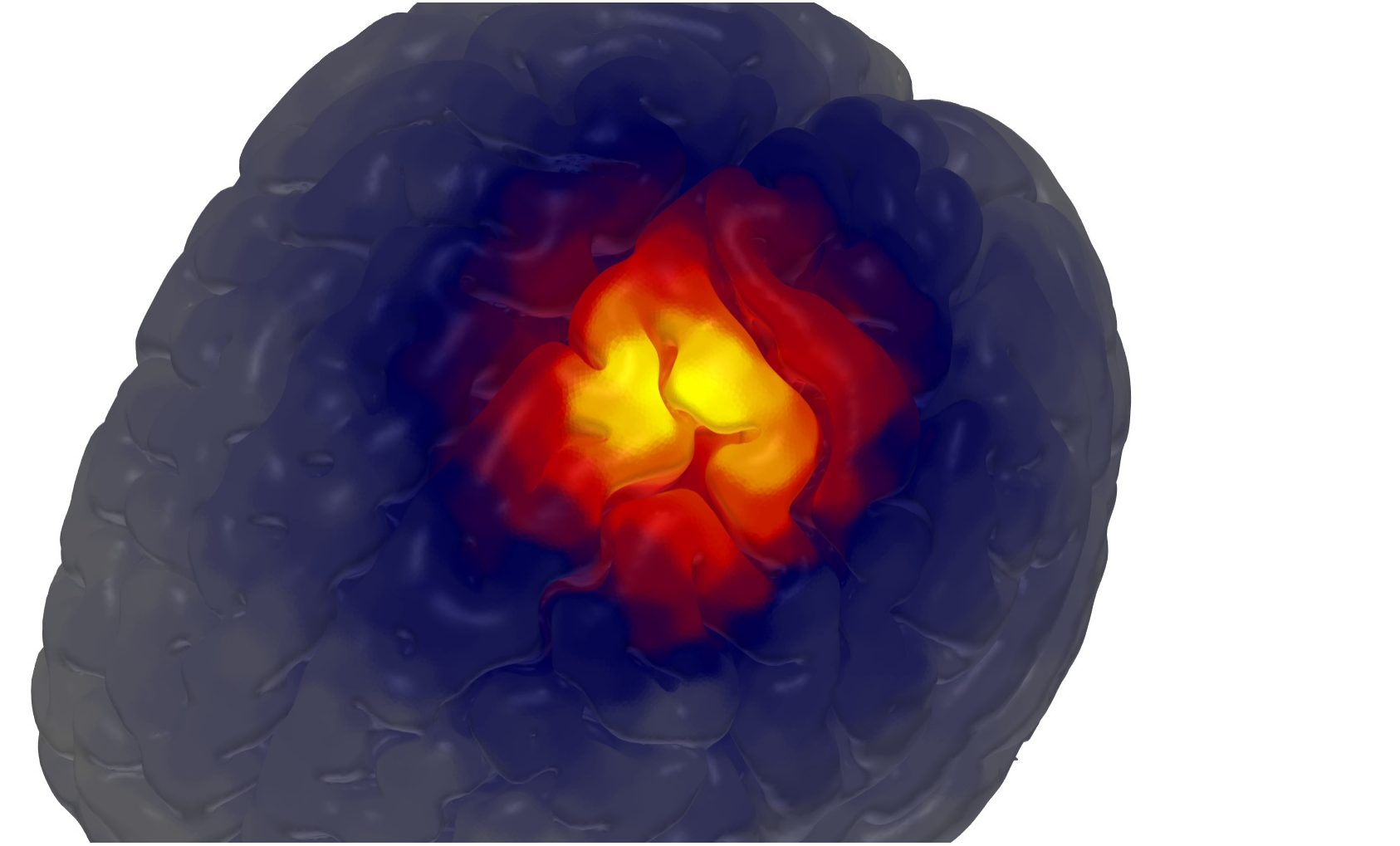}
        \end{minipage}%
        \begin{minipage}[b]{\imageWidth}
            \includegraphics[trim={6cm 2cm 6.5cm 2cm},clip,width=1\linewidth]{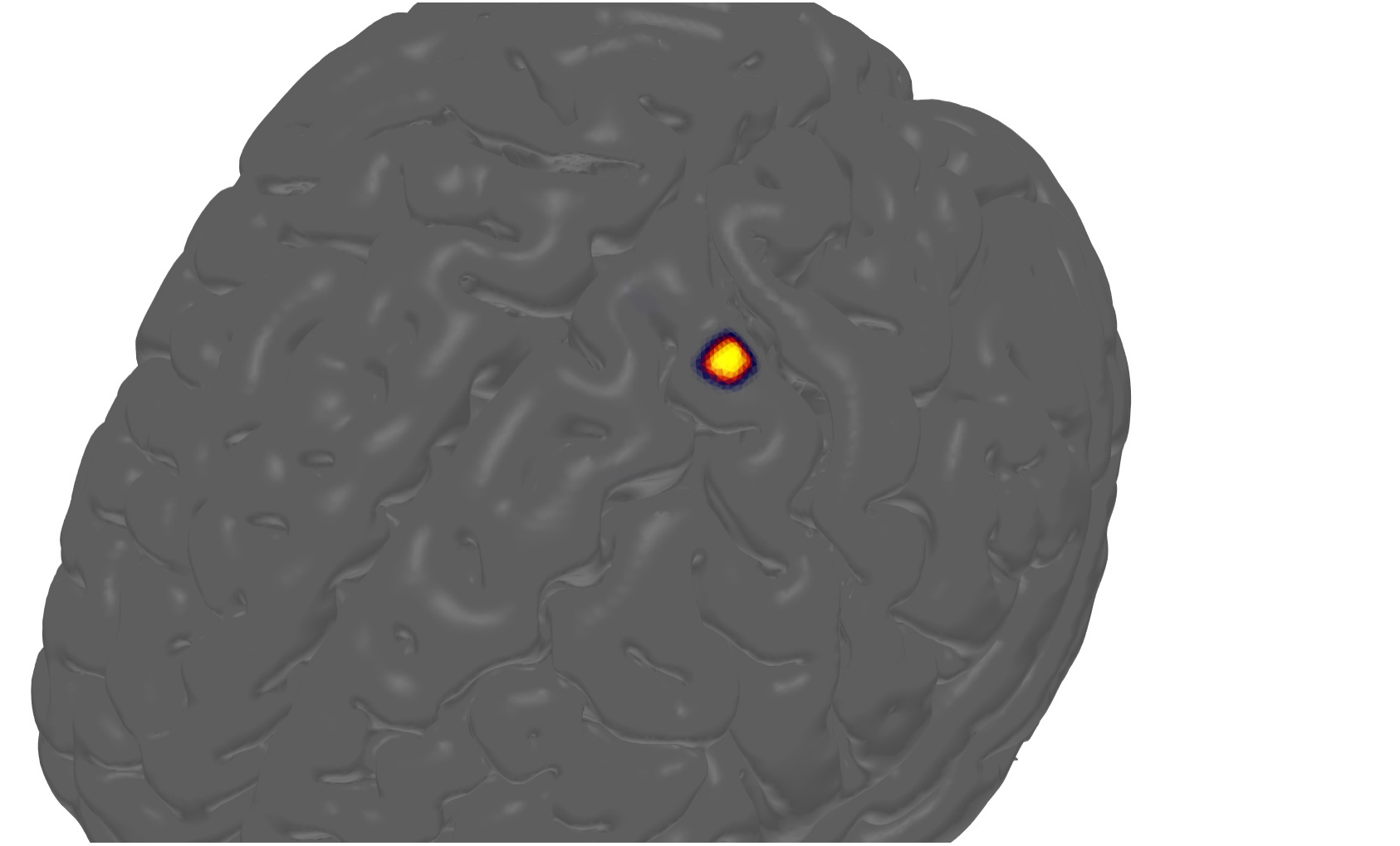}
        \end{minipage}%
        \begin{minipage}[b]{\imageWidth}
            \includegraphics[trim={6cm 2cm 6.5cm 2cm},clip,width=1\linewidth]{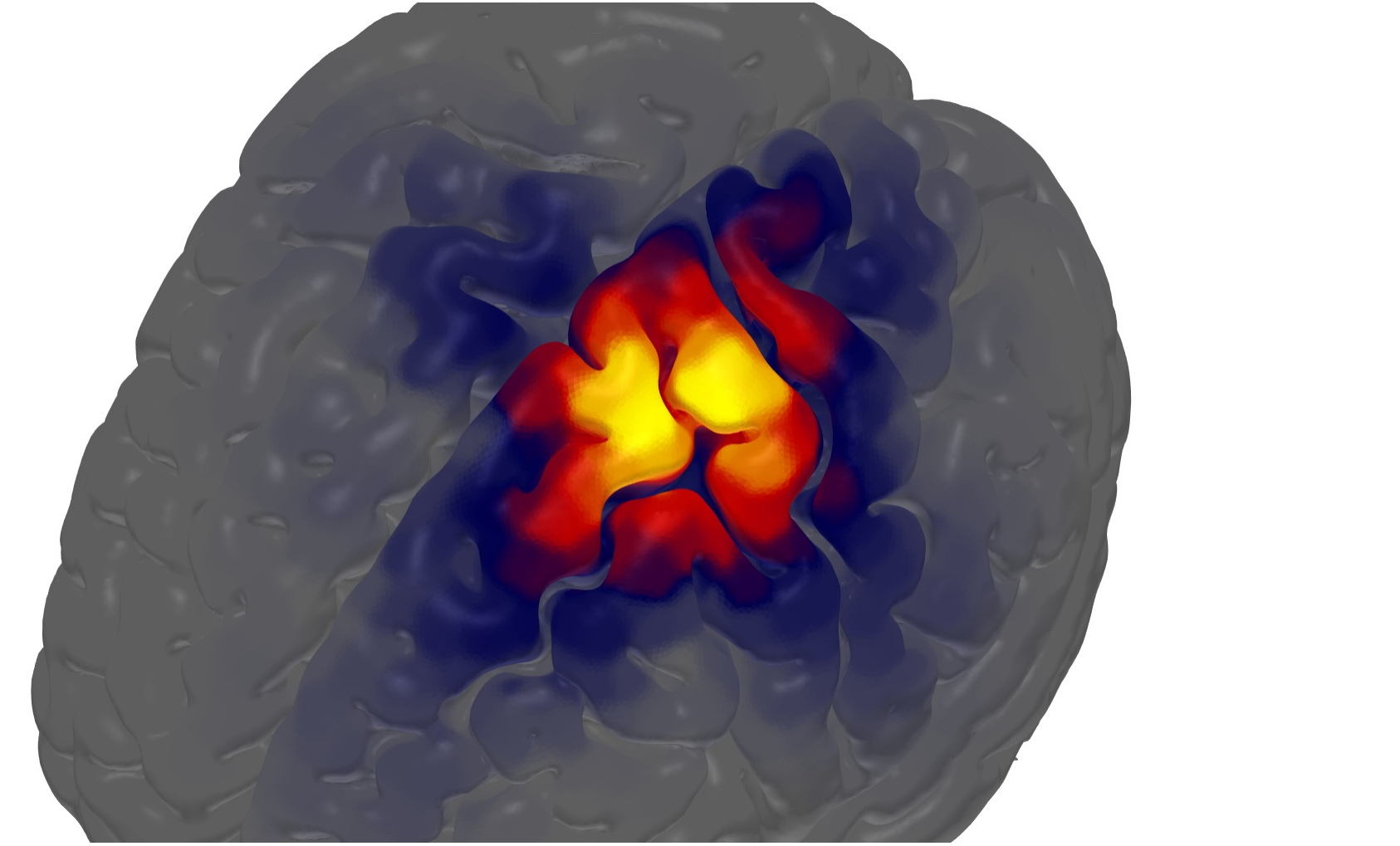}
        \end{minipage}%
        \begin{minipage}[b]{\imageWidth}
            \includegraphics[trim={6cm 2cm 6.5cm 2cm},clip,width=1\linewidth]{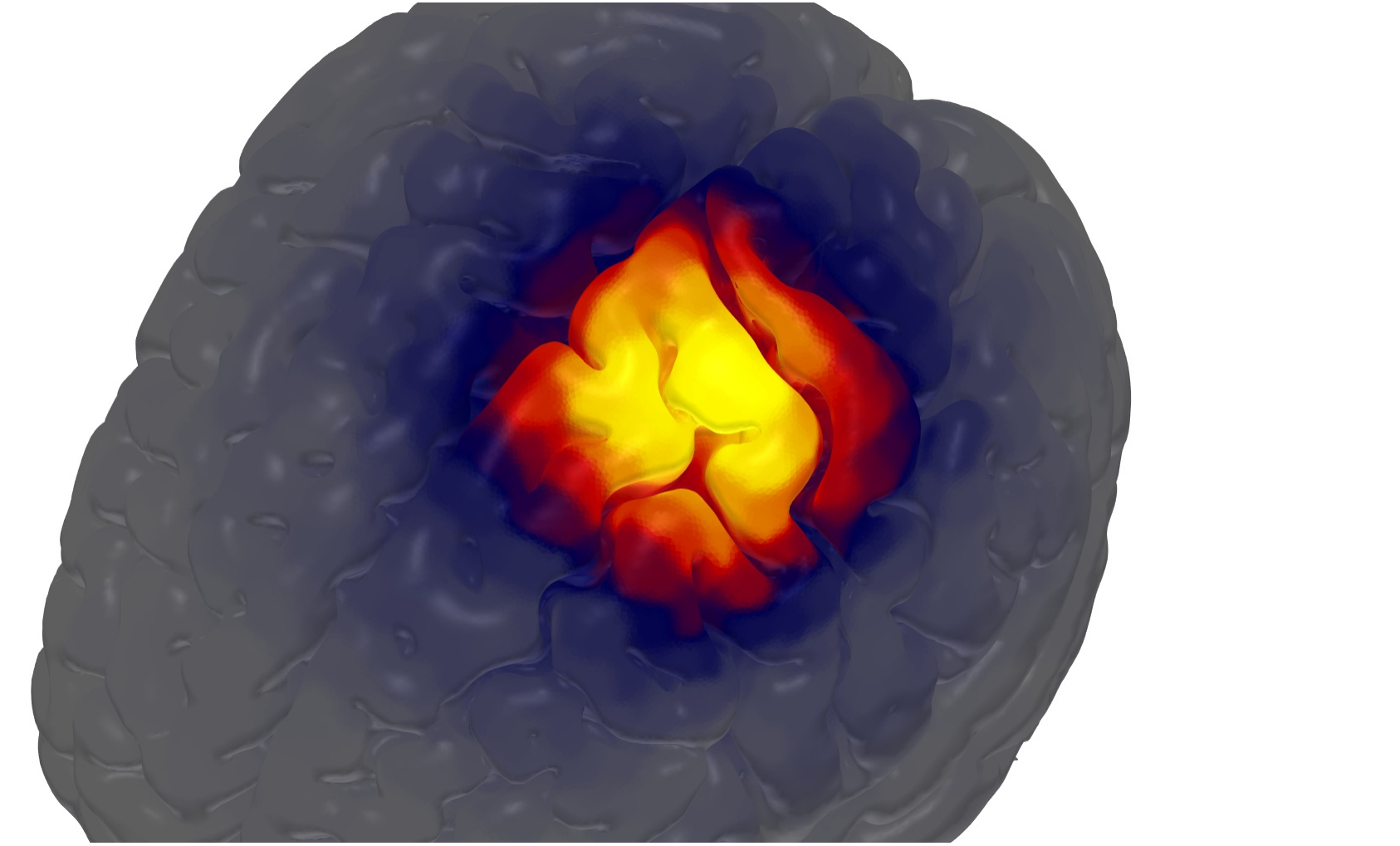}
        \end{minipage}

        \begin{minipage}[b]{\inverseMethodNameWidth}
            \rotatebox{90}{\hspace{0.1cm} $\Hdiv$}
        \end{minipage}%
        \begin{minipage}[b]{\imageWidth}
            \includegraphics[trim={6cm 2cm 6.5cm 2cm},clip,width=1\linewidth]{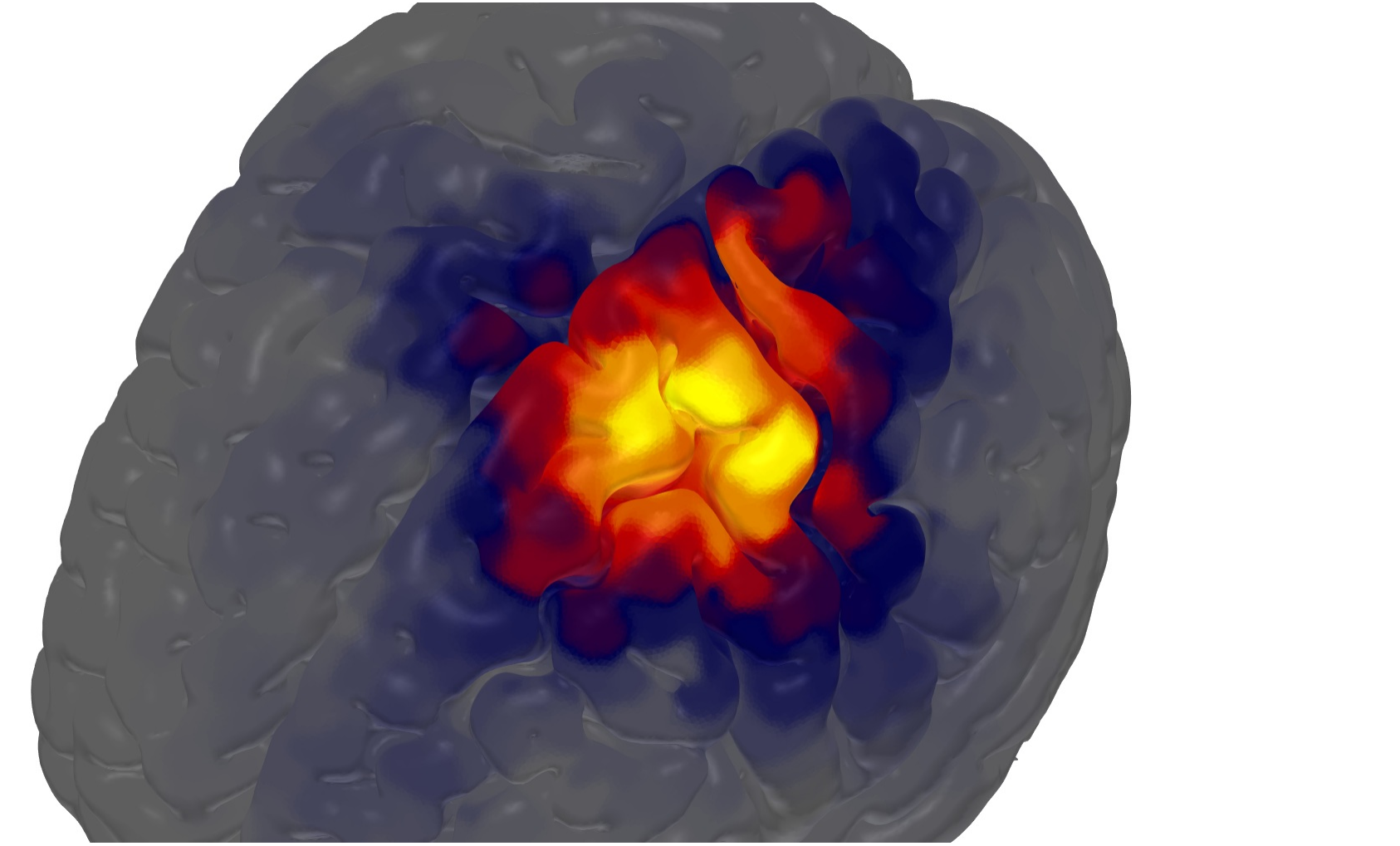}
        \end{minipage}%
        \begin{minipage}[b]{\imageWidth}
            \includegraphics[trim={6cm 2cm 6.5cm 2cm},clip,width=1\linewidth]{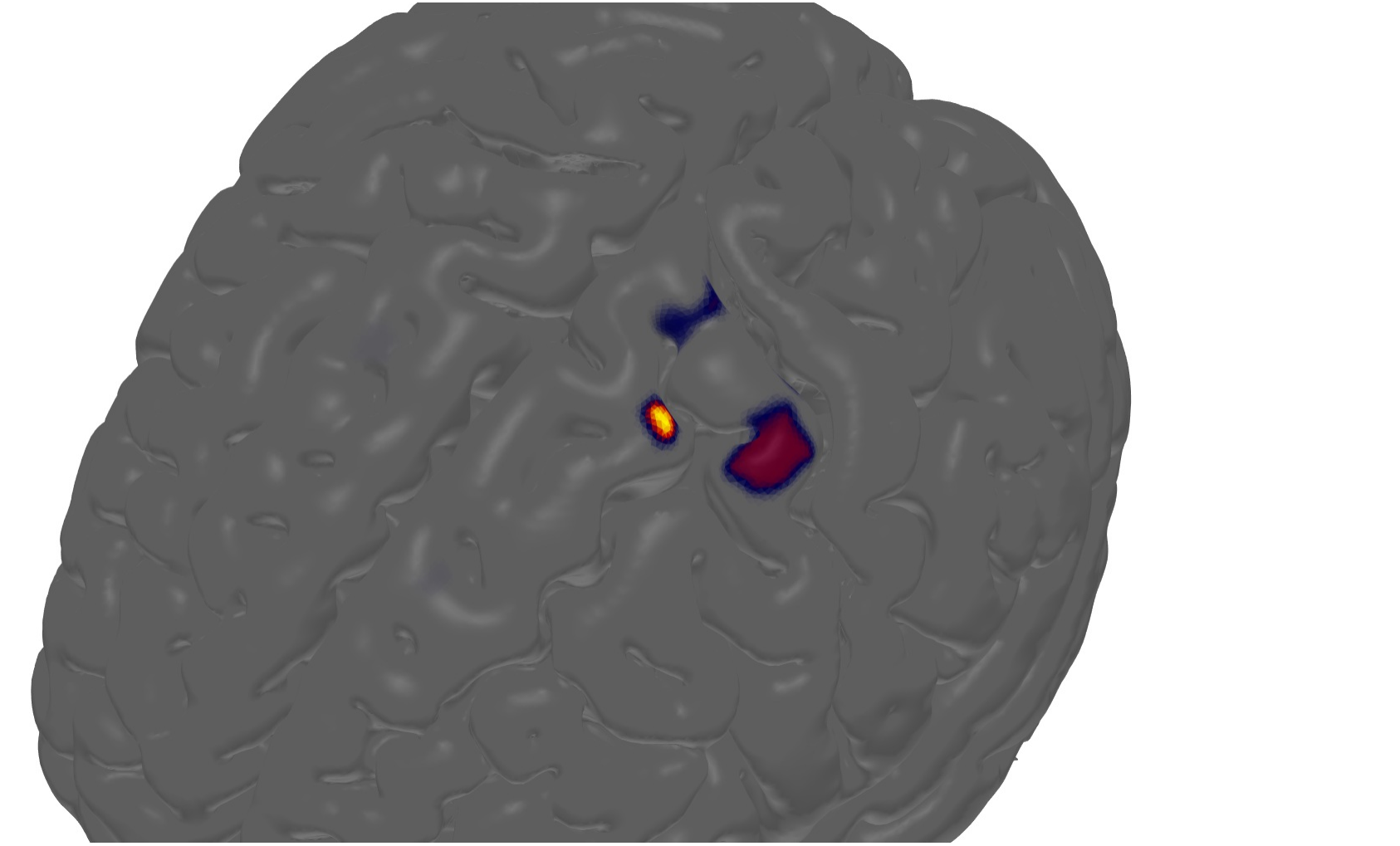}
        \end{minipage}%
        \begin{minipage}[b]{\imageWidth}
            \includegraphics[trim={6cm 2cm 6.5cm 2cm},clip,width=1\linewidth]{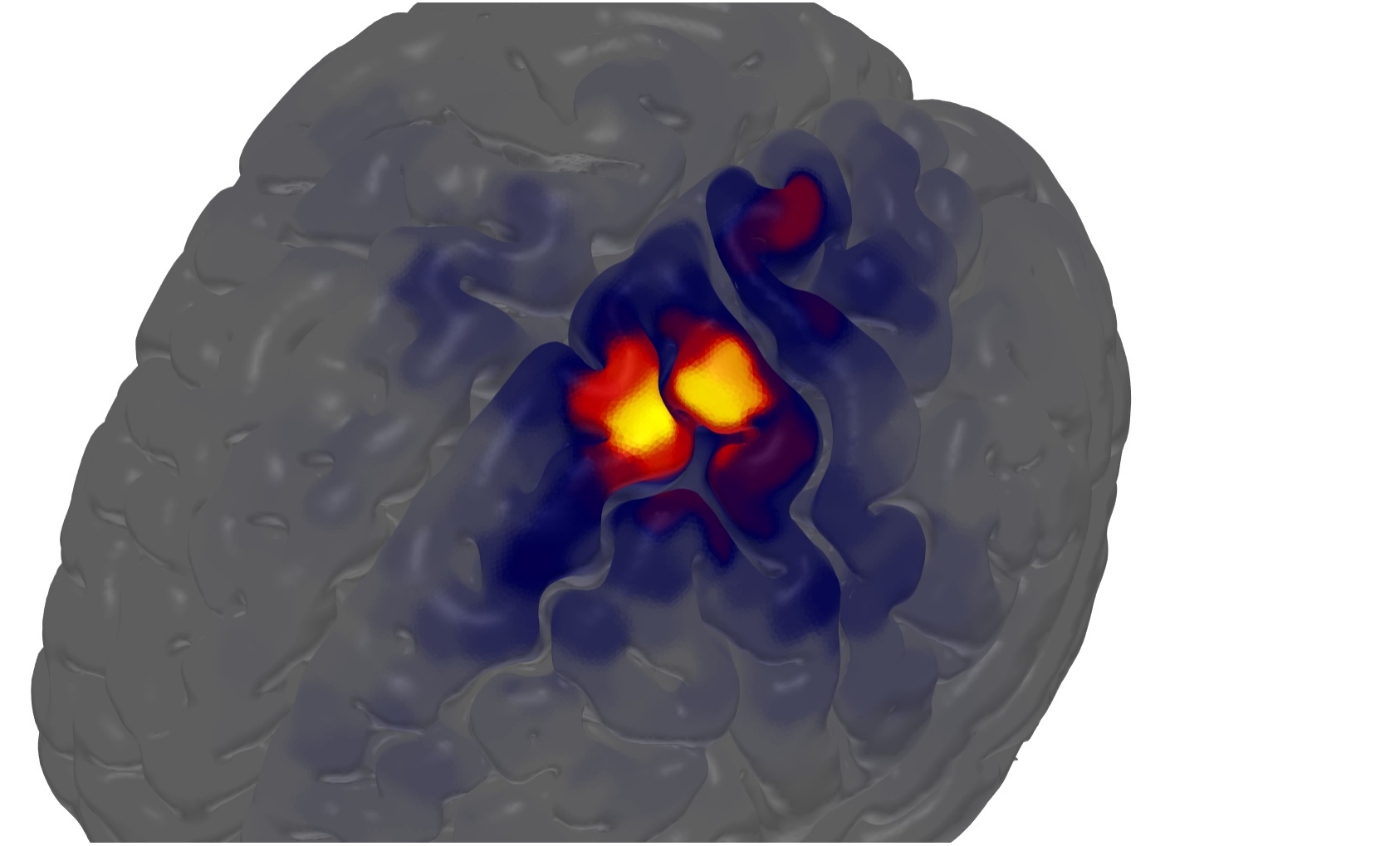}
        \end{minipage}%
        \begin{minipage}[b]{\imageWidth}
            \includegraphics[trim={6cm 2cm 6.5cm 2cm},clip,width=1\linewidth]{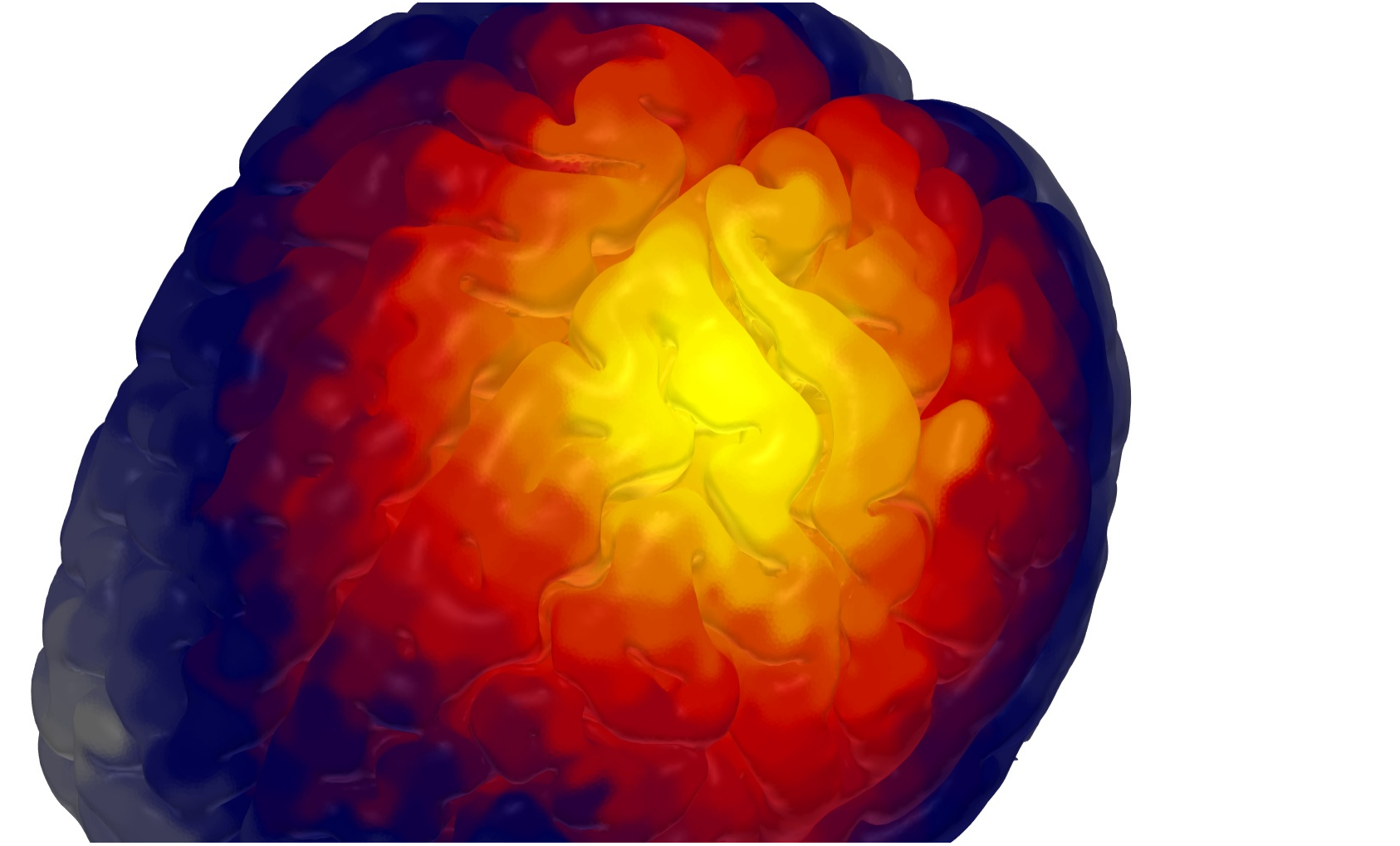}
        \end{minipage}

        \begin{minipage}[b]{\inverseMethodNameWidth}
            \rotatebox{90}{\hspace{0.1cm} Patch}
        \end{minipage}%
        \begin{minipage}[b]{\imageWidth}
            \includegraphics[trim={6cm 2cm 6.5cm 2cm},clip,width=1\linewidth]{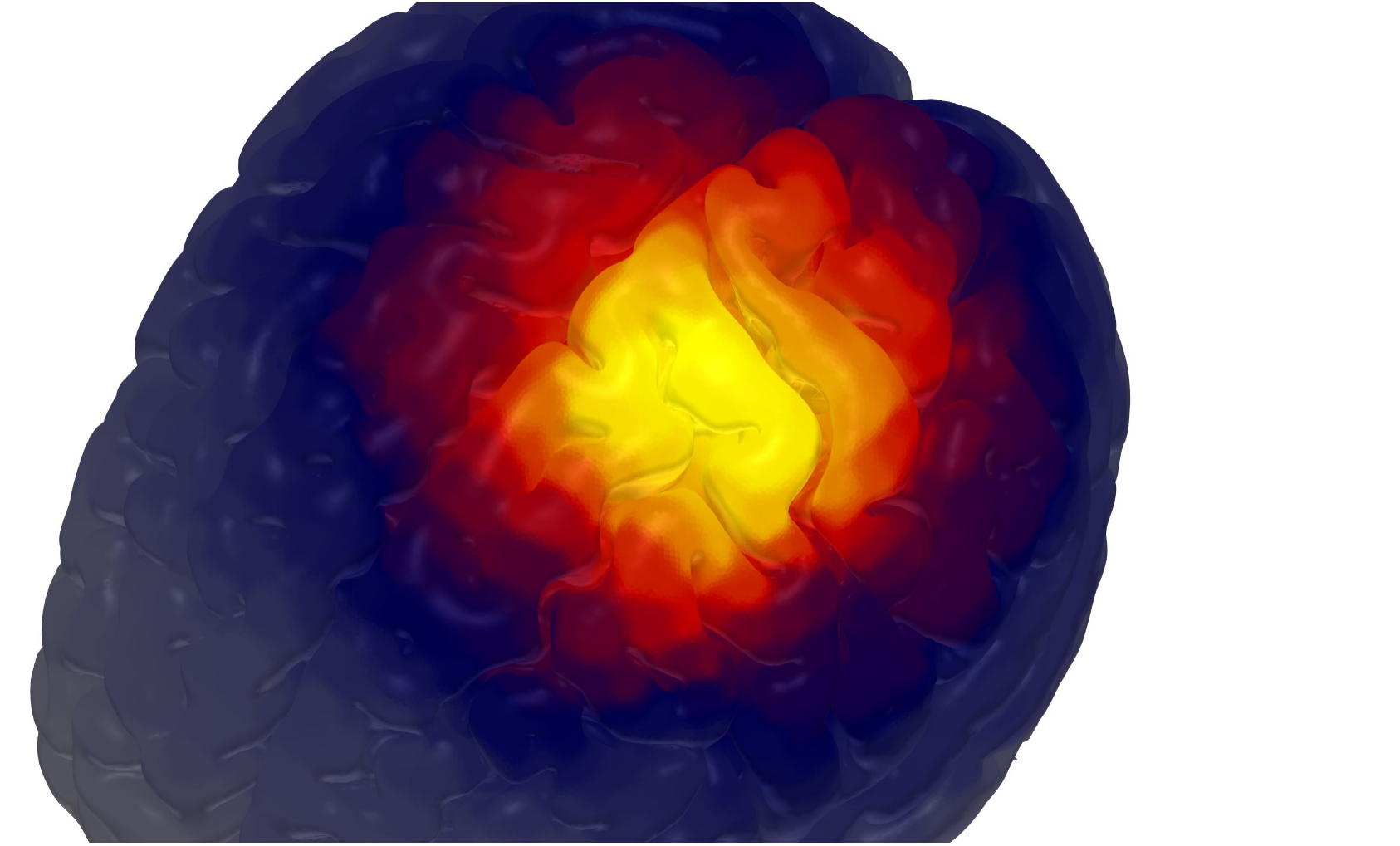}
        \end{minipage}%
        \begin{minipage}[b]{\imageWidth}
            \includegraphics[trim={6cm 2cm 6.5cm 2cm},clip,width=1\linewidth]{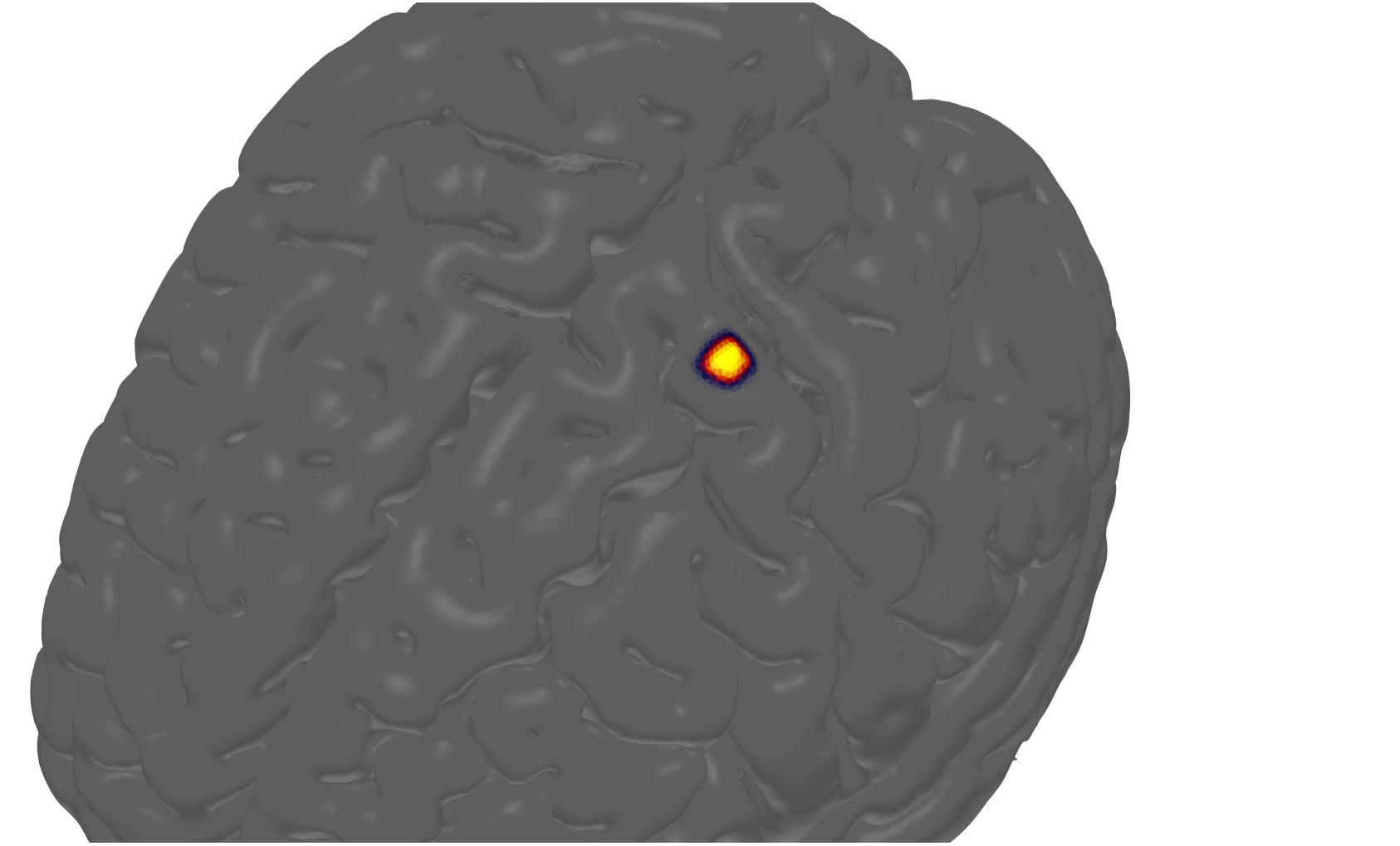}
        \end{minipage}%
        \begin{minipage}[b]{\imageWidth}
            \includegraphics[trim={6cm 2cm 6.5cm 2cm},clip,width=1\linewidth]{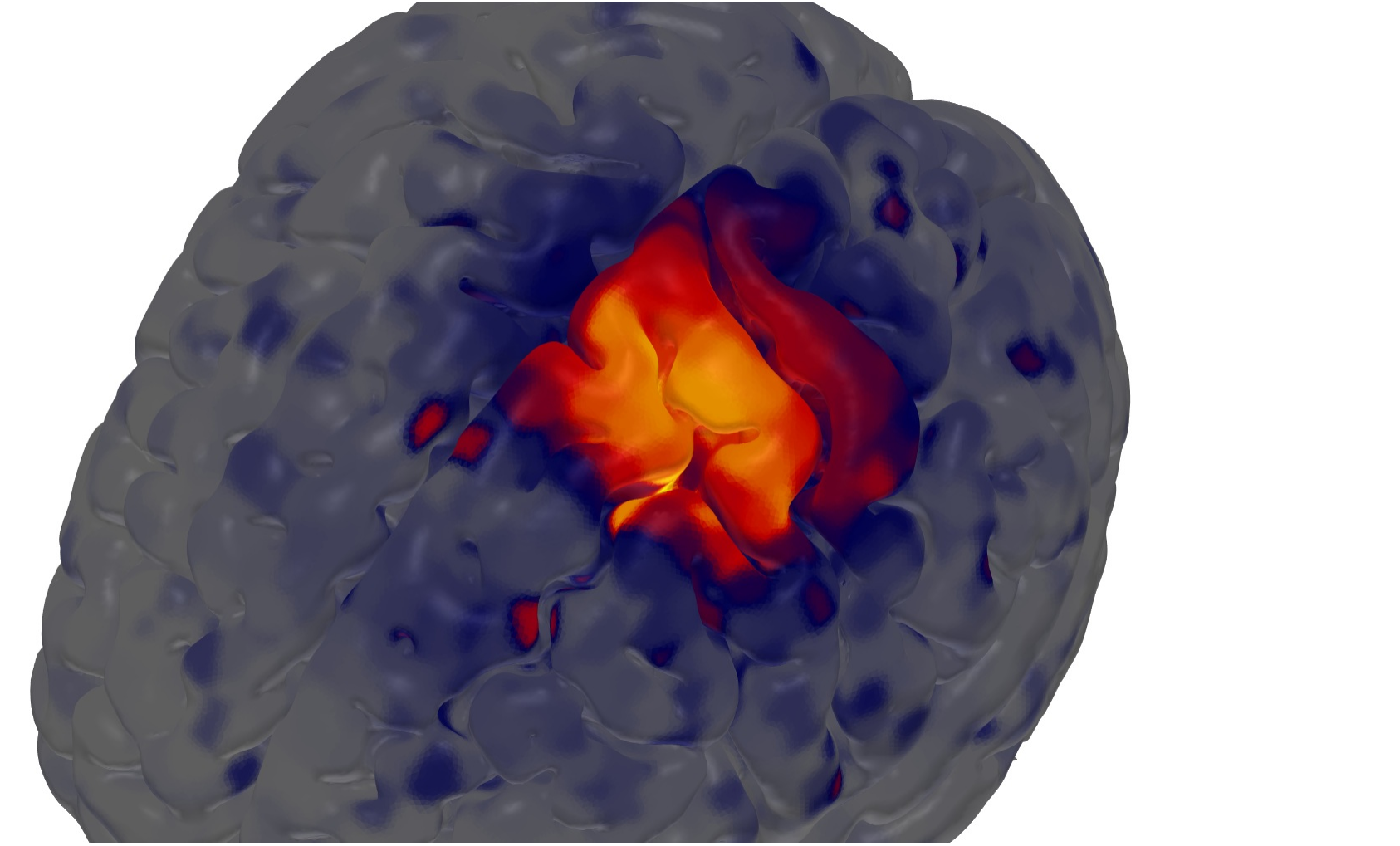}
        \end{minipage}%
        \begin{minipage}[b]{\imageWidth}
            \includegraphics[trim={6cm 2cm 6.5cm 2cm},clip,width=1\linewidth]{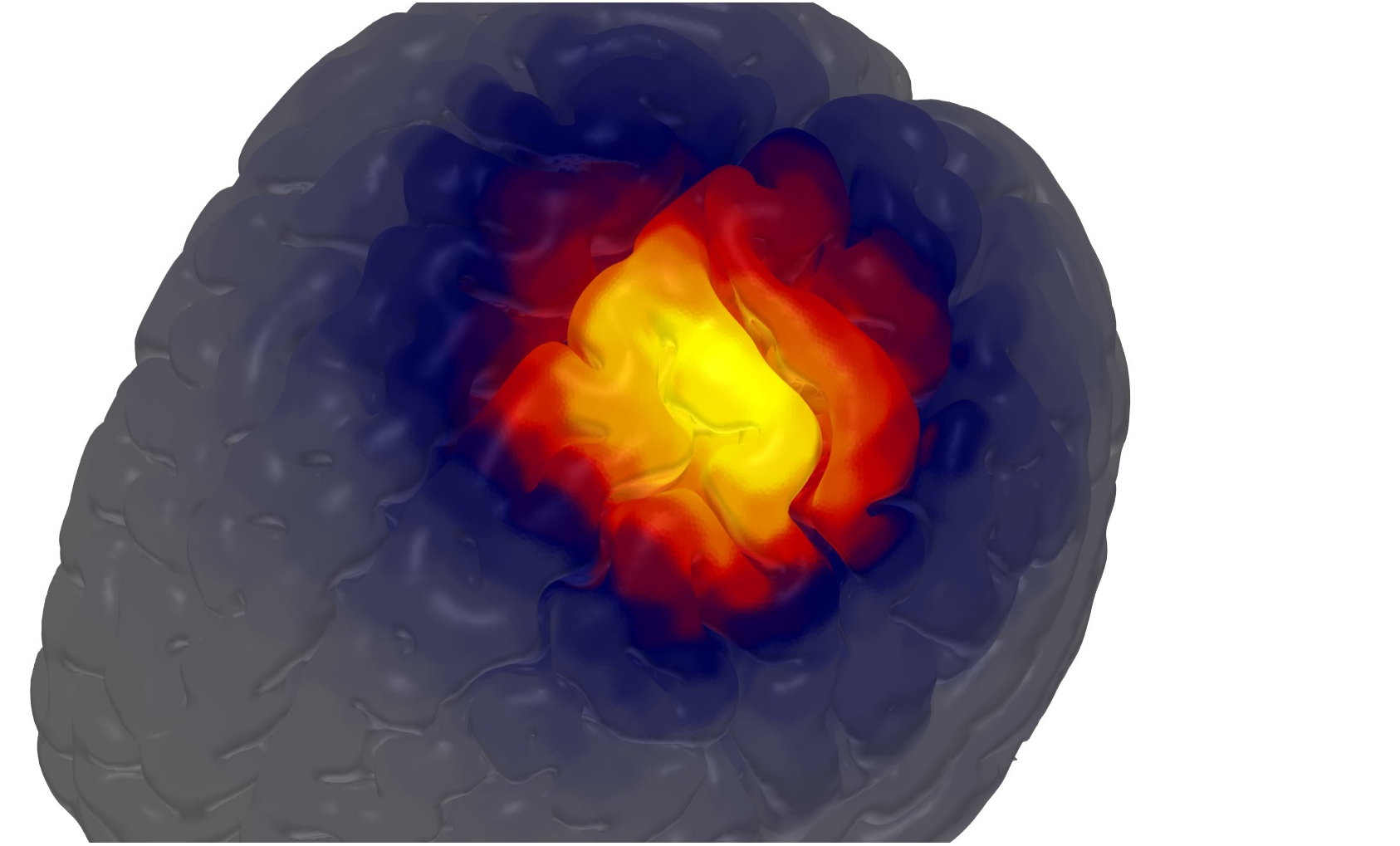}
        \end{minipage}

        \begin{minipage}[b]{\inverseMethodNameWidth}
            \rotatebox{90}{\hspace{0.1cm}Whitney} \rotatebox{90}{\hspace{0.08cm}difference}
        \end{minipage}%
        \begin{minipage}[b]{\imageWidth}
            \includegraphics[trim={6cm 2cm 6.5cm 2cm},clip,width=1\linewidth]{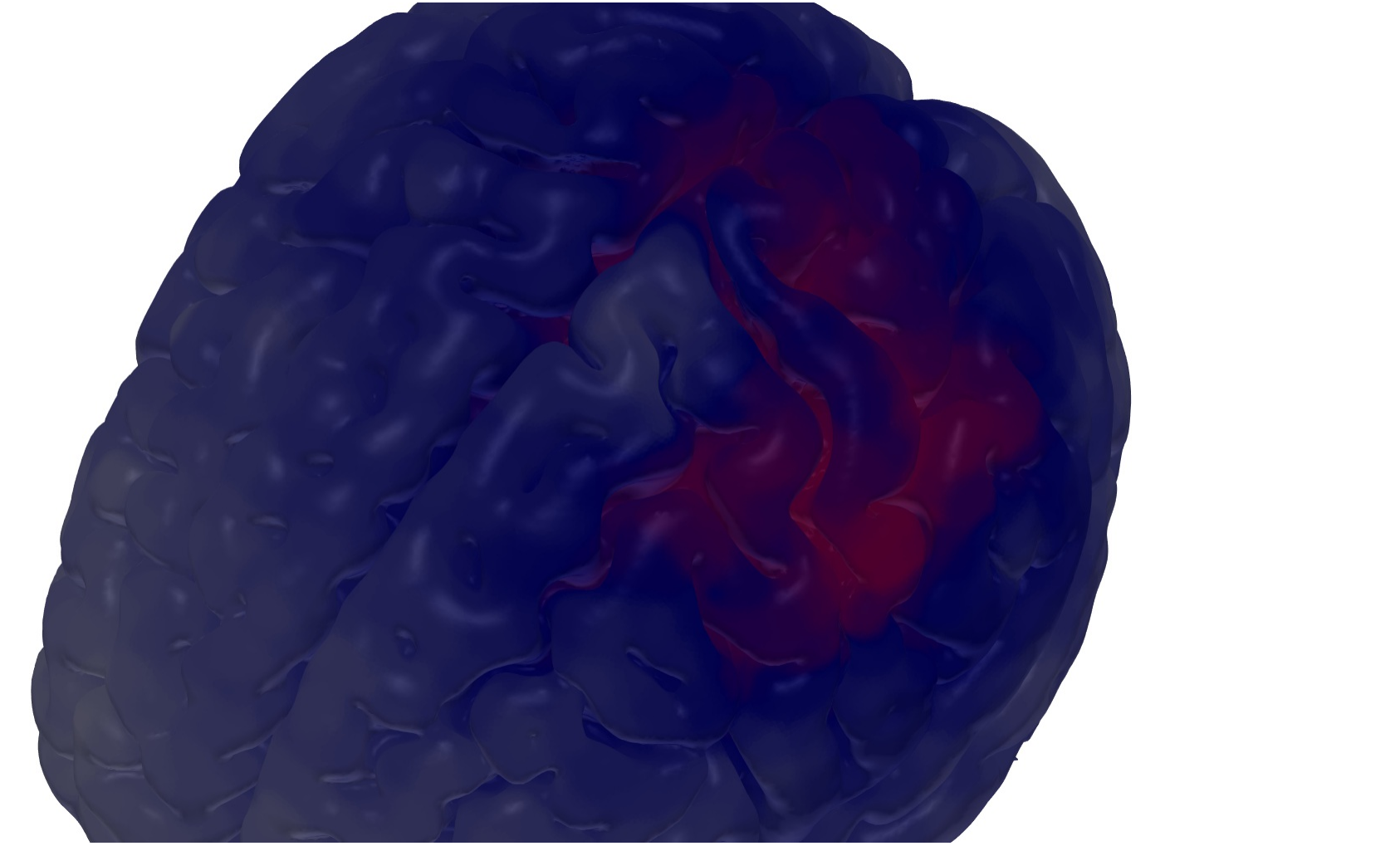}
        \end{minipage}%
        \begin{minipage}[b]{\imageWidth}
            \includegraphics[trim={6cm 2cm 6.5cm 2cm},clip,width=1\linewidth]{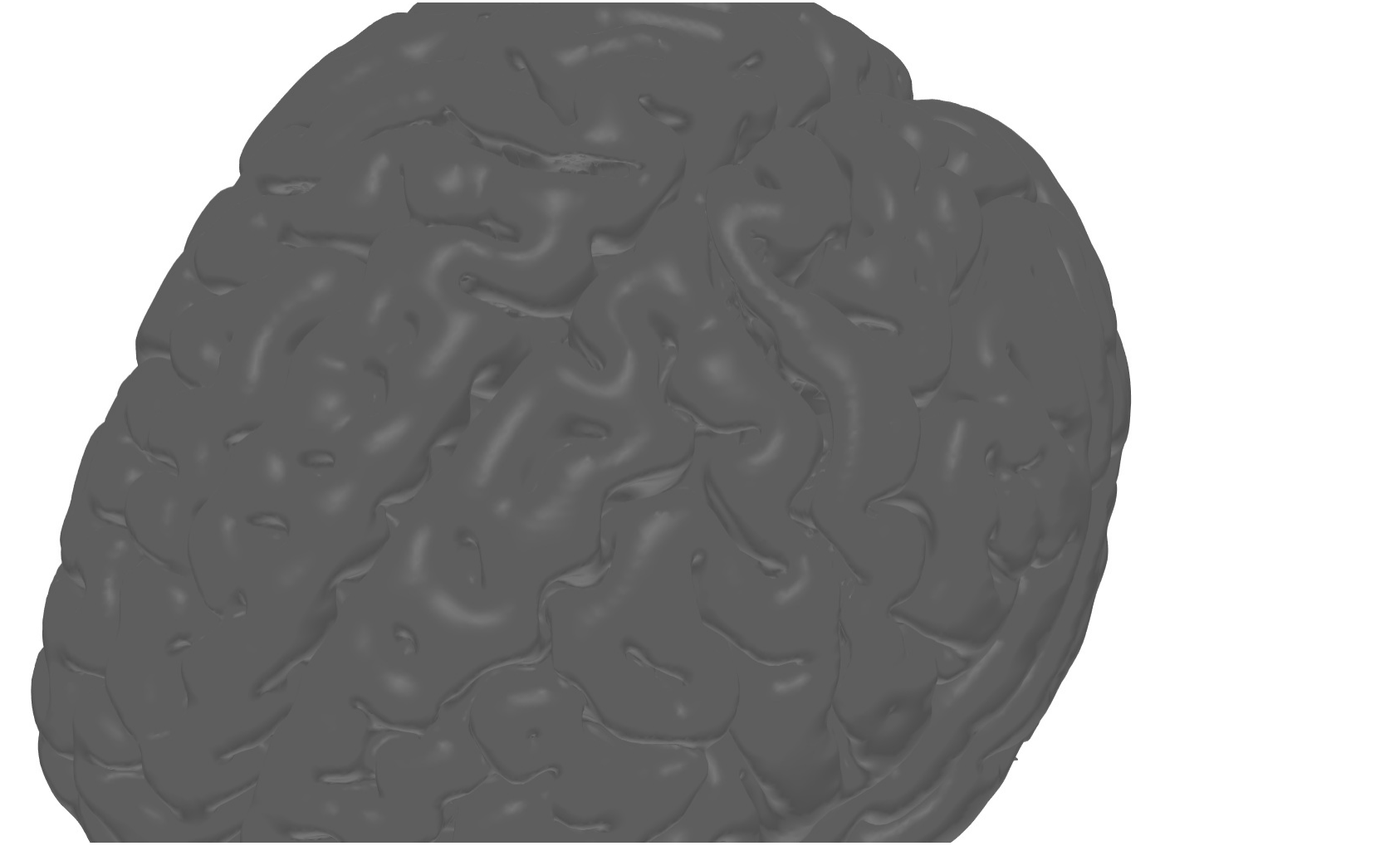}
        \end{minipage}%
        \begin{minipage}[b]{\imageWidth}
            \includegraphics[trim={6cm 2cm 6.5cm 2cm},clip,width=1\linewidth]{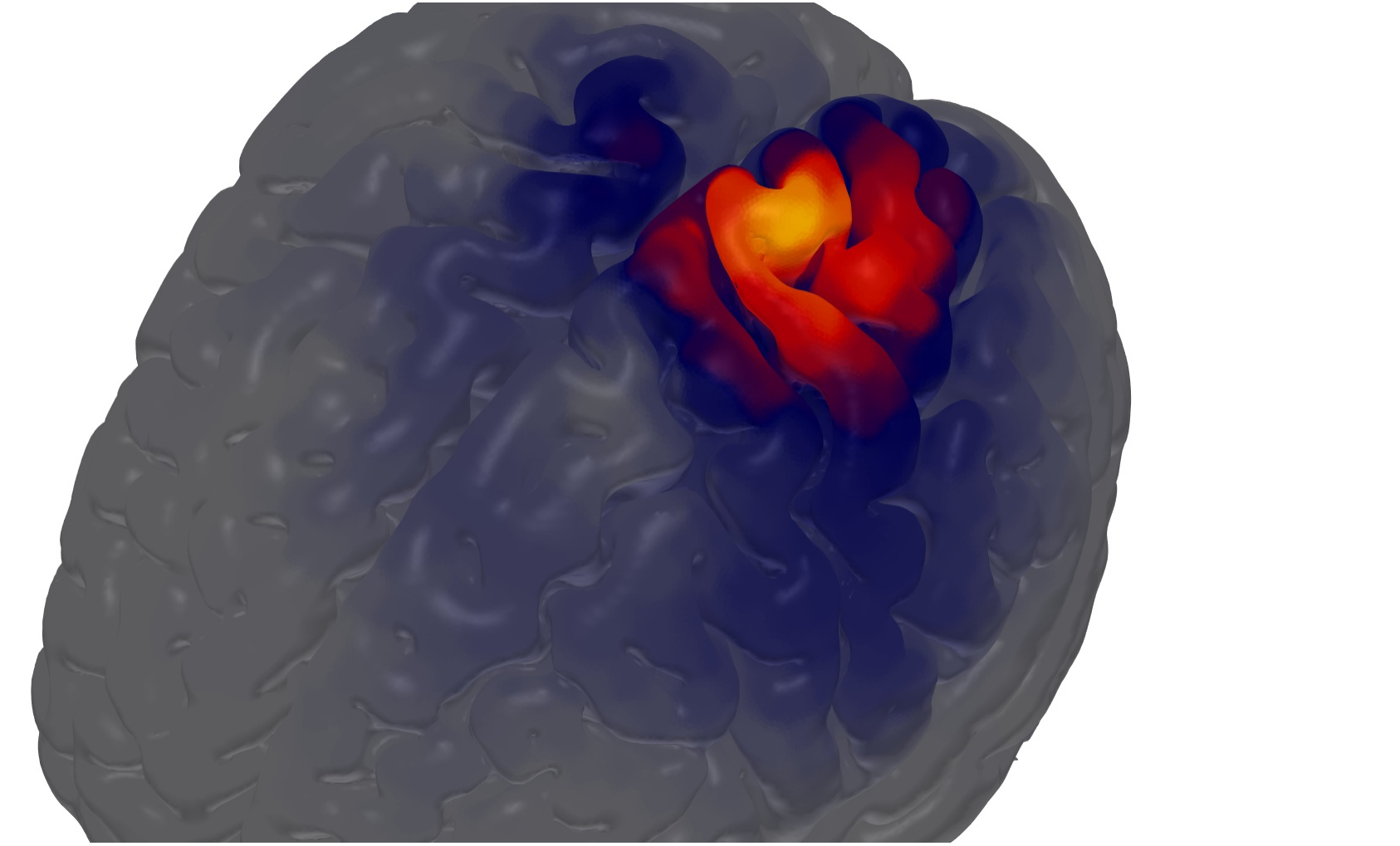}
        \end{minipage}%
        \begin{minipage}[b]{\imageWidth}
            \includegraphics[trim={6cm 2cm 6.5cm 2cm},clip,width=1\linewidth]{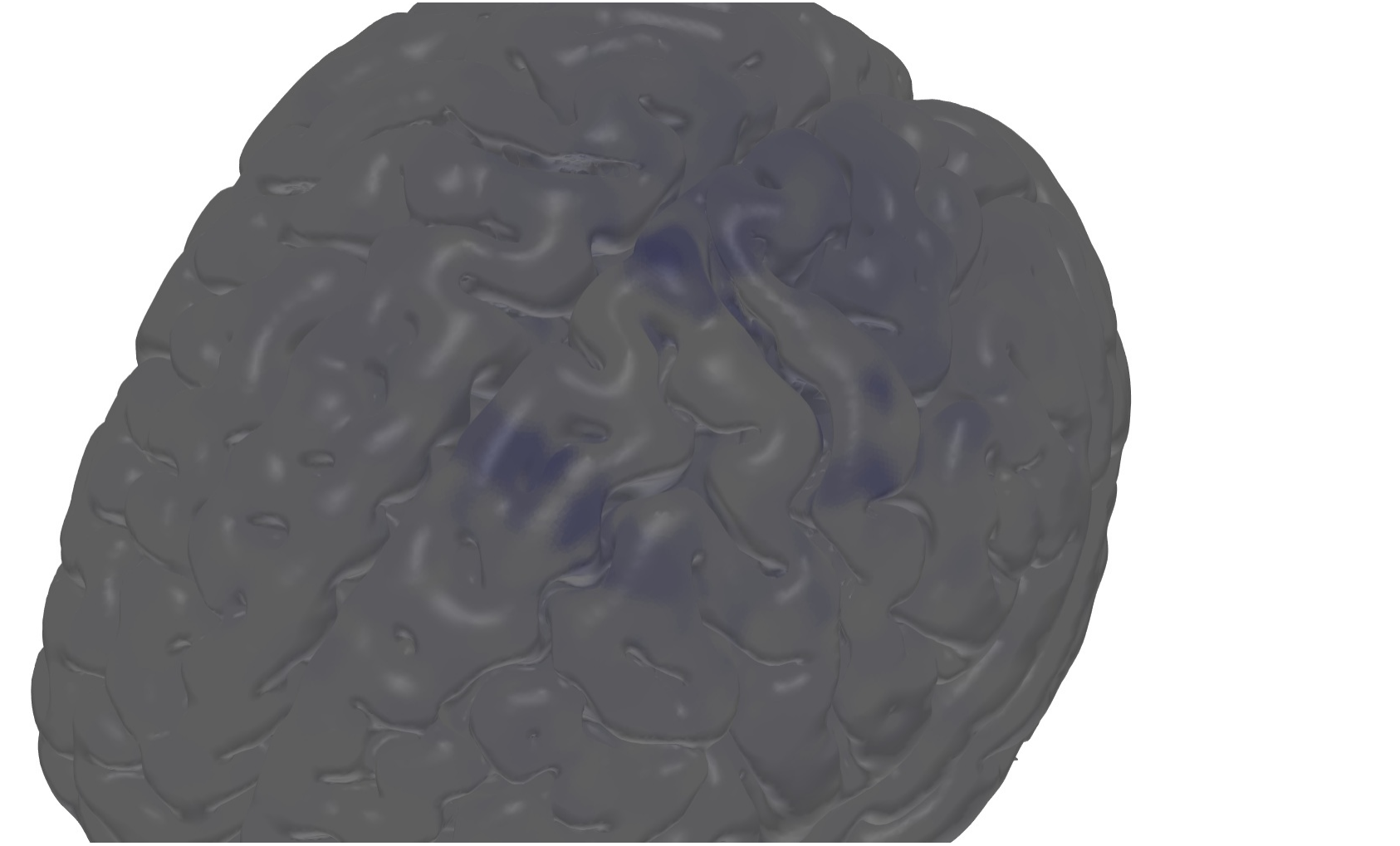}
        \end{minipage}

        \begin{minipage}[b]{\inverseMethodNameWidth}
            \rotatebox{90}{\hspace{0.2cm}$\Hdiv$} \rotatebox{90}{difference}
        \end{minipage}%
        \begin{minipage}[b]{\imageWidth}
            \includegraphics[trim={6cm 2cm 6.5cm 2cm},clip,width=1\linewidth]{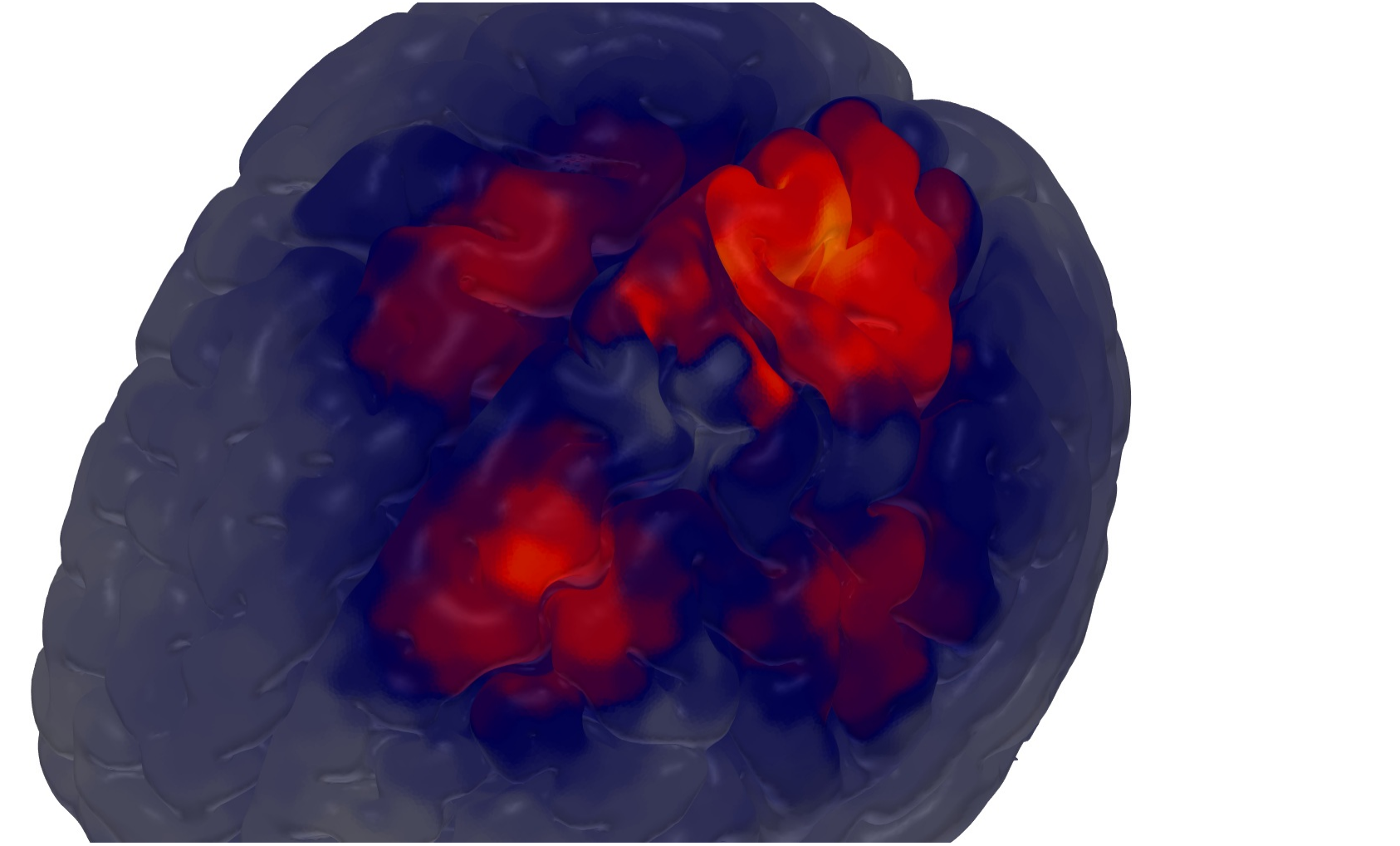}
        \end{minipage}%
        \begin{minipage}[b]{\imageWidth}
            \includegraphics[trim={6cm 2cm 6.5cm 2cm},clip,width=1\linewidth]{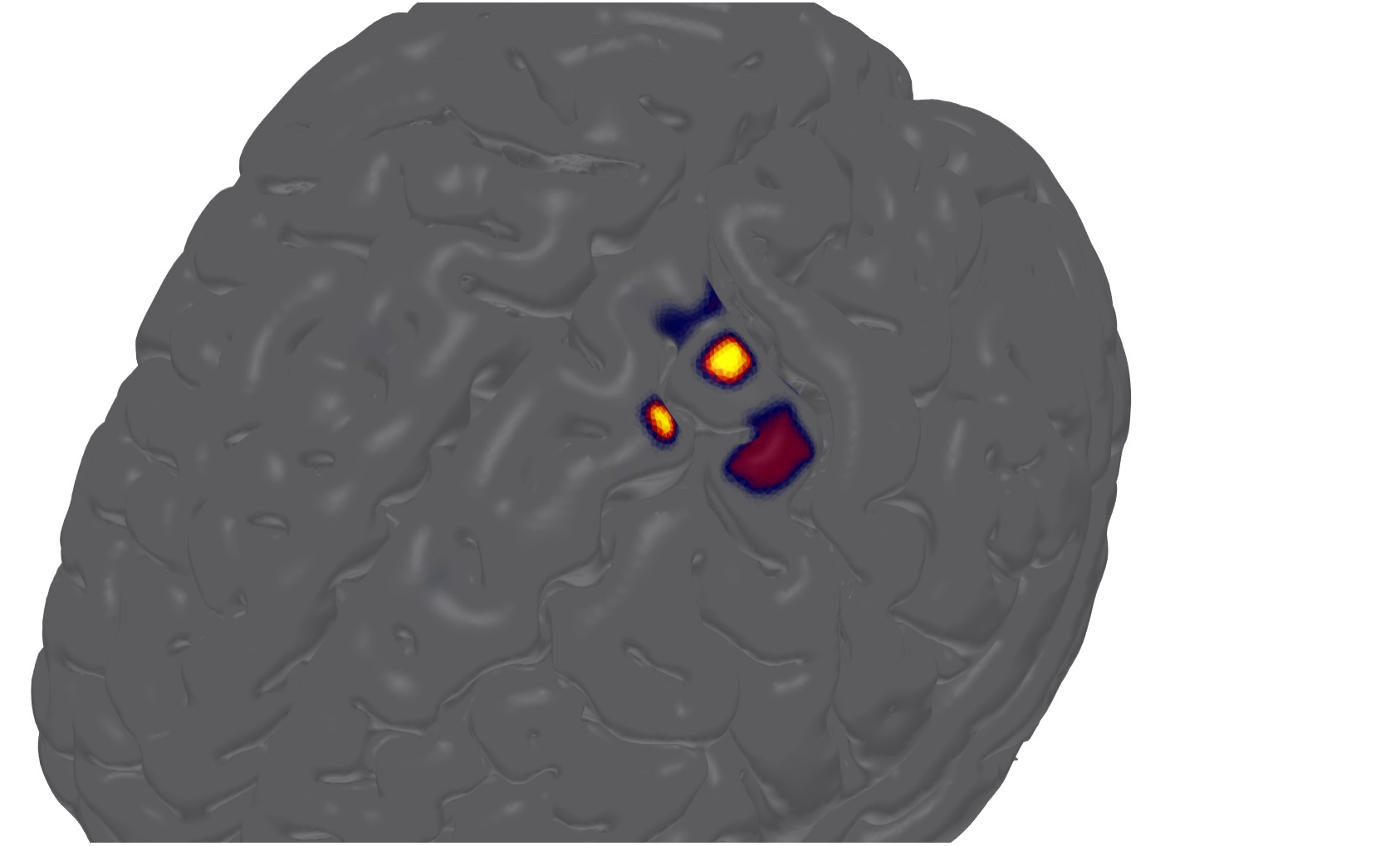}
        \end{minipage}%
        \begin{minipage}[b]{\imageWidth}
            \includegraphics[trim={6cm 2cm 6.5cm 2cm},clip,width=1\linewidth]{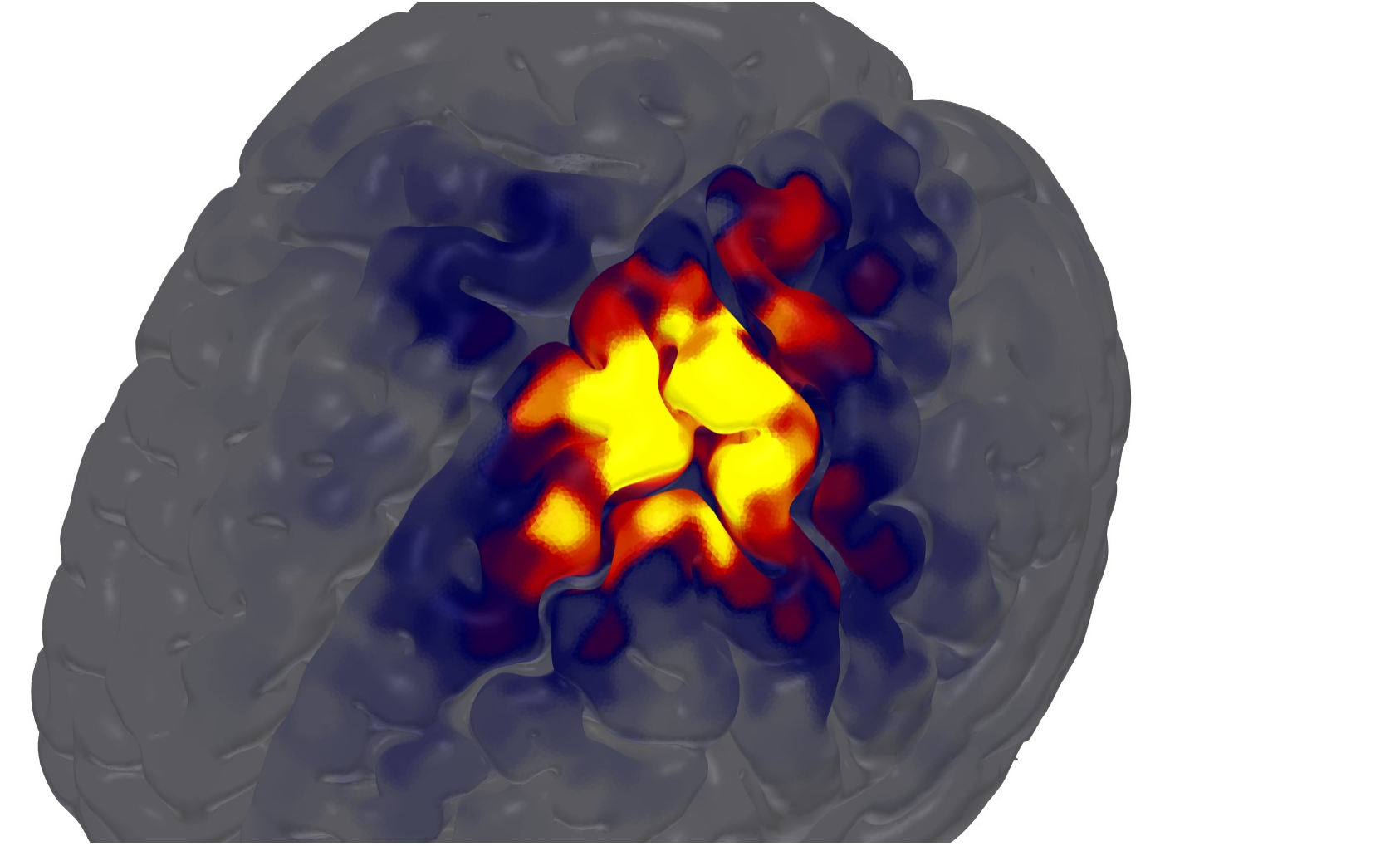}
        \end{minipage}%
        \begin{minipage}[b]{\imageWidth}
            \includegraphics[trim={6cm 2cm 6.5cm 2cm},clip,width=1\linewidth]{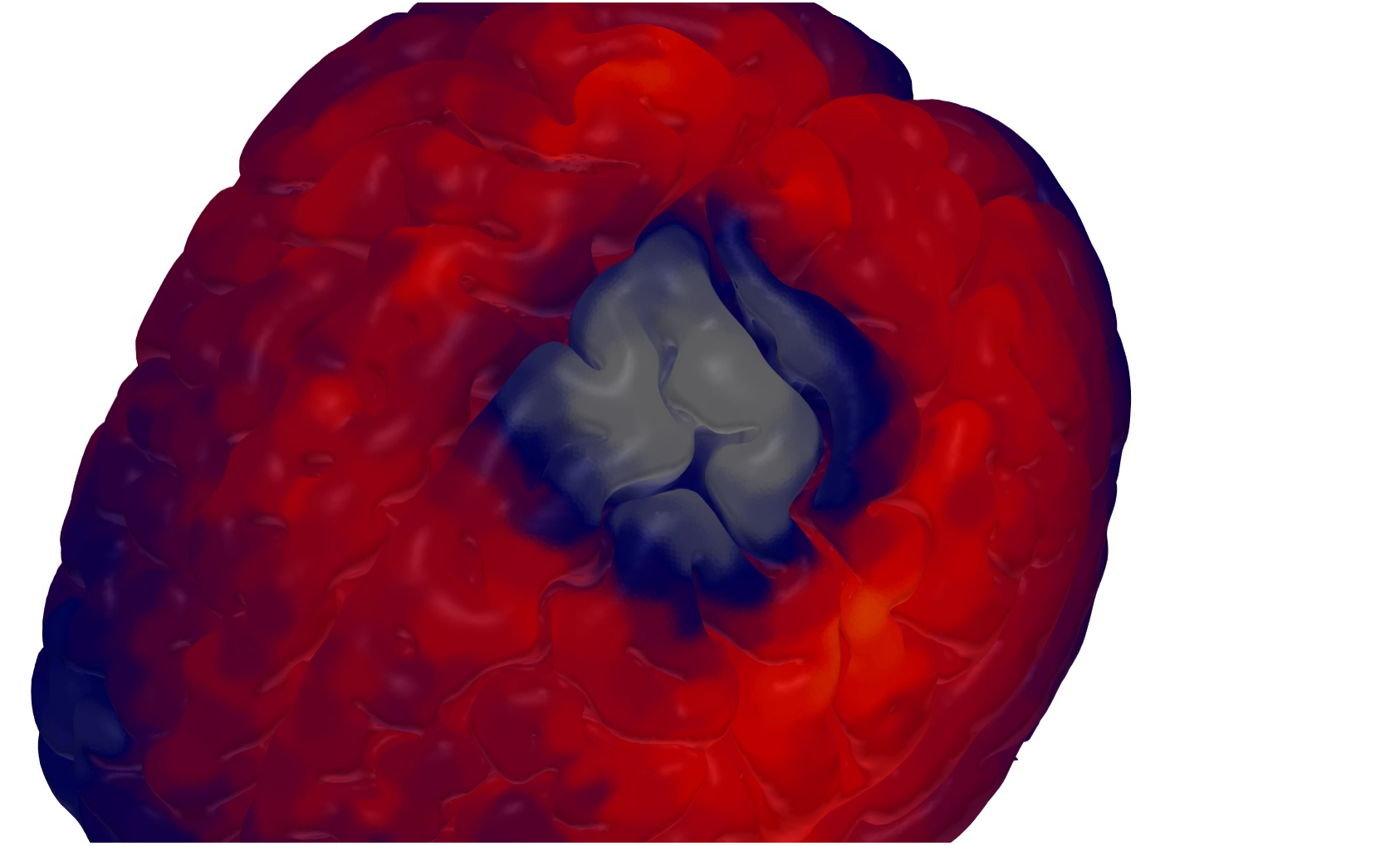}
        \end{minipage}

    \begin{minipage}[b]{\inverseMethodNameWidth}
            \rotatebox{90}{\hspace{0.1cm}Patch} \rotatebox{90}{difference}
        \end{minipage}%
        \begin{minipage}[b]{\imageWidth}
            \includegraphics[trim={6cm 2cm 6.5cm 2cm},clip,width=1\linewidth]{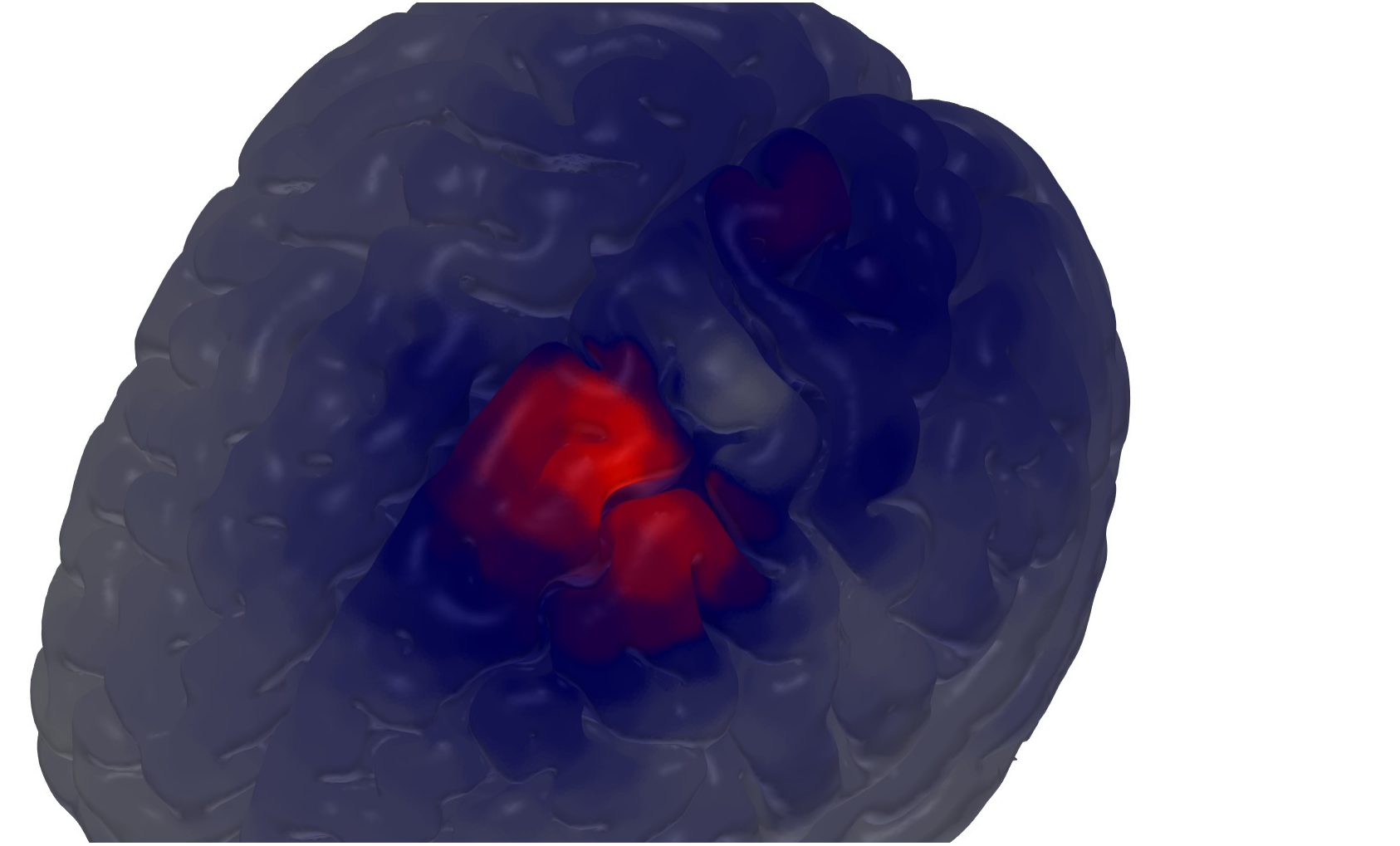}
        \end{minipage}%
        \begin{minipage}[b]{\imageWidth}
            \includegraphics[trim={6cm 2cm 6.5cm 2cm},clip,width=1\linewidth]{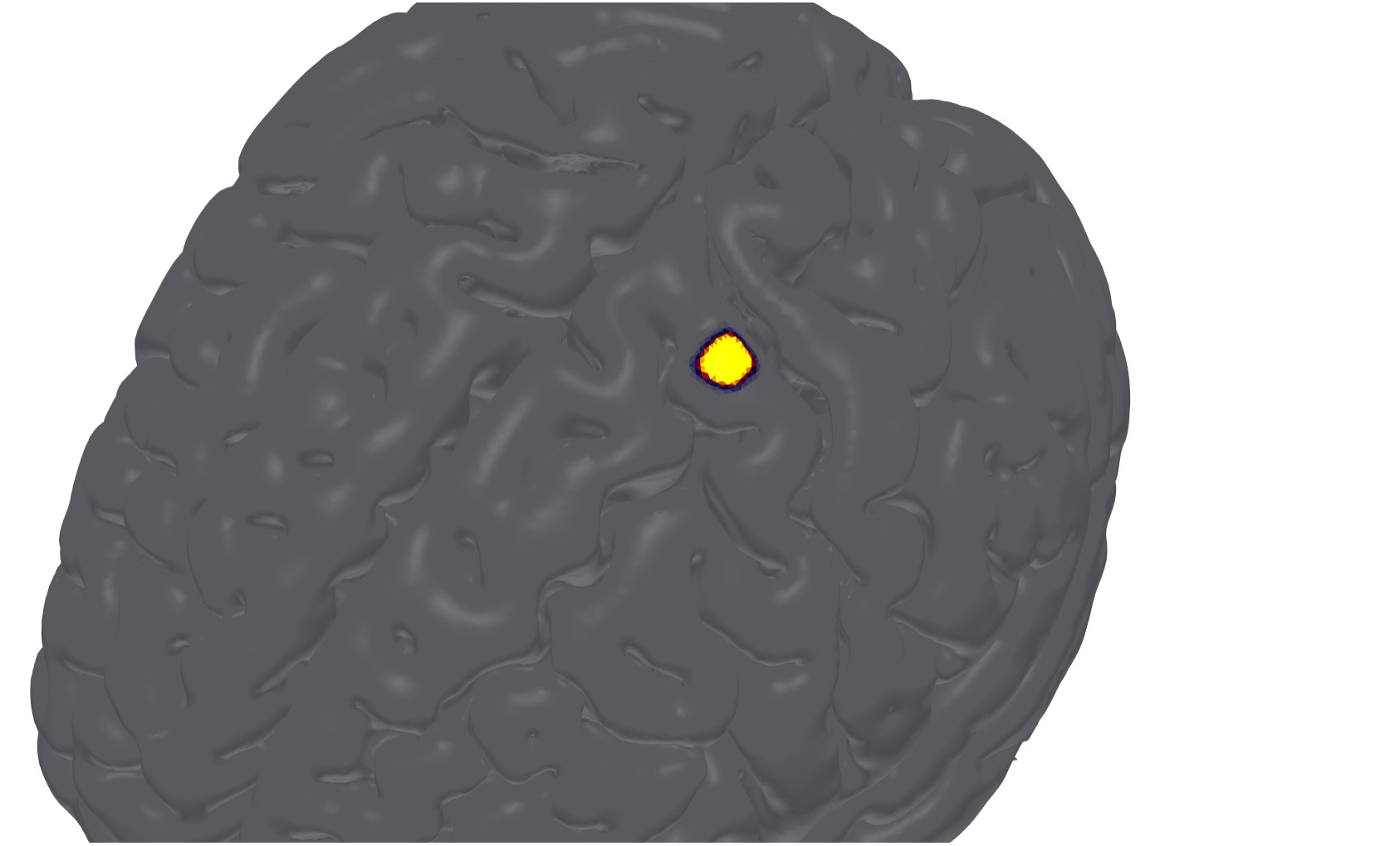}
        \end{minipage}%
        \begin{minipage}[b]{\imageWidth}
            \includegraphics[trim={6cm 2cm 6.5cm 2cm},clip,width=1\linewidth]{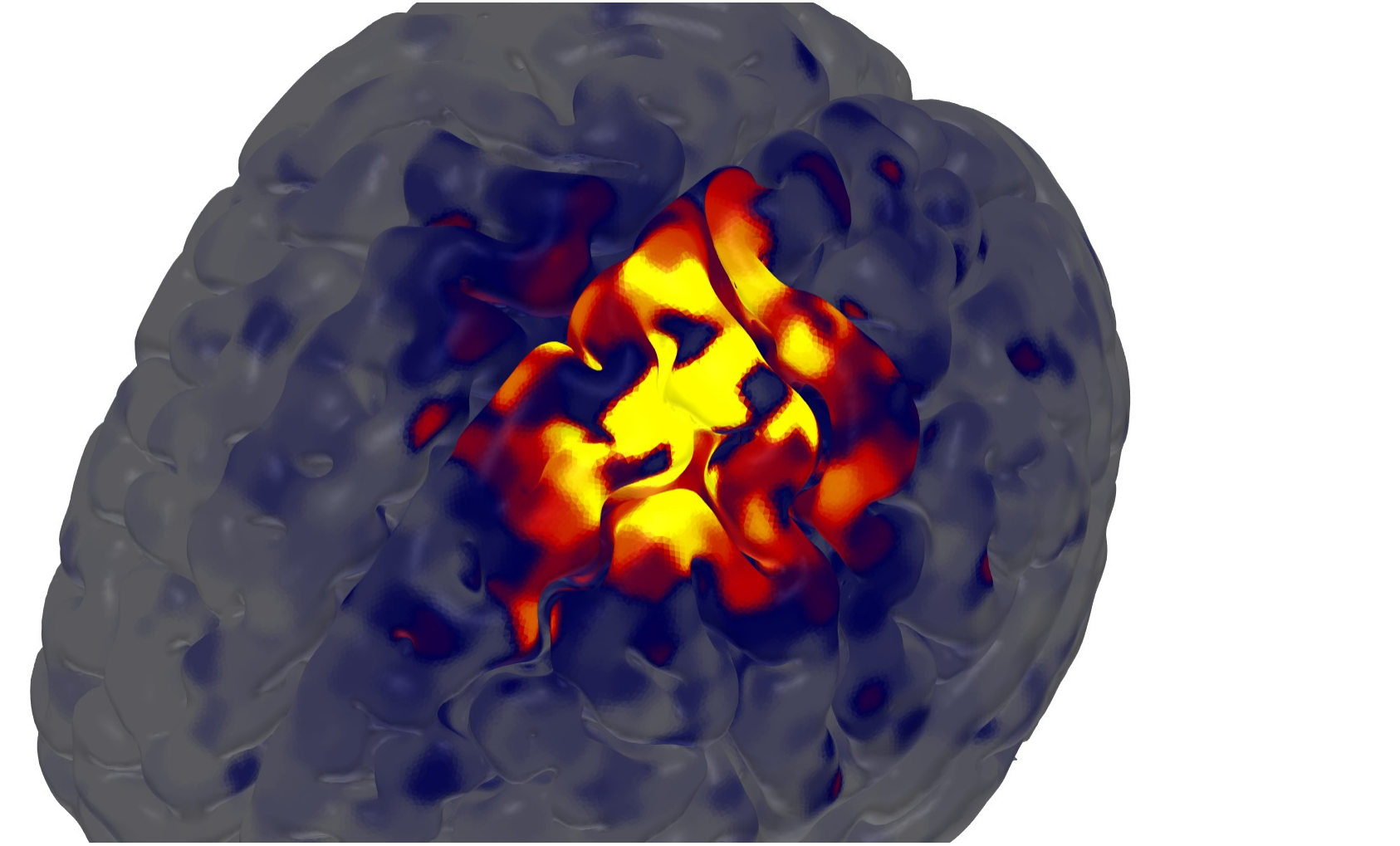}
        \end{minipage}%
        \begin{minipage}[b]{\imageWidth}
            \includegraphics[trim={6cm 2cm 6.5cm 2cm},clip,width=1\linewidth]{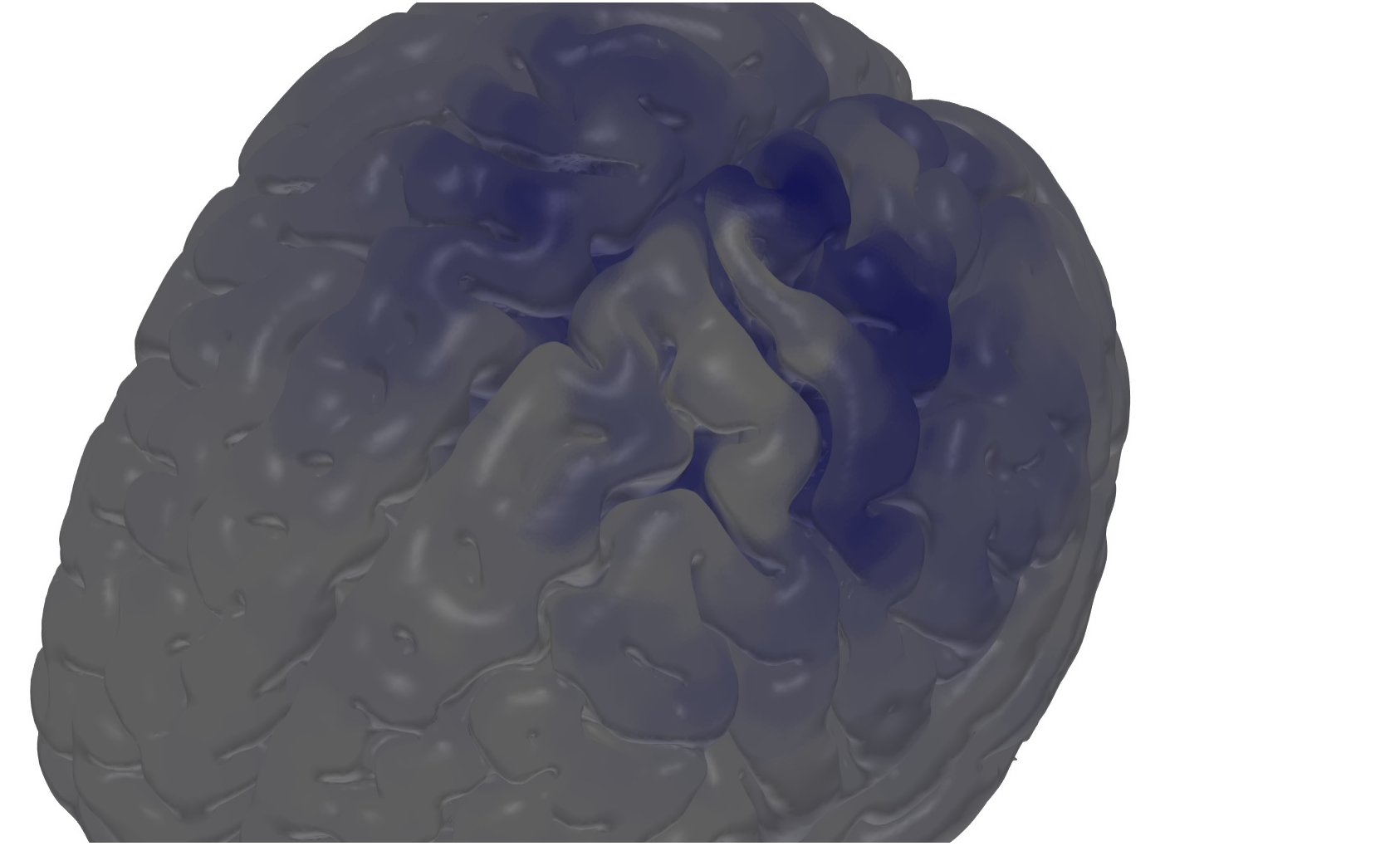}
        \end{minipage}
    \end{minipage}
    \hfill
    \begin{minipage}[b]{\colorbarWidth}
        \includegraphics[width=\linewidth]{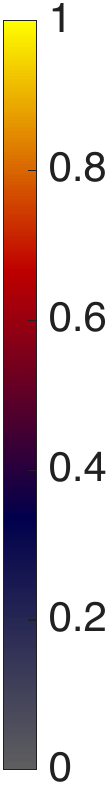}
    \end{minipage}
    \hfill
    \caption{Close up of reconstructions of a synthetic cortical dipole with sLORETA, SHAL1R, SKF, and DS using DUNEuro's Whitney basis function implementation (row 1), DUNEuro's Local subtraction (row 2), and Zeffiro Interface's $\Hdiv$ (row 3). Row 4 shows reconstructions of a patch-like source with Local subtraction as the interpolation method. Row 5 shows the difference between Whitney and Local subtraction reconstructions, while row 6 does the same for $\Hdiv$ and Local subtraction. Row 7 is the difference between a patch-like source reconstruction and focal Local subtraction. SHAL1R reconstructions show the true source location for Whitney and Local subtraction, and do not differ from each other. A SHAL1R reconstruction of a $\Hdiv$ source spreads the source across 3 positions near the true source position, resulting in an observable difference between it and a Local subtraction. The patch-like SHAL1R reconstruction difference behaves similarly to Whitney and Local subtraction, but is slightly more spread out.}
    \label{fig.erotuskuva}
\end{figure}

The estimated activity distribution of the cortical dipole source is presented in Figure~\ref{fig.erotuskuva} with the distributions from which the reference distribution is subtracted (last three rows). Here, we use the estimated distribution obtained with Local subtraction as the reference. The first three rows display the estimations obtained with Whitney-basis, Local subtraction, and $\Hdiv$, respectively. Additionally, the fourth row displays these results for a simulation of a patch-like source obtained by smoothing a high-density source space using \eqref{eq.gaussian}, with $\position$ restricted to the original focal source positions, while the window centers vary over all $\sourcePosition$. The source estimates are obtained using sLORETA, SHAL1R, SKF, and DS, as presented in the subsequent columns.

The sLORETA source estimation distributions in the first column show similar shapes for Whitney and Local subtraction, with the difference between the distributions appearing in the tail, as shown in the fifth row. The most concentrated, i.e., focal distribution is observed with $\Hdiv$, while the highest estimation values ($>0.9$), depicted in yellow, are slightly off from the location of the true source. The widest estimation distribution is perceived with the patch model.

Based on the estimation distribution, SHAL1R produces the most focal estimate, concentrated at the true source location, with Whitney and Local subtraction. With $\Hdiv$, we can observe activity around the true source location but not at the position of the true source. With the patch model, the estimate is slightly spread from a point, and the orientation of the reconstructed dipole is incorrect; thus, there is a large difference in the distribution shown in the last row.

SKF provides highly similar estimation distributions as sLORETA; however, the estimates are slightly more focal and irregular. Moreover, SKF estimation is highly focused on its 0.9 magnitude range. It is also less robust to changes in the chosen lead-field interpolation approach, which can be seen in the large errors when the reference distribution is subtracted from the distributions computed with Whitney, $\Hdiv$, and patch models. Unlike sLORETA, the most concentrated estimate is observed with the patch model, where the highest peak, depicted in yellow, is at the sulcus. The second-most-focal distribution appears with $\Hdiv$.

The goodness-of-fit value distributions of DS are similar to those of Whitney, Local subtraction, and the patch model, although the patch model's distribution is smoother than the other two. A significant spread for the highest values is obtained with $\Hdiv$. Out of these four inversion methods, the smallest differences near the actual source position are obtained with DS, although after a certain distance from the true source, SHAL1R reconstructions coincide better.

The results obtained with the $\Hdiv$ source model and the reconstruction differences relative to Local subtraction show that $\Hdiv$ yields wider reconstructions for sLORETA, SHAL1R, and DS. SKF, however, exhibits more focal estimation than the ones obtained with Local subtraction and the Whitney basis, while the reconstruction is spotty. The patch-like reconstruction behaves as expected, with all reconstructions spanning a larger volume than the focal lead field interpolation methods. With sLORETA and DS, the largest differences between the patch-like and focal Local subtraction approach occur away from the true source position, while with SHAL1R and SKF, the differences appear in the vicinity of the true source position. Surprising spottiness is evident in the SKF reconstruction, resembling that seen with $\Hdiv$. The effect is especially pronounced in the estimation-distribution differences between Patch and $\Hdiv$ and Local subtraction.

\subsection{Butterfly plots of synthetic cortical activity}\label{sec.butterfly.plots}

\begin{figure*}
    \centering
    \begin{minipage}{0.23\linewidth}
    \centering
    \footnotesize{Whitney}\vspace{0.1cm}
        \includegraphics[width=\linewidth]{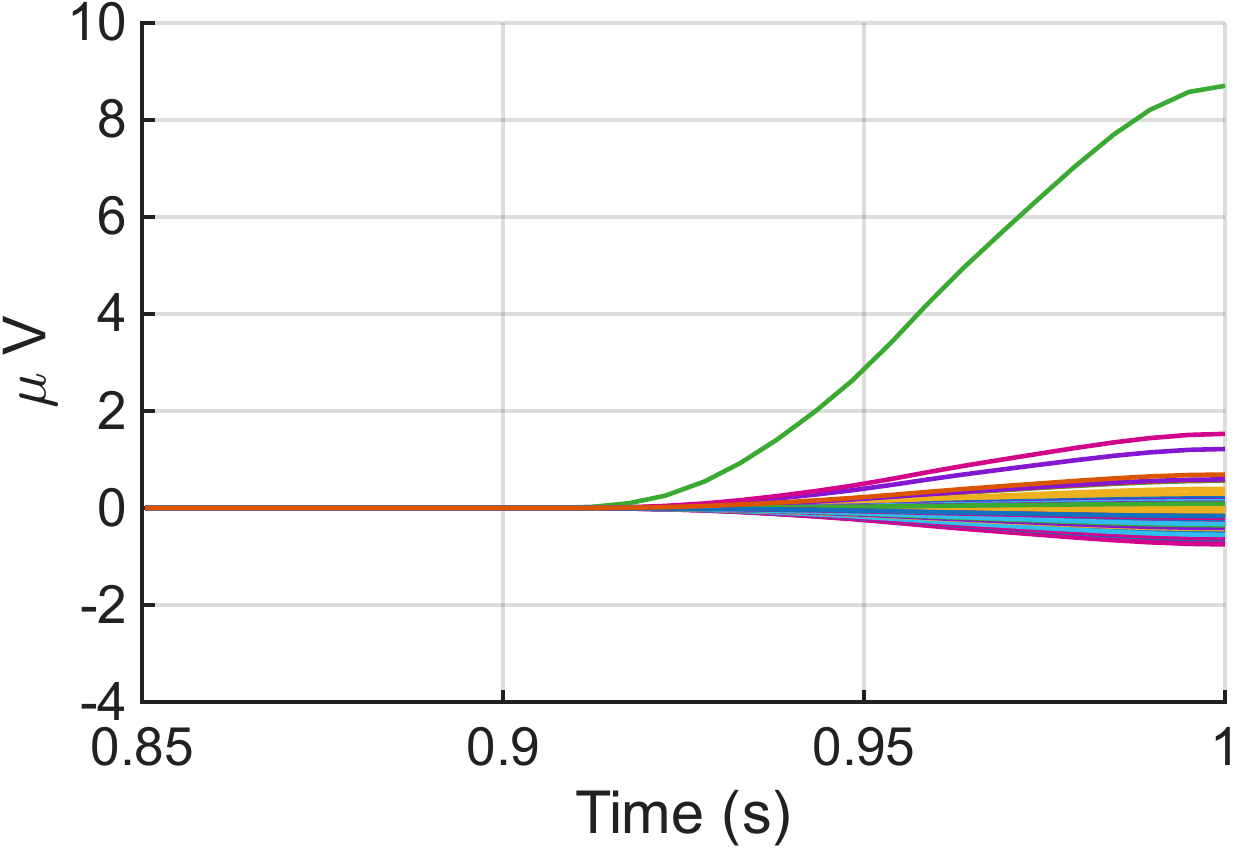}
    \end{minipage}\begin{minipage}{0.23\linewidth}
    \centering
    \footnotesize{Local subtraction}\vspace{0.1cm}
        \includegraphics[width=\linewidth]{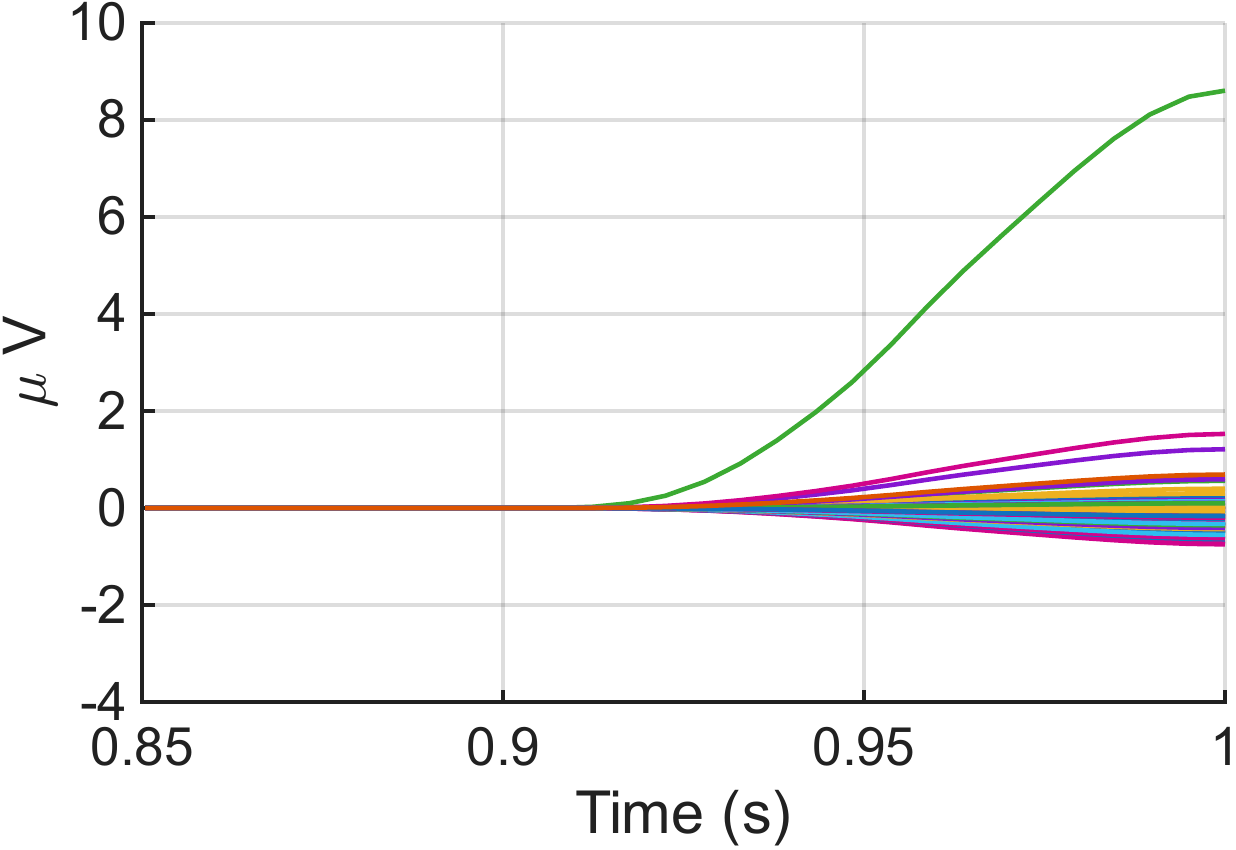}
    \end{minipage}\begin{minipage}{0.23\linewidth}
    \centering
    \footnotesize{$\Hdiv$}\vspace{0.1cm}
        \includegraphics[width=\linewidth]{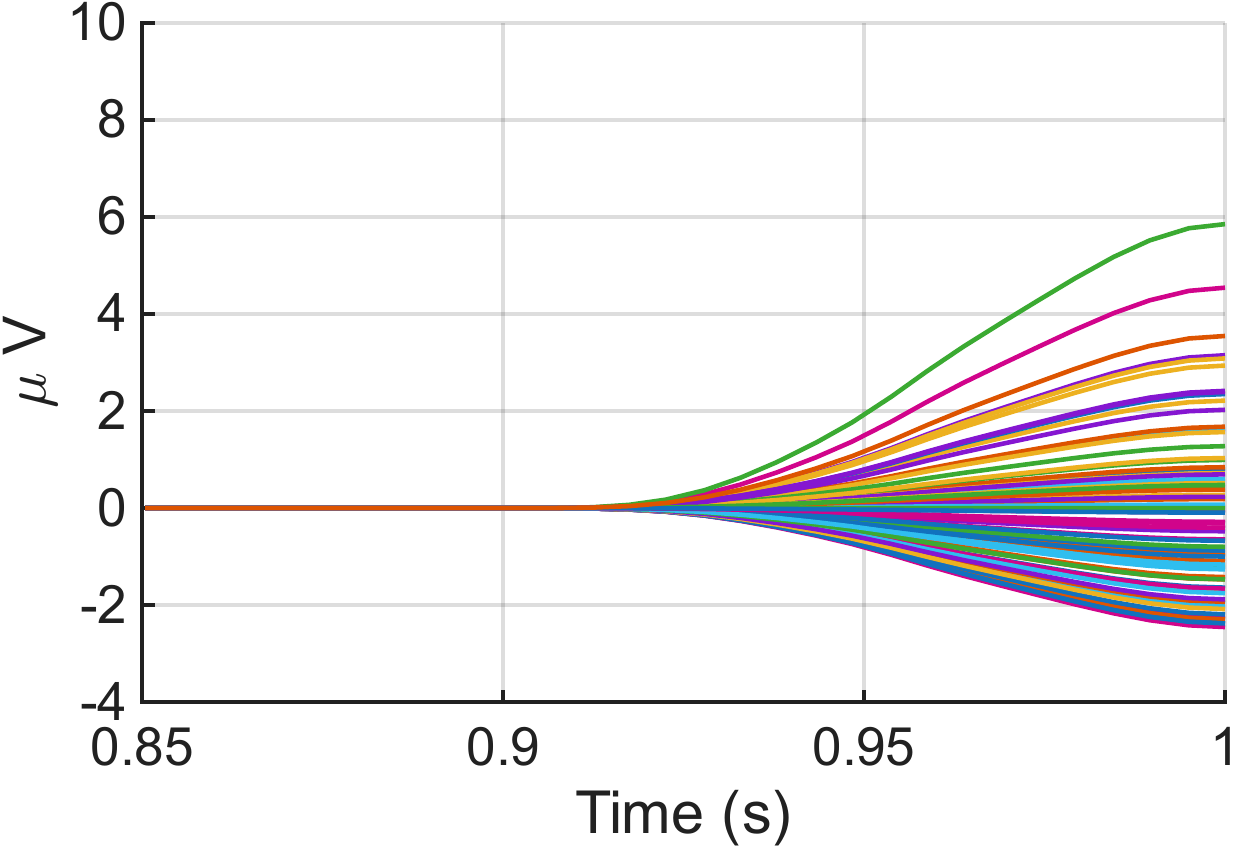}
    \end{minipage}\begin{minipage}{0.23\linewidth}
    \centering
    \footnotesize{Patch}\vspace{0.1cm}
        \includegraphics[width=\linewidth]{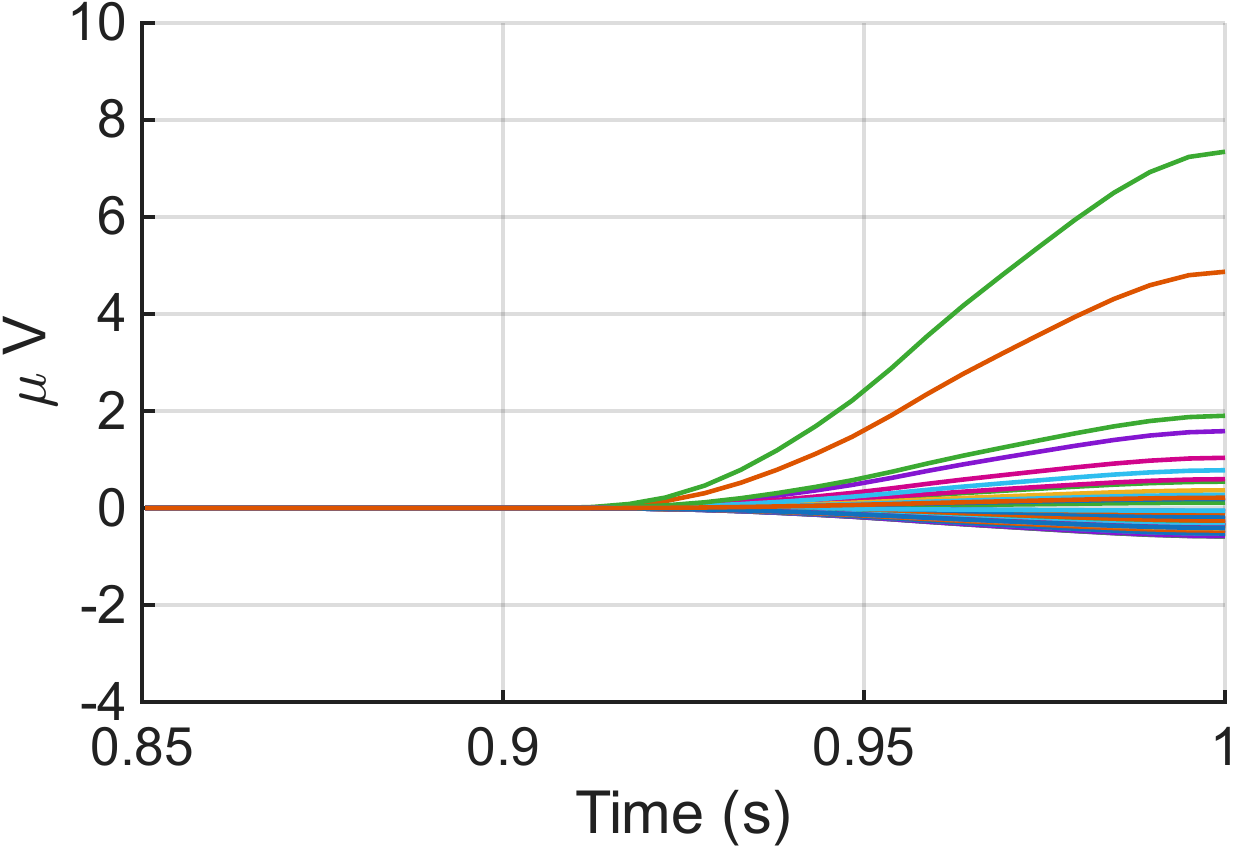}
    \end{minipage}
    \caption{Butterfly plots of the noiseless time series generated using anisotropic Whitney, Local subtraction, $\Hdiv$ basis functions, and patch source model.}
    \label{fig.butterfly}
\end{figure*}

Figure~\ref{fig.butterfly} presents butterfly plots of simulated, noiseless electrode potentials generated by a single cortical dipole source. Forward interpolation was performed using the Whitney, local subtraction, $\Hdiv$, and the patch model source configurations.

The butterfly plots in Figure~\ref{fig.butterfly} indicate that the Whitney and local subtraction source models yield highly similar electrode potential distributions. In contrast, the patch model produces potentials that are more broadly distributed among the electrodes, with a reduced maximum amplitude. The $\Hdiv$ model, on the other hand, shows pronounced sensitivity across all sensors, resulting in an even greater spatial spread of potentials than the patch model. The potentials are significantly more spread out than with the Patch model, where the notable correlation is only among sources no more than one centimeter apart.

\subsection{Source estimation at various depths from noiseless data}\label{sec.dbsp.noiseless}

In this experiment, we examined how well sLORETA, SHA\-L1R, and SKF can recover sources with correct depths measured from the inner skull surface, and how the spread of the estimates varies as a function of the true source depth. Whitney basis functions, Local subtraction, $\Hdiv$, and the patch model are used as the source models. Furthermore, we investigated how the patch sources can be estimated using Local subtraction. Theoretically, we can expect sLORETA and SHAL1R to have depth-unbiased estimates, which means that the regression lines computed from the graphs showing the estimated source depth as a function of true source depth should be \emph{identity lines} that run diagonally through the origin at a \num{45}\textdegree, yielding a slope of 1.

\begin{figure*}
    \centering
    \begin{minipage}{0.05\linewidth}
        \rotatebox{90}{Whitney}
    \end{minipage}\begin{minipage}{0.3\linewidth}
    \centering
    sLORETA \vspace{0.1cm}
    
        \includegraphics[width=0.98\linewidth]{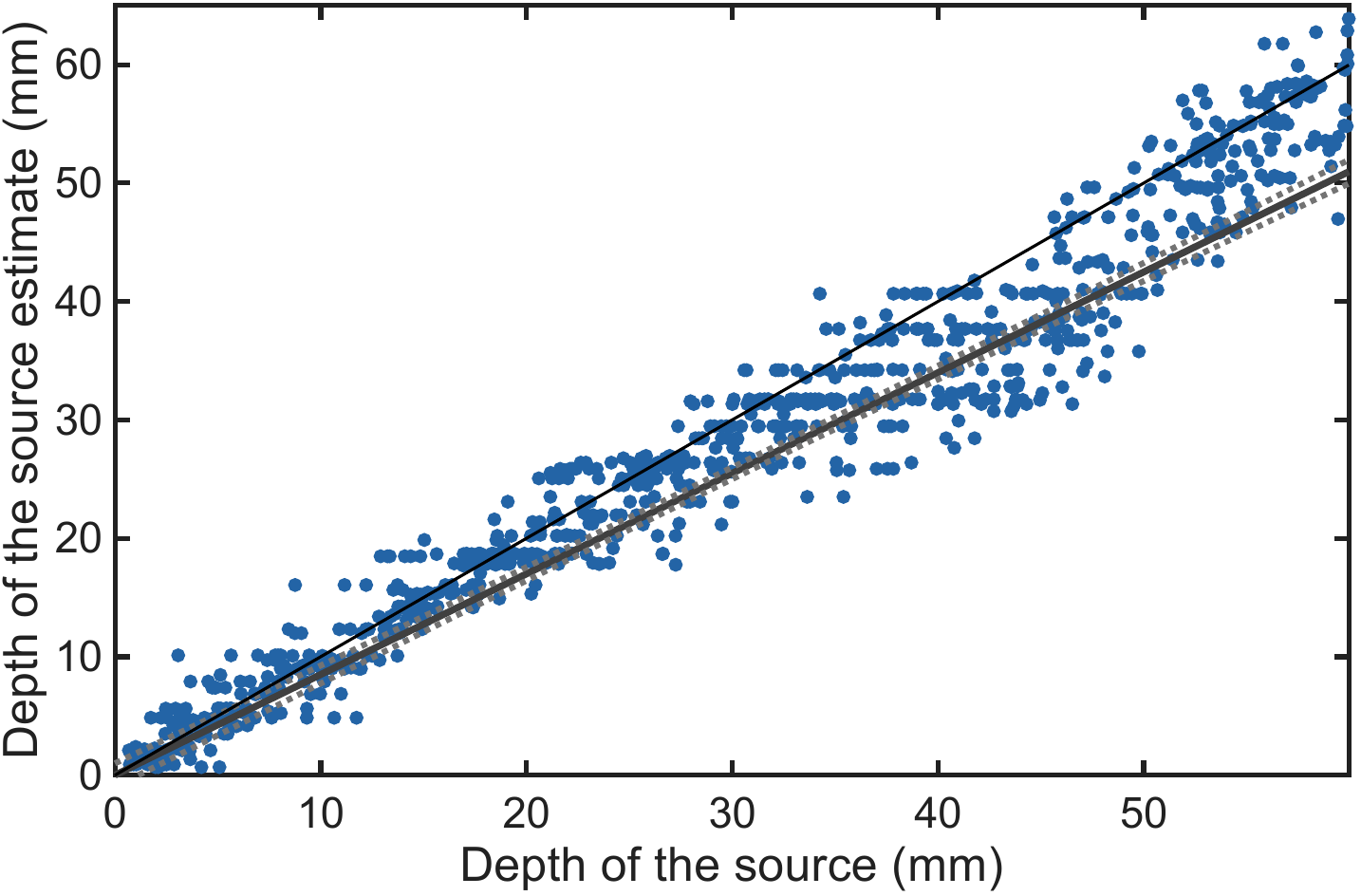}
    \end{minipage}\begin{minipage}{0.3\linewidth}
    \centering
    SHAL1R \vspace{0.1cm}
    
        \includegraphics[width=0.98\linewidth]{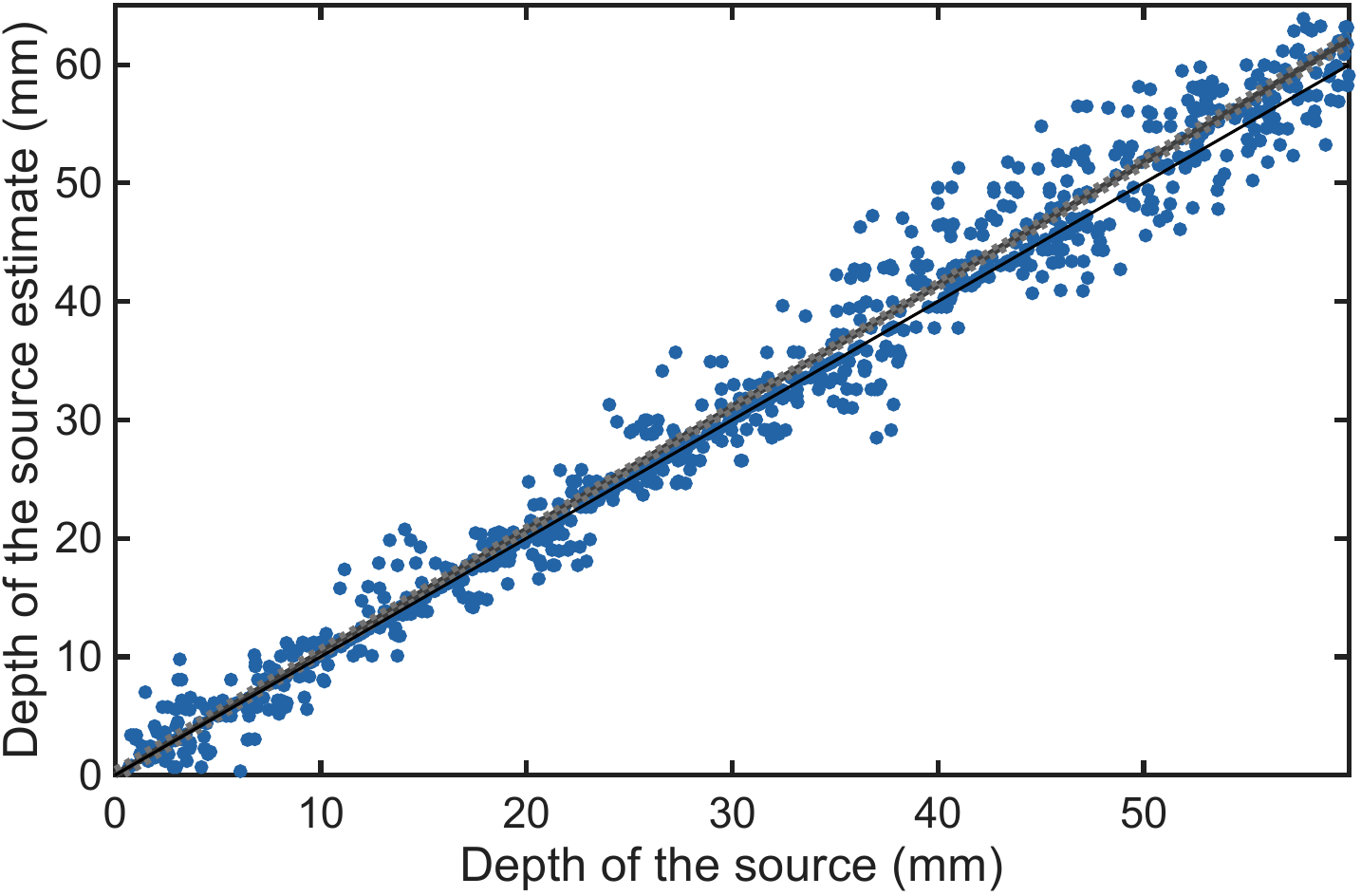}
    \end{minipage}\begin{minipage}{0.3\linewidth}
    \centering
    SKF \vspace{0.1cm}
    
        \includegraphics[width=0.98\linewidth]{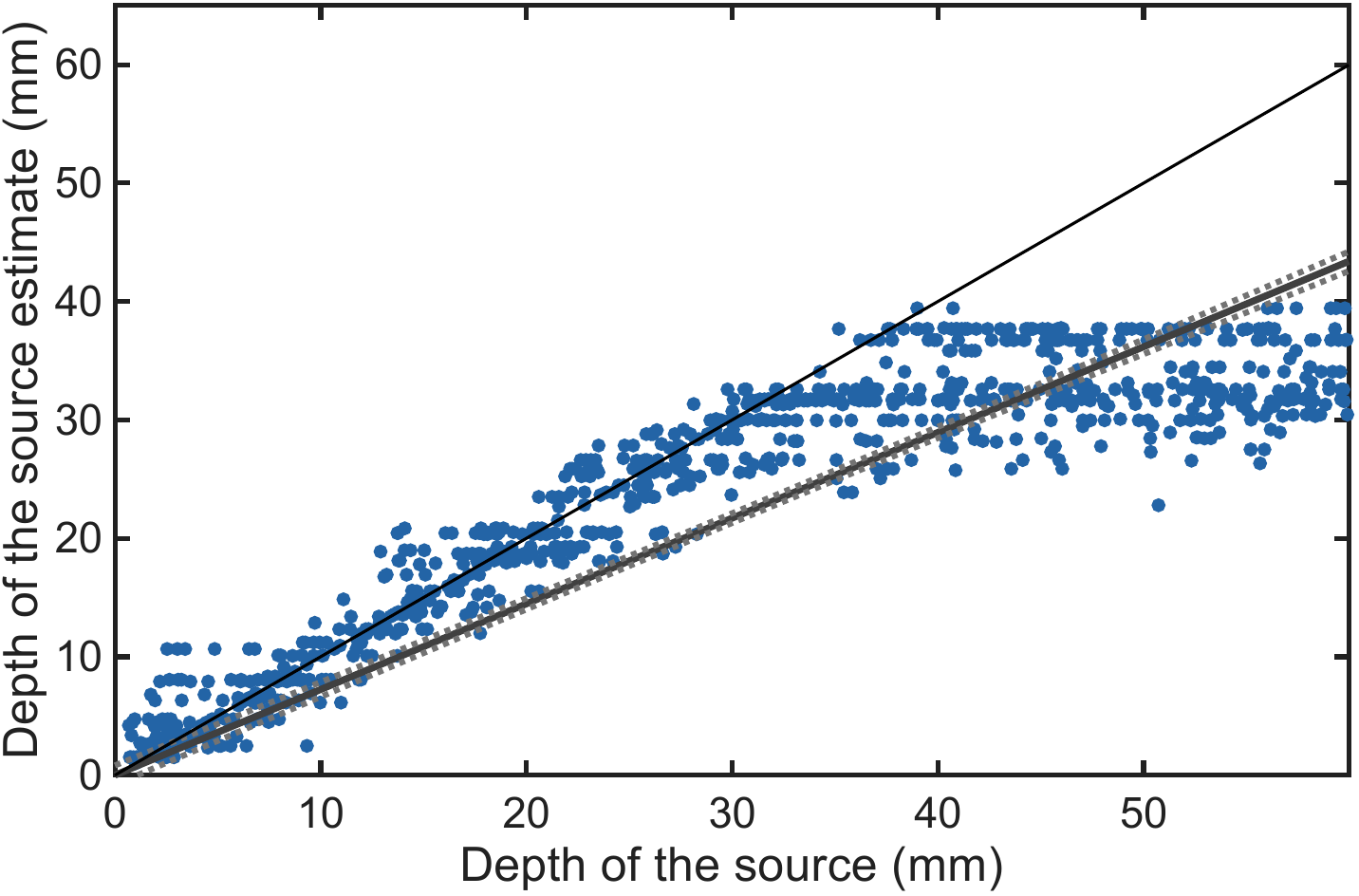}
    \end{minipage}

    \begin{minipage}{0.05\linewidth}
        \rotatebox{90}{$\Hdiv$}
    \end{minipage}\begin{minipage}{0.3\linewidth}
    \centering
        \includegraphics[width=0.98\linewidth]{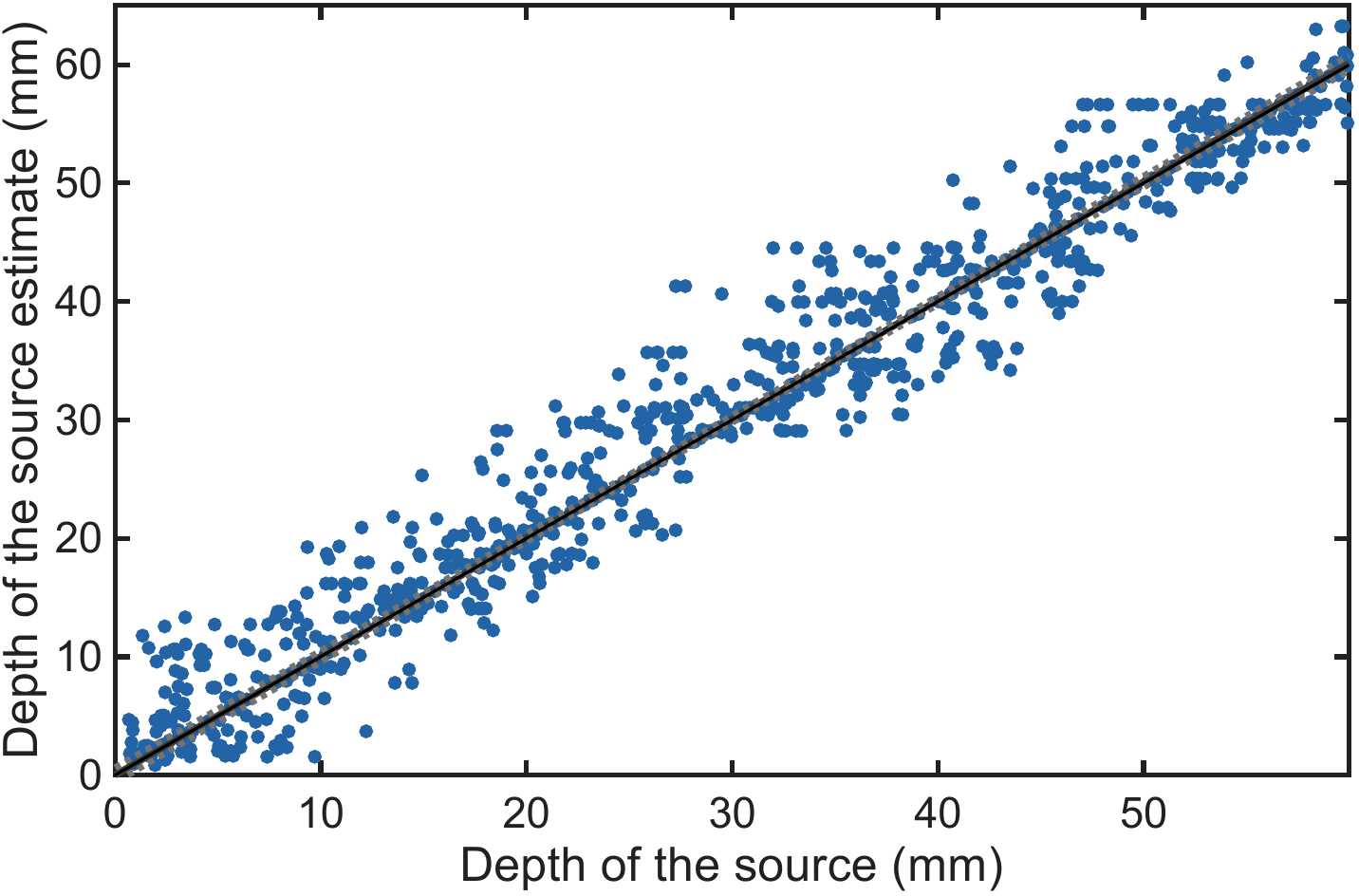}
    \end{minipage}\begin{minipage}{0.3\linewidth}
    \centering
        \includegraphics[width=0.98\linewidth]{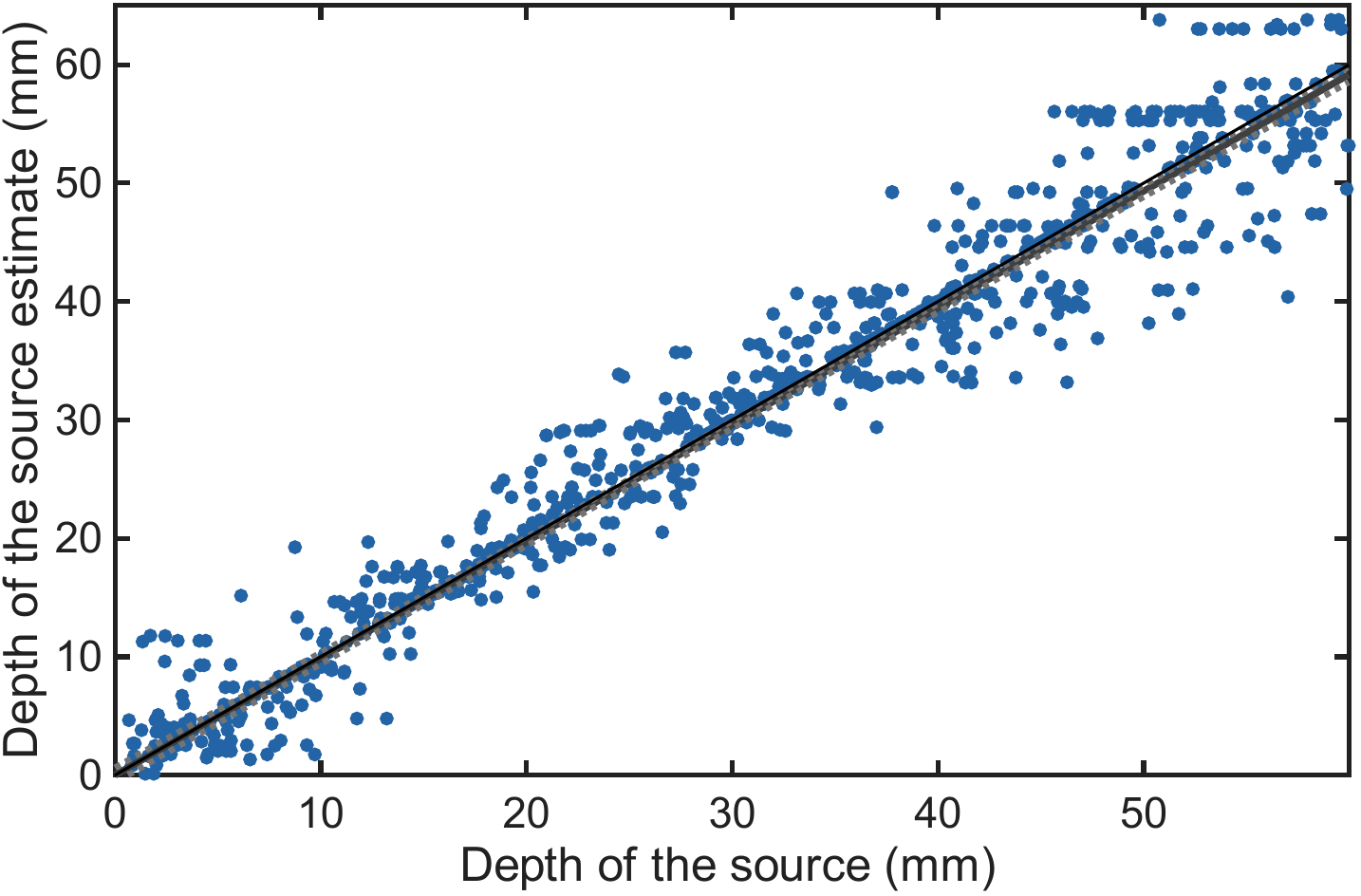}
    \end{minipage}\begin{minipage}{0.3\linewidth}
    \centering
        \includegraphics[width=0.98\linewidth]{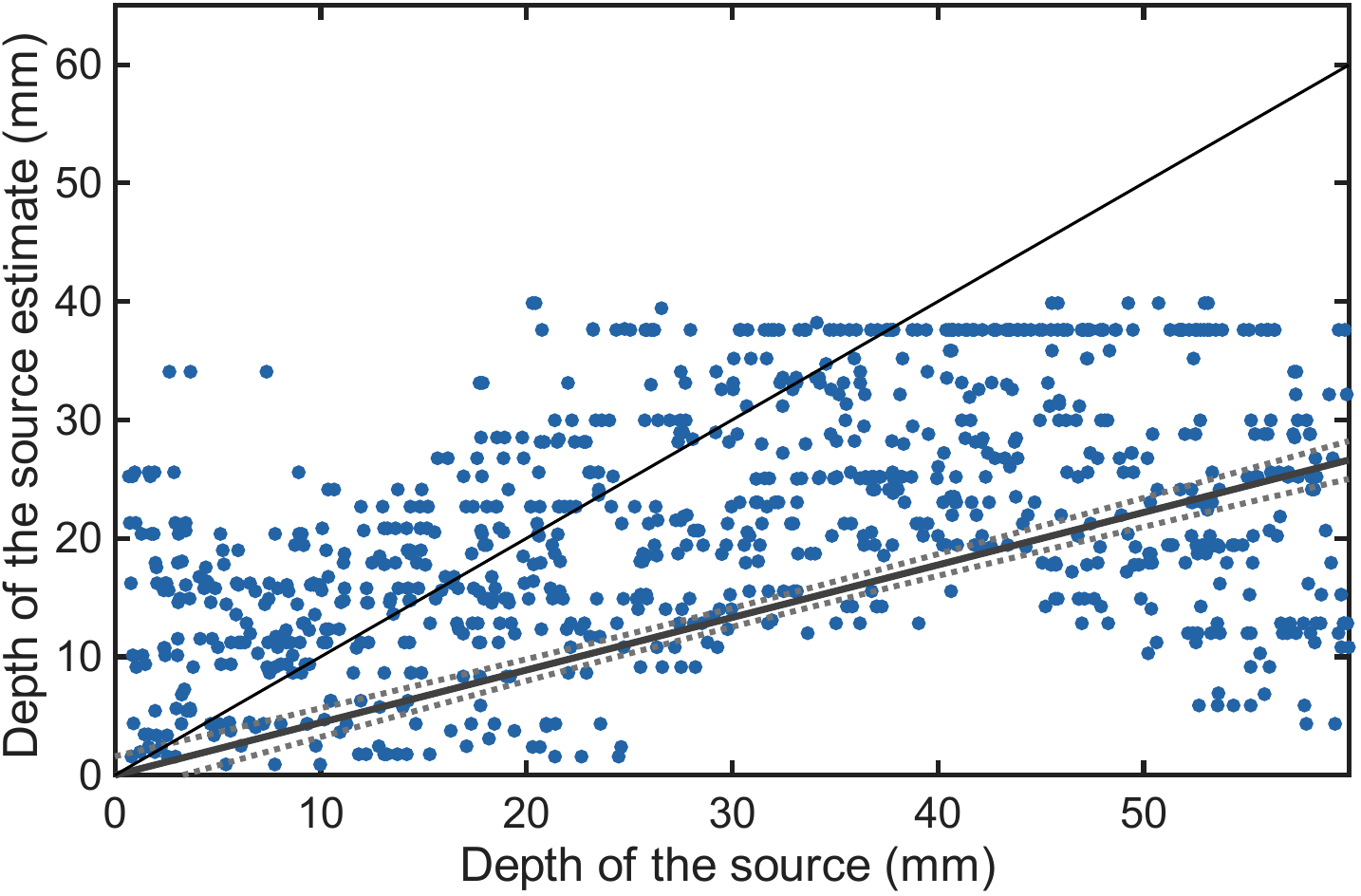}
    \end{minipage}

    \begin{minipage}{0.05\linewidth}
        \rotatebox{90}{Local subt.}
    \end{minipage}\begin{minipage}{0.3\linewidth}
    \centering
        \includegraphics[width=0.98\linewidth]{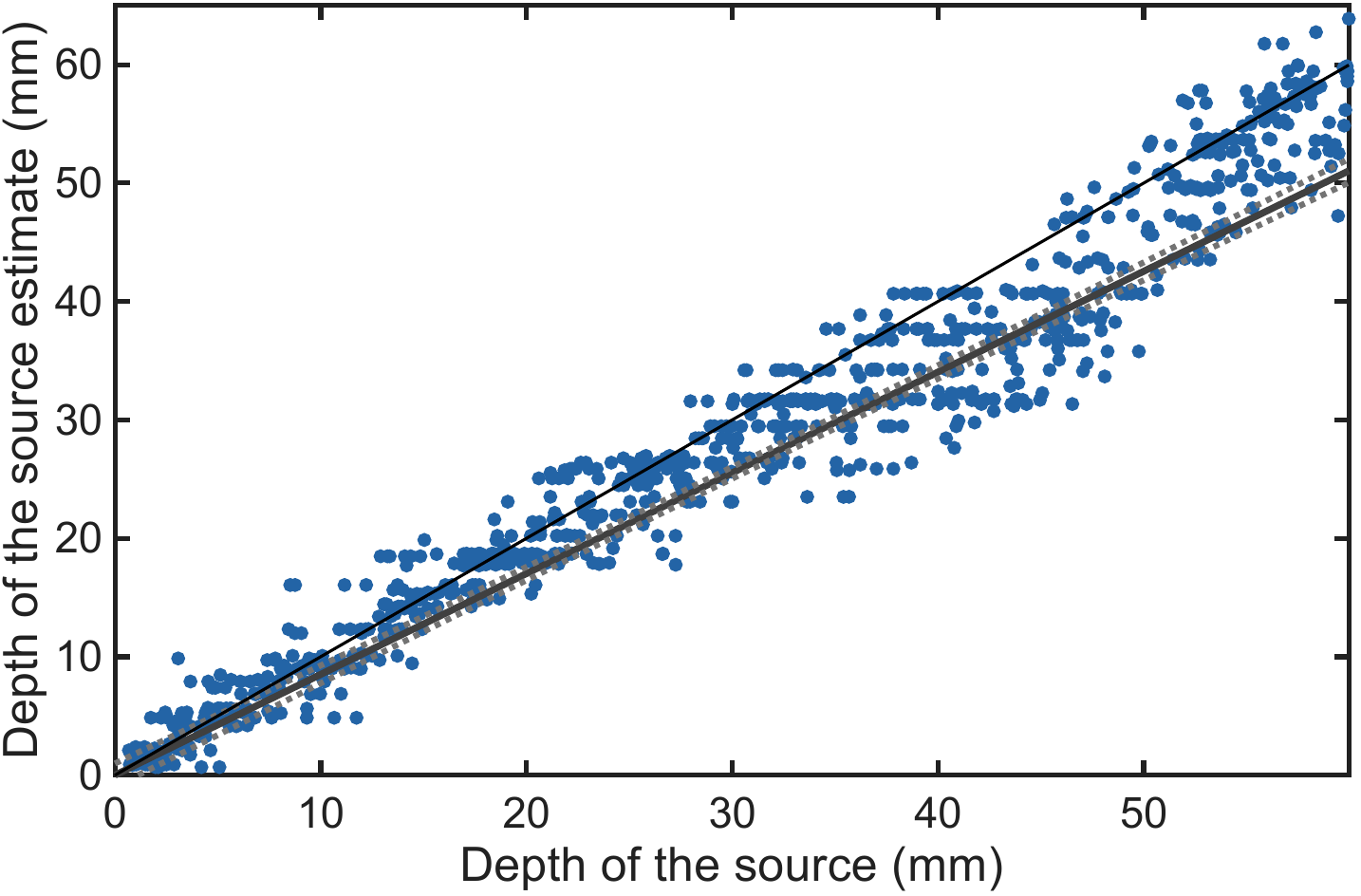}
    \end{minipage}\begin{minipage}{0.3\linewidth}
    \centering
        \includegraphics[width=0.98\linewidth]{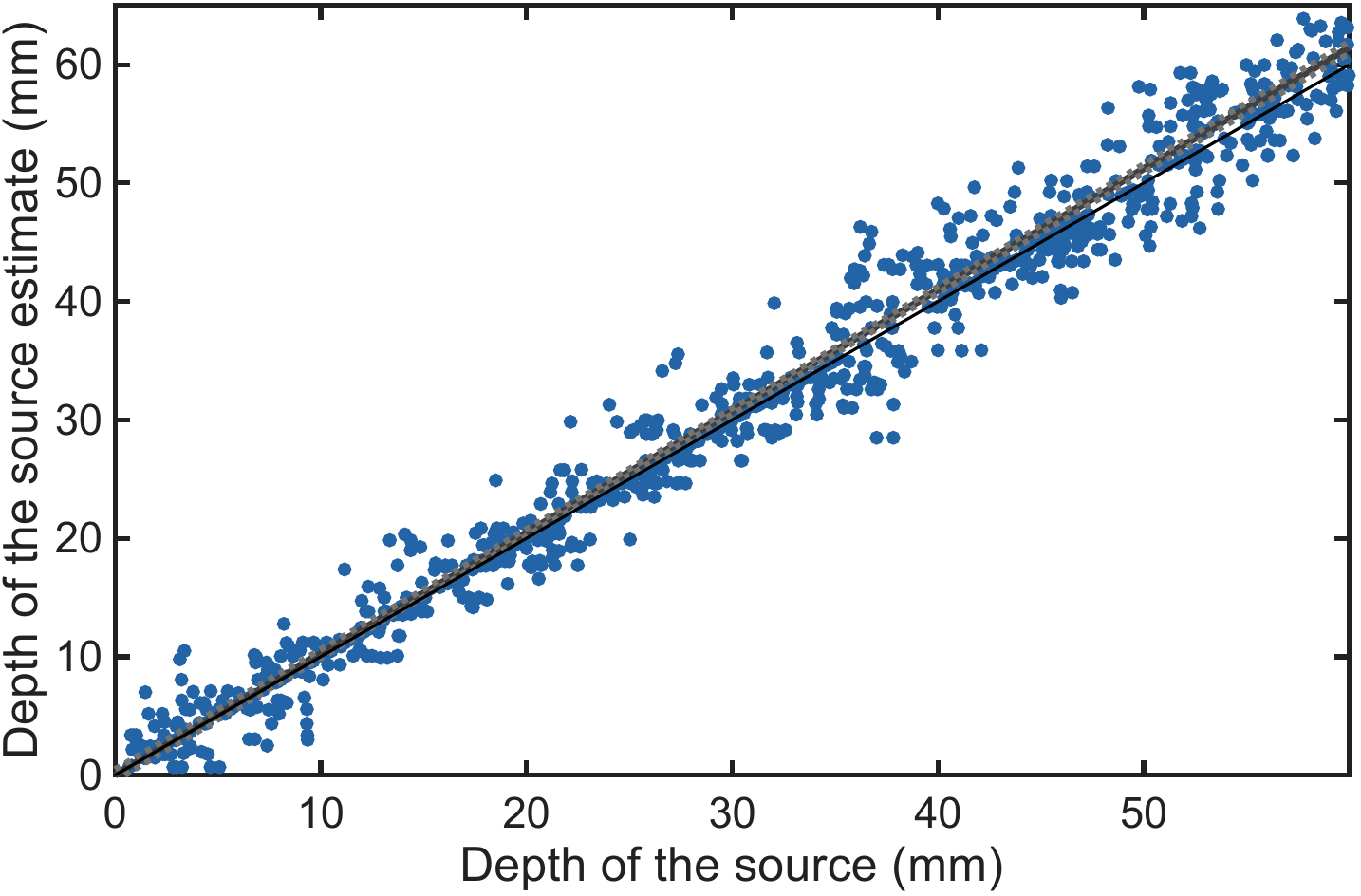}
    \end{minipage}\begin{minipage}{0.3\linewidth}
    \centering
        \includegraphics[width=0.98\linewidth]{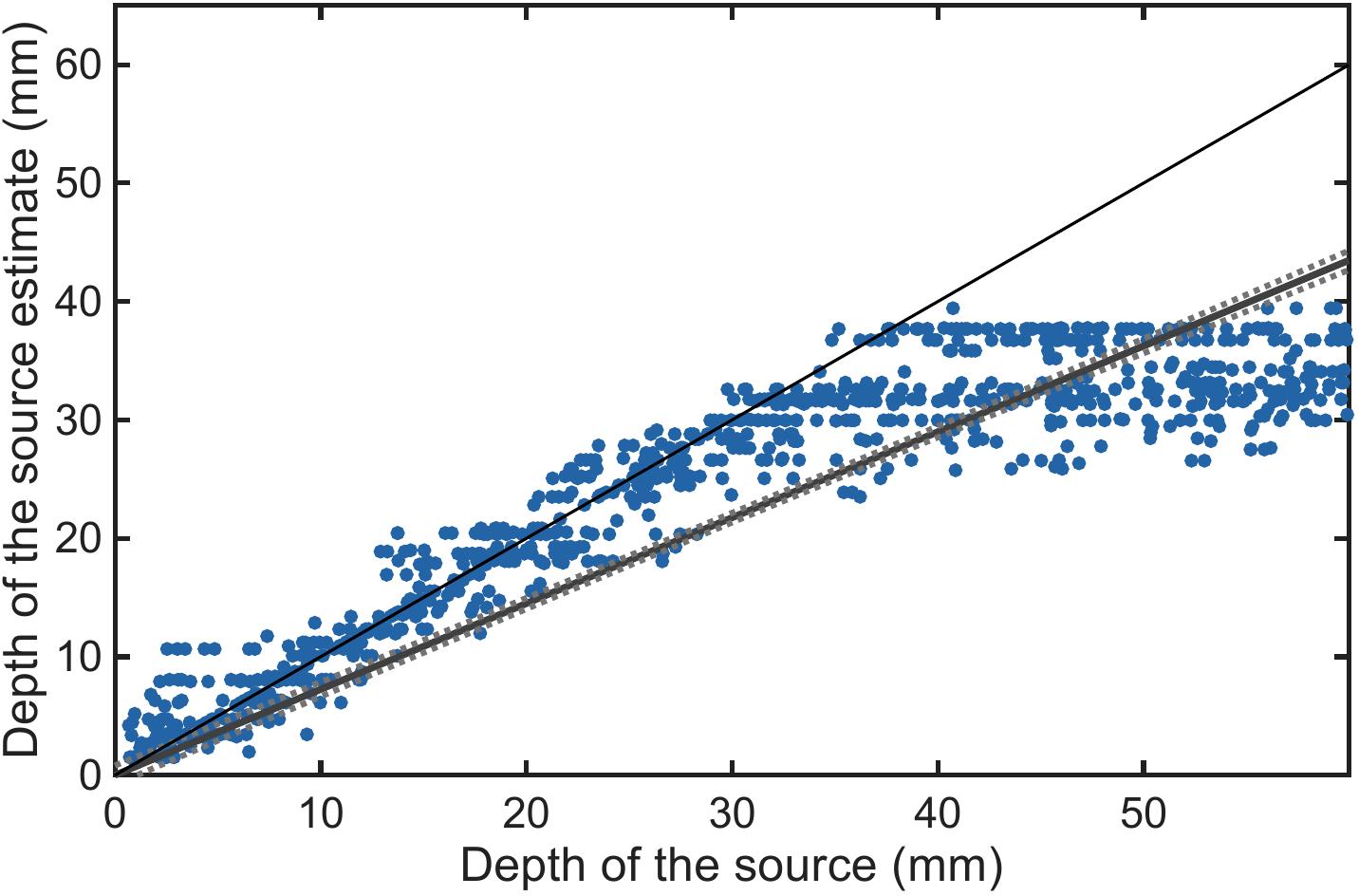}
    \end{minipage}

    \begin{minipage}{0.05\linewidth}
        \rotatebox{90}{Patch}
    \end{minipage}\begin{minipage}{0.3\linewidth}
    \centering
        \includegraphics[width=0.98\linewidth]{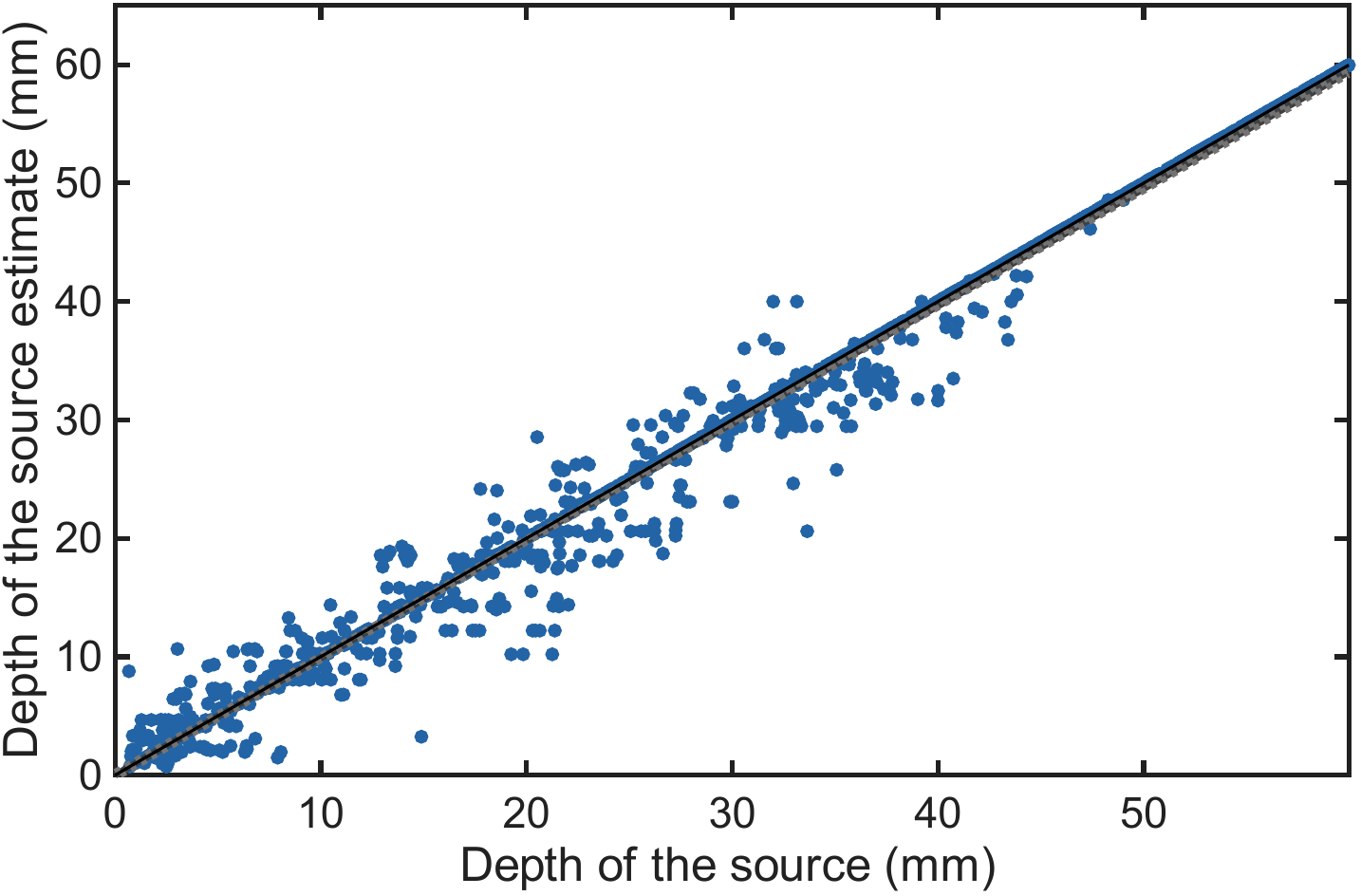}
    \end{minipage}\begin{minipage}{0.3\linewidth}
    \centering
        \includegraphics[width=0.98\linewidth]{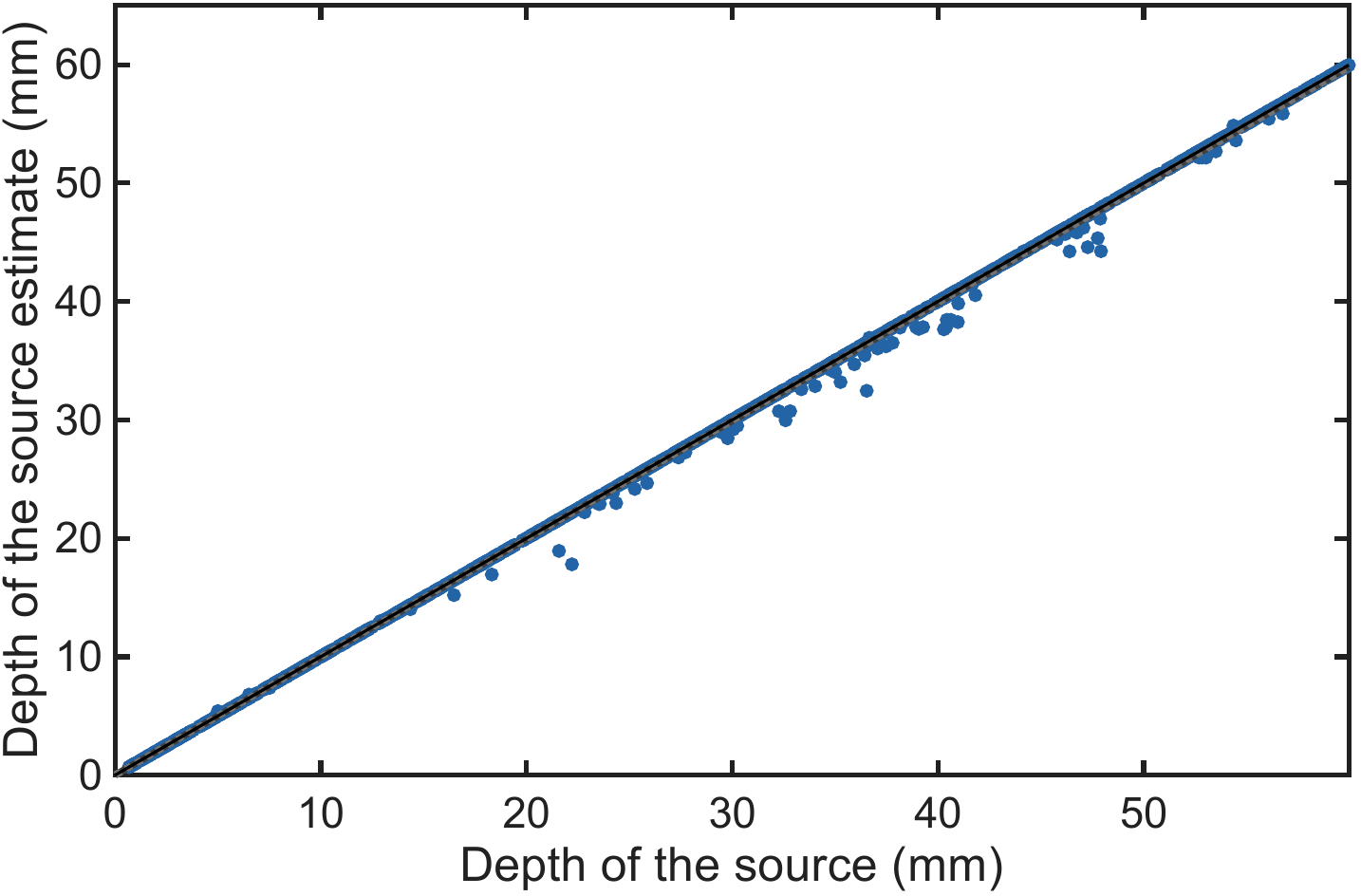}
    \end{minipage}\begin{minipage}{0.3\linewidth}
    \centering
        \includegraphics[width=0.98\linewidth]{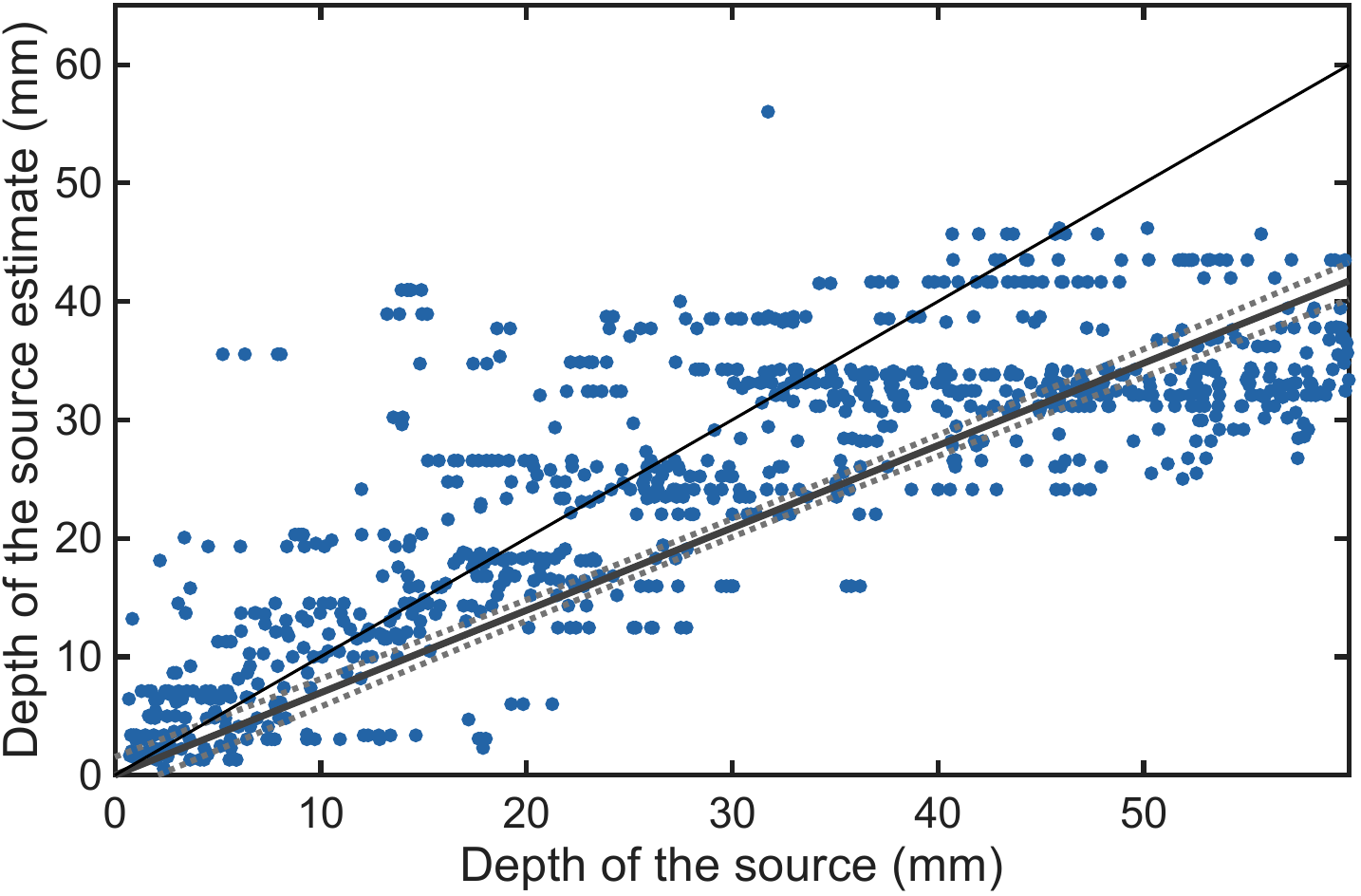}
    \end{minipage}

    \begin{minipage}{0.05\linewidth}
        \rotatebox{90}{Patch -- LS}
    \end{minipage}\begin{minipage}{0.3\linewidth}
    \centering
        \includegraphics[width=0.98\linewidth]{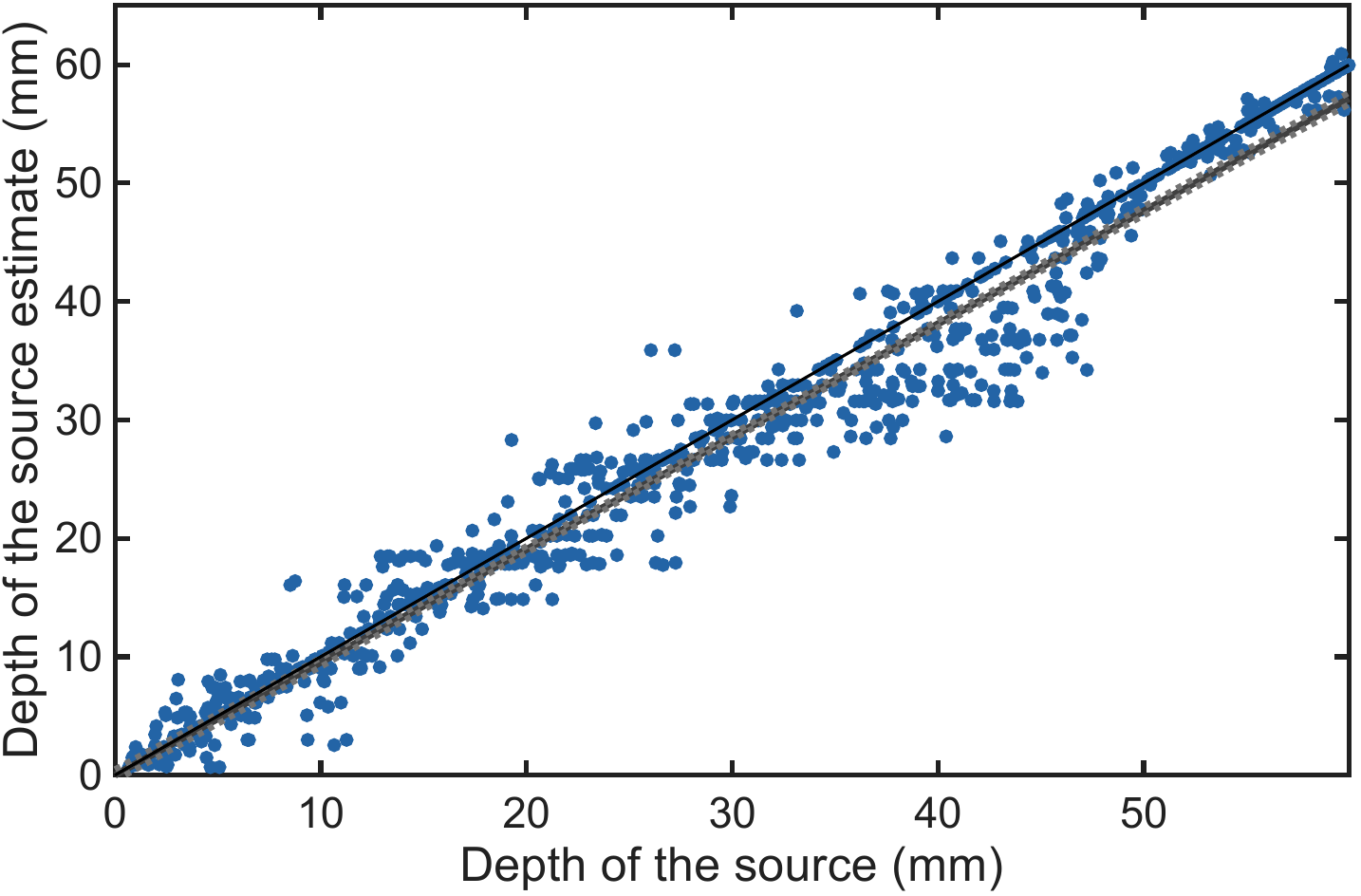}
    \end{minipage}\begin{minipage}{0.3\linewidth}
    \centering
        \includegraphics[width=0.98\linewidth]{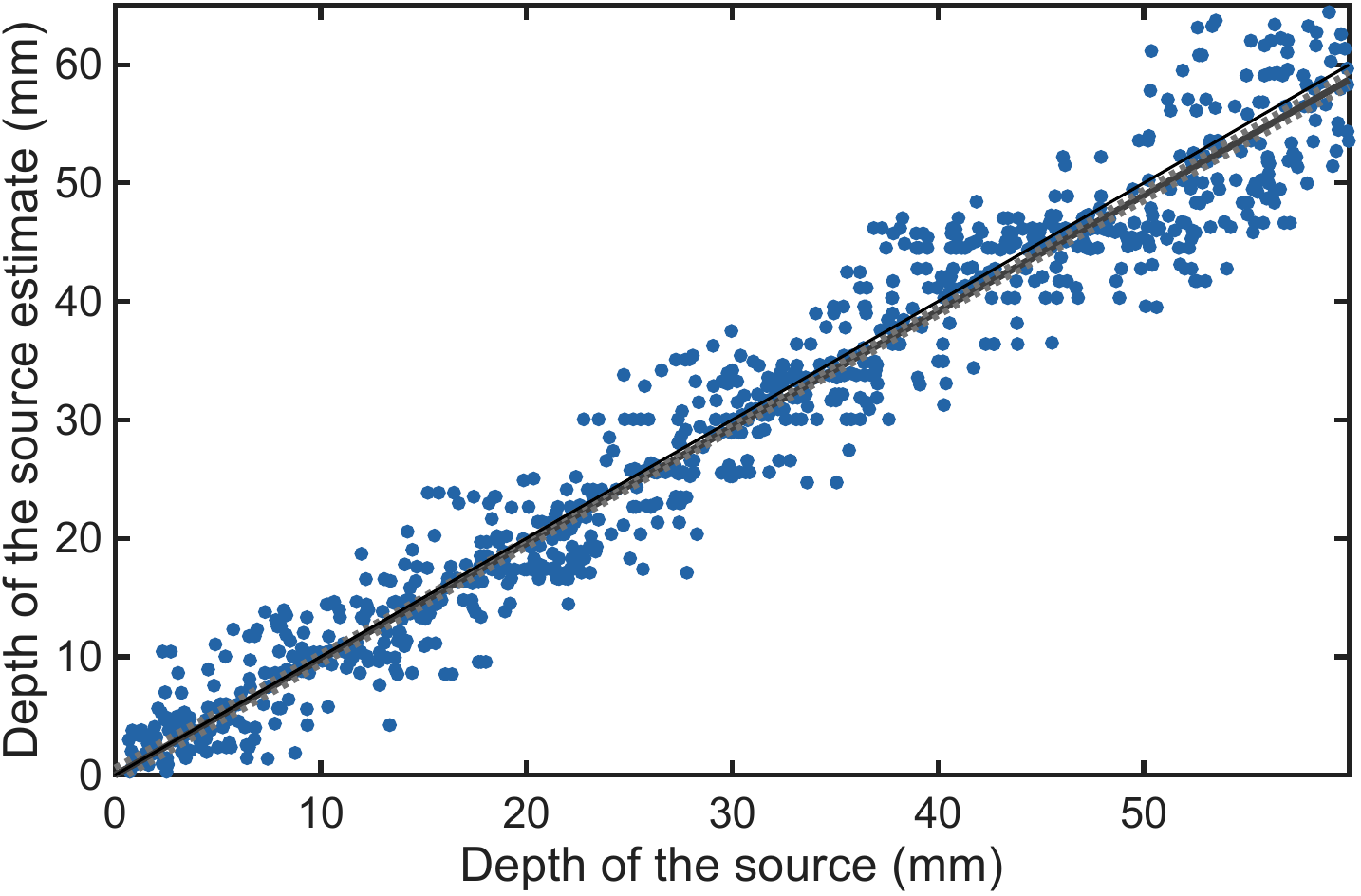}
    \end{minipage}\begin{minipage}{0.3\linewidth}
    \centering
        \includegraphics[width=0.98\linewidth]{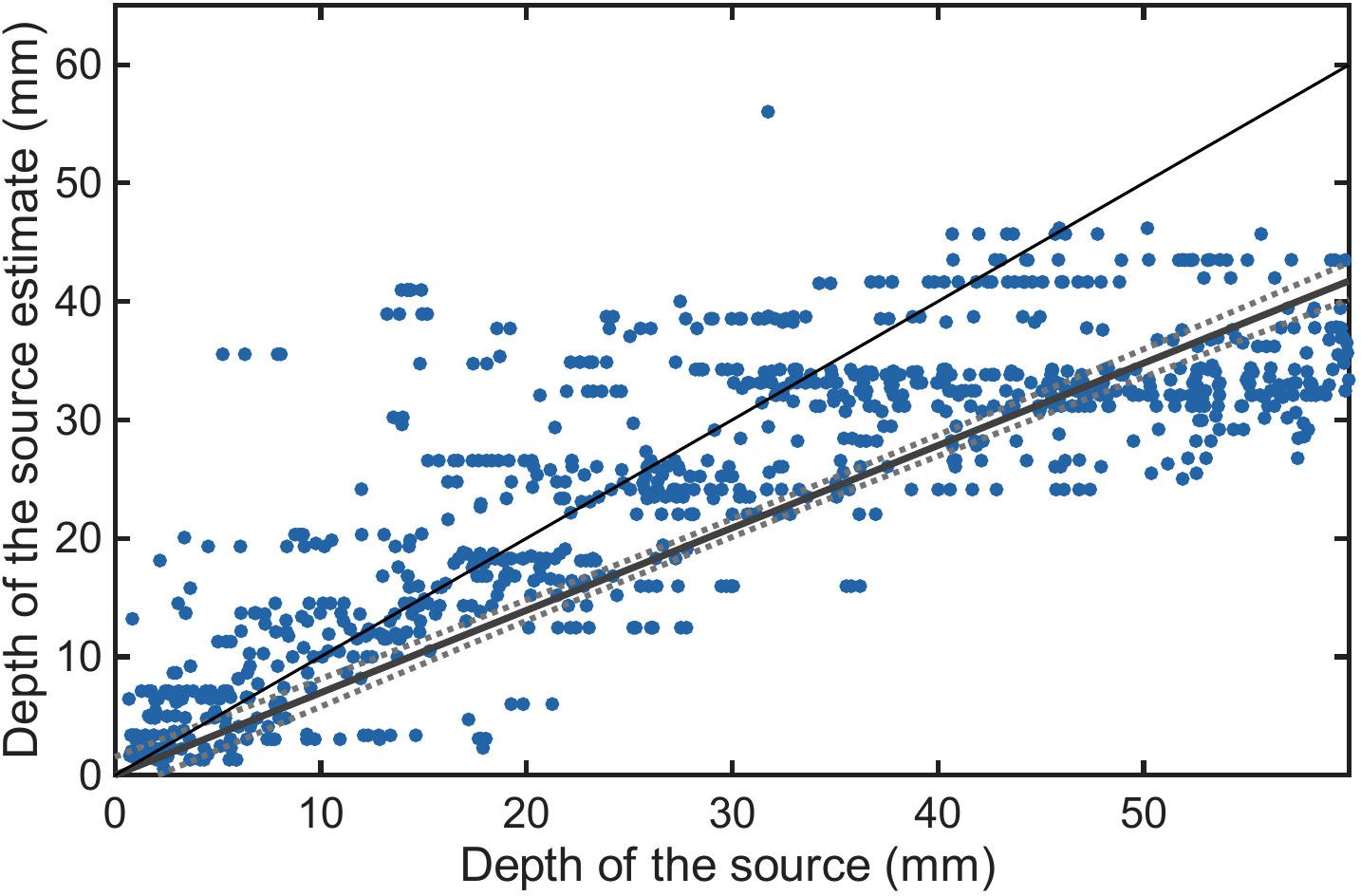}
    \end{minipage}
    
    \caption{Depth of true source plotted against the depth of estimation done by sLORETA and SHAL1R. The thin black line shows the optimal agreement between the true and estimated depth when the localization error is zero. The dark gray solid line displays the linear regression, and the dashed gray curves are the 95 \% confidence intervals.}
    \label{fig.duneuro.depthbias}
\end{figure*}

\begin{table*}[h!]
    \centering
    \begin{tabular}{c|c|c|c|c|c} \hline
    \multicolumn{6}{c}{Noiseless}\\ \hline
        \diagbox{Method}{Model} & Whitney & $\Hdiv$ & Local subt. & Patch & Patch--LS\\ \hline
            sLORETA & 0.85 & 1.00 & 0.85 & 0.99 & 0.95 \\
            SHAL1R & 1.03 & 0.99 & 1.02 & 1.00 & 0.98\\
            SKF & 0.72 & 0.44 & 0.72 & 0.70 & 0.70\\ \hline
            \multicolumn{6}{c}{\qty{15}{\decibel} signal-to-noise ratio}\\ \hline
                 & Whitney & $\Hdiv$ & Local subt. & Patch & Patch--LS\\ \hline
            sLORETA & 0.89 & 0.94 & 0.89 & 0.99 & 0.95 \\
            SHAL1R & 1.00 & 0.88 & 1.00 & 0.97 & 0.96 \\
            SKF & 0.70 & 0.58 & 0.70 & 0.71 & 0.69\\ \hline
    \end{tabular}
    \caption{Slopes of the regression lines plotted on the depth scatter plot in Figure \ref{fig.duneuro.depthbias}. Method is completely unbiased when the slope is exactly 1.}
    \label{tab.regressionslopes}
\end{table*}

By examining the depth scatter plots and the regression lines based on the scatter data in Figure~\ref{fig.duneuro.depthbias}, sLORETA admits to higher depth bias than SHAL1R with Whitney-basis and Local subtraction, because the slopes of both cases are 0.85 for sLORETA, while the SHAL1R's slopes of the regression lines are 1.03 for Whitney and 1.02 for Local subtraction, as presented in Table~\ref{tab.regressionslopes}. In both cases, Local subtraction, as a forward interpolation scheme, produces slightly better localization results than Whitney, which again shows a better agreement with the true source and depth of the estimate. A better overall agreement for these two inversion methods is obtained with the $\Hdiv$ source model, since the sLORETA slope is 1.00 and the SHAL1R slope is 0.99. SKF appears to behave differently from the two other methods, since a better slope is obtained with both Whitney and Local subtraction for value 0.72 compared to 0.44 obtained with $\Hdiv$. Based on the results, SKF appears to be depth unbiased up to \qty{30}{\milli\meter}.

The best agreement across all the methods is obtained with the patch-like model. The regression slopes for SHAL1R and sLORETA are 1.00 and 0.99, respectively. With SKF, the estimates are slightly more biased than those from the Whitney basis and the Local subtraction, with a slope of 0.70.

In the case where the true source is set to follow the patch model and the source is estimated using a local subtraction-based lead field, a worsening of depth agreement is observed with respect to the use of Local subtraction in both forward and inverse. is expected, and categorically, SHAL1R slope worsens to 0.98, and the deviation from the regression line increases, and SKF slope is 0.70. Contrarily, the depth bias of sLORETA reduces when inverting with Local subtraction compared to the case where the true source follows Local subtraction.

\begin{figure*}
    \centering
    \begin{minipage}{0.05\linewidth}
        \rotatebox{90}{Whitney}
    \end{minipage}\begin{minipage}{0.3\linewidth}
    \centering
    sLORETA \vspace{0.1cm}
    
        \includegraphics[width=0.98\linewidth]{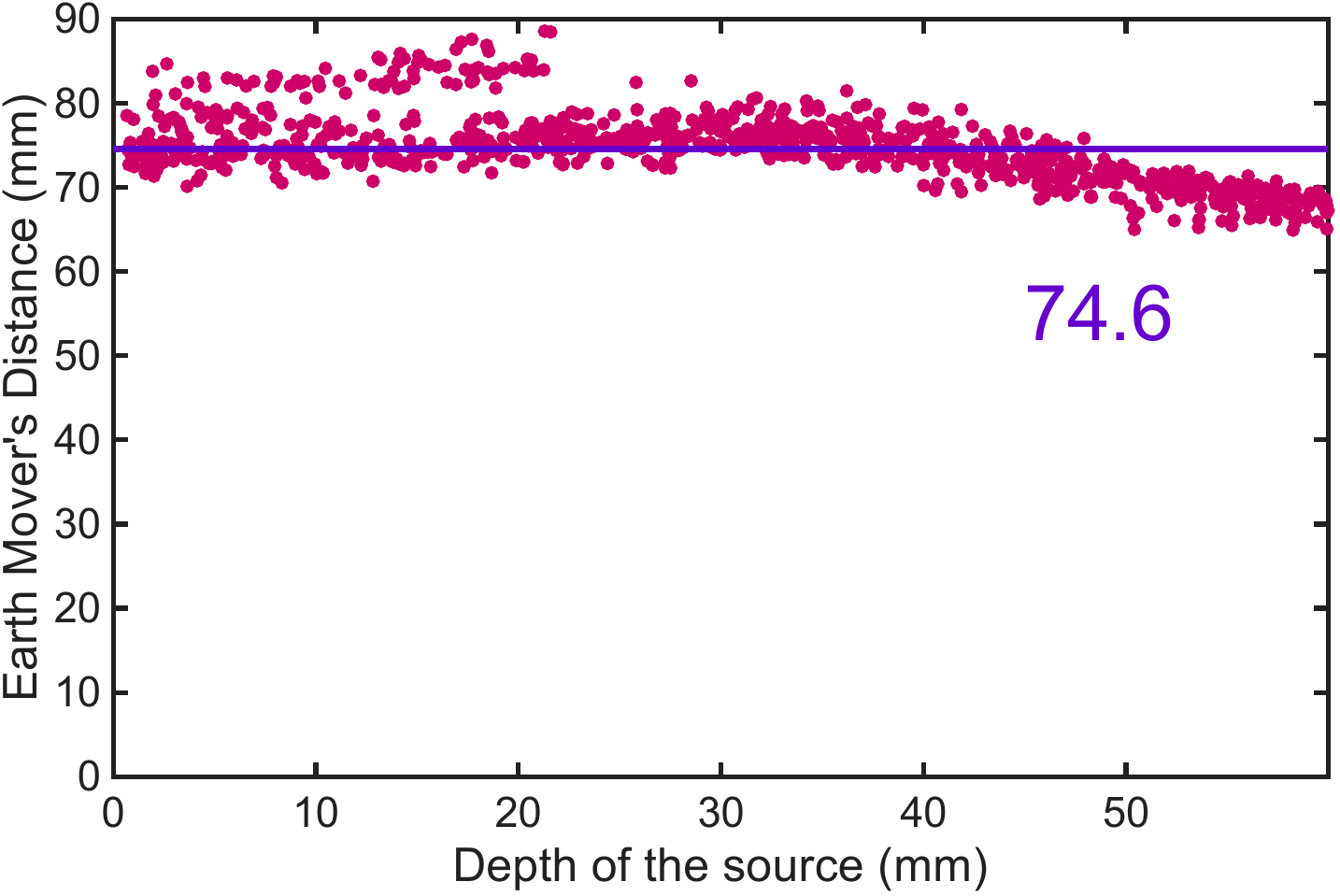}
    \end{minipage}\begin{minipage}{0.3\linewidth}
    \centering
    SHAL1R \vspace{0.1cm}
    
        \includegraphics[width=0.98\linewidth]{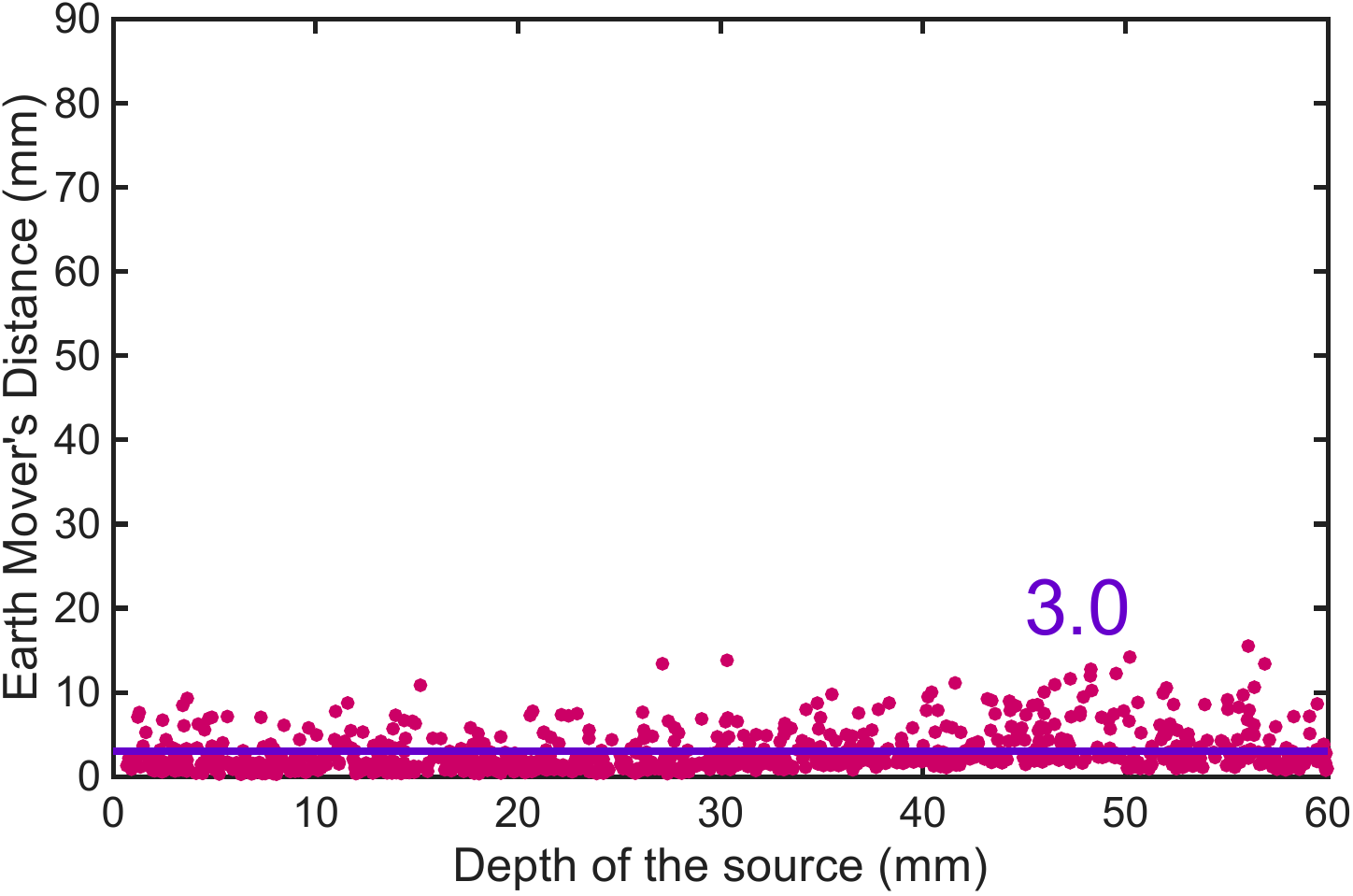}
    \end{minipage}\begin{minipage}{0.3\linewidth}
    \centering
    SKF \vspace{0.1cm}
    
        \includegraphics[width=0.98\linewidth]{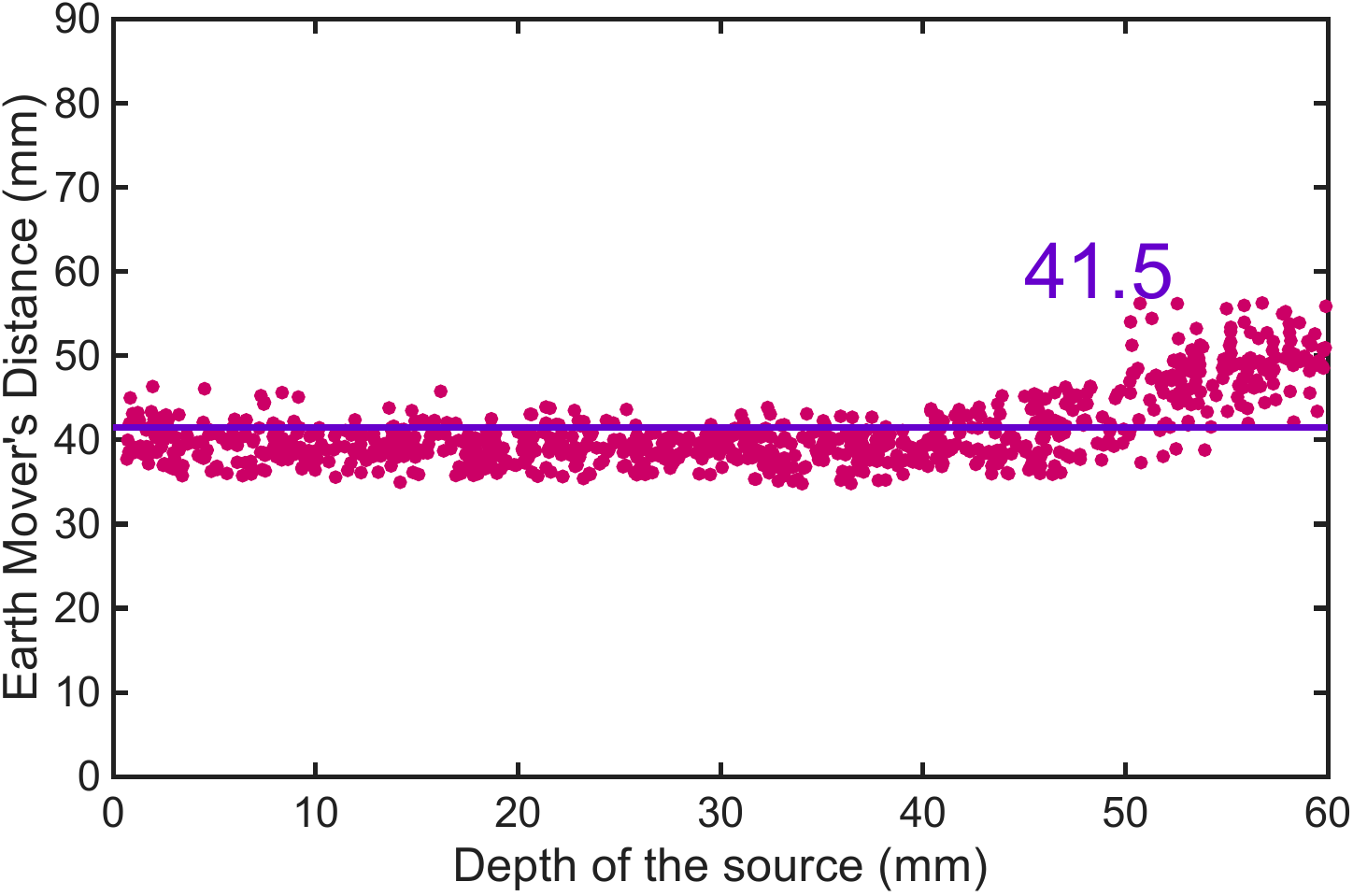}
    \end{minipage}

    \begin{minipage}{0.05\linewidth}
        \rotatebox{90}{$\Hdiv$}
    \end{minipage}\begin{minipage}{0.3\linewidth}
    \centering
        \includegraphics[width=0.98\linewidth]{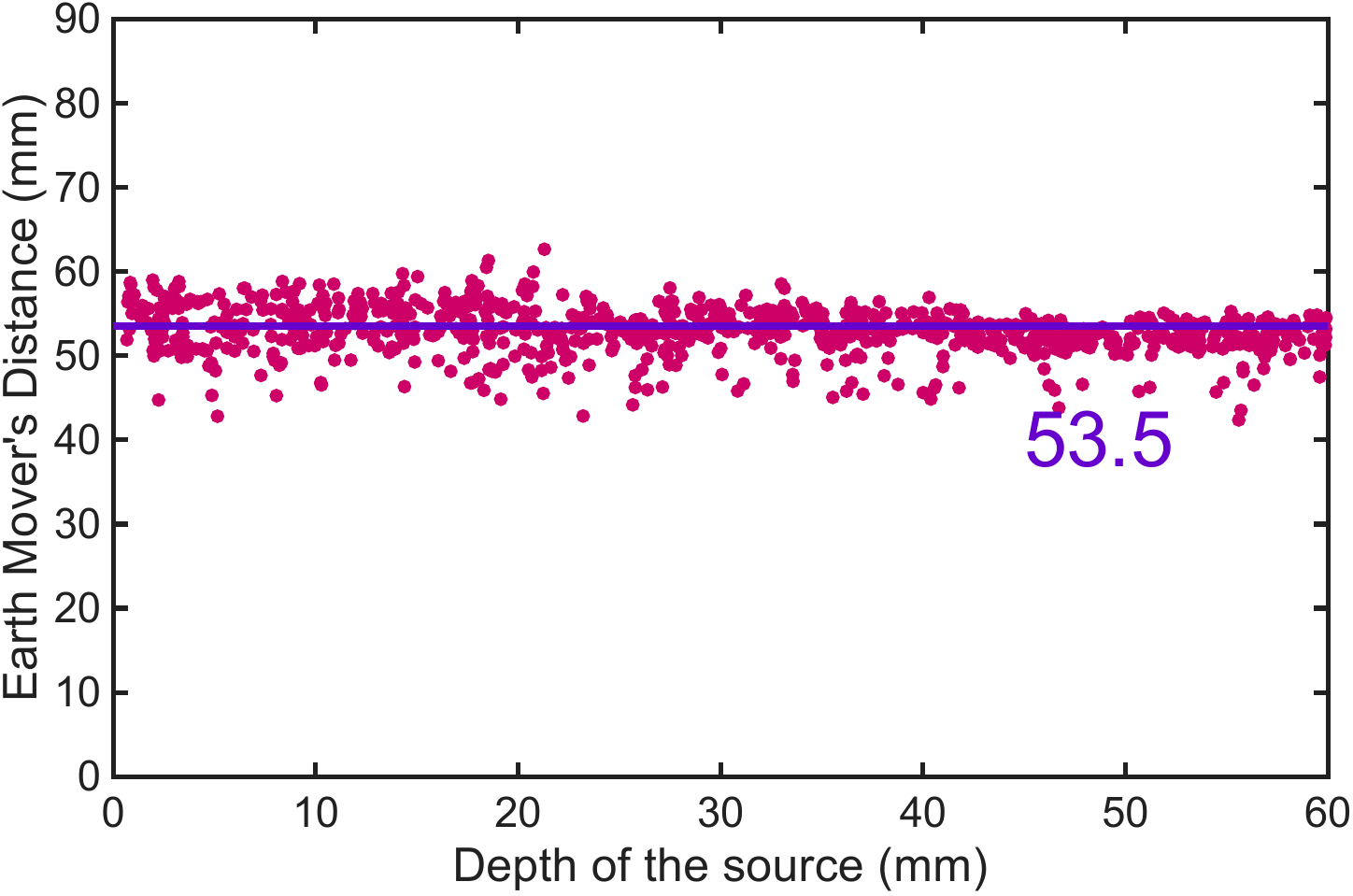}
    \end{minipage}\begin{minipage}{0.3\linewidth}
    \centering
        \includegraphics[width=0.98\linewidth]{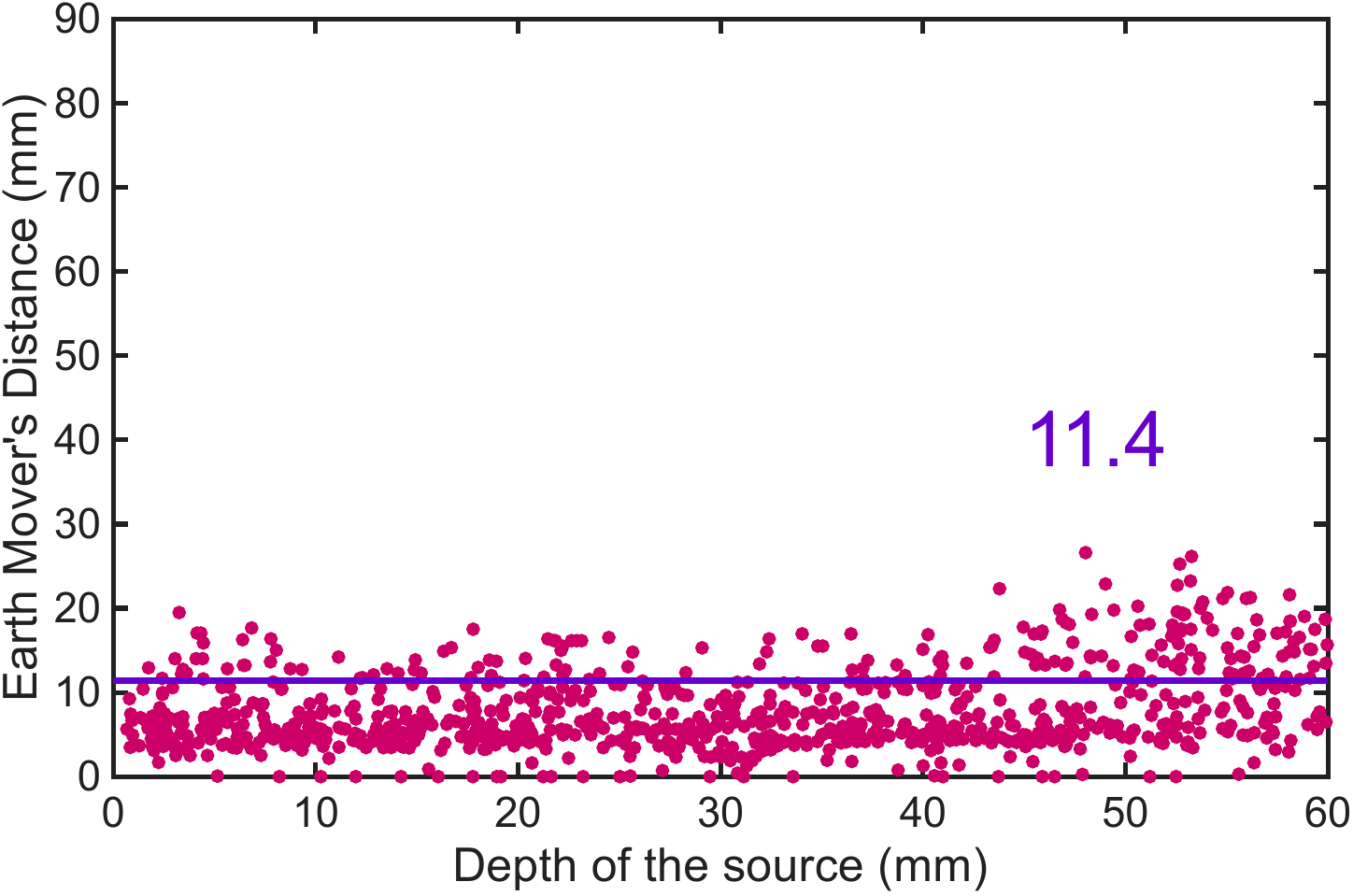}
    \end{minipage}\begin{minipage}{0.3\linewidth}
    \centering
        \includegraphics[width=0.98\linewidth]{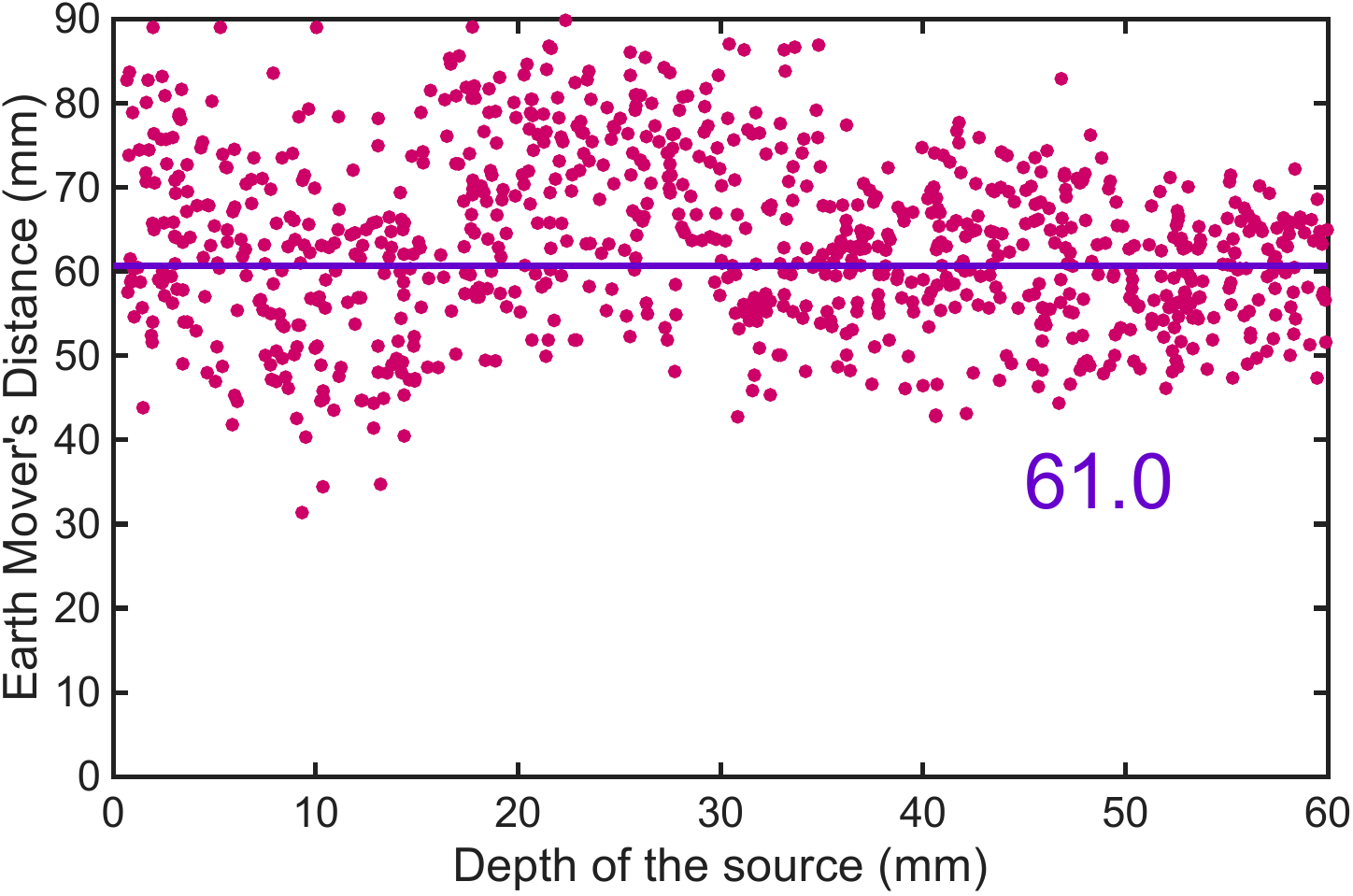}
    \end{minipage}

    \begin{minipage}{0.05\linewidth}
        \rotatebox{90}{Local subt.}
    \end{minipage}\begin{minipage}{0.3\linewidth}
    \centering
        \includegraphics[width=0.98\linewidth]{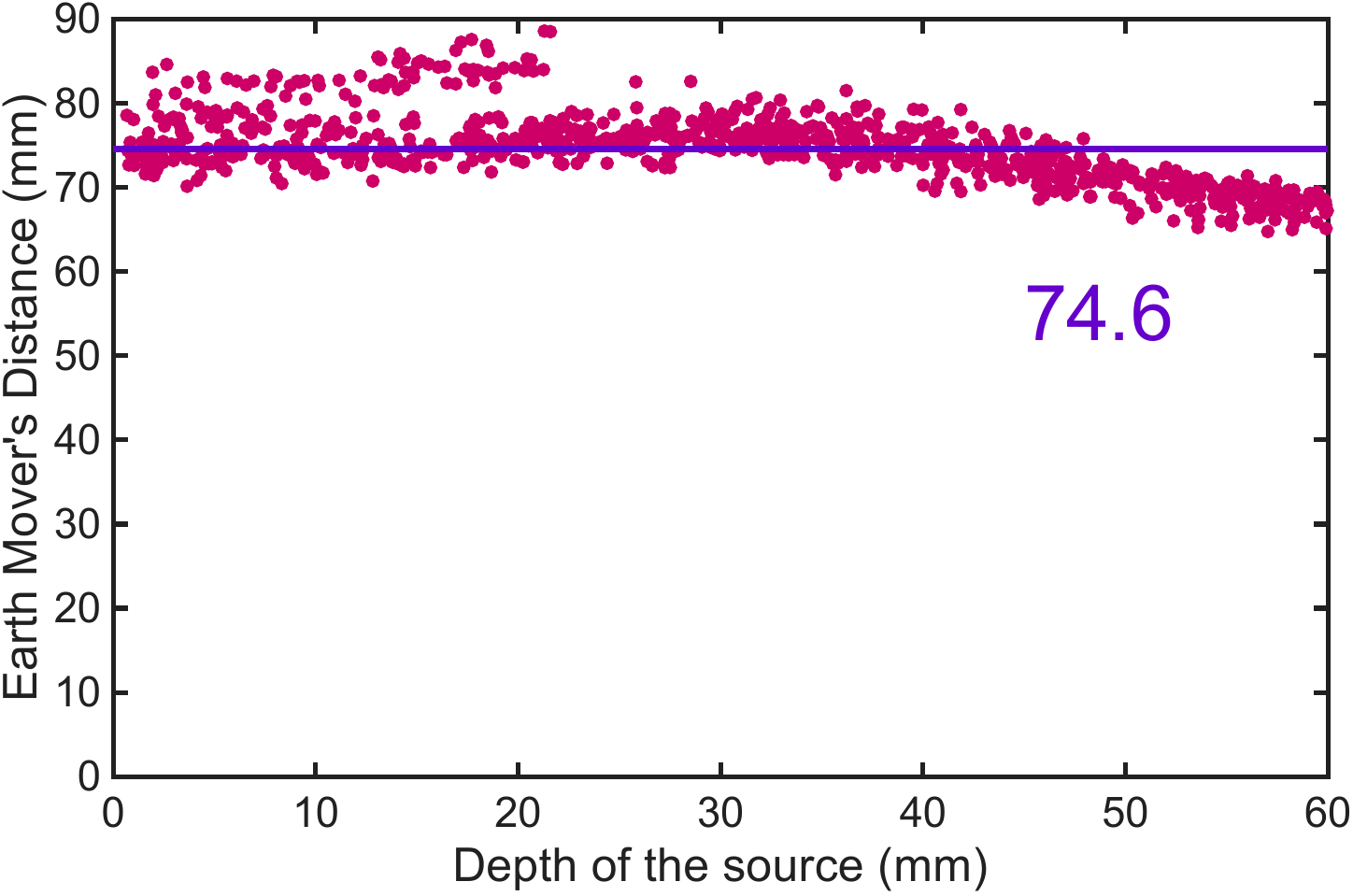}
    \end{minipage}\begin{minipage}{0.3\linewidth}
    \centering
        \includegraphics[width=0.98\linewidth]{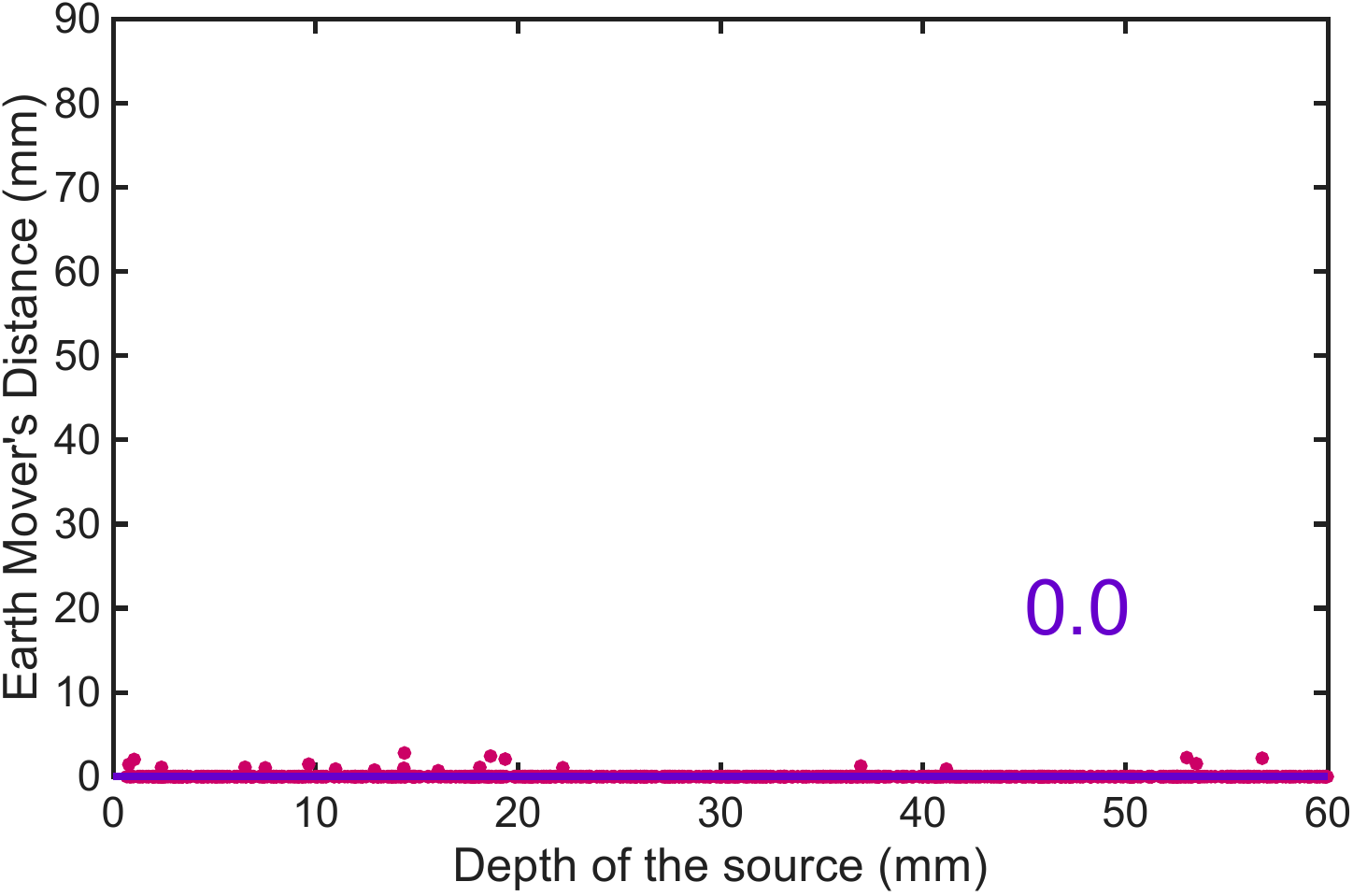}
    \end{minipage}\begin{minipage}{0.3\linewidth}
    \centering
        \includegraphics[width=0.98\linewidth]{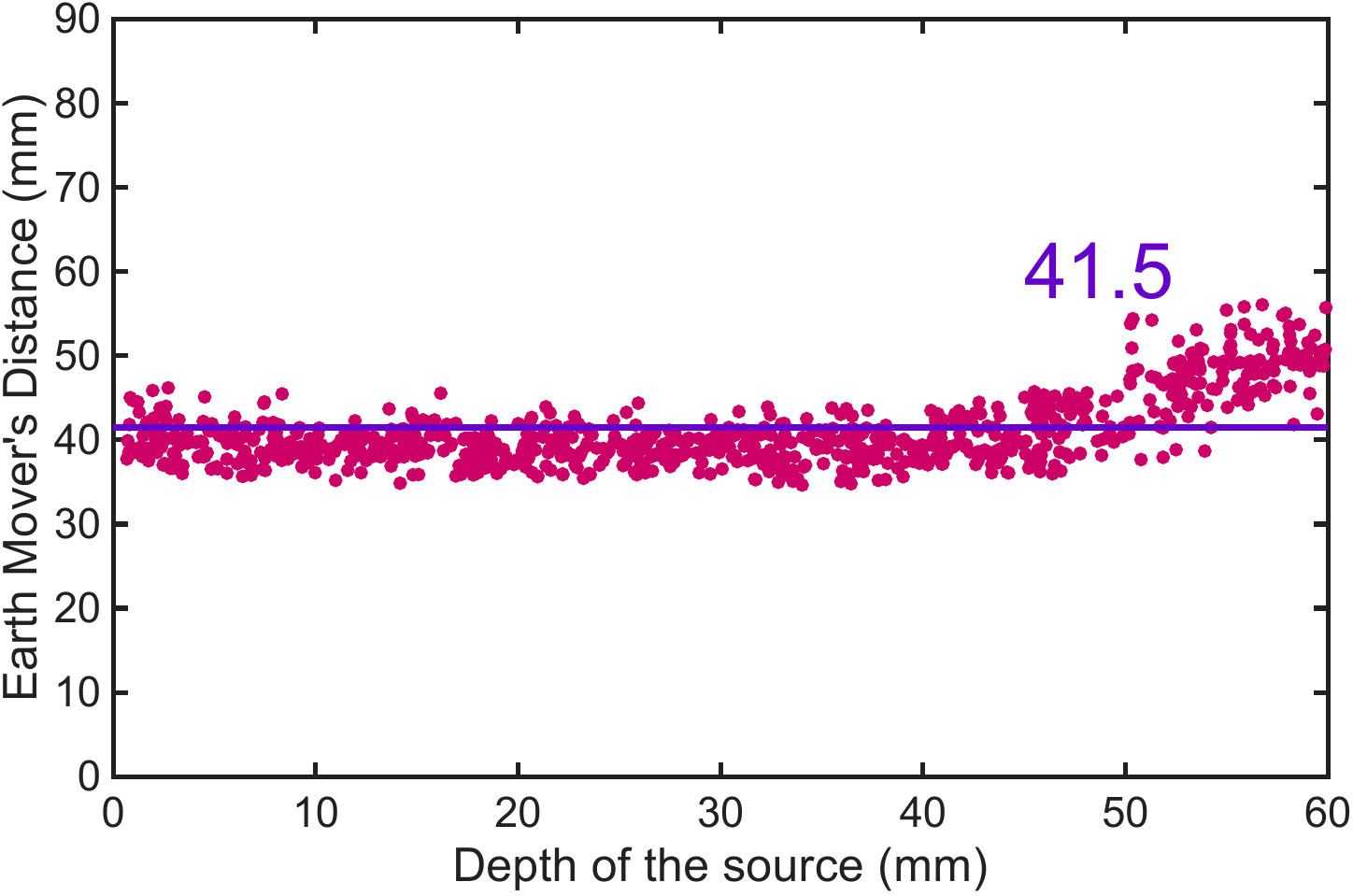}
    \end{minipage}

    \begin{minipage}{0.05\linewidth}
        \rotatebox{90}{Patch}
    \end{minipage}\begin{minipage}{0.3\linewidth}
    \centering
        \includegraphics[width=0.98\linewidth]{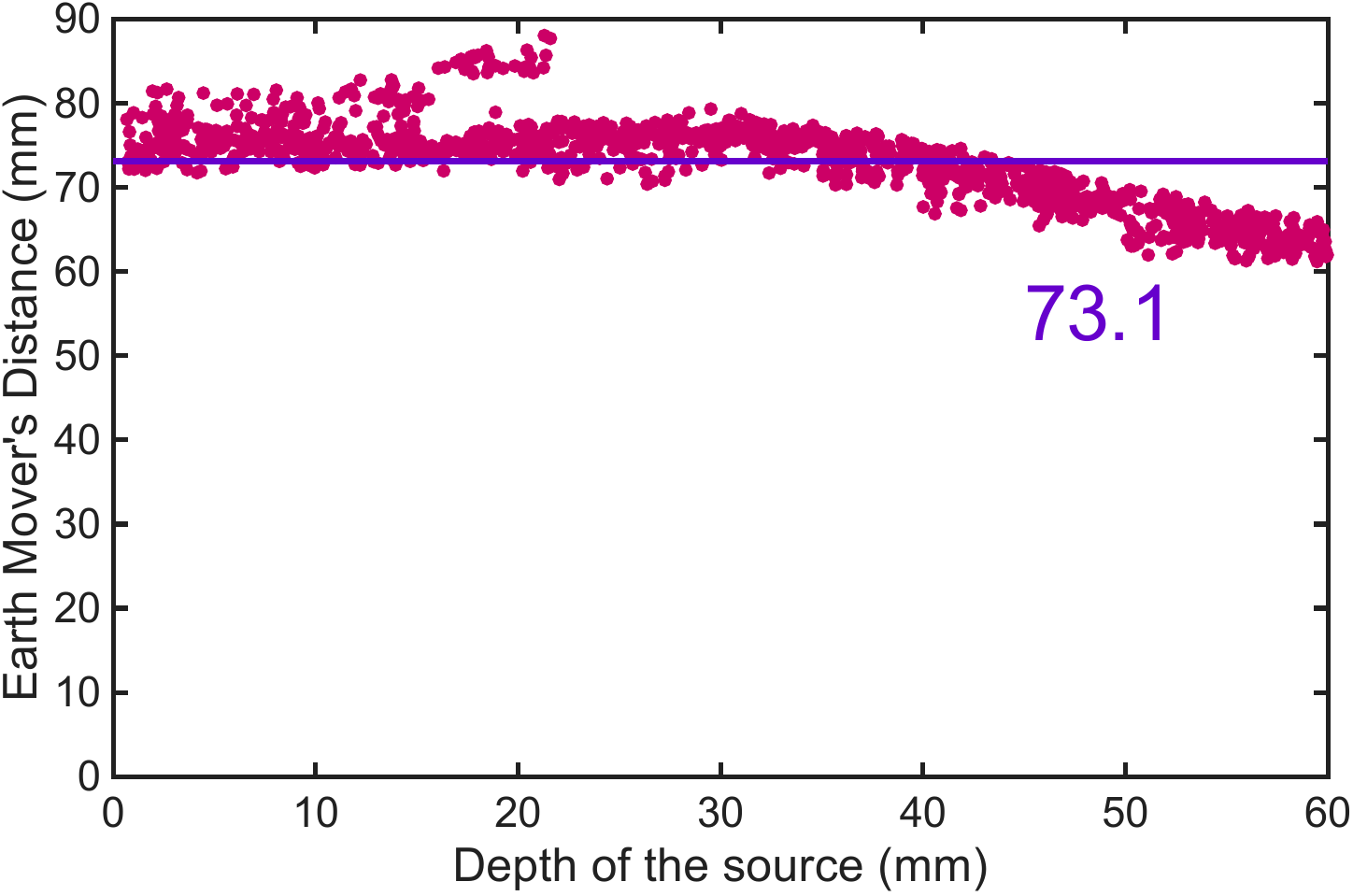}
    \end{minipage}\begin{minipage}{0.3\linewidth}
    \centering
        \includegraphics[width=0.98\linewidth]{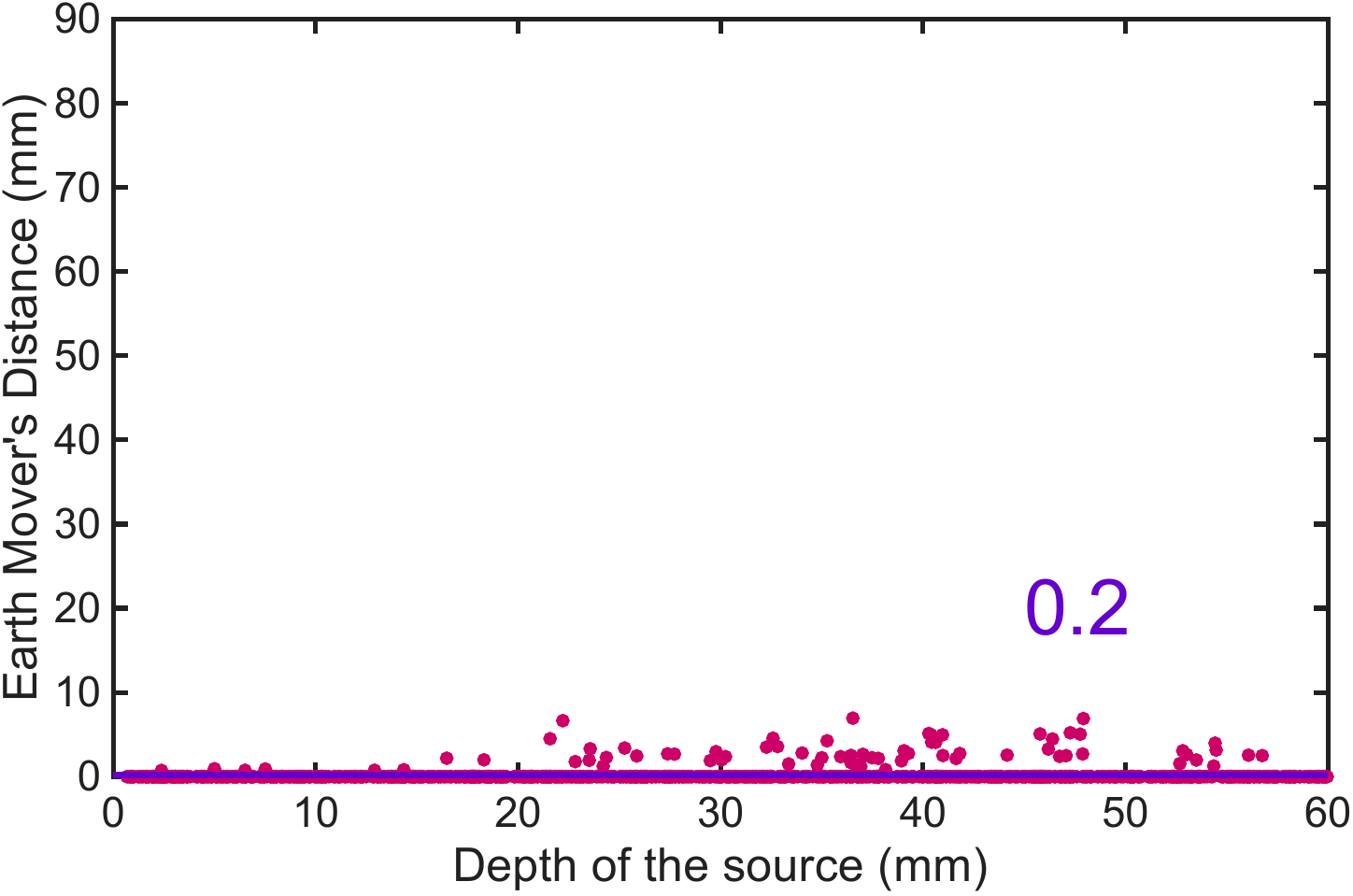}
    \end{minipage}\begin{minipage}{0.3\linewidth}
    \centering
        \includegraphics[width=0.98\linewidth]{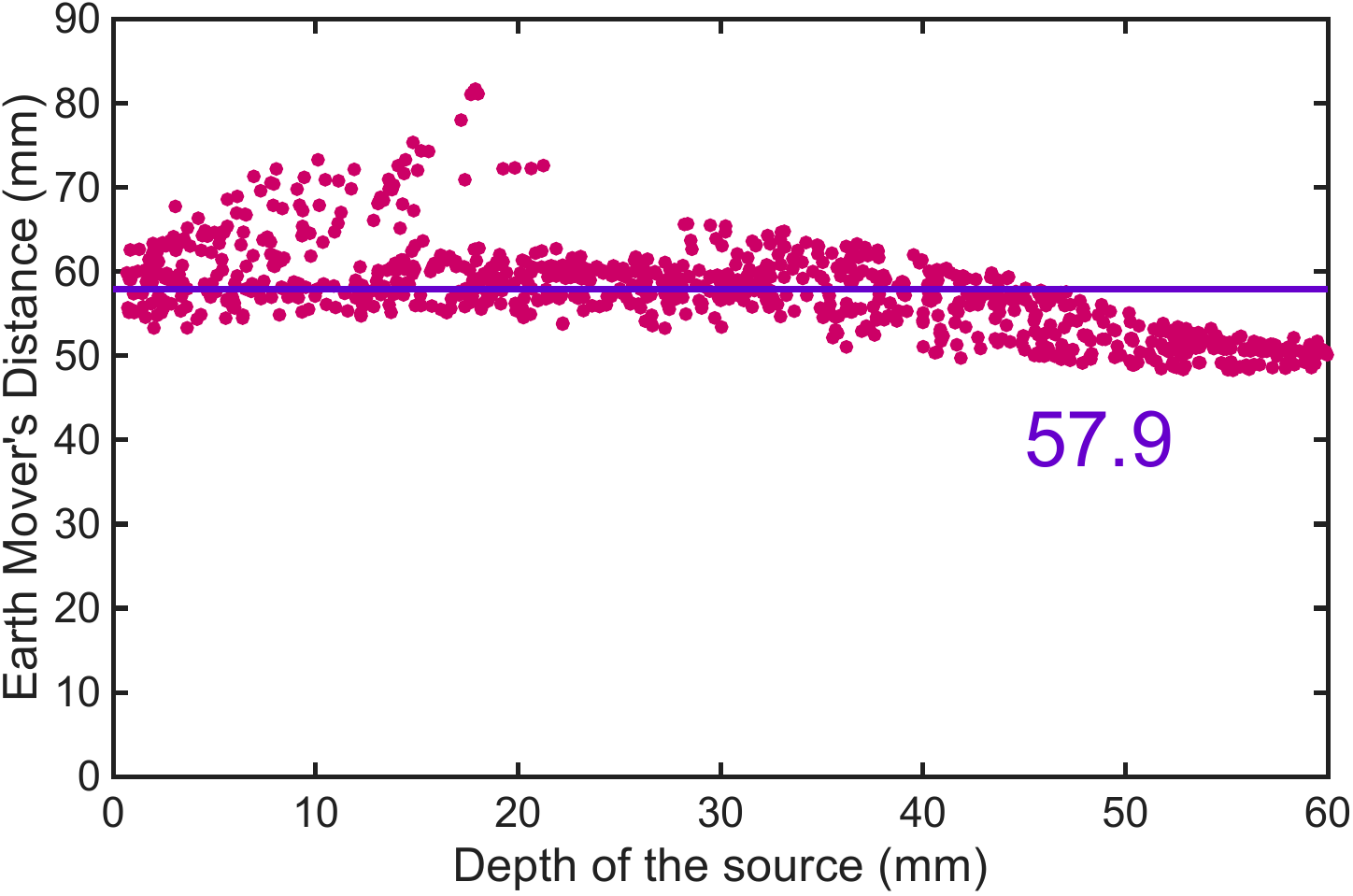}
    \end{minipage}

    \begin{minipage}{0.05\linewidth}
        \rotatebox{90}{Patch -- LS}
    \end{minipage}\begin{minipage}{0.3\linewidth}
    \centering
        \includegraphics[width=0.98\linewidth]{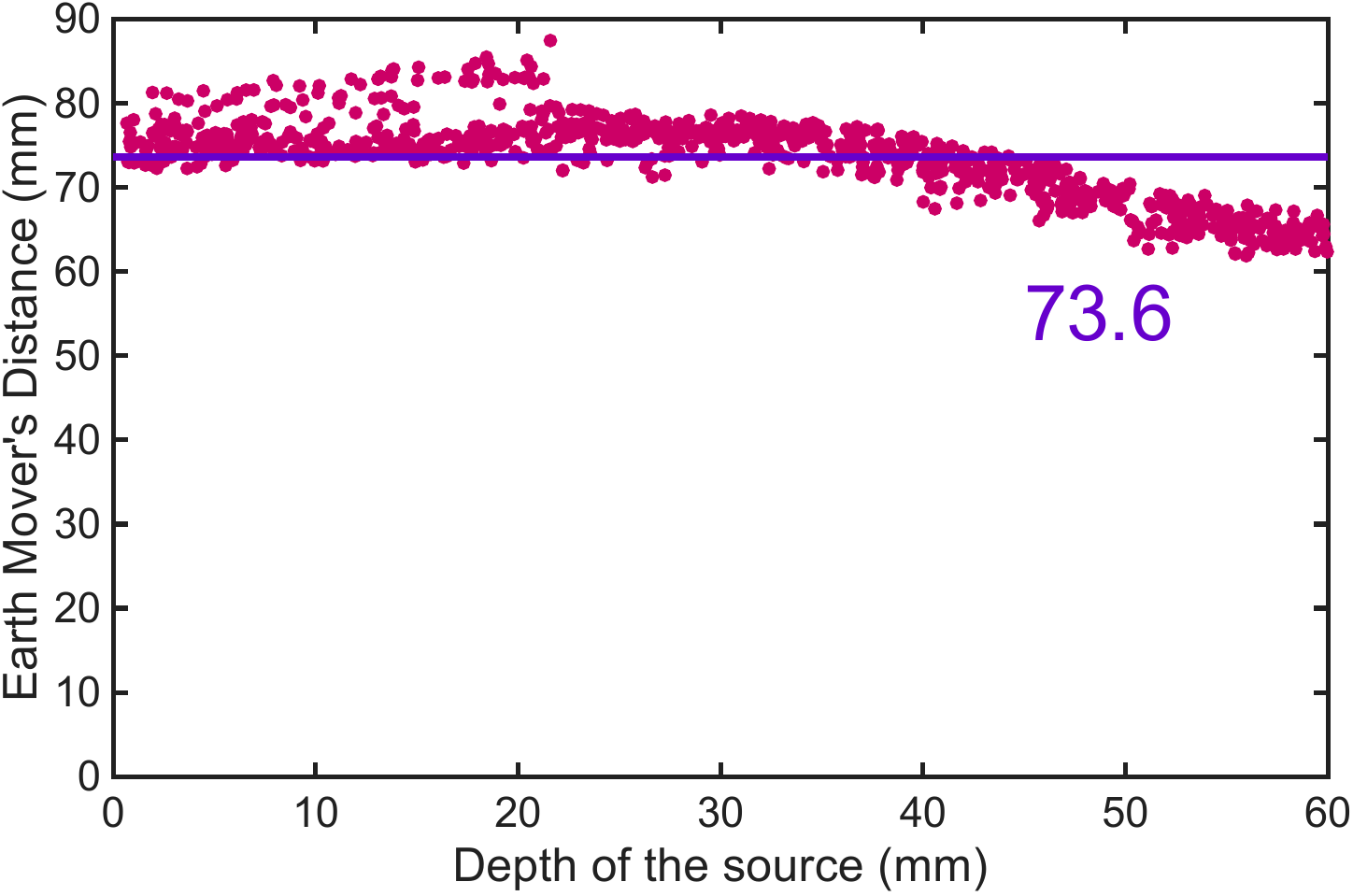}
    \end{minipage}\begin{minipage}{0.3\linewidth}
    \centering
        \includegraphics[width=0.98\linewidth]{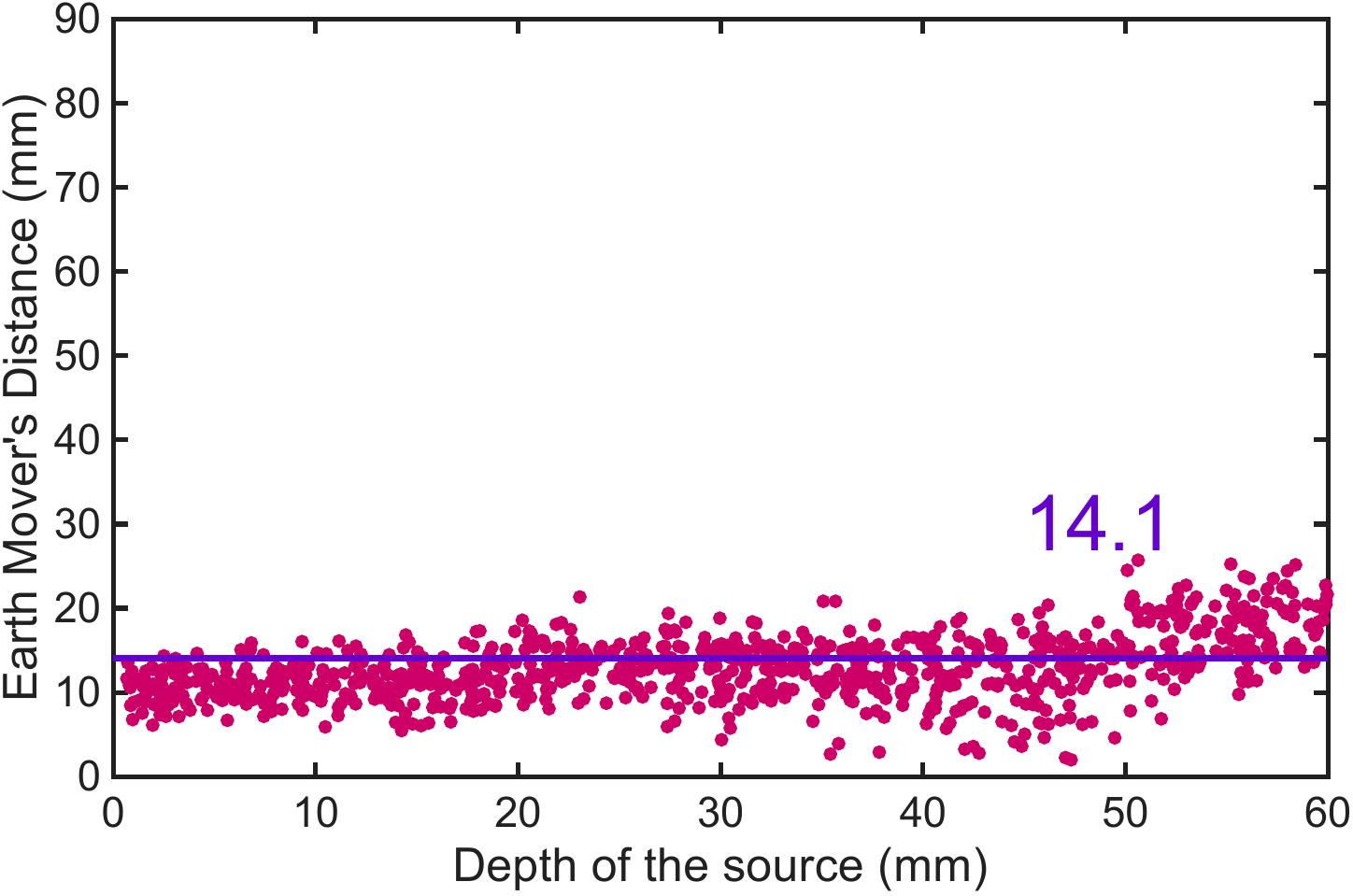}
    \end{minipage}\begin{minipage}{0.3\linewidth}
    \centering
        \includegraphics[width=0.98\linewidth]{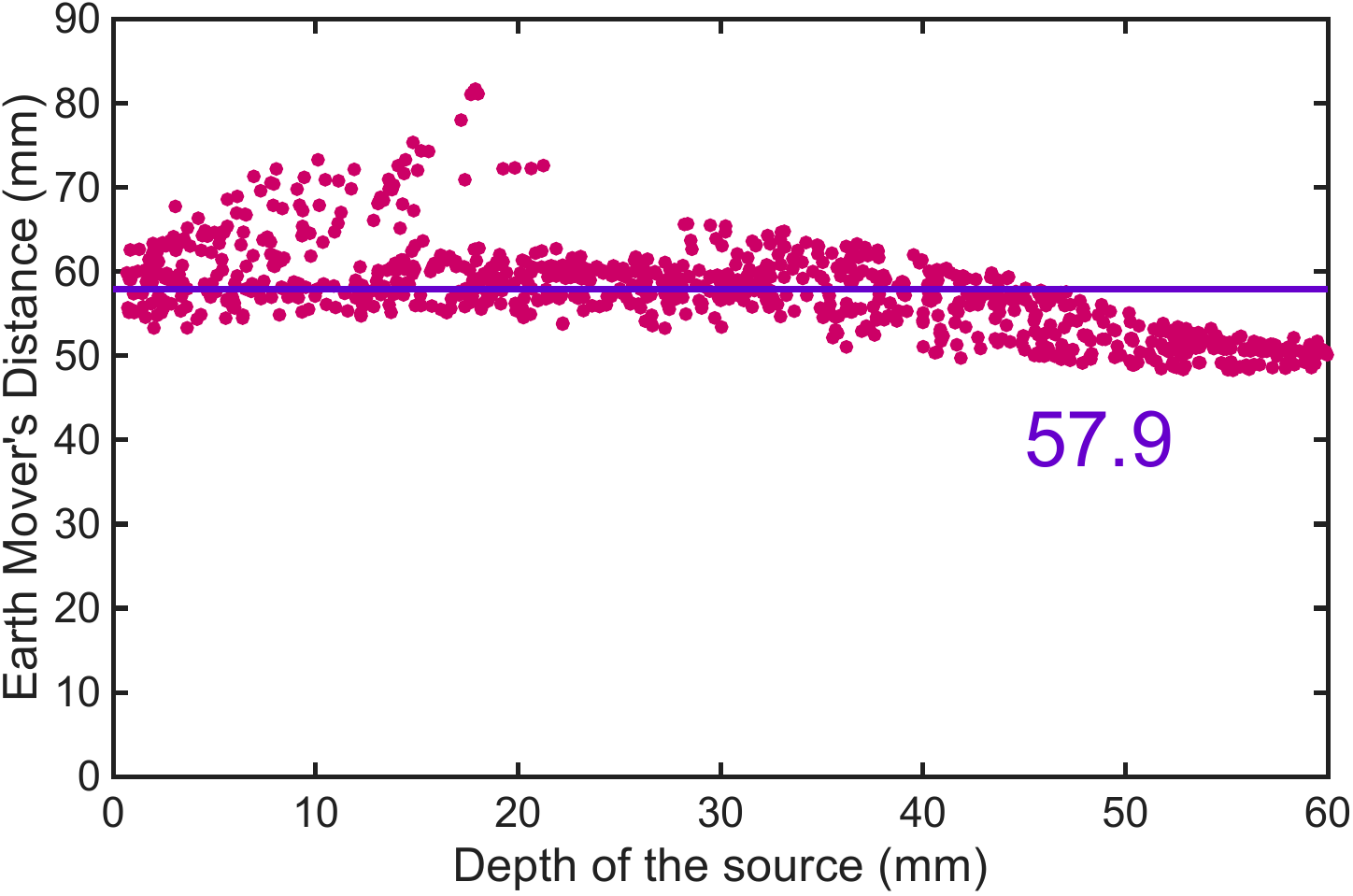}
    \end{minipage}
    
    \caption{Earth Mover's Distances for estimated sources. The violet horizontal line indicates the sample mean, which is also shown numerically on the graph.}
    \label{fig.duneuro.EMD}
\end{figure*}

Similarly to the depth estimation, Figure~\ref{fig.duneuro.EMD} shows a scatter plot of EMDs described by \eqref{eq.emd} between true and estimated sources at different distances from the inner skull surface. The results show high EMD values for sLORETA, which are highly similar across Whitney, Local subtraction, and the Patch model. The $\Hdiv$ model provides a surprisingly steady trend for EMD, which decreases only slightly for deep sources after \qty{30}{\milli\meter}.

The sparsity-promoting SHAL1R consistently produces the smallest EMD across methods and models. The mean EMDs of SHAL1R are \qty{3.0}{\milli\meter} for Whitney-basis, \qty{0.0}{\milli\meter} for Local subtraction, \qty{11.4}{\milli\meter} for $\Hdiv$, and \qty{0.2}{\milli\meter} for the patch model. The highest mean EMD is obtained when the patch source is estimated with Local subtraction. For Local subtraction, the estimate of SHAL1R is concentrated around a small region around the true source, resulting in near-zero EMD with almost no variance. The largest variance is with $\Hdiv$.

The EMDs of SKF can be grouped into three: Whitney and Local subtraction, the patch model and forward patch -- inverted Local subtraction, and $\Hdiv$. With Witney and Local subtraction, we observe the lowest mean EMD of \qty{41.5}{\milli\meter}. The trend is that EMD is steady and highly concentrated around \qty{40}{\milli\meter} up to the depth of \qty{50}{\milli\meter}, and beyond that, the EMD increases. With the patch modelled forward activity, the inversion model between the patch and Local subtraction does not impact EMD. We see a bifurcation of the main trend to steady \qty{57.9}{\milli\meter} mean value and increase to \qty{85}{\milli\meter} that last to \qty{30}{\milli\meter} depth. The steady trend around the mean stays at \qty{45}{\milli\meter} after which the EMD starts decreasing. The highest mean EMD of \qty{61.0}{\milli\meter} and the highest variety of the individual estimations are obtained with $\Hdiv$.

Results for sLORETA are highly similar to the ones for SKF, except that the EMD means are mostly higher: \qty{74.6}{\milli\meter} for Whitney and Local subtraction, \qty{73.1}{\milli\meter} for the patch model, and a slightly worse mean of \qty{73.6}{\milli\meter} when inverting the patch source with Local subtraction. The only exception is $\Hdiv$ with the lowest average EMD of \qty{53.5}{\milli\meter}, which is also the least variant. The same EMD trend observed with SKF and the patch model is now present with every model except $\Hdiv$. The differences are the lower variability at lower depths and a more pronounced drop in the deepest regime of \num{5}--\qty{6}{\centi\meter}.

\subsection{Source estimation at various depths from noisy data}\label{sec.noisy.signatures}

\begin{figure*}
    \centering
    \begin{minipage}{0.05\linewidth}
        \rotatebox{90}{Whitney}
    \end{minipage}\begin{minipage}{0.3\linewidth}
    \centering
    sLORETA \vspace{0.1cm}
    
        \includegraphics[width=0.98\linewidth]{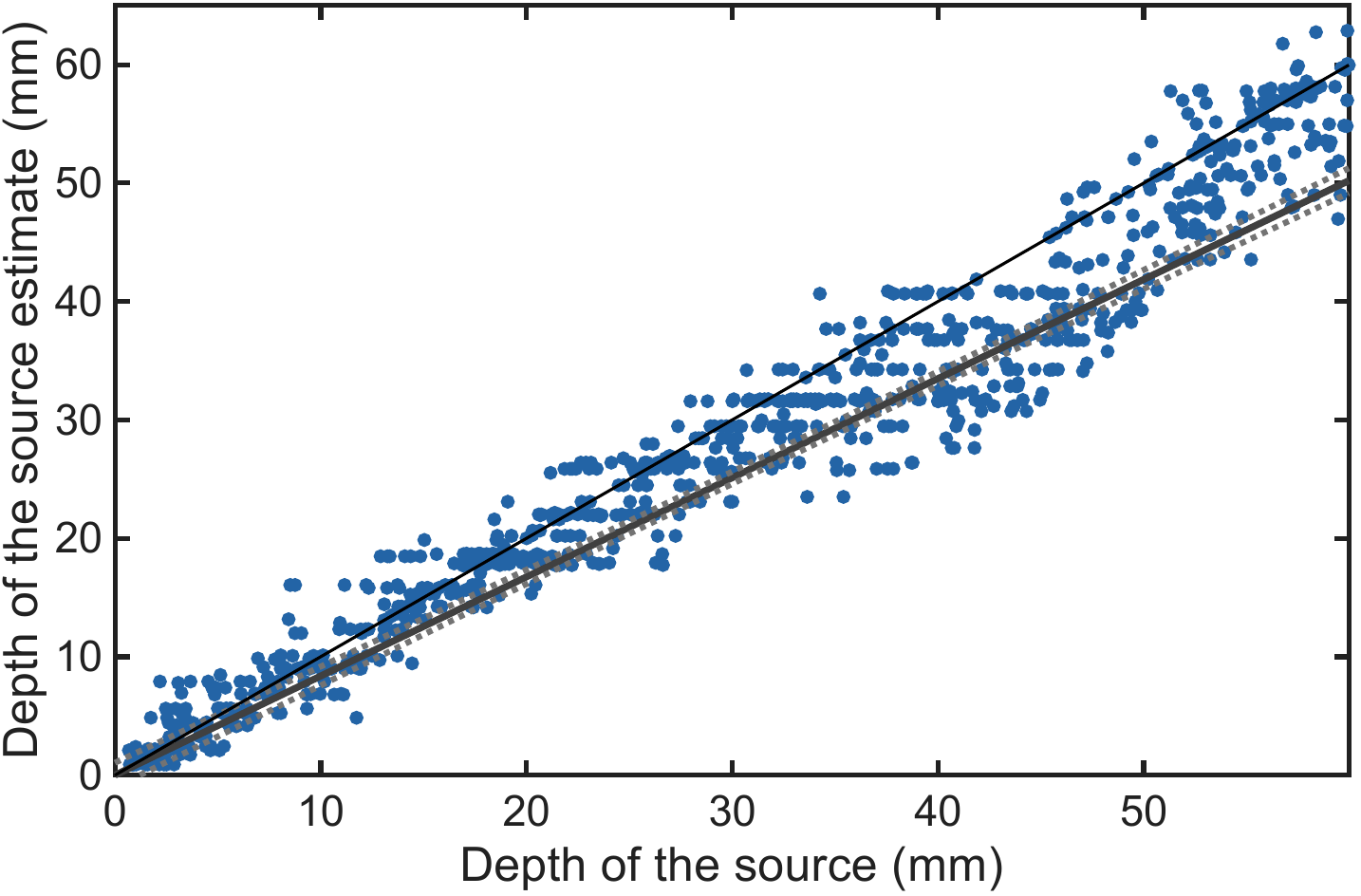}
    \end{minipage}\begin{minipage}{0.3\linewidth}
    \centering
    SHAL1R \vspace{0.1cm}
    
        \includegraphics[width=0.98\linewidth]{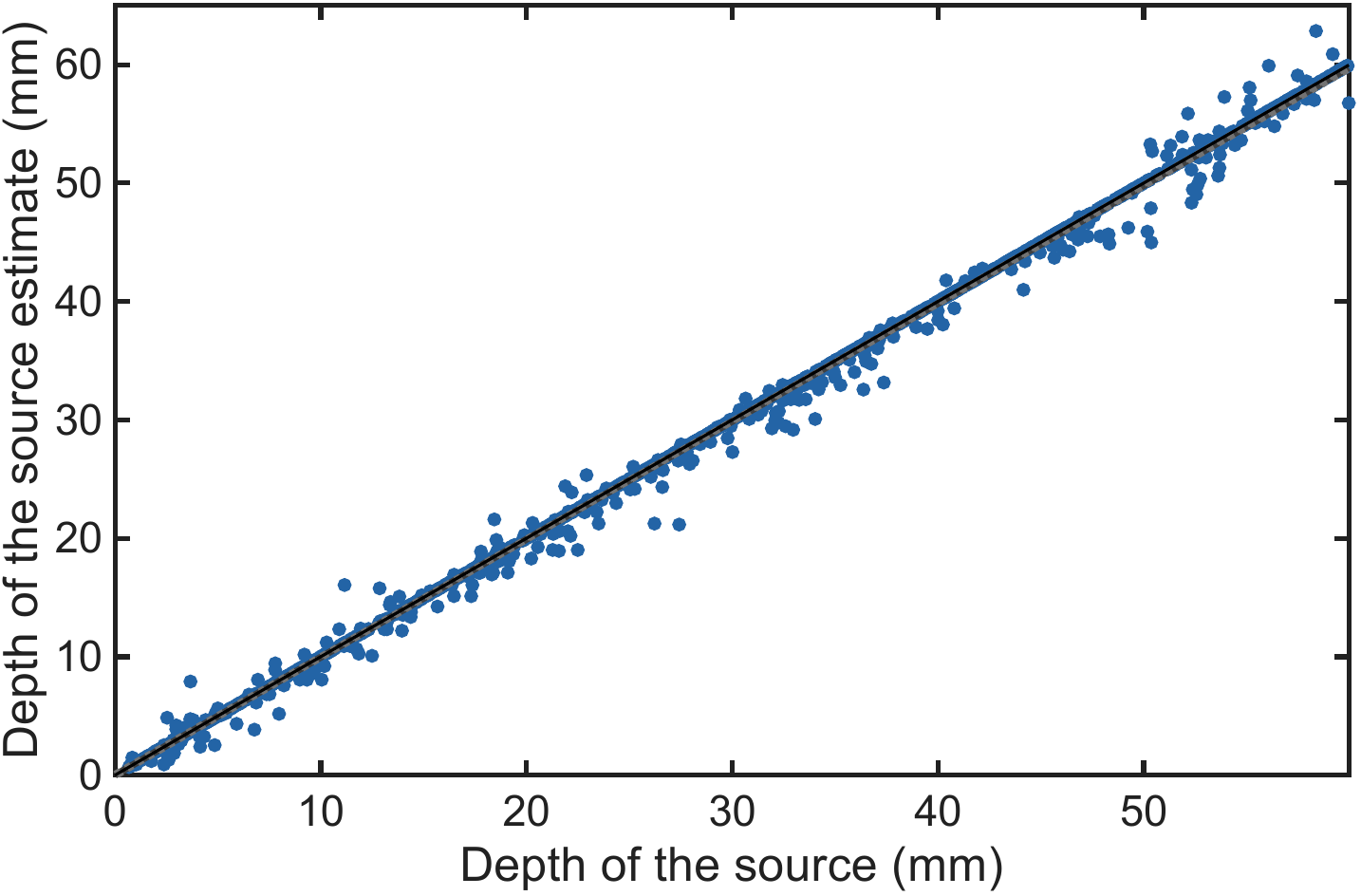}
    \end{minipage}\begin{minipage}{0.3\linewidth}
    \centering
    SKF \vspace{0.1cm}
    
        \includegraphics[width=0.98\linewidth]{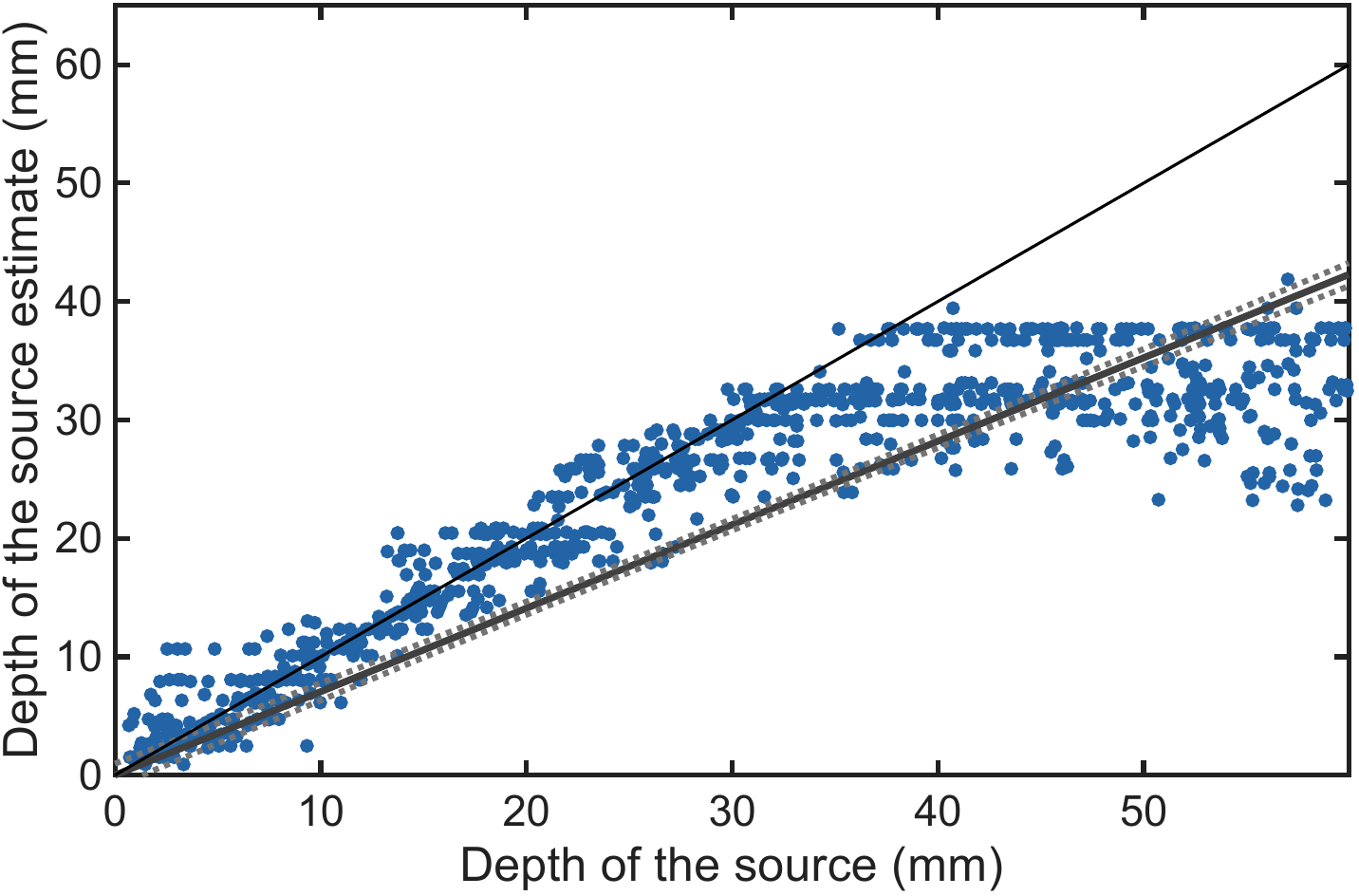}
    \end{minipage}

    \begin{minipage}{0.05\linewidth}
        \rotatebox{90}{$\Hdiv$}
    \end{minipage}\begin{minipage}{0.3\linewidth}
    \centering
        \includegraphics[width=0.98\linewidth]{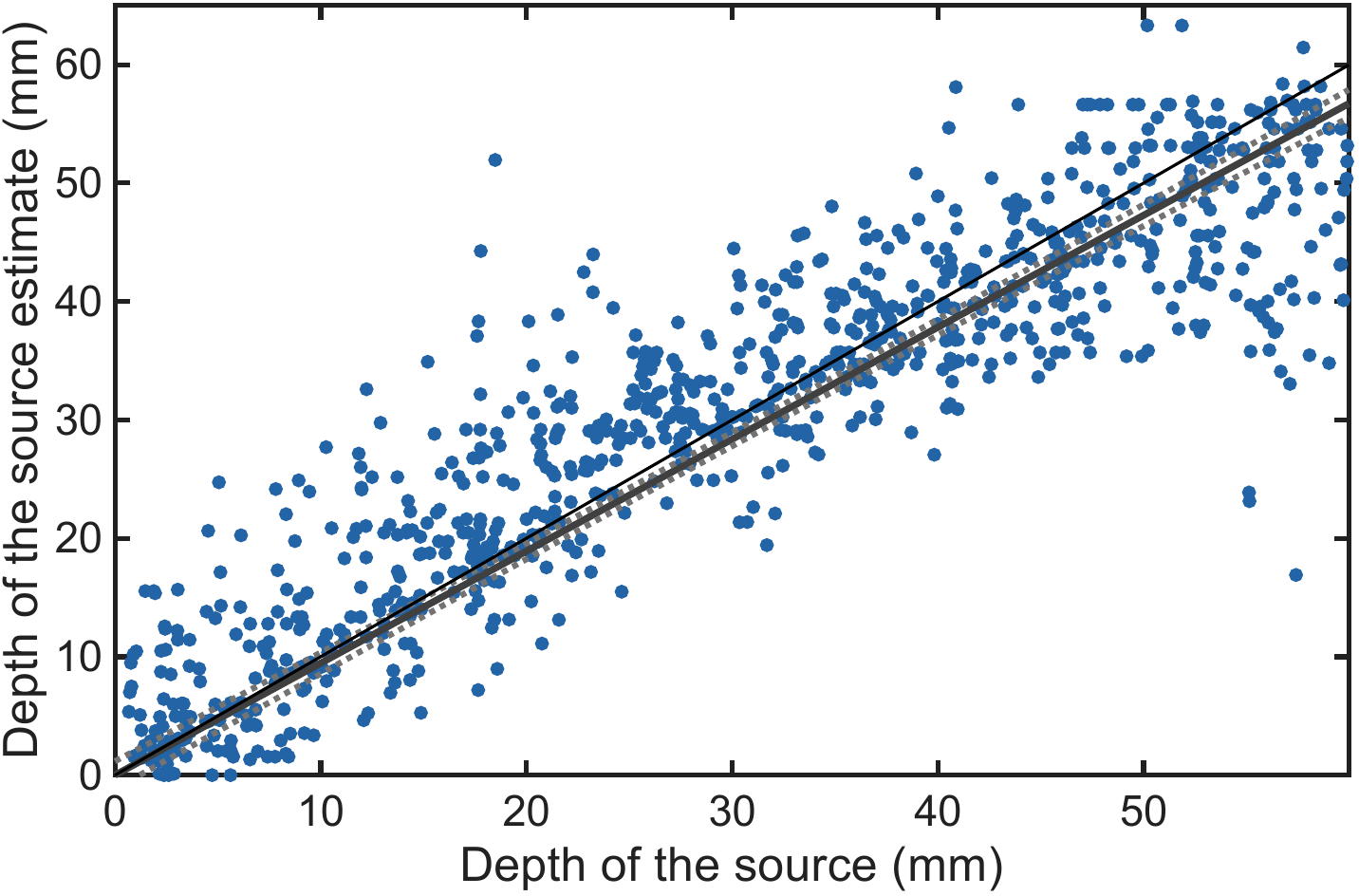}
    \end{minipage}\begin{minipage}{0.3\linewidth}
    \centering
        \includegraphics[width=0.98\linewidth]{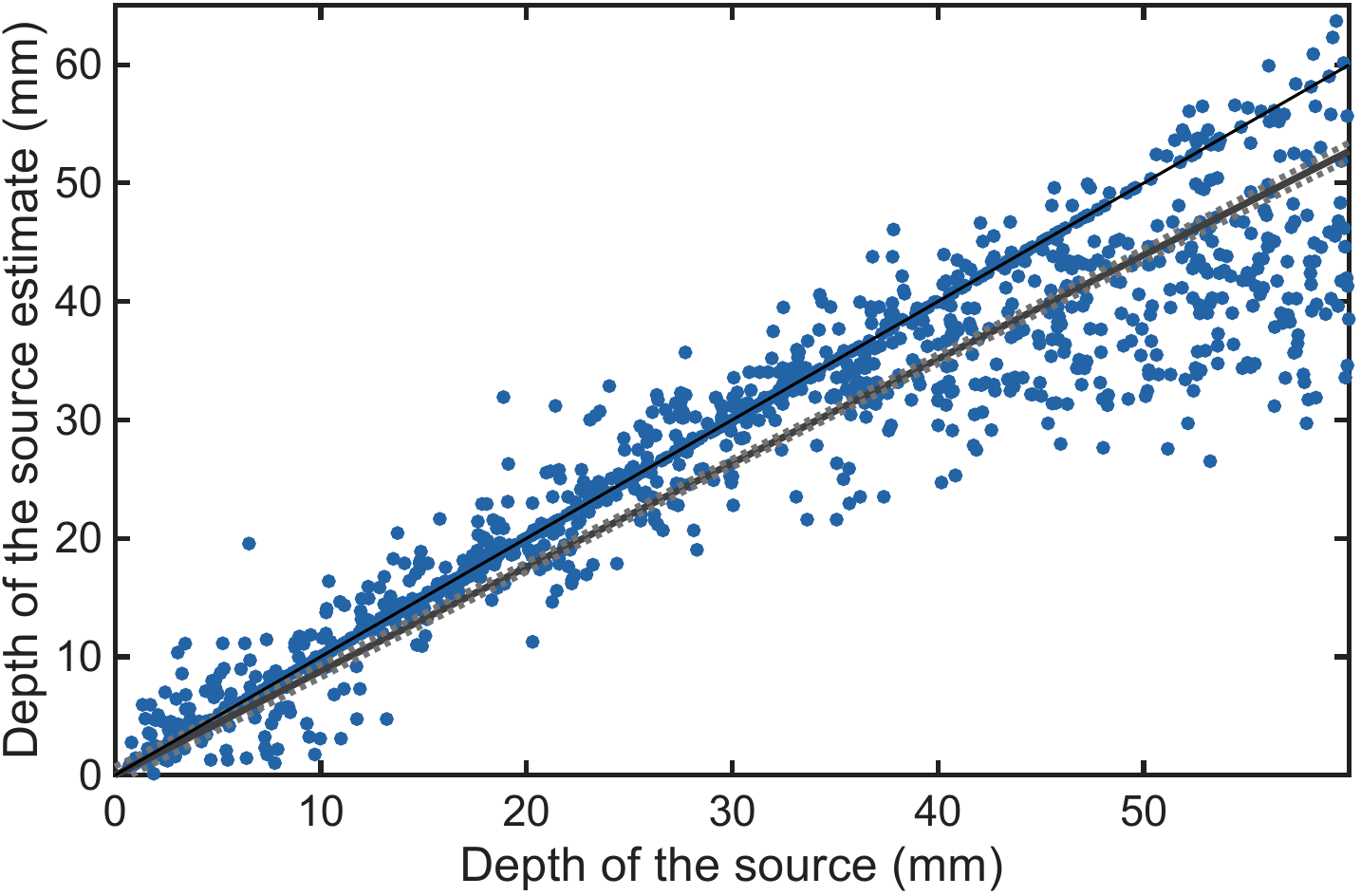}
    \end{minipage}\begin{minipage}{0.3\linewidth}
    \centering
        \includegraphics[width=0.98\linewidth]{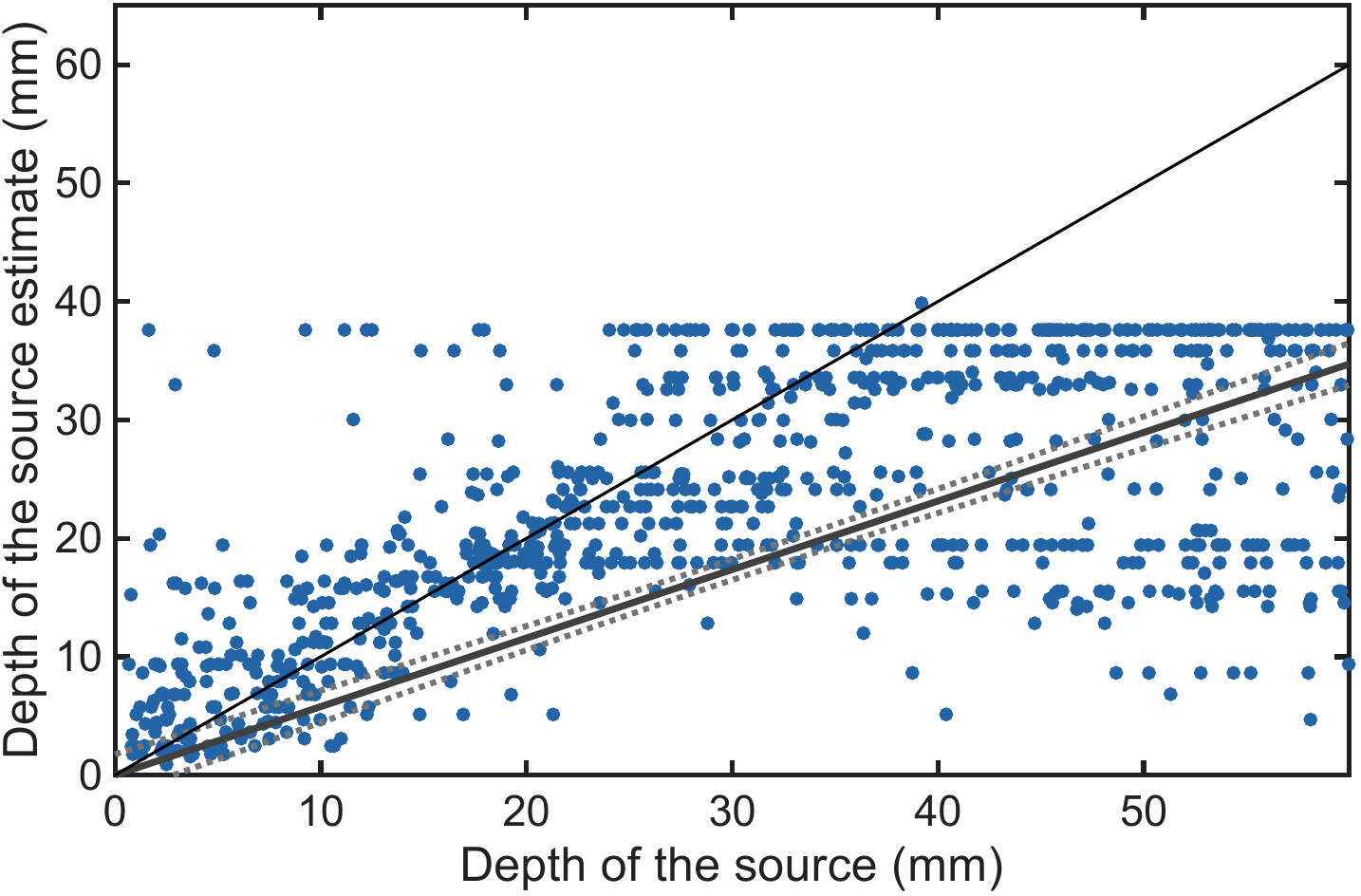}
    \end{minipage}

    \begin{minipage}{0.05\linewidth}
        \rotatebox{90}{Local subt.}
    \end{minipage}\begin{minipage}{0.3\linewidth}
    \centering
        \includegraphics[width=0.98\linewidth]{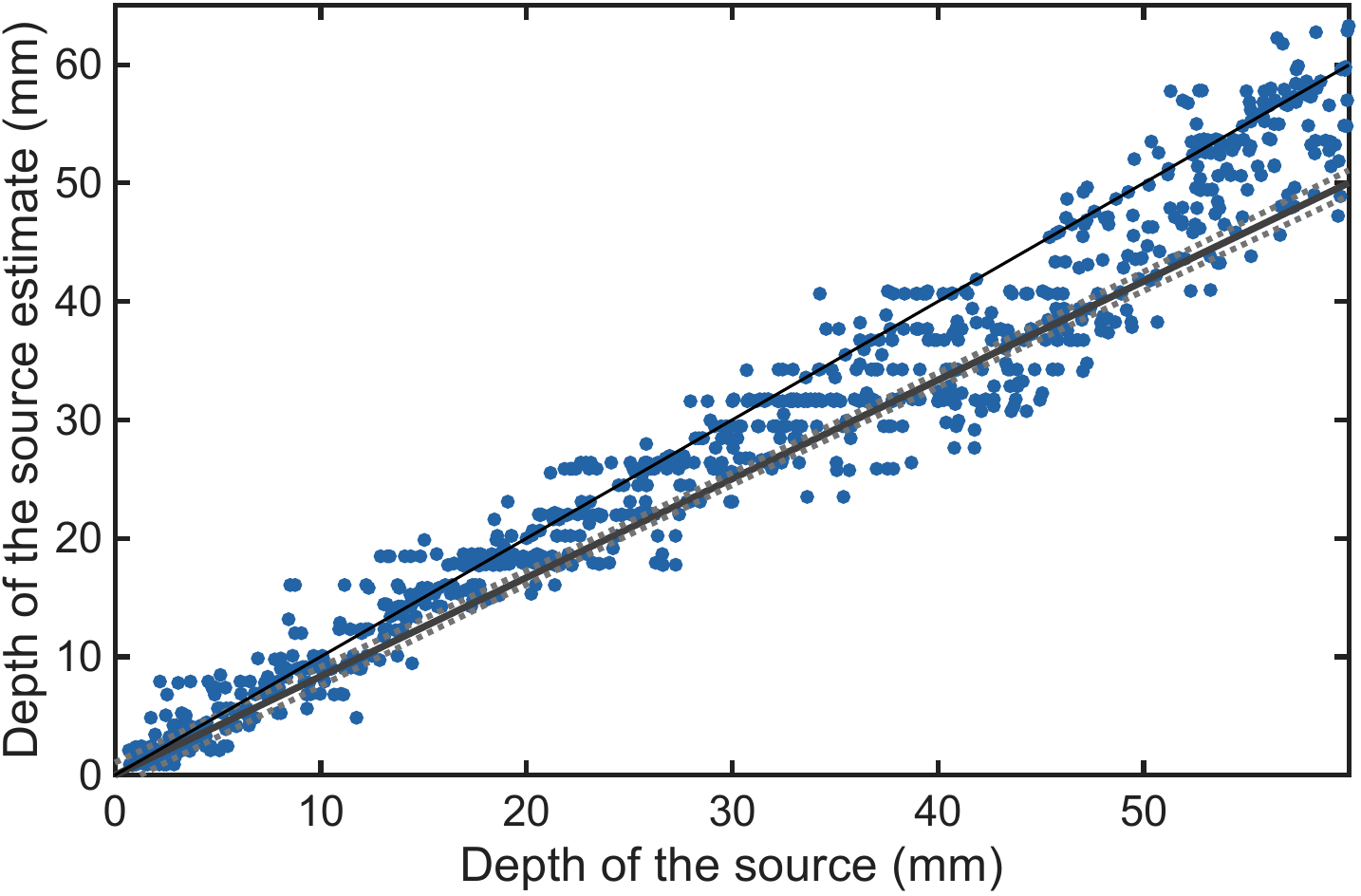}
    \end{minipage}\begin{minipage}{0.3\linewidth}
    \centering
        \includegraphics[width=0.98\linewidth]{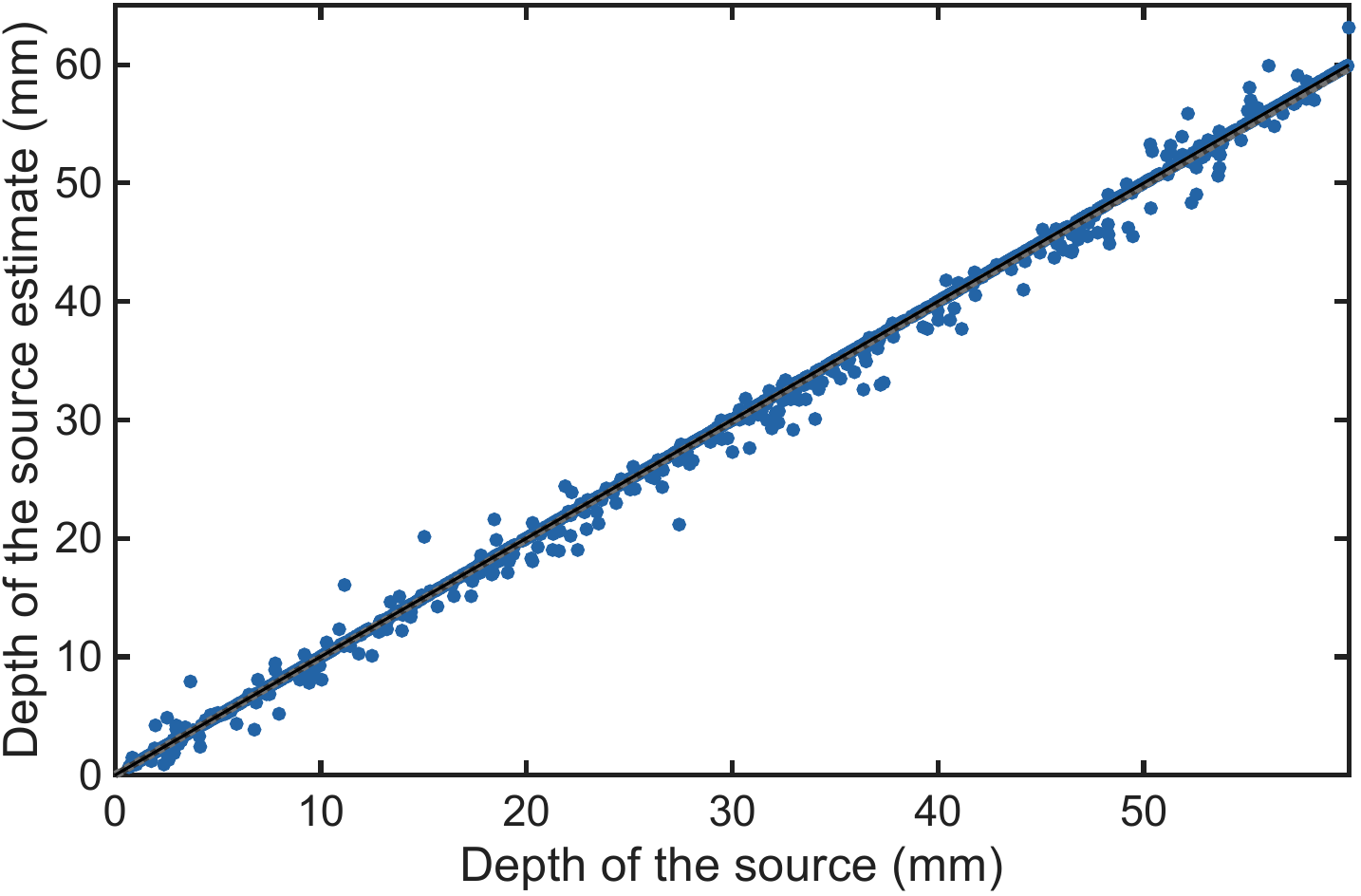}
    \end{minipage}\begin{minipage}{0.3\linewidth}
    \centering
        \includegraphics[width=0.98\linewidth]{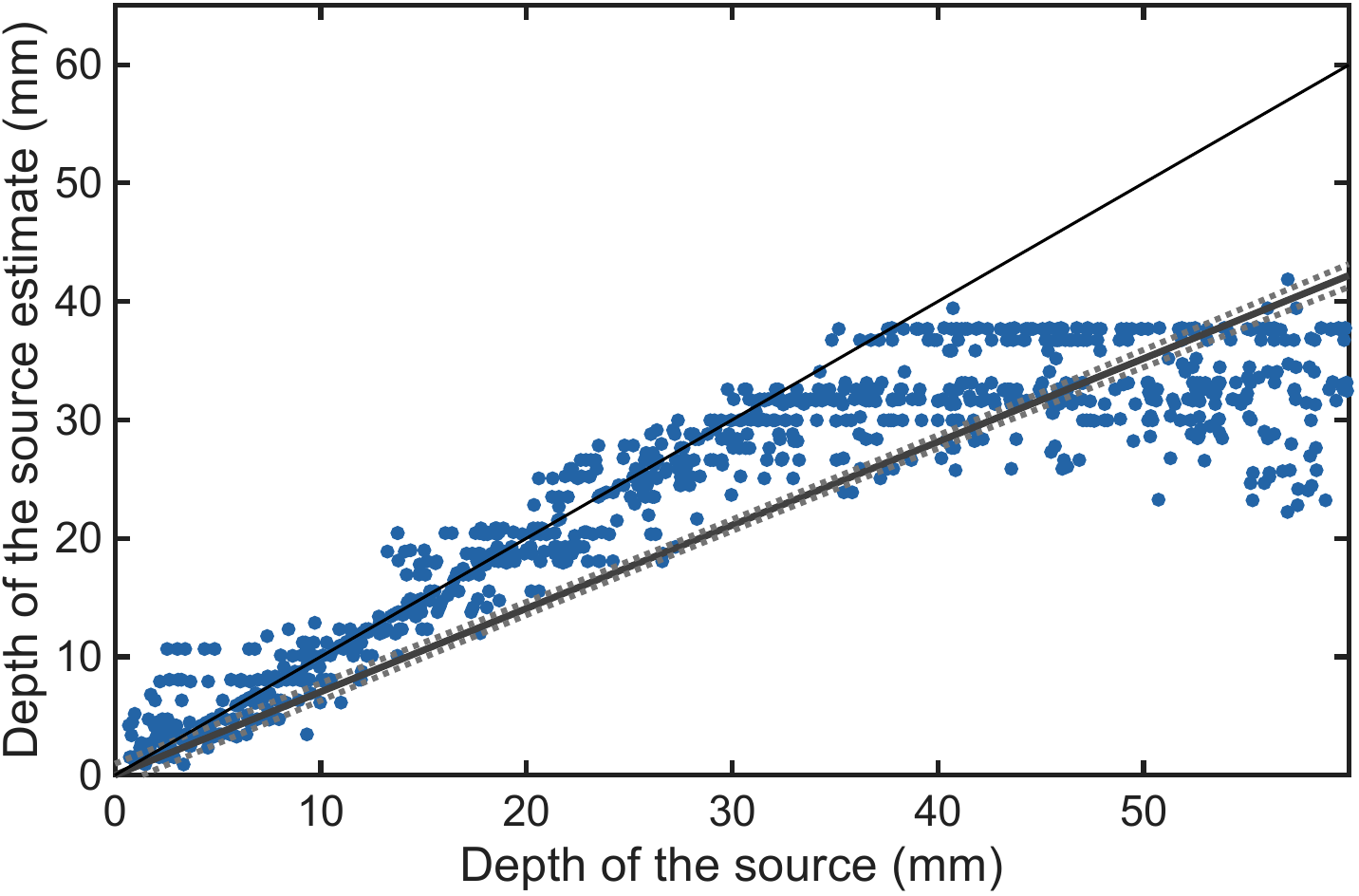}
    \end{minipage}

    \begin{minipage}{0.05\linewidth}
        \rotatebox{90}{Patch}
    \end{minipage}\begin{minipage}{0.3\linewidth}
    \centering
        \includegraphics[width=0.98\linewidth]{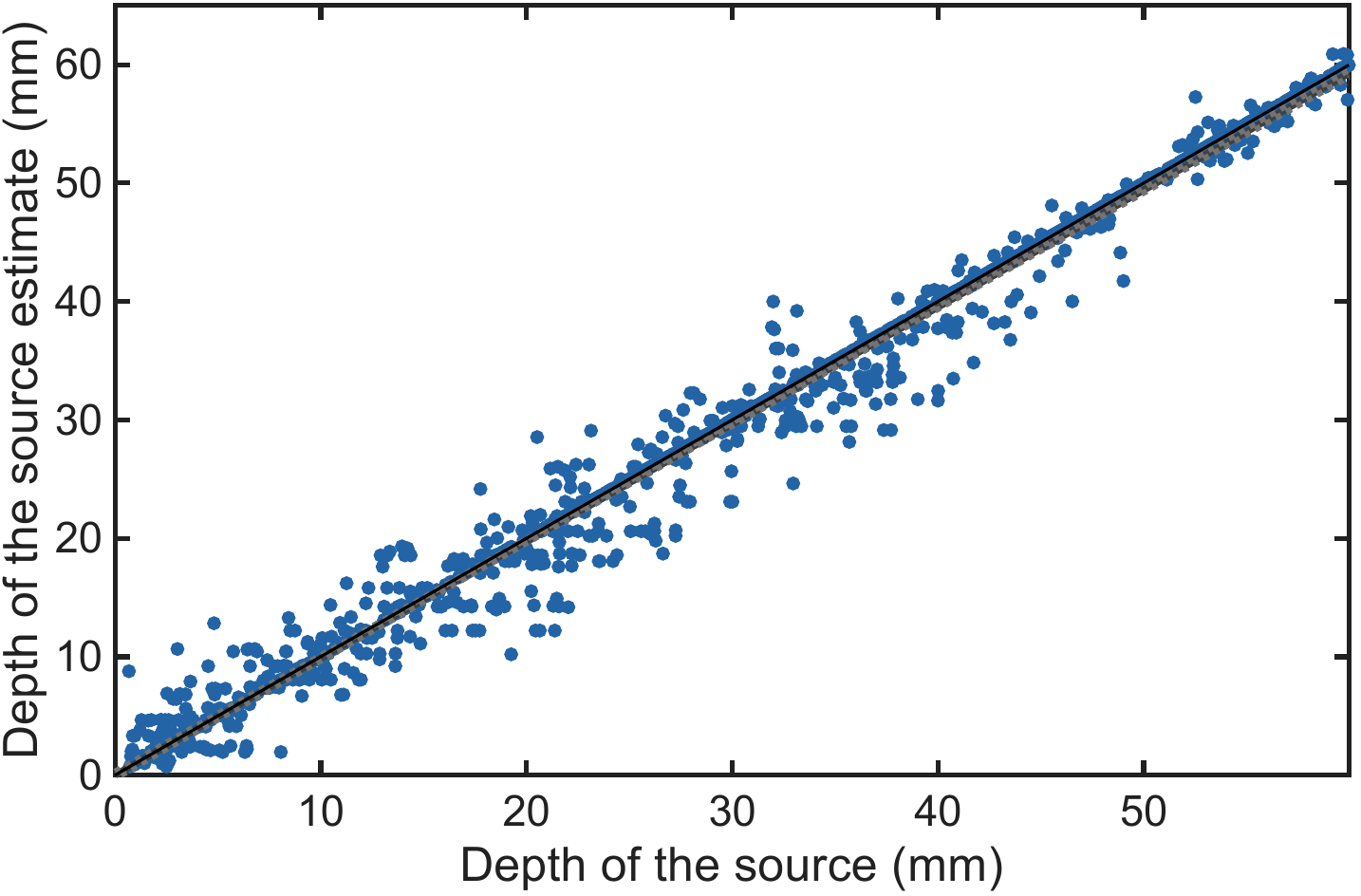}
    \end{minipage}\begin{minipage}{0.3\linewidth}
    \centering
        \includegraphics[width=0.98\linewidth]{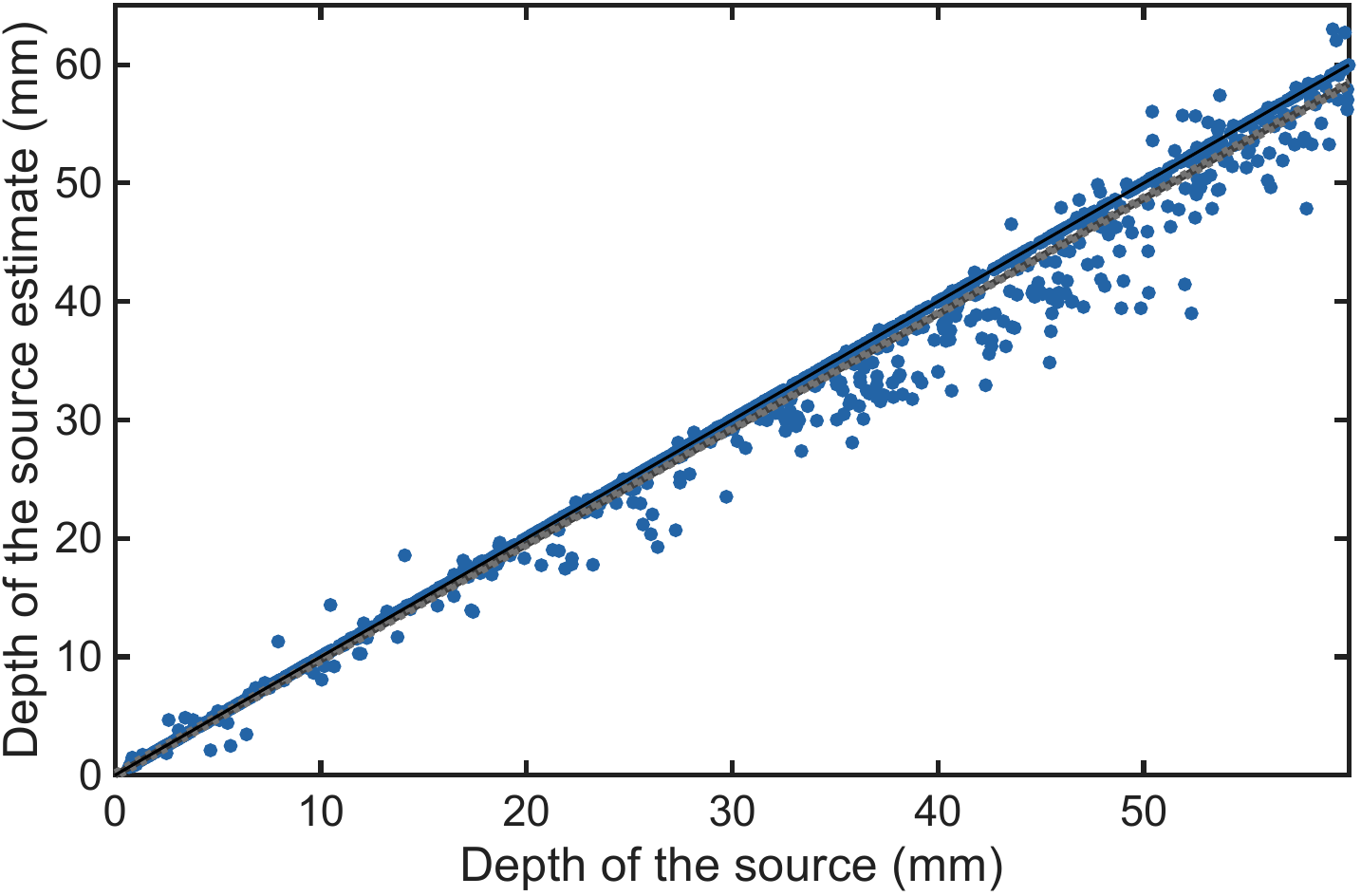}
    \end{minipage}\begin{minipage}{0.3\linewidth}
    \centering
        \includegraphics[width=0.98\linewidth]{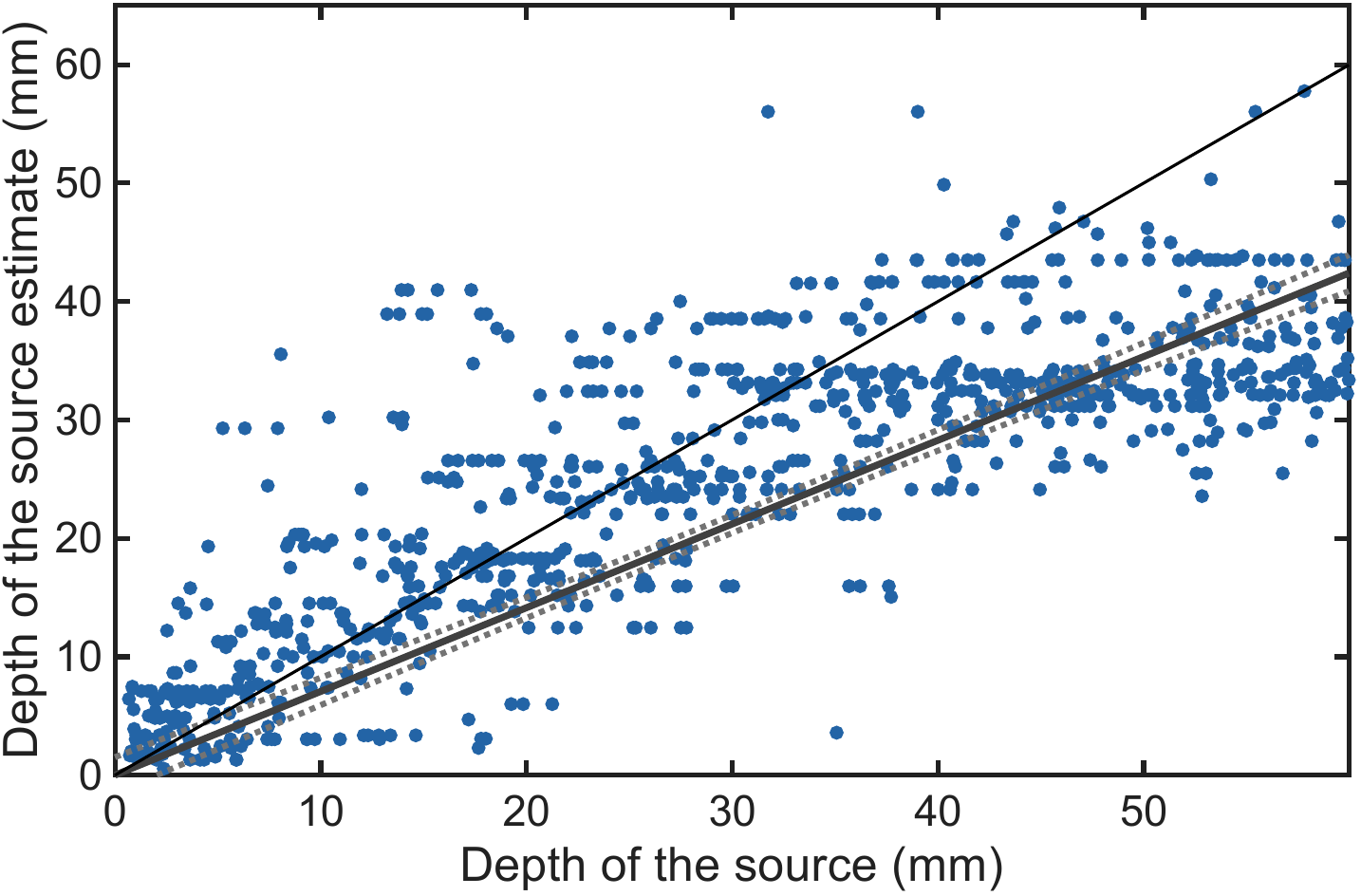}
    \end{minipage}

     \begin{minipage}{0.05\linewidth}
         \rotatebox{90}{Patch -- LS}
     \end{minipage}\begin{minipage}{0.3\linewidth}
     \centering
         \includegraphics[width=0.98\linewidth]{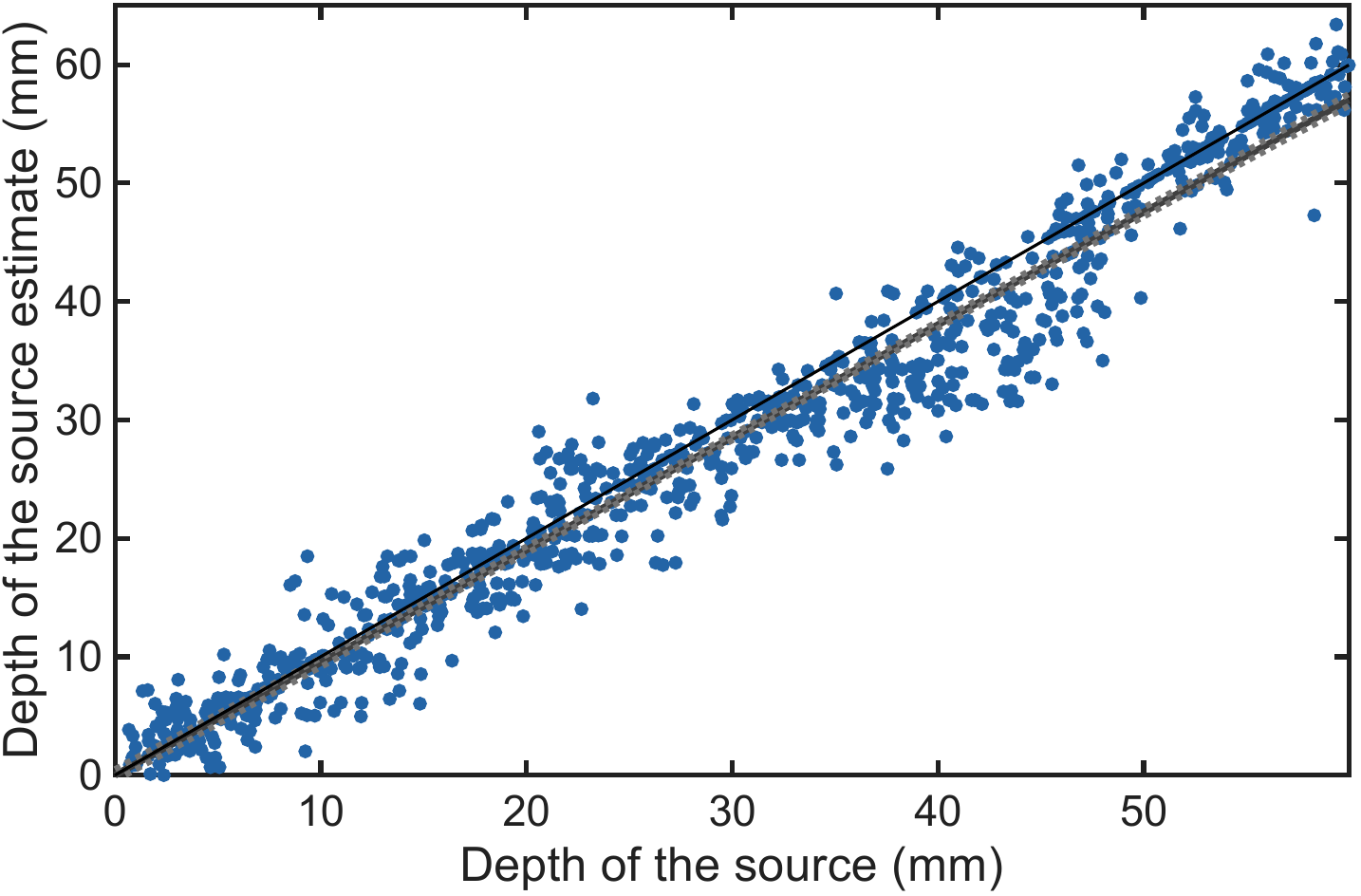}
     \end{minipage}\begin{minipage}{0.3\linewidth}
     \centering
         \includegraphics[width=0.98\linewidth]{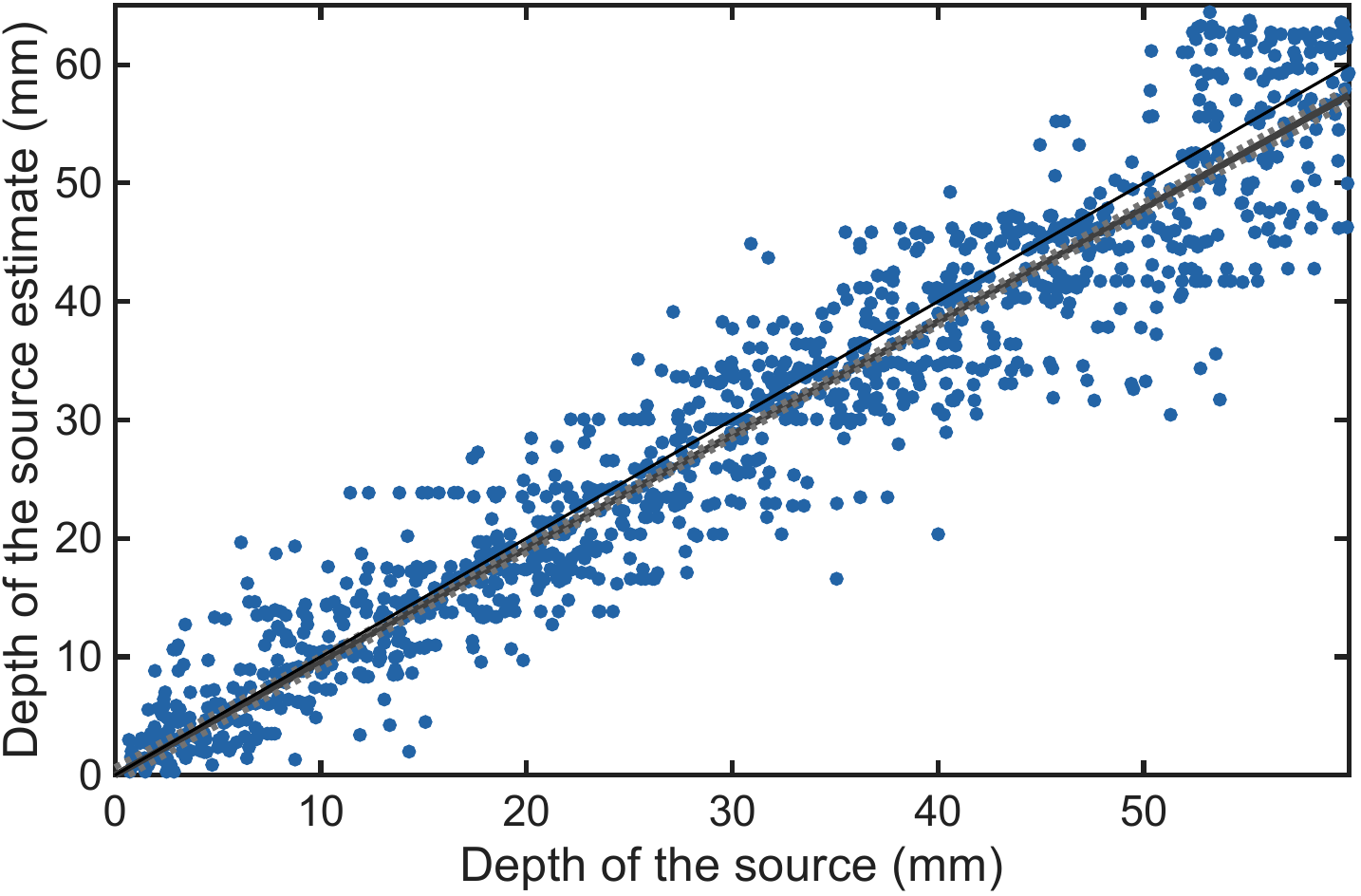}
     \end{minipage}\begin{minipage}{0.3\linewidth}
     \centering
         \includegraphics[width=0.98\linewidth]{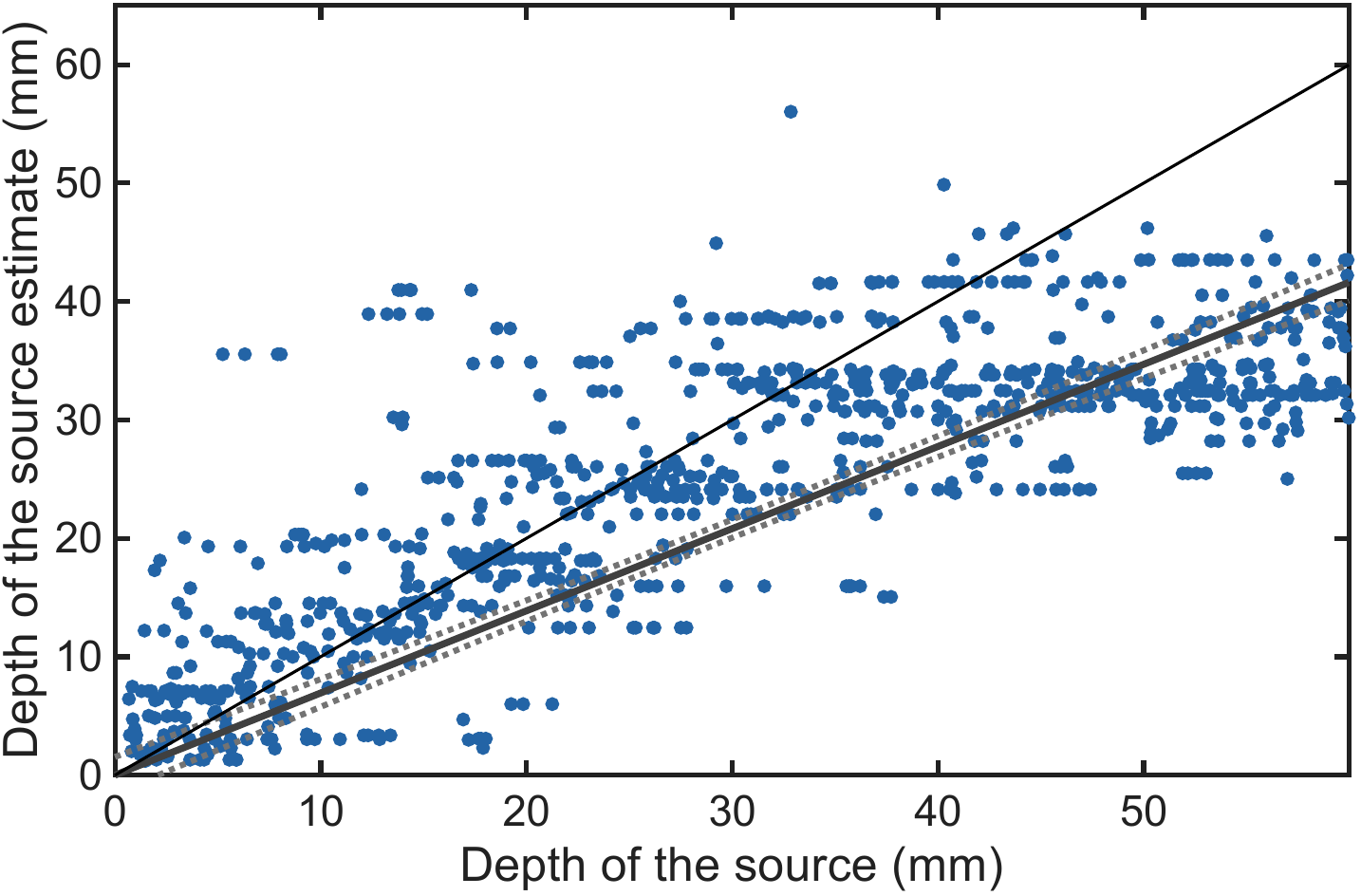}
     \end{minipage}
    
    \caption{Depth of the true source plotted against the depth estimated by sLORETA and SHAL1R at a signal-to-noise ratio (SNR) of \qty{15}{\decibel}. The thin black line shows the optimal agreement between the true and estimated depth when the localization error is zero. The dark gray solid line displays the linear regression, and the dashed gray curves are the 95 \% confidence intervals.}
    \label{fig.duneuro.depthbias15dB}
\end{figure*}

To assess how the depth bias of the various inverse methods would behave in a more realistic scenario, we also added Gau\ss ian noise to the simulated data yielding \qty{15}{\decibel} SNR, and then inverted the measurements with an equivalent noiseless $\leadFieldMatrix$. The resulting depth-bias scatter plots (DBSP) are seen in Figure~\ref{fig.duneuro.depthbias15dB}, while the exact slopes of the regression lines are found in Table~\ref{tab.regressionslopes}.

The added noise primarily causes greater deviations from ideal unbiased behavior, as seen in the increased spread of the results around the lines describing the ideal depth localization. Although the noise does not drastically alter the shapes of any of the point distributions compared to the noiseless case presented in Figure~\ref{fig.duneuro.depthbias}. One notable change, however, is with the interaction of Whitney and Local subtraction with SHAL1R: the agreement with the depth of the source and estimate is improved, resulting in a better coincidence of the ideal and regression line slopes, with the new values being 1.00 for both when compared to the previous \num{1.03} and \num{1.02}. In contrast, the deviation between $\Hdiv$ and the patch model increases. There is no significant impact for sLORETA and SKF with most of the source model. The regression slope of sLORETA decreases to 0.94 with $\Hdiv$, while the slope of SKF increases to 0.58 for the same source model as the most drastic change of them all.

\begin{figure*}
    \centering
    \begin{minipage}{0.05\linewidth}
        \rotatebox{90}{Whitney}
    \end{minipage}\begin{minipage}{0.3\linewidth}
    \centering
    sLORETA \vspace{0.1cm}
    
        \includegraphics[width=0.98\linewidth]{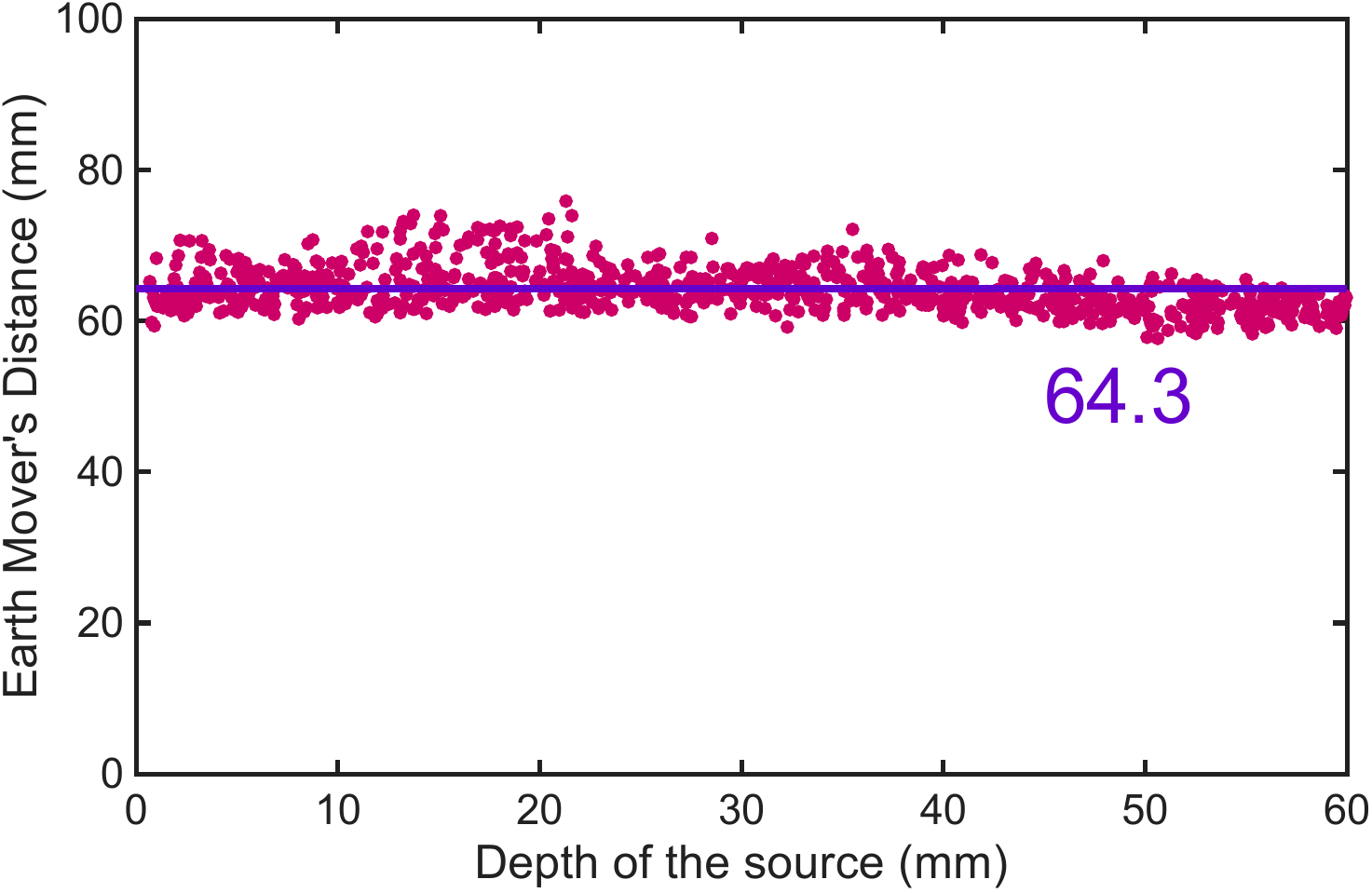}
    \end{minipage}\begin{minipage}{0.3\linewidth}
    \centering
    SHAL1R \vspace{0.1cm}
    
        \includegraphics[width=0.98\linewidth]{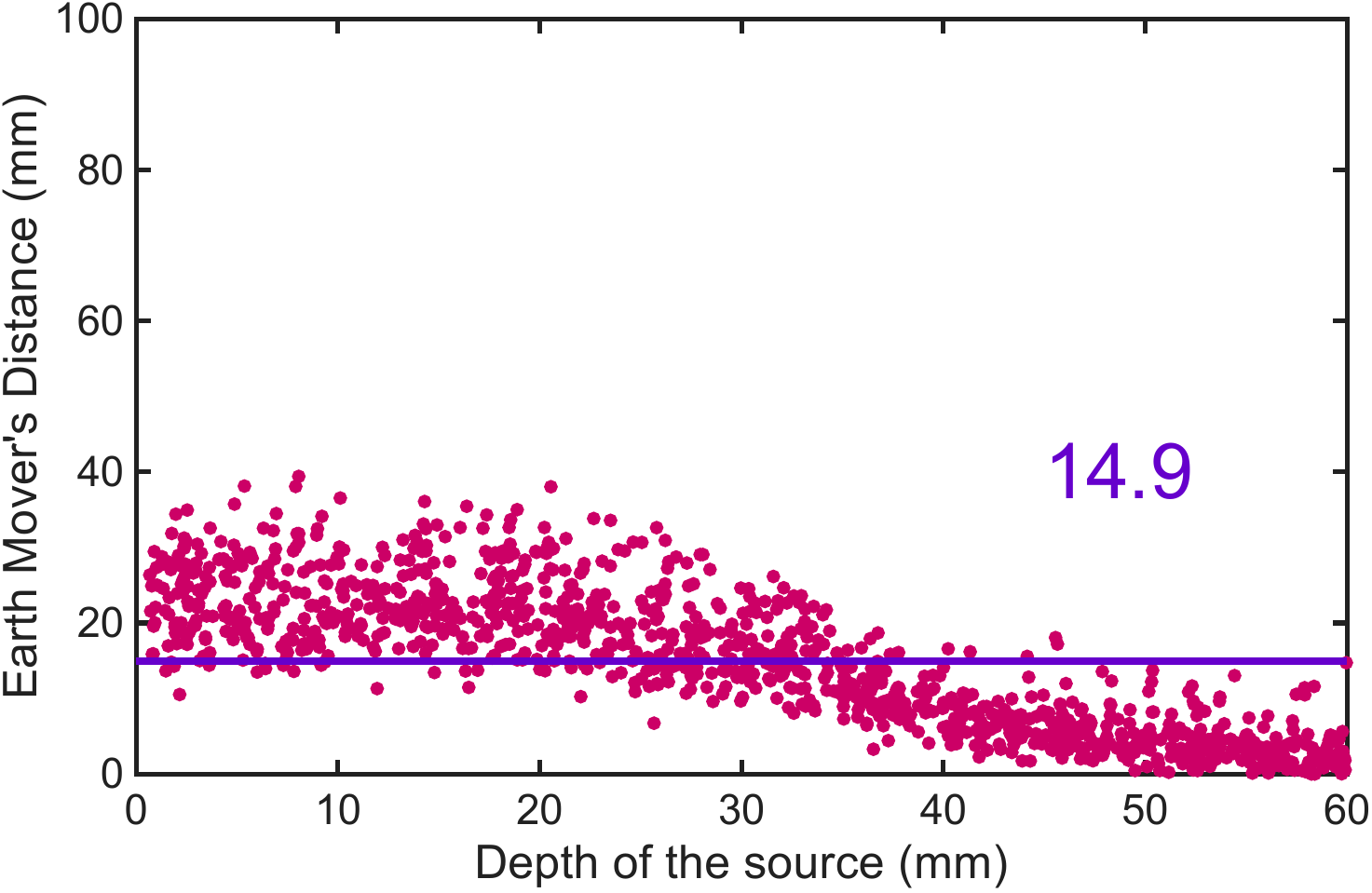}
    \end{minipage}\begin{minipage}{0.3\linewidth}
    \centering
    SKF \vspace{0.1cm}
    
        \includegraphics[width=0.98\linewidth]{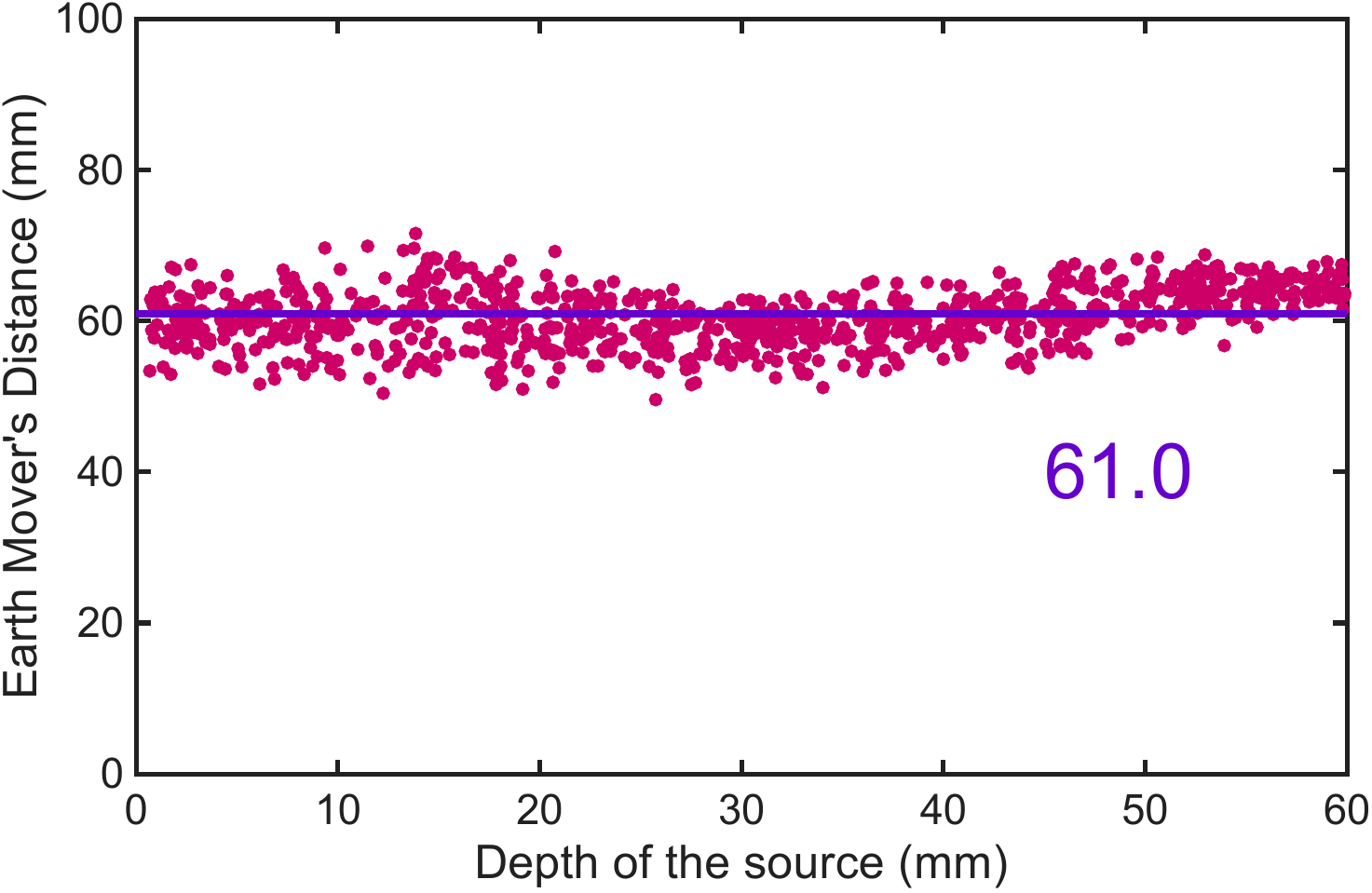}
    \end{minipage}

    \begin{minipage}{0.05\linewidth}
        \rotatebox{90}{$\Hdiv$}
    \end{minipage}\begin{minipage}{0.3\linewidth}
    \centering
        \includegraphics[width=0.98\linewidth]{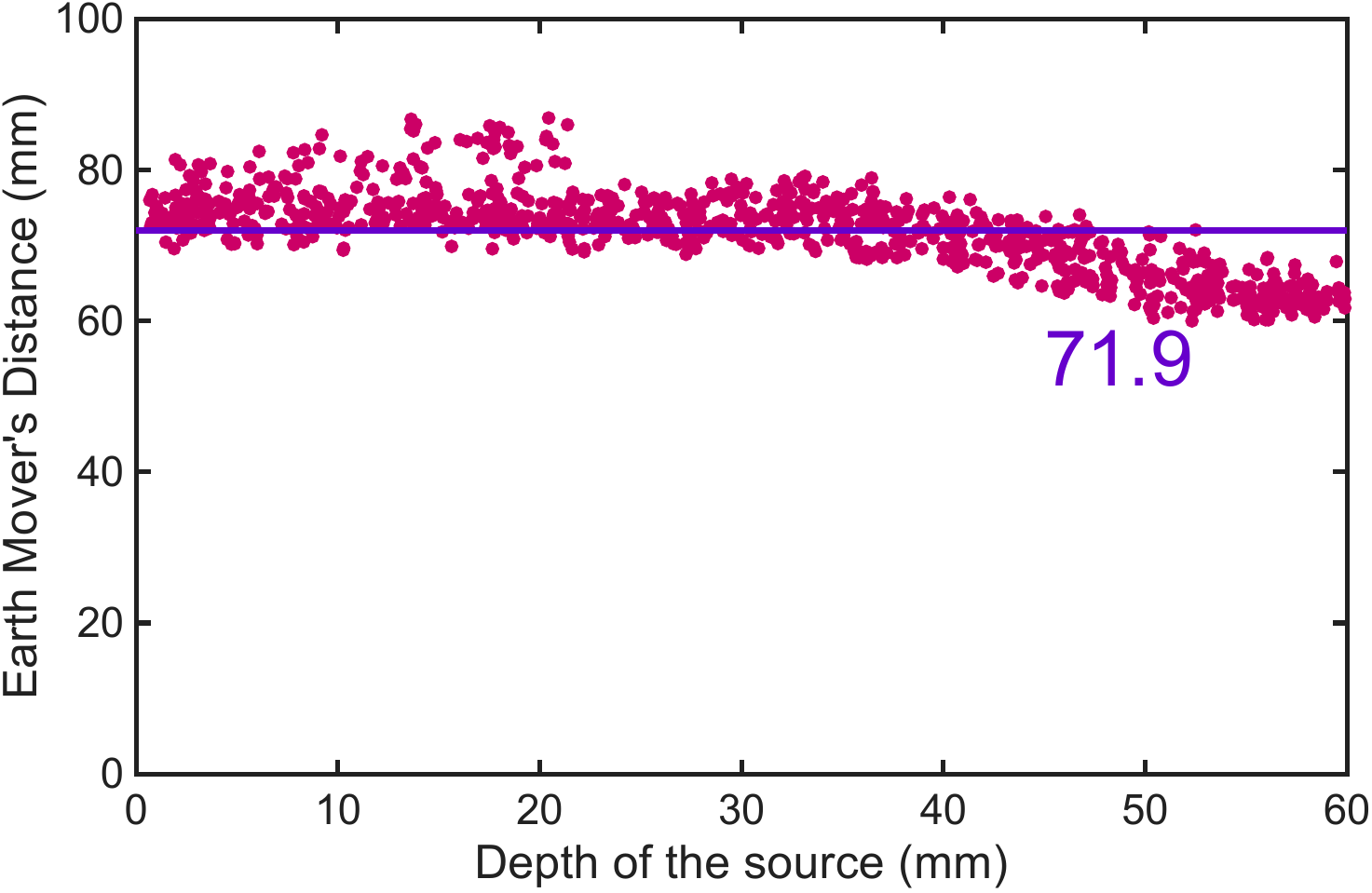}
    \end{minipage}\begin{minipage}{0.3\linewidth}
    \centering
        \includegraphics[width=0.98\linewidth]{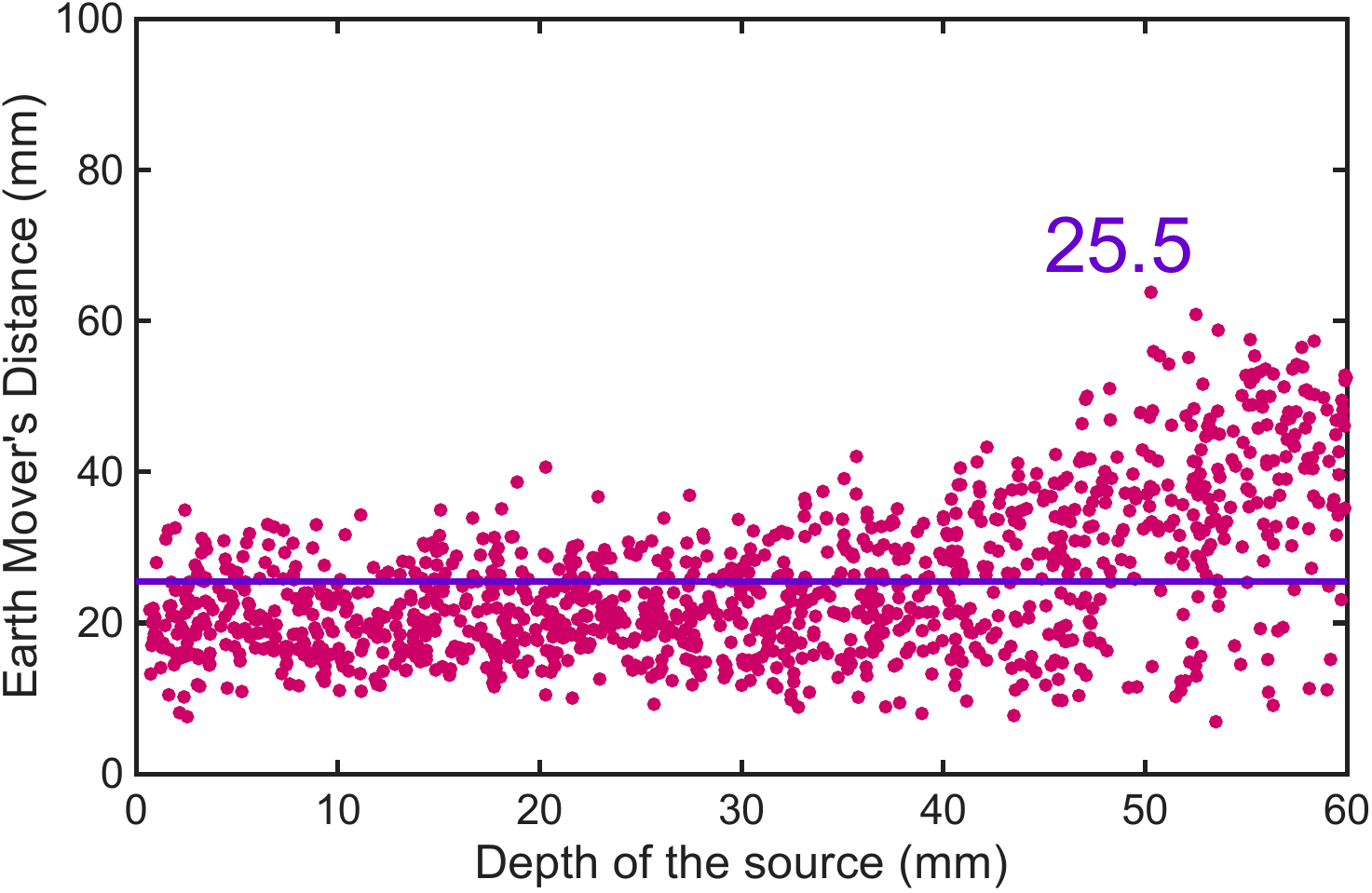}
    \end{minipage}\begin{minipage}{0.3\linewidth}
    \centering
        \includegraphics[width=0.98\linewidth]{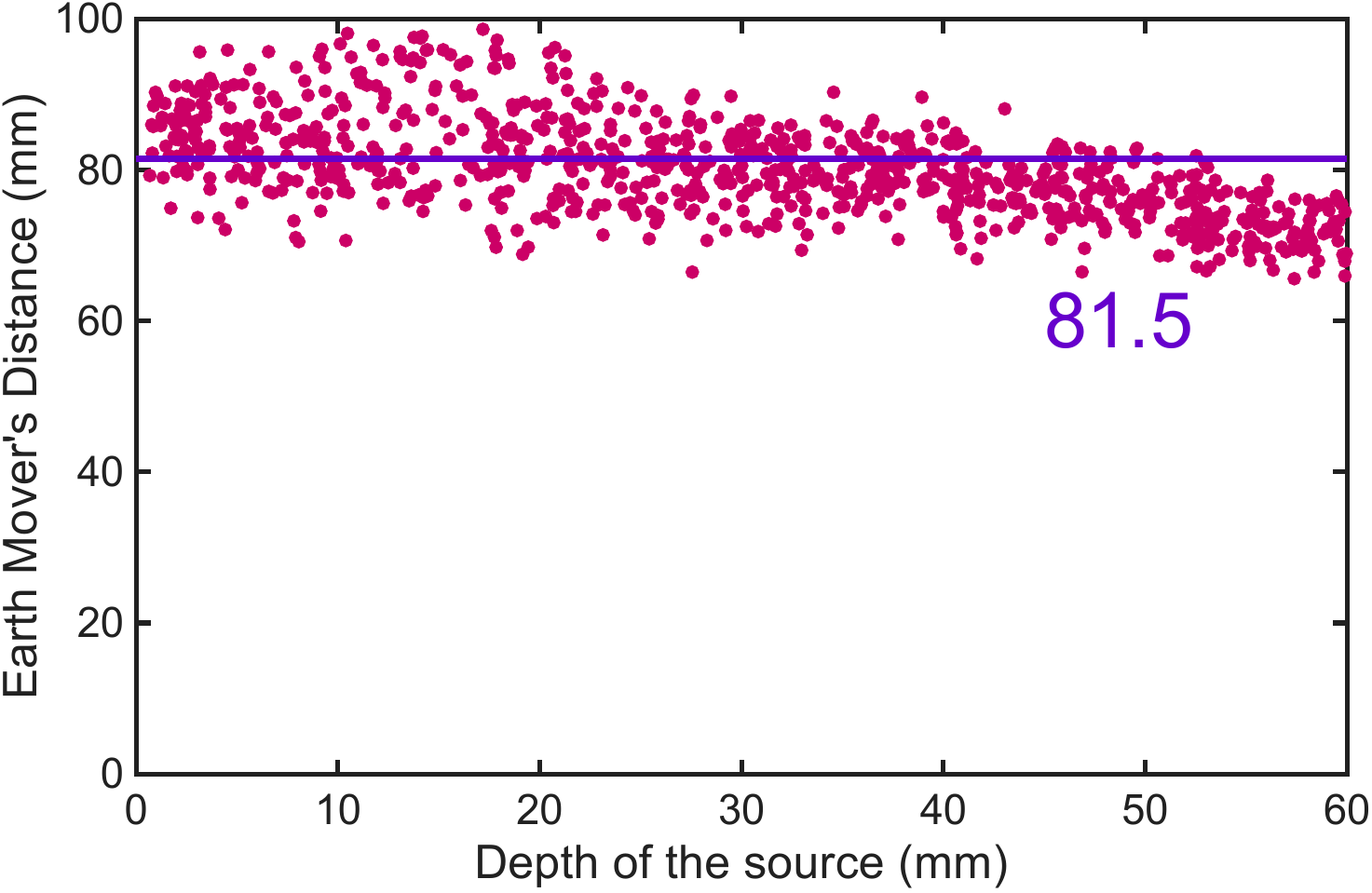}
    \end{minipage}

    \begin{minipage}{0.05\linewidth}
        \rotatebox{90}{Local subt.}
    \end{minipage}\begin{minipage}{0.3\linewidth}
    \centering
        \includegraphics[width=0.98\linewidth]{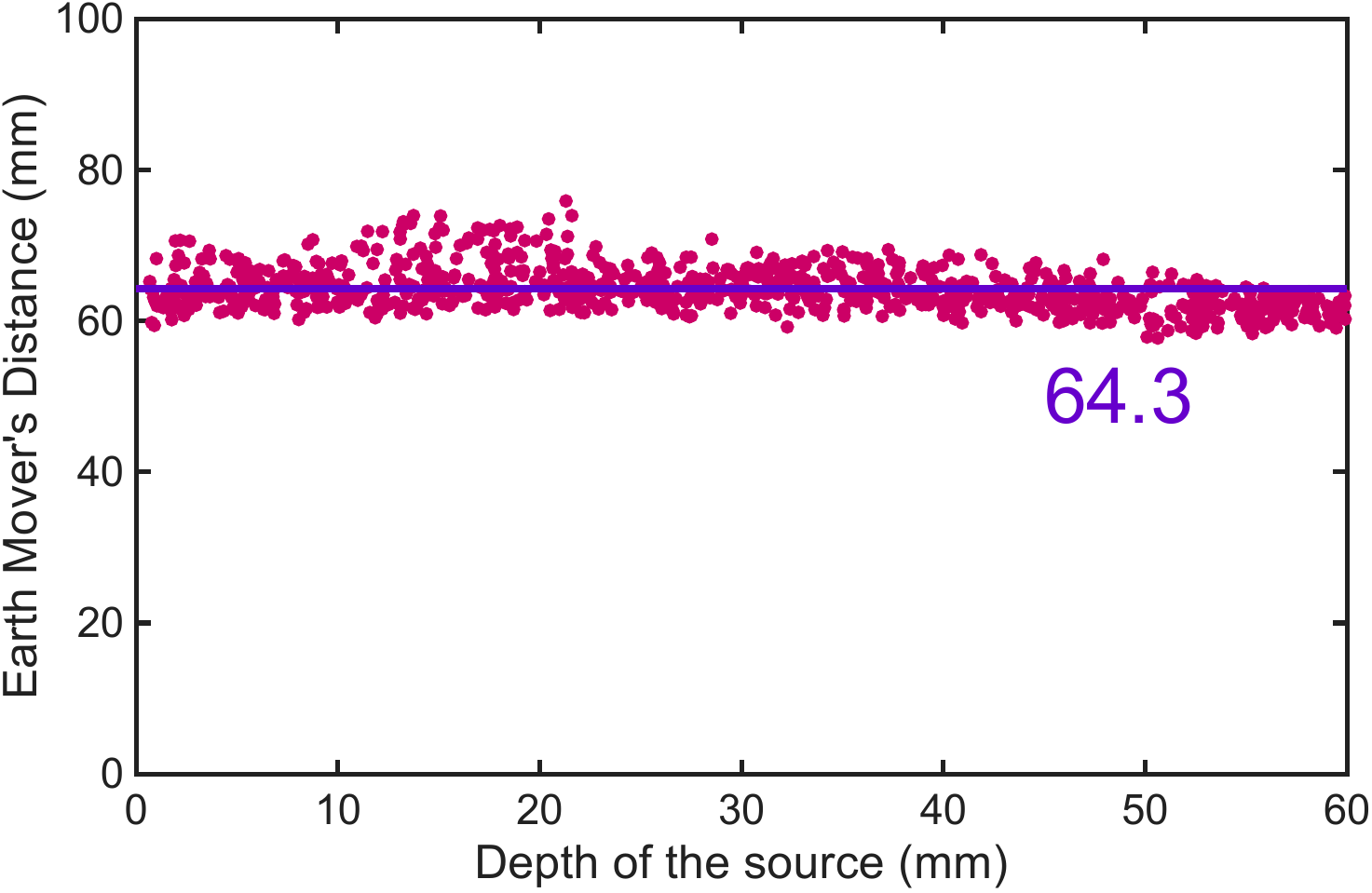}
    \end{minipage}\begin{minipage}{0.3\linewidth}
    \centering
        \includegraphics[width=0.98\linewidth]{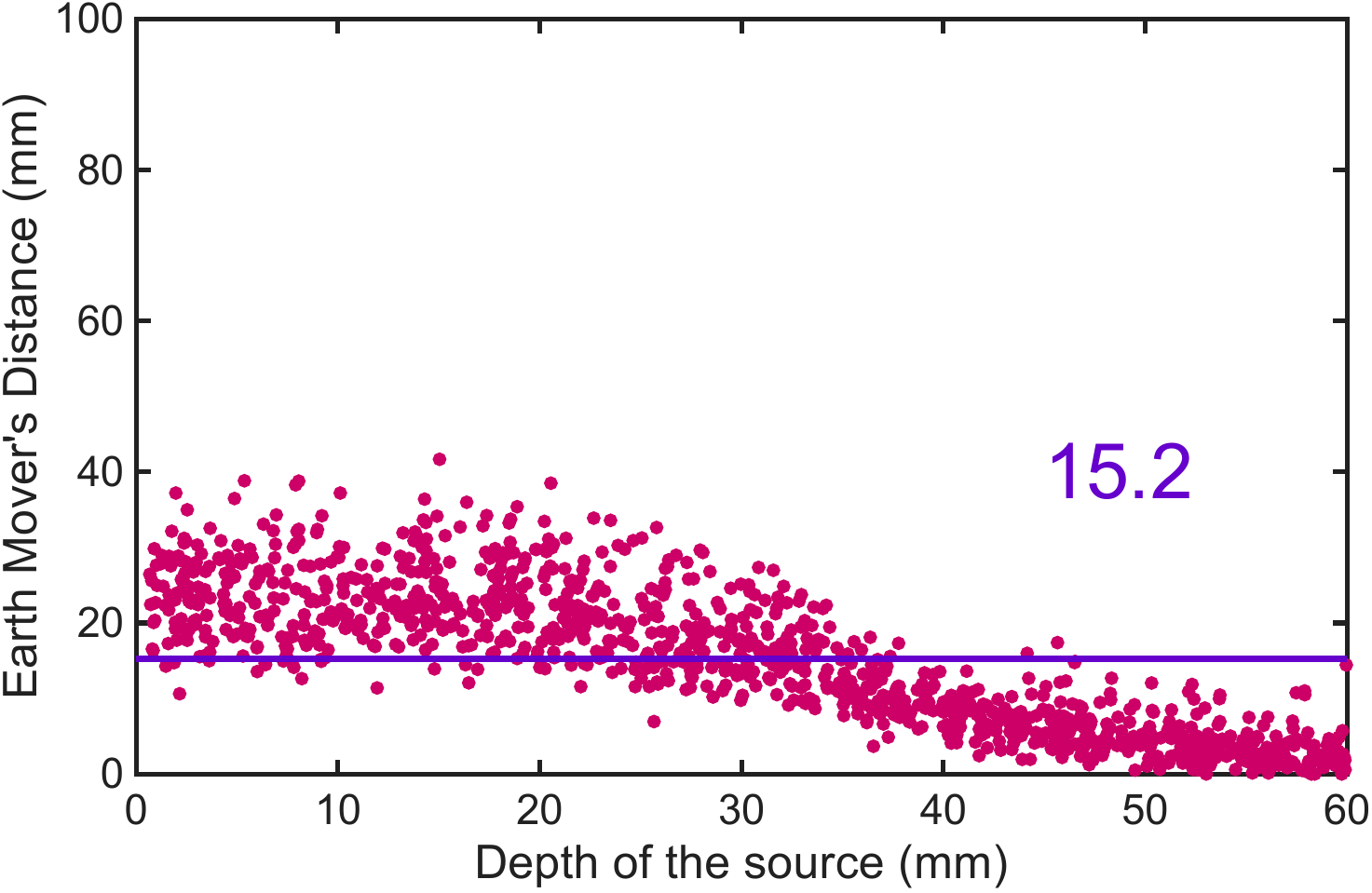}
    \end{minipage}\begin{minipage}{0.3\linewidth}
    \centering
        \includegraphics[width=0.98\linewidth]{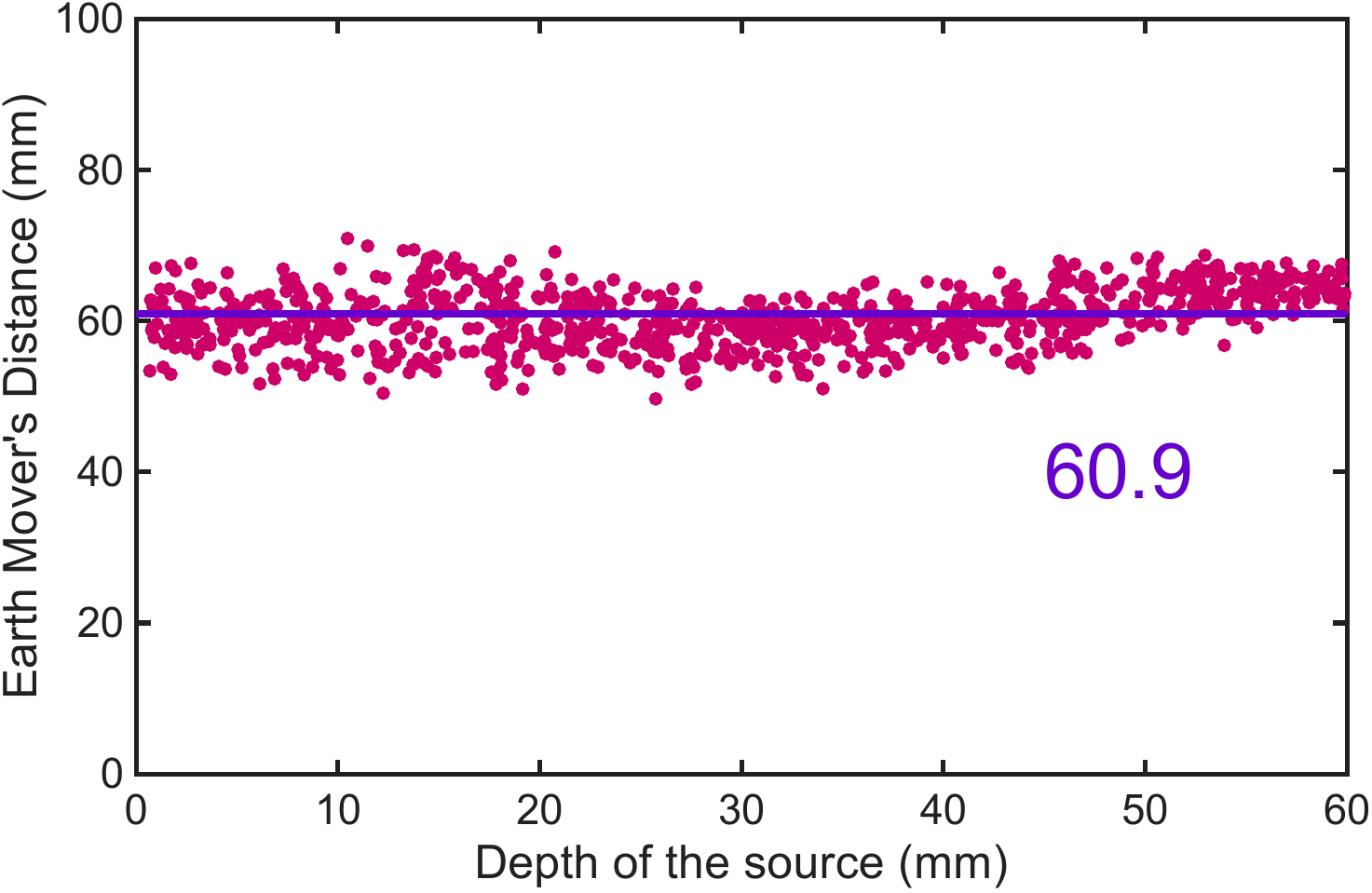}
    \end{minipage}

    \begin{minipage}{0.05\linewidth}
        \rotatebox{90}{Patch}
    \end{minipage}\begin{minipage}{0.3\linewidth}
    \centering
        \includegraphics[width=0.98\linewidth]{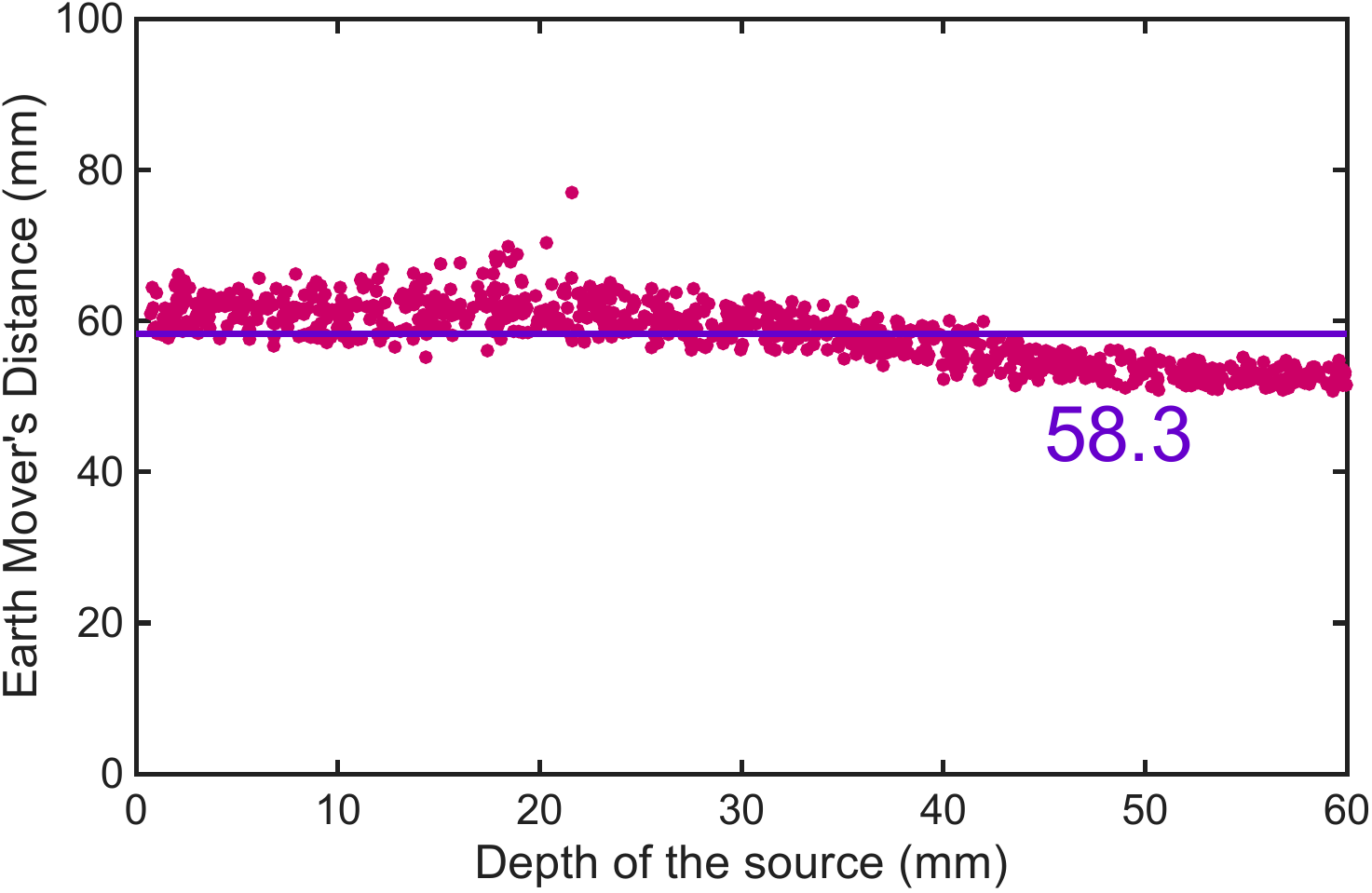}
    \end{minipage}\begin{minipage}{0.3\linewidth}
    \centering
        \includegraphics[width=0.98\linewidth]{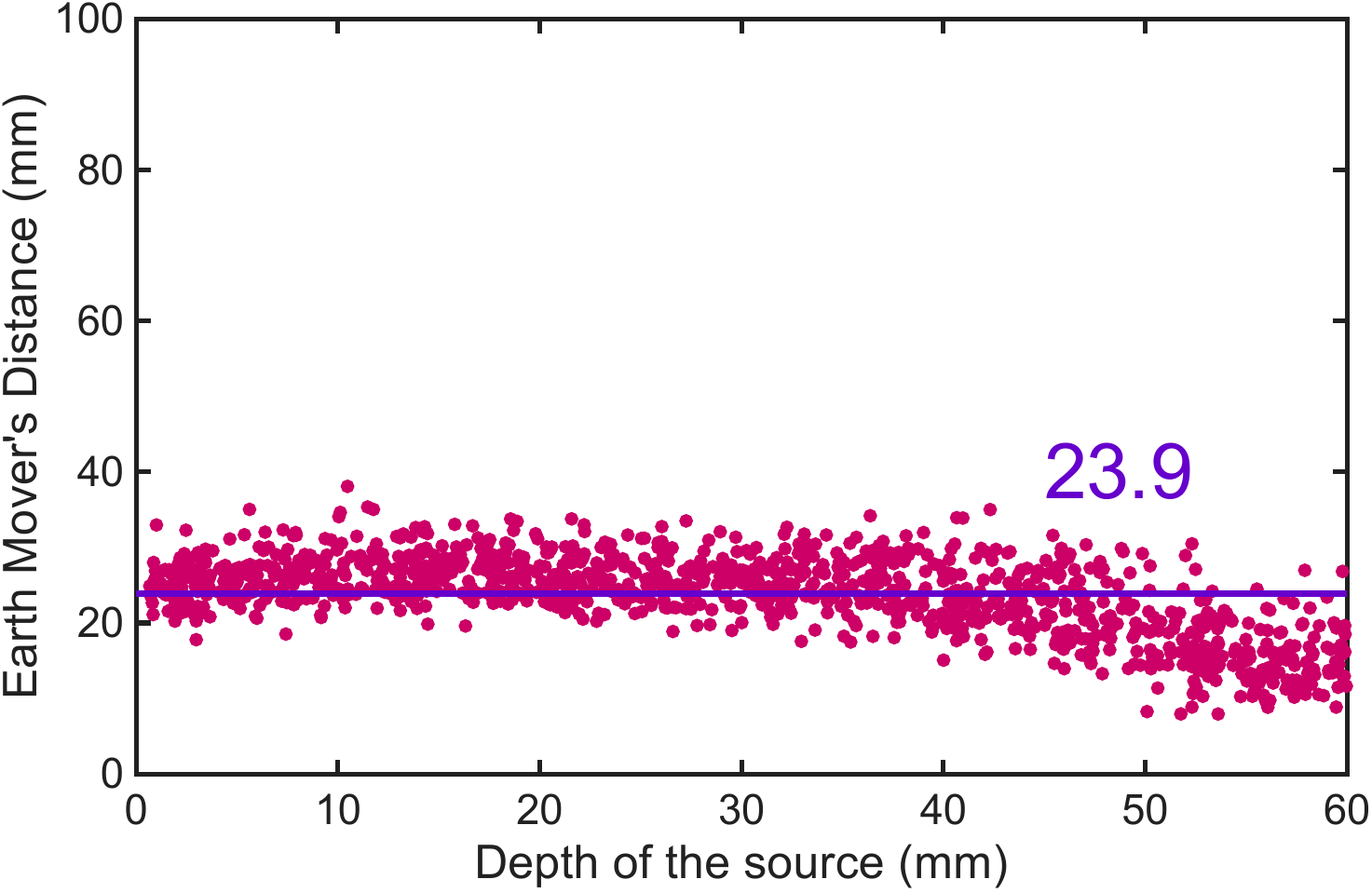}
    \end{minipage}\begin{minipage}{0.3\linewidth}
    \centering
        \includegraphics[width=0.98\linewidth]{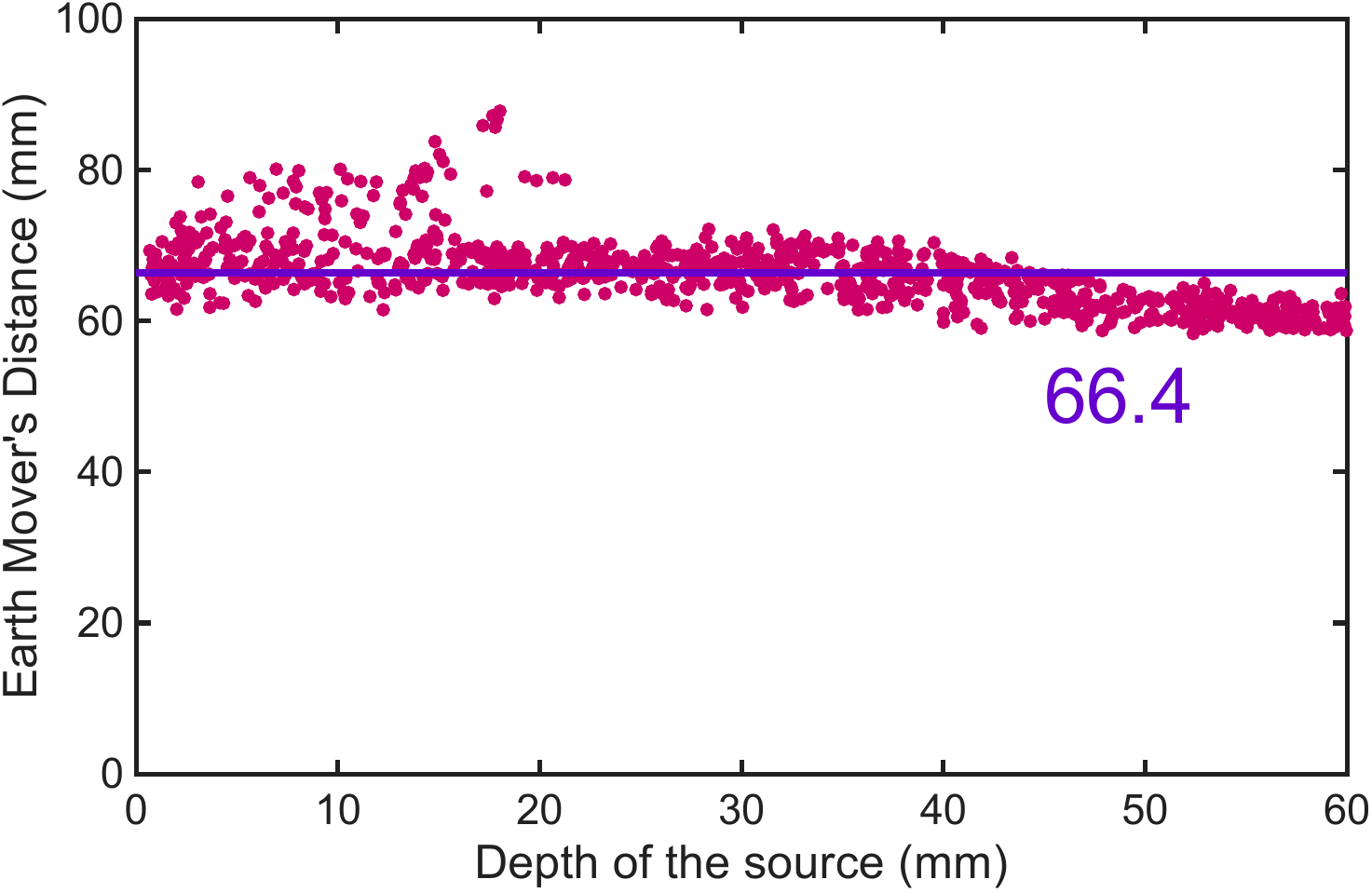}
    \end{minipage}

     \begin{minipage}{0.05\linewidth}
         \rotatebox{90}{Patch -- LS}
     \end{minipage}\begin{minipage}{0.3\linewidth}
     \centering
         \includegraphics[width=0.98\linewidth]{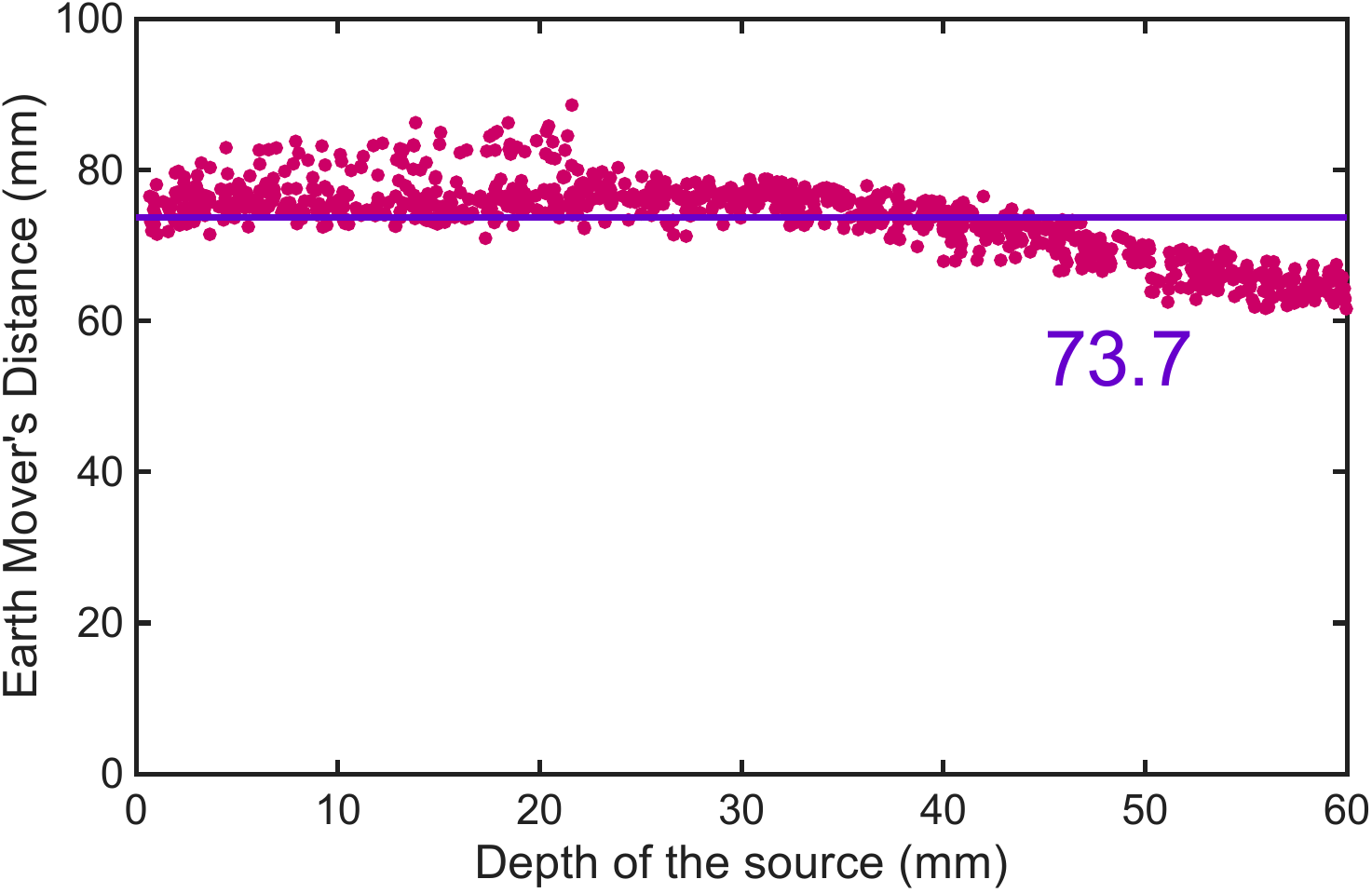}
     \end{minipage}\begin{minipage}{0.3\linewidth}
     \centering
         \includegraphics[width=0.98\linewidth]{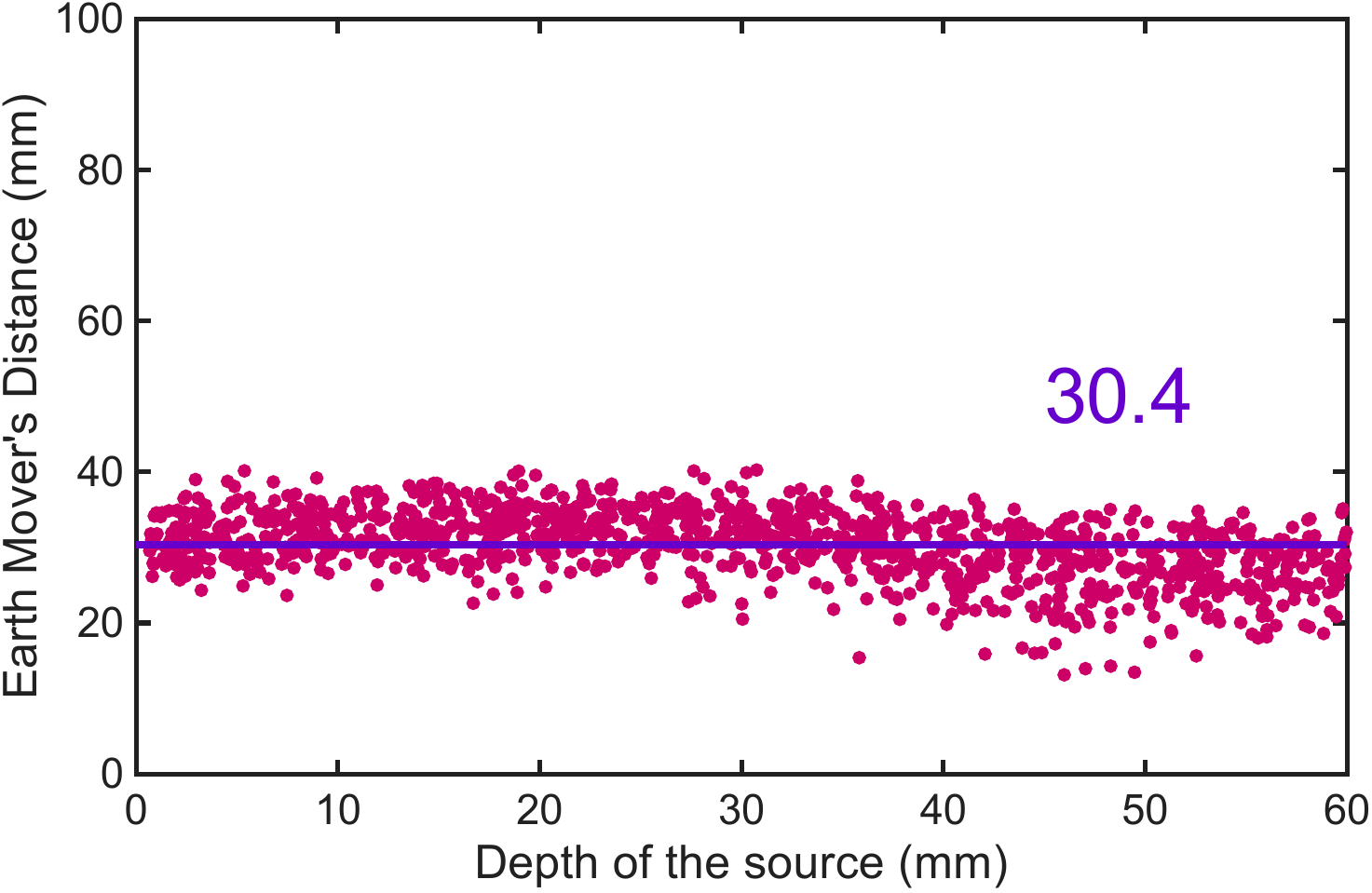}
     \end{minipage}\begin{minipage}{0.3\linewidth}
     \centering
         \includegraphics[width=0.98\linewidth]{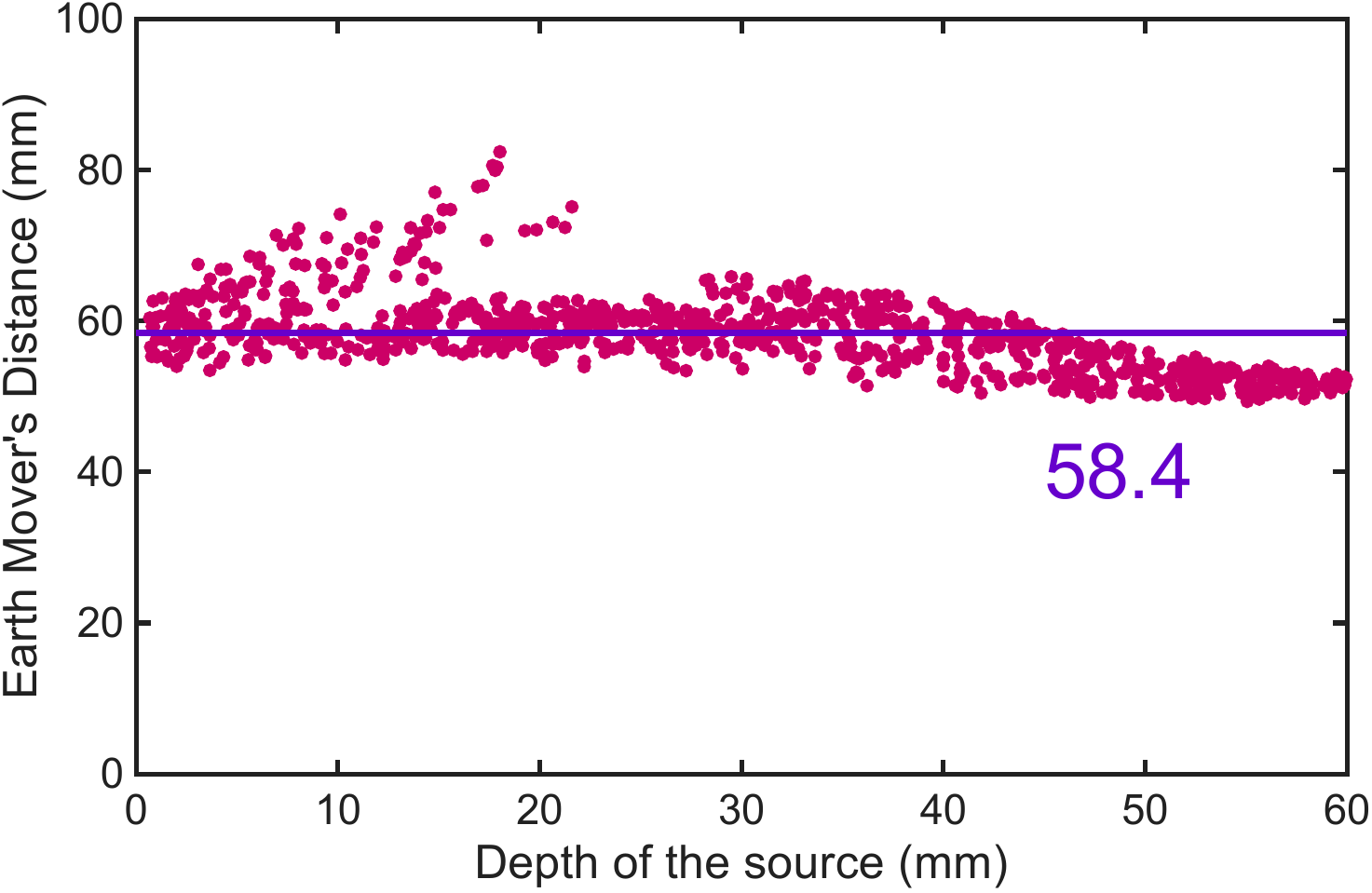}
     \end{minipage}
    
    \caption{Earth Mover's Distances for estimated sources when SNR is \qty{15}{\decibel}. The violet horizontal line indicates the sample mean, which is also shown numerically on the graph.}
    \label{fig.duneuro.EMD15dB}
\end{figure*}

To mirror the presentation of the localization depth bias in both the noiseless and \qty{15}{\decibel} SNR cases, Figure~\ref{fig.duneuro.EMD15dB} displays EMD's with different forward--inverse method pairs. The amount of noise added to the forward lead field or measurements is the same as in the case of DBSP.

In the case of EMD, a general trend with added noise is that the measure shows an average increase across all source depths for most forward--inverse combinations, while slightly increasing the vertical spread of the point clouds. Whitney and Local subtraction, in combination with SHAL1R, are especially affected, with shallower vertical values experiencing a jump of approximately \qty{15}{\milli\meter}. For Whitney-basis, the average is \qty{14.9}{\milli\meter}, and for Local subtraction, the mean is \qty{15.2}{\milli\meter}. Both Whitney and Local subtraction also show a downward trend in EMD as depth increases, approaching \num 0 at the greatest depth. There is again a notable exception to the rule however, just like there were in the case of the DBSP of Figure~\ref{fig.duneuro.depthbias15dB}: the combination of the patch lead field and SHAL1R sees a drastic decline in EMD across all depths, decreasing from approximately \qtyrange{70}{80}{\milli\meter} to \qtyrange{20}{35}{\milli\meter}, while the mean EMD is \qty{23.9}{\milli\meter}. The vertical spread of the EMD values also clearly decreases for cases such as $\Hdiv$ coupled with SKF, but with an increase in the lower bound of roughly \qty{30}{\milli\meter}, ignoring one outlier value near the \qty{20}{\milli\meter} mark at around \qty{15}{\milli\meter} depth in the noiseless case. Otherwise, the results of SKF show that the overall trends are converging towards a flat line compared to the noiseless case. The same can be observed for sLORETA. The sLORETA column is also mostly unaffected by the added noise, discounting the case of $\Hdiv$ + sLORETA, whic goes from a roughly flat \qty{55}{\milli\meter} distribution to a spread where the shallower values start at approximately \qty{80}{\milli\meter} and decrease towards the \qty{60}{\milli\meter} mark at the greatest depth. With sLORETA, a clear decrease of the EMD mean is observed except when patch sources are inverted with Local subtraction (mean of \qty{73.7}{\milli\meter}) causes almost no change, and the mean is increased to \qty{71.9}{\milli\meter} with $\Hdiv$. There is no consistent indication of which source model yields the most concise estimates based on the average EMD: the patch model is the optimal choice for sLORETA (\qty{58.3}{\milli\meter}), Whitney-basis for SHAL1R (\qty{14.9}{\milli\meter}), and Local subtraction for SKF (\qty{60.9}{\milli\meter}). However, it should be noted that with SHAL1R and SKF, the difference between the averages for Whitney-basis and Local subtraction is within fractions of millimeters.

\section{Discussion}\label{sec.discussion}

In this paper, we studied how the choice of the source model can affect the quality of single-source estimation, both quantitatively and qualitatively. The study was conducted in the finite element method settings. Two software implementations dedicated to the neuroelectromagnetic forward problem, DUNEuro \cite{schrader2021duneuro}
and Zeffiro Interface \cite{he2020zeffiro} were used with a software-package-connecting pipeline that enables a fair comparison. The compared source models are face-intersecting and edgewise Whitney basis functions and Local subtraction \cite{holtershinken-etal-2025} provided by DUNEuro, and $\Hdiv$ implemented in Zeffiro Interface. The visual smoothness and spread of an estimated cortical source were examined, along with depth estimates and Earth Mover's Distance (EMD), for sources at various depths measured from the upper hemisphere of the inner skull layer, with inversion methods aimed at reducing depth bias.

    Zeffiro Interface \cite{he2020zeffiro} and DUNEuro \cite{schrader2021duneuro} were chosen as the lead field construction utilities for this study due to their demonstrated ability to construct accurate FEM-based EEG forward mappings, which is a requirement for accurate volumetric modelling of tissue conductivity. Boundary-element-method-based solvers such as the one found in OpenMEEG \cite{gramfort2010openmeeg} would have been less suitable for this task, as they only accurately model tissues at the tissue type boundaries. We also chose to interface directly with these tools instead of using sofware such as Brainstorm \cite{tadel2011brainstorm} as an unnecessary intermediary interface.
    The lead fields produced by DUNEuro match well with the expected outcome of a reciprocally produced lead field \cite{rush1969reciprocity,weinstein2000reciprocity,vanrumste2004reciprocity}, at least when an active-element-centroid-based source positioning scheme is used in the transfer matrix interpolation step: the field is strong near the superficial electrodes and shows a smoothly diminishing sensitivity pattern when moving towards deeper compartments. The $\Hdiv$ implementation of Zeffiro does not match with the expectation as clearly, as its field is more spread out the entire cortical surface, resembling the smoothed patch-like lead field constructed with DUNEuro. The $\Hdiv$ lead field also contained a few drastic outliers in field strength that needed to be filtered out of the final solution using a \qty{95}{\percent} quantile cap on the field strength, to gain access to the underlying lead field sensitivity pattern.
    Yet another reason for the use of Zeffiro Interface was its MATLAB namespace \texttt{+inverse}, which contains implementations of the utilized inverse methods.

We explored the spread and smoothness of cortical source reconstructions with sLORETA, SHAL1R, SKF, and DS; and measured the depth bias and EMD of the standardized methodologies sLORETA, SHAL1R, and SKF. The depth bias experiment was conducted with both noiseless and noisy synthetic data of \qty{15}{\decibel} SNR.
    Based on the results, while the characteristics of inversion methods vary within a given source model, the Whitney basis and Local subtraction yield almost identical behavior for a given inverse method.

sLORETA and SHAL1R exhibit a nearly biasless trend in the noiseless case, whereas the same is true only for the Whitney basis and Local subtraction when noise is added. SHAL1R is an alternative to sLORETA that explicitly assumes a sparse activity distribution \cite{lahtinen2024shalpr}, which is appropriate since the activity originates from a point source.

The surrogate patch-like source model created via Local subtraction produced unbiased results with both inversion approaches. The overall largest bias is obtained with $\Hdiv$ in the presence of noise. A notable exception for highly non-biased estimations of sLORETA and SHAL1R is observed with SKF. Given that SKF (the spatiotemporally standardized, post-weighted Kalman filter) is a recursive estimator that relies on prior state information and the given evolution model, the parametrization of the evolution model may affect the biasing effect. Moreover, the presence of strong inter-source correlations in the forward model --- particularly with $\Hdiv$ and the patch model --- noticeably influenced the outcomes. These interdependencies arise by design in the patch model due to Gau\ss ian smoothing and appear inherently in the $\Hdiv$ framework. Notably, the SKF approach exhibited a notable deviation in depth estimates at approximately \qty{4}{\centi\meter}, a feature absent from the other standardized methods. This deviation can be attributed to the mismatch between the model's prior and evolution assumptions and the true dipolar point source: SKF is tailored for spatially distributed activity, so it underestimates isolated deep sources, treating their weak signals as noise. When the ground truth was adjusted to match the patch model and inversion was performed using local subtraction, the estimation outcomes mirrored those from the patch model itself. This suggests that SKF may be less effective at resolving patch-type sources without introducing depth bias. In contrast, sLORETA and SHAL1R behaved more consistently: sLORETA achieved higher accuracy for patch-like sources, while SHAL1R, which enforces spatial sparsity, demonstrated slightly reduced performance.

The uniqueness of the surface potentials caused by a single source is the requirement for unbiased estimation \cite{lahtinen-standardization-2024}. Due to the high agreement of depths of the source and its estimate, in addition to the nearly zero EMD for SHAL1R when Local subtraction is used, indicates that the Local subtraction provides the best opportunity for unbiased estimates across the source models compared. This property is most likely realistic as reflected in studies on reconstructing deep brain activity from human EEG recordings; EEG's capability to distinguish deep sources is evident, but the main challenge is in recovering those weak signals from noisy data \cite{KrishnaswamyPavitra2017Seeo,Fahimi2020ECoGvdEEGdeeppossible,Rezaei2021}. Therefore, the ambiguity should come from the modelled noise, not from the lead field in numerical settings.

In contrast to the previous case, the worst EMD for SHAL1R is obtained with the patch model and $\Hdiv$. The cortical reconstructions with SHAL1R and DS indicate that this source model yields a region of nearly identical sources at the cortical level; hence, the difficulty of localizing the source is greater than with the Whitney basis and Local subtraction. High EMD value for both sLORETA and SHAL1R indicates high ambiguity in source location at every depth. Moreover, the results show that the negative total impact of the noise was largest with $\Hdiv$ and the patch model, highlighting the suboptimal source distinguishability.

    The anomaly of decreasing EMD for deep sources observed mainly with Whitney basis and Local subtraction could be explained by local perturbations of the sub-lead fields in deep regions. If we draw a conclusion about how the lead field should behave by analyzing analytically solvable spherical conductivity models, three principles can be drawn: (1) the relative perturbation can never be zero; (2) the variety of relative difference values should shrink along with depth; (3) the spatial distance between compared sub-lead fields should correlate positively with the variability. In this regard, we can conclude that all the methods have problems with short-scale variability, as the values are mostly effectively zero, probably due to cumulative floating-point errors. We could also conclude that the Whitney basis can not quite hit the mark for principles (2) and (3) due to its heavy tail for higher relative differences at depths \qtyrange{0}{30}{\milli\meter}. In contrast, Local subtraction follows principles (2) and (3). For both of these methods provided by DUNeuro, the behaviour at high depth coincides with the expectations drawn from the local properties of the analytical solution in a sphere. This sudden shift in local behaviour could explain the value-dropping phenomena witnessed with EMD.

Comparing the performances of sLORETA and SHAL1R in reconstructing patch-modelled sources, we can see the possible advantage of a hypothetical source model that maps point sources to patch sources. Based on the results, the accuracy of sparsity-promoting solvers to estimate point sources could be leveraged to estimate patch-type
sources. Considering that some accuracy is inevitably lost when we estimate patch sources using a point-source model.

\section{Limitations of the study}\label{sec.limitations}
The numerical analysis of this study is limited to single-source cases. On the other hand, given the study's focus, the controllability of the one-source scenario and the theoretical background on unbiased estimation allow us to design meaningful test scenarios that reveal the source-model distinctions highlighted in the paper. Furthermore, standardized methods as well as many filtering methods, such as DS, are designed to recover a single source, posing unavoidable challenges for potential multisource analysis.

The other limitation concerns the indifference of the inversion results regarding the Whitney basis and Local subtraction: The high-density mesh used for this study makes the FEM solution highly accurate. Thus, the observed differences are due to differences in convergent forward solutions. Considering the fact that computing a FEM solution on a highly dense mesh is time-consuming, it would have been valuable information to know how mesh density affects the Whitney basis and Local subtraction, which could then have shown a difference between the said source models. If either of these methods were less dependent on mesh density, it would have solidified its superiority. Although the final comparison of the source models should be explored using multiple clinically measured EEG recordings to assess the realism of the point-source assumption and the modeling of deep sources. Addressing the shortcomings mentioned will streamline future work.

\section{Conclusions}\label{sec.conclusions}

    Of the inverse methods utilized in this study, SHAL1R and DS presuppose a very focal source. The best results for these methods are obtained with Local subtraction; thus, it is deemed to provide the best model for point sources among those seen in this study, in terms of estimation spread, DBSP, and EMD.

    In the case of a patch-like source, sLORETA appears to yield less biased results, as indicated by the depth-bias scatter plots and generally lower EMD. Furthermore, adding noise does not increase the method's bias, though it does narrow the performance gap relative to other source models, such as the Whitney basis. Maintaining a consistent bias level appears to be a characteristic of SKF as well, however, its overall results are leaning towards Local subtraction and Whitney basis.

    The results indicate clearly that Local subtraction, as a point source model, has a slightly lower ability to estimate patch sources than point sources. Instead, a patch-like model for inversion is favorable when the true source is patch-distributed.

\section{Author contributions}\label{sec.author.contributions}

Santtu Söderholm produced the original focal ECD EEG lead fields $\leadFieldMatrix$ using DUNEuro and Zeffiro Interface and was one of the main writers of the work. Joonas Lahtinen produced patch-like lead fields and inverse reconstructions of brain activity using the various inverse methods mentioned in this work and also wrote and improved many parts of the text. Sampsa Pursiainen provided the ICBM 2009a head model and the Zeffiro $\Hdiv$ lead field routines, participated in writing the introduction and abstract, and worked in a supervisory role throughout the writing process.

\section{Copyright notice}\label{sec.copyright}

This work is an extended version of the publication \cite{soderholm2026forwardinverse}, copyright \copyright\ IEEE 2026. Reused results from the publication concerning lead field norms are contained in \Cref{fig.duneuro.L.vecnorm,fig.zeffiro.L.vecnorm}. Regarding investigated source models and the inverse methods sLORETA, SHAL1R, and DS, reused results are contained in \Cref{fig.erotuskuva,fig.duneuro.depthbias,fig.duneuro.EMD}. Permission to publish this extended version has been granted by IEEE.

\printcredits


\interlinepenalty=10000

\printbibliography

@article{
    holtershinken-etal-2025,
    author = {Höltershinken, Malte B. and Lange, Pia and Erdbrügger, Tim and Buschermöhle, Yvonne and Wallois, Fabrice and Buyx, Alena and Pursiainen, Sampsa and Vorwerk, Johannes and Engwer, Christian and Wolters, Carsten H.},
    title = {{The Local Subtraction Approach for EEG and MEG Forward Modeling}},
    journal = {SIAM Journal on Scientific Computing},
    volume = {47},
    number = {1},
    pages = {B160-B189},
    year = {2025},
    doi = {10.1137/23M1582874},
}

@article{
    neugebauer-etal-2022,
    author = {Neugebauer, Frank and Antonakakis, Marios and Unnwongse, Kanjana and Parpaley, Yaroslav and Wellmer, Jörg and Rampp, Stefan and Wolters, Carsten H.},
    title = {{Validating EEG, MEG and Combined MEG and EEG Beamforming for an Estimation of the Epileptogenic Zone in Focal Cortical Dysplasia}},
    journal = {Brain Sciences},
    volume = {12},
    year = {2022},
    number = {1},
    article-number = {114},
    pubmedid = {35053857},
    issn = {2076-3425},
    doi = {10.3390/brainsci12010114},
}

@book{
    Knosche-Haueisen-2022,
	author = {Thomas R. Knösche and Jens Haueisen},
	title = {{{EEG/MEG Source Reconstruction}}},
	subtitle = {{Texbook for Electro- and Magnetoencephalography}},
	publisher = {Springer},
	year = {2022},
	isbn = {978-3-030-74918-7},
}

@article{
    lahtinen-standardization-2024,
    doi = {10.1088/1361-6420/ad5f53},
    year = {2024},
    month = {7},
    publisher = {IOP Publishing},
    volume = {40},
    number = {9},
    pages = {095002},
    author = {Lahtinen, Joonas},
    title = {{On bias and its reduction via standardization in discretized electromagnetic source localization problems}},
    journal = {Inverse Problems},
}

@article{
    Pursiainen-etal-2016,
    doi = {10.1088/0031-9155/61/24/8502},
    year = {2016},
    month = {nov},
    publisher = {IOP Publishing},
    volume = {61},
    number = {24},
    pages = {8502},
    author = {Pursiainen, S and Vorwerk, J and Wolters, C H},
    title = {{Electroencephalography (EEG) forward modeling via H(div) finite element sources with focal interpolation}},
    journal = {Physics in Medicine \& Biology},
}

@article{
    wolters2004transfer,
    author = {Wolters, Carsten and Grasedyck, Lars and Hackbusch, Wolfgang},
    year = {2004},
    month = {08},
    pages = {1099-1116},
    title = {{Efficient Computation of lead field bases and influence matrix for the FEM-based EEG and MEG inverse problem}},
    volume = {20},
    journal = {Inverse Problems},
    doi = {10.1088/0266-5611/20/4/007},
}

@article{
    wolters-etal-2007,
    author = {Wolters, C. H. and Köstler, H. and Möller, C. and Härdtlein, J. and Grasedyck, L. and Hackbusch, W.},
    title = {{Numerical Mathematics of the Subtraction Method for the Modeling of a Current Dipole in EEG Source Reconstruction Using Finite Element Head Models}},
    journal = {SIAM Journal on Scientific Computing},
    volume = {30},
    number = {1},
    pages = {24-45},
    year = {2008},
    doi = {10.1137/060659053},
}

@book{
    sekihara-etal-2008,
    author = {Kensuke Sekihara and Srikantan S Nagarajan},
    title = {{Adaptive Spatial Filters for Electromagnetic Brain Imaging}},
    publisher = {Springer},
    year = {2008},
    isbn = {978-3-540-79369-4},
    issn = {1864-5763},
}

@article{hamalainen1993magnetoencephalography,
  author  = {Matti S. Hämäläinen and Riitta Hari and Risto J. Ilmoniemi and Jukka Knuutila and Olli V. Lounasmaa},
  title = {{Magnetoencephalography---theory, instrumentation, and applications to noninvasive studies of the working human brain}},
  journal = {Reviews of Modern Physics},
  volume  = {65},
  number  = {2},
  pages   = {413--497},
  year    = {1993}
}

@article{pascual2002sloreta,
  author  = {Roberto D. Pascual-Marqui},
  title = {{Standardized low-resolution brain electromagnetic tomography (sLORETA)}},
  journal = {Methods and Findings in Experimental and Clinical Pharmacology},
  volume  = {24},
  pages   = {5--12},
  year    = {2002}
}

@article{miinalainen2019realistic,
  author  = {Tuuli Miinalainen and Atena Rezaei and Defne Us and Andreas Nü{\ss}ing and Christian Engwer and Carsten H. Wolters and Sampsa Pursiainen},
  title = {{A realistic, accurate and fast source modeling approach for the EEG forward problem}},
  journal = {NeuroImage},
  volume  = {184},
  pages   = {56--67},
  year    = {2019}
}

@article{he2020zeffiro,
  author  = {Q. He and A. Rezaei and S. Pursiainen},
  title = {{Zeffiro interface for electromagnetic brain imaging}},
  journal = {Neuroinformatics},
  volume  = {18},
  number  = {2},
  pages   = {237--250},
  year    = {2020}
}

@article{schrader2021duneuro,
  author  = {S. Schrader and A. Nü{\ss}ing and C. Engwer and J. Vorwerk and C. H. Wolters},
  title = {{DUNEuro: a software toolbox for forward modeling in bioelectromagnetism}},
  journal = {PLOS ONE},
  volume  = {16},
  number  = {6},
  pages   = {e0252431},
  year    = {2021}
}

@article{lahtinen2024shalpr,
  author  = {Joonas Lahtinen and Alexandra Koulouri and Stefan Rampp and Jörg Wellmer and Carsten Wolters and Sampsa Pursiainen},
  title = {{Standardized hierarchical adaptive Lp regression for noise robust focal epilepsy source reconstructions}},
  journal = {Clinical Neurophysiology},
  volume  = {159},
  pages   = {24--40},
  year    = {2024}
}

@article{holtershinken2025local,
  author  = {Malte B. Höltershinken and Pia Lange and Tim Erdbrügger and Yvonne Buschermöhle and Fabrice Wallois and Alena Buyx and Sampsa Pursiainen and Johannes Vorwerk and Christian Engwer and Carsten H. Wolters},
  title = {{The local subtraction approach for EEG and MEG forward modeling}},
  journal = {SIAM Journal on Scientific Computing},
  volume  = {47},
  number  = {1},
  pages   = {B160--B189},
  year    = {2025}
}

@article{soderholm2024effects,
  title = {{The effects of peeling on finite element method-based EEG source reconstruction}},
  author={Söderholm, Santtu and Lahtinen, Joonas and Wolters, Carsten H and Pursiainen, Sampsa},
  journal={Biomedical Signal Processing and Control},
  volume={89},
  pages={105695},
  year={2024},
  publisher={Elsevier}
}

@article{cuffin2001realistic,
  author = {B. N. Cuffin and D. L. Schomer and J. R. Ives and H. Blume},
  title = {{Experimental tests of EEG source localization accuracy in realistically shaped head models}},
  journal = {Clinical Neurophysiology},
  volume = {112},
  number = {12},
  pages = {2288--2292},
  year = {2001}
}

@article{cuffin2001spherical,
  author = {B. N. Cuffin and D. L. Schomer and J. R. Ives and H. Blume},
  title = {{Experimental tests of EEG source localization accuracy in spherical head models}},
  journal = {Clinical Neurophysiology},
  volume = {112},
  number = {1},
  pages = {46--51},
  year = {2001}
}

@article{lahtinen2024standardized,
  title = {{Standardized Kalman filtering for dynamical source localization of concurrent subcortical and cortical brain activity}},
  author={Lahtinen, Joonas and Ronni, Paavo and Subramaniyam, Narayan Puthanmadam and Koulouri, Alexandra and Wolters, Carsten and Pursiainen, Sampsa},
  journal={Clinical Neurophysiology},
  volume={168},
  pages={15--24},
  year={2024},
  publisher={Elsevier}
}

@article{
    fuchs-1998,
    title = {{Improving source reconstructions by combining bioelectric and biomagnetic data}},
    journal = {Electroencephalography and Clinical Neurophysiology},
    volume = {107},
    number = {2},
    pages = {93-111},
    year = {1998},
    issn = {0013-4694},
    doi = {https://doi.org/10.1016/S0013-4694(98)00046-7},
    author = {Manfred Fuchs and Michael Wagner and Hans-Aloys Wischmann and Thomas Köhler and Annette Theißen and Ralf Drenckhahn and Helmut Buchner},
}

@INPROCEEDINGS{
    rubner-etal-1998-emd,
    author={Rubner, Y. and Tomasi, C. and Guibas, L.J.},
    booktitle = {{Sixth International Conference on Computer Vision (IEEE Cat. No.98CH36271)}},
    title = {{A metric for distributions with applications to image databases}},
    year={1998},
    volume={},
    number={},
    pages={59-66},
    doi={10.1109/ICCV.1998.710701}
}

@article{
    fonov-2011,
    title = {{Unbiased average age-appropriate atlases for pediatric studies}},
    journal = {NeuroImage},
    volume = {54},
    number = {1},
    pages = {313-327},
    year = {2011},
    issn = {1053-8119},
    doi = {10.1016/j.neuroimage.2010.07.033},
    author = {Vladimir Fonov and Alan C. Evans and Kelly Botteron and C. Robert Almli and Robert C. McKinstry and D. Louis Collins}
}

@article{
    fonov-2009,
    title = {{Unbiased nonlinear average age-appropriate brain templates from birth to adulthood}},
    journal = {NeuroImage},
    volume = {47},
    pages = {S102},
    year = {2009},
    note = {Organization for Human Brain Mapping 2009 Annual Meeting},
    issn = {1053-8119},
    doi = {https://doi.org/10.1016/S1053-8119(09)70884-5},
    author = {VS Fonov and AC Evans and RC McKinstry and CR Almli and DL Collins}
}

@online {
    icbm-website,
    author = {{NeuroImaging and Surgical Technologies Lab}},
    title = {{ICBM 152 Nonlinear atlases}},
    date = {2009},
    url = {https://nist.mni.mcgill.ca/icbm-152-nonlinear-atlases-2009/},
    urldate = {2026-04-16},
}

@article{
    vergani-2024,
    author = {Vergani, Alberto Arturo},
    title = {{Hans Berger (1873–1941): the German psychiatrist who recorded the first electrical brain signal in humans 100 years ago}},
    journal = {Advances in Physiology Education},
    volume = {48},
    number = {4},
    pages = {878-881},
    year = {2024},
    doi = {10.1152/advan.00119.2024},
    note ={PMID: 39236103}
}

@article{
  missey-etal-2026,
  title = {{Non-invasive temporal interference stimulation of the hippocampus suppresses epileptic biomarkers in patients with Epilepsy: biophysical differences between kilohertz and amplitude modulated stimulation}},
  journal = {Brain Stimulation},
  volume = {19},
  number = {1},
  pages = {102981},
  year = {2026},
  issn = {1935-861X},
  doi = {10.1016/j.brs.2025.11.008},
  author = {Florian Missey and Emma Acerbo and Adam S. Dickey and Jan Trajlinek and Ondřej Studnička and Claudia Lubrano and Mariane {de Araújo e Silva} and Evan Brady and Vit Všianský and Johanna Szabo and Irena Dolezalova and Daniel Fabo and Martin Pail and Claire-Anne Gutekunst and Rosanna Migliore and Michele Migliore and Stanislas Lagarde and Romain Carron and Fariba Karimi and Raul Castillo Astorga and Antonino M. Cassara and Niels Kuster and Esra Neufeld and Fabrice Bartolomei and Nigel P. Pedersen and Robert E. Gross and Viktor Jirsa and Daniel L. Drane and Milan Brázdil and Adam Williamson},
}

@article{
    chaddad-etal-2023,
    title = {{Electroencephalography Signal Processing: A Comprehensive Review and Analysis of Methods and Techniques.}},
    author = {Chaddad, Ahmad and Wu, Yihang and Kateb, Reem and Bouridane, Ahmed},
    journal = {Sensors (Basel, Switzerland)},
    year = {2023},
    volume = {23},
    number = {14},
    doi = {10.3390/s23146434},
    issn = {1424-8220 (Electronic)},
    pmid = {37514728}
}

@book{
    hari-and-puce-2023,
    author = {Hari, Riitta and Puce, Aina},
    title = {{MEG - EEG Primer}},
    publisher = {Oxford University Press},
    year = {2023},
    month = {09},
    isbn = {9780197542187},
    doi = {10.1093/med/9780197542187.001.0001},
}

@article{
    clark-and-plonsey-1966,
    title = {{A Mathematical Evaluation of the Core Conductor Model}},
    journal = {Biophysical Journal},
    volume = {6},
    number = {1},
    pages = {95-112},
    year = {1966},
    issn = {0006-3495},
    doi = {https://doi.org/10.1016/S0006-3495(66)86642-0},
    author = {John Clark and Robert Plonsey}
}

@article{
    scherg-1985,
    title = {{Two bilateral sources of the late AEP as identified by a spatio-temporal dipole model}},
    journal = {Electroencephalography and Clinical Neurophysiology/Evoked Potentials Section},
    volume = {62},
    number = {1},
    pages = {32-44},
    year = {1985},
    issn = {0168-5597},
    doi = {10.1016/0168-5597(85)90033-4},
    author = {M Scherg and D {Von Cramon}}
}

@ARTICLE{
    ary-etal-1981,
    author={Ary, James P. and Klein, Stanley A. and Fender, Derek H.},
    journal={IEEE Transactions on Biomedical Engineering},
    title = {{Location of Sources of Evoked Scalp Potentials: Corrections for Skull and Scalp Thicknesses}},
    year={1981},
    volume={BME-28},
    number={6},
    pages={447-452},
    doi={10.1109/TBME.1981.324817}
}

@ARTICLE{
    demunck-1992,
    author={de Munck, J.C.},
    journal={IEEE Transactions on Biomedical Engineering},
    title = {{A linear discretization of the volume conductor boundary integral equation using analytically integrated elements (electrophysiology application)}},
    year={1992},
    volume={39},
    number={9},
    pages={986-990},
    doi={10.1109/10.256433}
}

@ARTICLE{
    schlitt-1995,
    author={Schlitt, H.A. and Heller, L. and Aaron, R. and Best, E. and Ranken, D.M.},
    journal={IEEE Transactions on Biomedical Engineering},
    title = {{Evaluation of boundary element methods for the EEG forward problem: effect of linear interpolation}},
    year={1995},
    volume={42},
    number={1},
    pages={52-58},
    doi={10.1109/10.362919}
}

@article{
    wolters-etal-2006,
    title = {{Influence of tissue conductivity anisotropy on EEG/MEG field and return current computation in a realistic head model: A simulation and visualization study using high-resolution finite element modeling}},
    journal = {NeuroImage},
    volume = {30},
    number = {3},
    pages = {813-826},
    year = {2006},
    issn = {1053-8119},
    doi = {10.1016/j.neuroimage.2005.10.014},
    author = {C.H. Wolters and A. Anwander and X. Tricoche and D. Weinstein and M.A. Koch and R.S. MacLeod}
}

@book{kaipio2006statistical,
	title = {{Statistical and computational inverse problems}},
	author = {Kaipio, Jari and Somersalo, Erkki},
	volume = {160},
	year = {2006},
	publisher = {Springer Science \& Business Media},
}

@article{
    hanke-2011,
    author = {Hanke, Martin and Harrach, Bastian and Hyvönen, Nuutti},
    title = {{Justification of Point Electrode Models in Electrical Impedance Tomography}},
    journal = {Mathematical Models and Methods in Applied Sciences},
    volume = {21},
    number = {06},
    pages = {1395-1413},
    year = {2011},
    doi = {10.1142/S0218202511005362}
}

@ARTICLE{
    agsten2018electrodes,
    author={Pursiainen, Sampsa and Agsten, Britte and Wagner, Sven and Wolters, Carsten H.},
    journal={IEEE Transactions on Neural Systems and Rehabilitation Engineering},
    title = {{Advanced Boundary Electrode Modeling for tES and Parallel tES/EEG}},
    year={2018},
    volume={26},
    number={1},
    pages={37-44},
    doi={10.1109/TNSRE.2017.2748930}
}

@article{
    gross2023reciprocity,
    title = {{Bioelectromagnetism in Human Brain Research: New Applications, New Questions.}},
    author = {Gross, Joachim and Junghöfer, Markus and Wolters, Carsten},
    journal = {The Neuroscientist : a review journal bringing neurobiology, neurology and psychiatry},
    year = {2023},
    volume = {29},
    number = {1},
    pages = {62-77},
    doi = {10.1177/10738584211054742},
    issn = {1089-4098 (Electronic)},
    pmid = {34873945}
}

@article{
    Pascual2007sqrtm,
    author = {Pascual-Marqui, Roberto D},
    address = {Ithaca},
    copyright = {Notwithstanding the ProQuest Terms and conditions, you may use this content in accordance with the associated terms available at http://arxiv.org/abs/0710.3341.},
    issn = {2331-8422},
    journal = {arXiv.org},
    language = {eng},
    publisher = {Cornell University Library, arXiv.org},
    title = {{Discrete, {3D} distributed, linear imaging methods of electric neuronal activity. {P}art 1: exact, zero error localization}},
    year = {2007},
}

@article{
    geselowitz-1967,
    title = {{On Bioelectric Potentials in an Inhomogeneous Volume Conductor}},
    journal = {Biophysical Journal},
    volume = {7},
    number = {1},
    pages = {1-11},
    year = {1967},
    issn = {0006-3495},
    doi = {https://doi.org/10.1016/S0006-3495(67)86571-8},
    author = {David B. Geselowitz}
}

@article{
    elvetun2025depthbias,
    title = {{Weighted sparsity regularization for solving the inverse EEG problem: A case study}},
    journal = {Biomedical Signal Processing and Control},
    volume = {107},
    pages = {107673},
    year = {2025},
    issn = {1746-8094},
    doi = {10.1016/j.bspc.2025.107673},
    author = {Ole Løseth Elvetun and Niranjana Sudheer}
}

@article{
	Lahtinen2023,
	author = {Lahtinen, Joonas and Moura, Fernando and Samavaki, Maryam and Siltanen, Samuli and Pursiainen, Sampsa},
	address = {England},
	copyright = {2023 The Author(s). Published by IOP Publishing Ltd},
	issn = {1741-2560},
	journal = {Journal of neural engineering},
	language = {eng},
	number = {2},
	pages = {26005-},
	publisher = {IOP Publishing},
	title = {{In silico study of the effects of cerebral circulation on source localization using a dynamical anatomical atlas of the human head}},
	volume = {20},
	year = {2023},
}

@article{
	KrishnaswamyPavitra2017Seeo,
	author = {Krishnaswamy, Pavitra and Obregon-Henao, Gabriel and Ahveninen, Jyrki and Khan, Sheraz and Babadi, Behtash and Iglesias, Juan Eugenio and Hämäläinen, Matti S. and Purdon, Patrick L.},
	address = {United States},
	title = {{Sparsity enables estimation of both subcortical and cortical activity from MEG and EEG}},
	volume = {114},
	year = {2017},
	copyright = {Volumes 1–89 and 106–114, copyright as a collective work only; author(s) retains copyright to individual articles},
	issn = {0027-8424},
	publisher = {National Academy of Sciences},
	journal = {PROCEEDINGS OF THE NATIONAL ACADEMY OF SCIENCES OF THE UNITED STATES OF AMERICA},
	language = {eng},
	number = {48},
	pages = {E10465-E10474},
}

@article{
	Rezaei2021,
	author = {Rezaei, Atena and Lahtinen, Joonas and Neugebauer, Frank and Antonakakis, Marios and Piastra, Maria Carla and Koulouri, Alexandra and Wolters, Carsten H. and Pursiainen, Sampsa},
	address = {United States},
	title = {{Reconstructing subcortical and cortical somatosensory activity via the RAMUS inverse source analysis technique using median nerve SEP data}},
	volume = {245},
	year = {2021},
	copyright = {2021 The Authors},
	issn = {1053-8119},
	publisher = {Elsevier Inc},
	journal = {NeuroImage (Orlando, Fla.)},
	language = {eng},
	pages = {118726-},
}

@article{
	Fahimi2020ECoGvdEEGdeeppossible,
	author = {Fahimi Hnazaee, Mansoureh and Wittevrongel, Benjamin and Khachatryan, Elvira and Libert, Arno and Carrette, Evelien and Dauwe, Ine and Meurs, Alfred and Boon, Paul and Van Roost, Dirk and Van Hulle, Marc M.},
	address = {United States},
	title = {{Localization of deep brain activity with scalp and subdural EEG}},
	volume = {223},
	year = {2020},
	copyright = {2020},
	issn = {1053-8119},
	publisher = {Elsevier Inc},
	journal = {NeuroImage (Orlando, Fla.)},
	language = {eng},
	pages = {117344-},
}

@article{
	SSLOFO2005Liu,
	language = {eng},
	number = {10},
	pages = {1681-1691},
	publisher = {IEEE},
	author = {Hesheng Liu and Schimpf, P.H and Guoya Dong and Xiaorong Gao and Fusheng Yang and Shangkai Gao},
	address = {United States},
	copyright = {COPYRIGHT 2005 Institute of Electrical and Electronics Engineers, Inc.},
	issn = {0018-9294},
	journal = {{IEEE} Trans. Biomed. Eng.},
	title = {{Standardized shrinking {LORETA-FOCUSS} ({SSLOFO}): a new algorithm for spatio-temporal {EEG} source reconstruction}},
	volume = {52},
	year = {2005},
}

@article{
    vorwerk2019multipoles,
    title = {{{The multipole approach for EEG forward modeling using the finite element method}}},
    journal = {NeuroImage},
    volume = {201},
    pages = {116039},
    year = {2019},
    issn = {1053-8119},
    doi = {https://doi.org/10.1016/j.neuroimage.2019.116039},
    author = {Johannes Vorwerk and Anne Hanrath and Carsten H. Wolters and Lars Grasedyck}
}

@article{
    buchner1997stvenant,
    title = {{Inverse localization of electric dipole current sources in finite element models of the human head}},
    journal = {Electroencephalography and Clinical Neurophysiology},
    volume = {102},
    number = {4},
    pages = {267-278},
    year = {1997},
    issn = {0013-4694},
    doi = {https://doi.org/10.1016/S0013-4694(96)95698-9},
    author = {Helmut Buchner and Gunter Knoll and Manfred Fuchs and Adrian Rienäcker and Rainer Beckmann and Michael Wagner and Jiri Silny and Jörg Pesch}
}

@article{
    soderholm2026interference,
    title = {{A complete-electrode-model-based forward approach for transcranial temporal interference stimulation with linearization: A numerical simulation study}},
    journal = {Biomedical Signal Processing and Control},
    volume = {124},
    pages = {110624},
    year = {2026},
    issn = {1746-8094},
    doi = {https://doi.org/10.1016/j.bspc.2026.110624},
    author = {Santtu Söderholm and Maryam Samavaki and Sampsa Pursiainen},
}

@Article{
    hämäläinen1994mne,
    author={Hämäläinen, M. S. and Ilmoniemi, R. J.},
    title = {{Interpreting magnetic fields of the brain: minimum norm estimates}},
    journal={Medical {\&} Biological Engineering {\&} Computing},
    year={1994},
    month={1},
    day={01},
    volume={32},
    number={1},
    pages={35-42},
    issn={1741-0444},
    doi={10.1007/BF02512476},
    url={https://doi.org/10.1007/BF02512476}
}

@article{
    lin2006dsm,
    author = {Lin, Fa-Hsuan and Belliveau, John W. and Dale, Anders M. and Hämäläinen, Matti S.},
    title = {{Distributed current estimates using cortical orientation constraints}},
    journal = {Human Brain Mapping},
    volume = {27},
    number = {1},
    pages = {1-13},
    doi = {10.1002/hbm.20155},
    eprint = {https://onlinelibrary.wiley.com/doi/pdf/10.1002/hbm.20155},
    year = {2006}
}

@article{
    lin2021dsm,
    title = {{Distributed source modeling of intracranial stereoelectro-encephalographic measurements}},
    journal = {NeuroImage},
    volume = {230},
    pages = {117746},
    year = {2021},
    issn = {1053-8119},
    doi = {10.1016/j.neuroimage.2021.117746},
    author = {Fa-Hsuan Lin and Hsin-Ju Lee and Jyrki Ahveninen and Iiro P. Jääskeläinen and Hsiang-Yu Yu and Cheng-Chia Lee and Chien-Chen Chou and Wen-Jui Kuo}
}

@article{deGooijer2013,
author = {de Gooijer-van de Groep, Karin L. and Leijten, Frans S.S. and Ferrier, Cyrille H. and Huiskamp, Geertjan J.M.},
address = {New York, NY},
issn = {1065-9471},
journal = {Human brain mapping},
language = {eng},
number = {9},
pages = {2032-2044},
publisher = {Blackwell Publishing Ltd},
title = {{Inverse modeling in magnetic source imaging: Comparison of {MUSIC}, {SAM}(g2), and s{LORETA} to interictal intracranial {EEG}}},
volume = {34},
year = {2013},
}

@article{Coito2019,
author = {Coito, Ana and Biethahn, Silke and Tepperberg, Janina and Carboni, Margherita and Roelcke, Ulrich and Seeck, Margitta and Mierlo, Pieter and Gschwind, Markus and Vulliemoz, Serge},
address = {United States},
issn = {2470-9239},
journal = {Epilepsia open},
language = {eng},
number = {2},
pages = {281-292},
publisher = {John Wiley & Sons, Inc},
title = {{Interictal epileptogenic zone localization in patients with focal epilepsy using electric source imaging and directed functional connectivity from low‐density {EEG}}},
volume = {4},
year = {2019},
}

@article{Mouthaan2019,
author = {Mouthaan, Brian E. and Rados, Matea and Boon, Paul and Carrette, Evelien and Diehl, Beate and Jung, Julien and Kimiskidis, Vasilios and Kobulashvili, Teia and Kuchukhidze, Giorgi and Larsson, Pål G. and Leitinger, Markus and Ryvlin, Philippe and Rugg-Gunn, Fergus and Seeck, Margitta and Vulliémoz, Serge and Huiskamp, Geertjan and Leijten, Frans S.S. and Van Eijsden, Pieter and Trinka, Eugen and Braun, Kees P.J.},
address = {Netherlands},
copyright = {2019 International Federation of Clinical Neurophysiology},
issn = {1388-2457},
journal = {Clinical neurophysiology},
language = {eng},
number = {5},
organization = {on behalf of the E-PILEPSY consortium},
pages = {845-855},
publisher = {Elsevier B.V},
title = {{Diagnostic accuracy of interictal source imaging in presurgical epilepsy evaluation: A systematic review from the E-PILEPSY consortium}},
volume = {130},
year = {2019},
}

@article{Diamond2023,
    author = {Diamond, Joshua M and Withers, C Price and Chapeton, Julio I and Rahman, Shareena and Inati, Sara K and Zaghloul, Kareem A},
    title = "{Interictal discharges in the human brain are travelling waves arising from an epileptogenic source}",
    journal = {Brain},
    volume = {146},
    number = {5},
    pages = {1903-1915},
    year = {2023},
    month = {02},
    issn = {0006-8950},
    doi = {10.1093/brain/awad015},
}

@article{RineyCatherine2012,
author = {Riney, Catherine J and Chong, William K and Clark, Chris A and Cross, J. Helen},
address = {Ireland},
copyright = {Elsevier Ireland Ltd},
issn = {0720-048X},
journal = {Eur. J. Radiol.},
language = {eng},
number = {6},
pages = {1299-1305},
publisher = {Elsevier Ireland Ltd},
title = {{Voxel based morphometry of FLAIR MRI in children with intractable focal epilepsy: Implications for surgical intervention}},
volume = {81},
year = {2012},
}

@article{Aydin2017,
author = {Aydin, {\"U}. and Rampp, S. and Wollbrink, A. and Kugel, H. and Cho, J. -H. and Knösche, T. R. and Grova, C. and Wellmer, J. and Wolters, C. H.},
address = {New York},
copyright = {The Author(s) 2017},
issn = {0896-0267},
journal = {Brain topography},
language = {eng},
number = {4},
pages = {417-433},
publisher = {Springer US},
title = {{Zoomed {MRI} Guided by Combined {EEG}/{MEG} Source Analysis: A Multimodal Approach for Optimizing Presurgical Epilepsy Work-up and its Application in a Multi-focal Epilepsy Patient Case Study}},
volume = {30},
year = {2017},
}

@book{
    brenner2008fembook,
    author = {
        Susanne C. Brenner
        and
        L. Ridgway Scott
    },
    title = {{{The Mathematical Theory of Finite Element Methods}}},
    publisher = {Springer},
    year = {2008},
    doi = {10.1007/978-0-387-75934-0},
    
}

@article{
    fischl2012freesurfer,
    title = {{FreeSurfer}},
    journal = {NeuroImage},
    volume = {62},
    number = {2},
    pages = {774-781},
    year = {2012},
    note = {20 YEARS OF fMRI},
    issn = {1053-8119},
    doi = {10.1016/j.neuroimage.2012.01.021},
    author = {Bruce Fischl},
}

@article{
    puonti2020simnibs,
    title = {{Accurate and robust whole-head segmentation from magnetic resonance images for individualized head modeling}},
    journal = {NeuroImage},
    volume = {219},
    pages = {117044},
    year = {2020},
    issn = {1053-8119},
    doi = {https://doi.org/10.1016/j.neuroimage.2020.117044},
    url = {https://www.sciencedirect.com/science/article/pii/S1053811920305309},
    author = {Oula Puonti and Koen {Van Leemput} and Guilherme B. Saturnino and Hartwig R. Siebner and Kristoffer H. Madsen and Axel Thielscher},
}

@article{
    galazprieto2023zeffiromesh,
    doi = {10.1371/journal.pone.0290715},
    author = {Galaz Prieto, Fernando AND Lahtinen, Joonas AND Samavaki, Maryam AND Pursiainen, Sampsa},
    journal = {PLOS ONE},
    publisher = {Public Library of Science},
    title = {{Multi-compartment head modeling in EEG: Unstructured boundary-fitted tetra meshing with subcortical structures}},
    year = {2023},
    month = {09},
    volume = {18},
    pages = {1-25},
    number = {9},
}

@book{
    peterson2022basisfunctions,
    author = {Andrew F Peterson},
    title = {{Mapped vector basis functions for electromagnetic integral equations}},
    publisher = {Springer},
    year = 2022,
    isbn = {978-3-031-01686-8},
}

@article{LuckaFelix2012HBif,
    author = {Lucka, Felix and Pursiainen, Sampsa and Burger, Martin and Wolters, Carsten H.},
    title = {{Hierarchical Bayesian inference for the EEG inverse problem using realistic FE head models: Depth localization and source separation for focal primary currents}}, 
    language = {eng},
    number = {4},
    pages = {1364-1382},
    publisher = {Elsevier Inc},
    volume = {61},
    year = {2012},
    address = {United States},
    copyright = {2012 Elsevier Inc.},
    issn = {1053-8119},
    journal = {NeuroImage (Orlando, Fla.)}
}

@ARTICLE{
    bauer2015whitney,
    author={Bauer, Martin and Pursiainen, Sampsa and Vorwerk, Johannes and Köstler, Harald and Wolters, Carsten H.},
    journal={IEEE Transactions on Biomedical Engineering}, 
    title={{Comparison Study for Whitney (Raviart–Thomas)-Type Source Models in Finite-Element-Method-Based EEG Forward Modeling}}, 
    year={2015},
    volume={62},
    number={11},
    pages={2648-2656},
    doi={10.1109/TBME.2015.2439282},
}

@Article{
    weinstein2000reciprocity,
    author={Weinstein, David
    and Zhukov, Leonid
    and Johnson, Chris},
    title={Lead-field Bases for Electroencephalography Source Imaging},
    journal={Annals of Biomedical Engineering},
    year={2000},
    month={Sep},
    day={01},
    volume={28},
    number={9},
    pages={1059-1065},
    issn={1573-9686},
    doi={10.1114/1.1310220},
}

@INPROCEEDINGS{
	soderholm2026forwardinverse,
	AUTHOR={Santtu {Söderholm} and Joonas Lahtinen and Sampsa Pursiainen},
	TITLE="{Forward-Inverse} Interplay in {FEM-Based} {EEG} Source Imaging: Distributional Signatures of Advanced Source Models and Inverse Solvers",
	BOOKTITLE="2026 IEEE International Conference on Metrology for eXtended Reality, Artificial Intelligence and Neural Engineering (MetroXRAINE) (IEEE MetroXRAINE 2026)",
	ADDRESS="Chemnitz, Germany",
	PAGES=6,
	DAYS=19,
	MONTH=oct,
	YEAR=2026
}

@article {
    gramfort2010openmeeg,
    author={Gramfort, Alexandre and Papadopoulo, Th{\'e}odore and Olivi, Emmanuel and Clerc, Maureen},
    title={{OpenMEEG: opensource software for quasistatic bioelectromagnetics}},
    journal={BioMedical Engineering OnLine},
    year={2010},
    month={Sep},
    day={06},
    volume={9},
    number={1},
    pages={45},
    issn={1475-925X},
    doi={10.1186/1475-925X-9-45},
}

@article{
    tadel2011brainstorm,
    author = {Tadel, François and Baillet, Sylvain and Mosher, John C. and Pantazis, Dimitrios and Leahy, Richard M.},
    title = {Brainstorm: A User-Friendly Application for MEG/EEG Analysis},
    journal = {Computational Intelligence and Neuroscience},
    volume = {2011},
    number = {1},
    pages = {879716},
    doi = {https://doi.org/10.1155/2011/879716},
    year = {2011}
}

@article{
  rush1969reciprocity,
  author={Rush, Stanley and Driscoll, Daniel A.},
  journal={IEEE Transactions on Biomedical Engineering},
  title={{EEG Electrode Sensitivity-An Application of Reciprocity}},
  year={1969},
  volume={BME-16},
  number={1},
  pages={15-22},
  doi={10.1109/TBME.1969.4502598}
}

@article{
    vanrumste2004reciprocity,
    title={{The Validation of the Finite Difference Method and Reciprocity for Solving the Inverse Problem in EEG Dipole Source Analysis}},
    author={Bart Vanrumste and Gert Van Hoey and Rik Van de Walle and Michel D'Have and Ignace Lemahieu and Paul Boon},
    journal={Brain Topography},
    year={2004},
    volume={14},
    pages={83-92},
}

\end{document}